\documentclass[11pt,a4paper]{article}
\usepackage{amsthm,amsfonts,amsmath,amscd,amssymb}
\usepackage[english]{babel}
\usepackage[T1]{fontenc}
\usepackage[utf8]{inputenc}
\usepackage[
  margin=2cm,includefoot,footskip=30pt]{geometry}
\usepackage{lmodern}
\usepackage{bm}
\usepackage{booktabs}
\usepackage{array}
\usepackage{graphics,graphicx}
\usepackage[font=small]{caption}
\usepackage{cite}
\usepackage[dvipsnames]{xcolor}
\usepackage{tikz-cd}
\usetikzlibrary{decorations.markings,angles}
\usepackage{stmaryrd} \SetSymbolFont{stmry}{bold}{U}{stmry}{m}{n}
\usepackage{mathrsfs}
\usepackage{extarrows}
\usepackage{mathtools}
\usepackage{tensor}
\usepackage{tabularx}
\usepackage{enumerate}
\usepackage{adjustbox}
\usepackage[mathcal]{eucal}
\usepackage{bbold}
\usepackage{bm}
\usepackage{tabularx}
\usepackage{framed}

\newcolumntype{P}[1]{>{\centering\arraybackslash}p{#1}}
\newcolumntype{M}[1]{>{\centering\arraybackslash}m{#1}}

\usepackage[strict]{changepage}
\usepackage{tensor}
\usepackage[all]{xy}
\usepackage{ulem}

\usepackage{fancybox}  
\usepackage[breakable]{tcolorbox}  
\usepackage{graphicx}
\usepackage{tcolorbox}  

\usepackage{verbatim}

\usepackage{hyperref}
\hypersetup{
	colorlinks,
	urlcolor={red!55!black},
	citecolor={green!60!black},
	linkcolor={red!55!black}
}

\mathchardef\mhyphen="2D
\def\on{\operatorname}

\makeatletter
\providecommand{\leftsquigarrow}{%
  \mathrel{\mathpalette\reflect@squig\relax}%
}
\newcommand{\reflect@squig}[2]{%
  \reflectbox{$\m@th#1\rightsquigarrow$}%
}
\makeatother

\definecolor{ao}{rgb}{0.0, 0.5, 0.0}

\usepackage{mdframed}
\newtheorem{ftheo}{Theorem}

\newenvironment{fthm}
  {\begin{mdframed}[innertopmargin = 3pt, innerbottommargin=3pt,skipabove=5pt,skipbelow=5pt,linewidth=0.25pt,nobreak=true,align=center]\begin{ftheo}}
  {\end{ftheo}\end{mdframed}}

\usepackage[capitalize, noabbrev, nosort]{cleveref}
\newtheorem{theorem}{Theorem}[section]
\newtheorem{lemma}[theorem]{Lemma}
\newtheorem{conjecture}[theorem]{Conjecture}
\newtheorem{proposition}[theorem]{Proposition}
\newtheorem{corollary}[theorem]{Corollary}

\theoremstyle{definition}
\newtheorem{construction}[theorem]{Construction}
\newtheorem{definition}[theorem]{Definition}
\newtheorem{notation}[theorem]{Notation}
\newtheorem{remark}[theorem]{Remark}
\newtheorem{example}[theorem]{Example}

\newcommand\noloc{%
  \nobreak
  \mspace{6mu plus 1mu}
  {:}
  \nonscript\mkern-\thinmuskip
  \mathpunct{}
  \mspace{2mu}
}

\newcommand{\rgraph}{{\bf G}}
\newcommand{\Wgraph}{{\bf G_{\mathfrak{w}}}}
\newcommand{\gw}{{\bf G_{\mathfrak{w}}}}
\newcommand{\simp}{{\bm L}}

\newcommand{\lsimp}{\mathbb{S}}
\newcommand{\rsimp}{\reflectbox{$\mathbb{S}$}}
\newcommand{\lsilt}{\mathbb{P}}
\newcommand{\rsilt}{\reflectbox{$\mathbb{P}$}}
\newcommand{\cusilt}{\reflectbox{$\mathbb{P}$}^\uparrow}
\newcommand{\cdsilt}{\reflectbox{$\mathbb{P}$}^\downarrow}
\newcommand{\susilt}{\mathbb{P}^\uparrow}
\newcommand{\sdsilt}{\mathbb{P}^\downarrow}
\newcommand{\ssimp}{\sS}
\newcommand{\csimp}{\backS}
\newcommand{\backZ}{\reflectbox{$\mathbb{Z}$}}
\newcommand{\zssimp}{\mathbb{Z}}
\newcommand{\zcsimp}{\backZ}

\newcommand{\C}{\mathcal{C}}
\newcommand{\D}{\mathcal{D}}
\newcommand{\E}{\mathcal{E}}
\newcommand{\F}{\mathcal{F}}
\newcommand{\V}{\mathcal{V}}
\newcommand{\M}{\mathcal{M}}
\newcommand{\N}{\mathcal{N}}

\newcommand{\cC}{\mathbb{C}}
\newcommand{\Z}{\mathbb{Z}}

\newcommand{\glsec}{R\Gamma}
\newcommand{\losec}{\mathcal{L}}
\newcommand{\glsecF}{\glsec(\Wgraph,\F_\w)}

\newcommand{\Glsec}[2]{\glsec(#1,#2)}
\newcommand{\Losec}[2]{\losec(#1,#2)}

\newcommand{\g}{\gamma}
\newcommand{\w}{\mathfrak{w}}

\newcommand{\st}[1]{{\Delta_#1}}
\newcommand{\cst}[1]{{\nabla_#1}}
\newcommand{\evinf}{\mbox{ev}_{\vinf}}

\renewcommand{\a}{\alpha}
\renewcommand{\b}{\beta}
\newcommand{\lr}{\longrightarrow}
\newcommand{\s}{\sigma}
\newcommand{\sS}{\mathbb{S}}
\newcommand{\backS}{\reflectbox{$\mathbb{S}$}}

\newcommand{\sse}{\subseteq}
\newcommand{\wsch}{{\F_{\mathfrak{w}}}}
\newcommand{\tw}[1]{{T_{S_{#1}}}}
\newcommand{\itw}[1]{{T^{-1}_{S_{#1}}}}

\newcommand{\vd}{\mathsf{v}}
\newcommand{\ed}{\mathsf{e}}

\newcommand{\sM}{\mathcal{M}}
\newcommand{\SE}{\mathcal{E}}

\newcommand{\la}{\lambda}
\newcommand{\e}{\varepsilon}

\newcommand{\dd}{\partial}

\newcommand{\Br}{\text{Br}}
\newcommand{\SL}{\mathcal{L}}

\newcommand{\vinf}{v_{-\infty}}

\newcommand{\einftop}{e_{\infty}}
\newcommand{\bw}{(\beta,{\bf a})}
\newcommand{\ba}{\beta,{\bf a}}
\newcommand{\lbw}{\SL_{\beta,{\bf a}}}
\newcommand{\cbw}{\mathcal{C}_{\beta,{\bf a}}}
\newcommand{\Lbw}{L_{\beta,{\bf a}}}
\newcommand{\Mor}{\on{Mor}}
\newcommand{\LC}{\bm L^\w}
\renewcommand{\k}{\kappa}

\newcommand{\stds}{\mathbb{\Delta}}

\newcommand{\fib}{\operatorname{fib}}
\newcommand{\cof}{\operatorname{cof}}
\newcommand{\mor}{\operatorname{Mor}}

\title{Categorical Lusztig cycles and weave schobers}
\author{Roger Casals and Merlin Christ}
\date{}

\begin{document}
\maketitle

\begin{abstract}
We establish the foundations of categorical weave calculus, developing the diagrammatic calculus of weaves and braid varieties within the study of Calabi-Yau triangulated categories and cluster tilting theory. This is achieved by associating a perverse sheaf of triangulated categories to each Demazure weave. A central contribution is the construction and study of the categorical Lusztig cycles and their duals, which we show form simple-minded and silting collections in the category of global sections of such a sheaf of categories. These categorical collections are built using the diagrammatics of weaves and we study their behavior under changes of weaves. For instance, we show that they undergo tilts under weave mutations. En route, we develop the study of categorical weighted braid words, as canonical rigid filtered dg modules over derived preprojective algebras, and the categorical incarnation of the tropical Lusztig rules, as a gluing mechanism for such filtered objects.
Appendix A also contains homological results, providing a novel construction of simple-minded and silting collections from full exceptional collections, and characterizing when these arise from a highest weight structure on an abelian category.
\end{abstract}

\setcounter{tocdepth}{3}
\tableofcontents

\section{Introduction}\label{sec:intro}

The object of this article is to establish the foundations of categorical weave calculus. This allows for many of the applications of weave calculus and braid varieties to be elevated to the study of relative Calabi-Yau categories and cluster tilting theory. We first associate certain perverse schobers to Demazure weaves, which we refer to as weave schobers. A central contribution is then the construction and study of the categorical Lusztig cycles, which are objects in the triangulated category of global section of such a weave schober and generalize the combinatorial Lusztig weave cycles. We show that the categorical Lusztig cycles of a weave schober form a simple-minded collection, are invariant under weave equivalence, and undergo a simple tilt under weave mutation. En route, we also introduce and develop the study of categorical weighted braid words, as canonical filtered objects in the derived category of the preprojective algebra, the categorical incarnation of the tropical Lusztig rules, as a gluing mechanism for such filtered objects, and the relation between mutations of exceptional collections and mutations of their simple-minded collections.\\

We also establish corresponding results for the dual categorical Lusztig cycles, which form a silting collection Koszul dual to the categorical Lusztig cycles, and construct explicit weave realizations. In a nutshell, we provide a complete foundational set of results enhancing weave calculus and the study of braid varieties to the context of 3-Calabi--Yau categories and cluster tilting theory. In addition, the results on weave schobers established in this article have a number of applications, some of which will be further explored in subsequent work. These include the construction of additive categorifications for the cluster algebras on braid varieties, and a form of homological mirror symmetry between the Fukaya-Seidel categories associated to a weave and the non-commutative crepant resolutions associated to weave quivers with potentials.


\subsection{Scientific context}\label{ssec:scientific_context} 
\begin{center}
	\begin{figure}[ht!]
		\centering
		\includegraphics[scale=1]{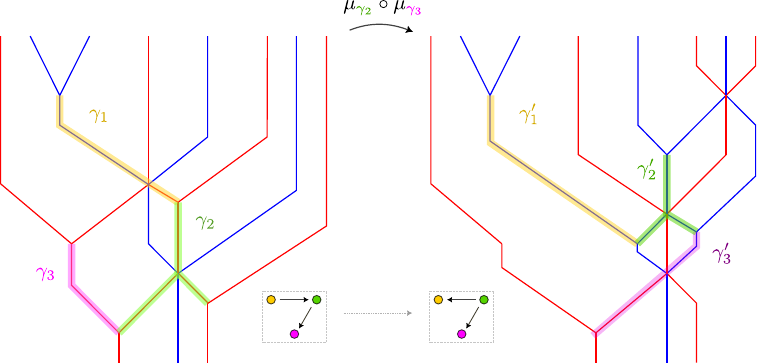}
		\caption{Two instances of Demazure weaves $\w,\w':\beta\lr\beta'$ from $\beta=s_2s_1(s_1s_2)^3\in\Br^+_3$ to $\beta'=s_2s_1s_2$, where blue weave edges are labeled by $s_1\in S_3$ and red weave edges by $s_2\in S_3$. (Left) The mutable Lusztig cycles $\gamma_1,\gamma_2,\gamma_3$ for $\w$ are highlighted in yellow, green and pink respectively, and their intersection quiver is displayed in the square box at the lower right of $\w$.  (Right) The weave $\w'$ obtained from $\w$ by performing weave mutations along the Lusztig cycles $\gamma_3$ and $\gamma_2$, in this order. Its mutable Lusztig cycles $\gamma'_1,\gamma'_2,\gamma'_3$ are also highlighted and their intersection quiver is displayed to the lower left of $\w'$.}
		\label{fig:weave_intro_example2}
	\end{figure}
\end{center}

Weaves were introduced in \cite{CasalsZaslow} within the context of symplectic topology. They have since been developed in different directions, including ones within combinatorial and algebraic frameworks, cf.~e.g.~\cite{CGGJ24,CGGLSS25}. In this context, the study of weaves can be understood as a planar diagrammatic calculus, with certain moves \& mutations between weaves, and the combinatorics of their cycles and associated functions on braid varieties. See \cref{fig:weave_intro_example2} for two instances of weaves related by weave mutations, and their Lusztig cycles. The basic results on the diagrammatics of weaves, including their moves \& mutations, and the Lusztig cycles and their intersection theory, are summarized in \cite[Section 4]{CGGLSS25}, see also \cite{CGGJ24,CasalsWeng22,CasalsZaslow}. This collection of results and diagrammatic techniques on weaves, which is centered around the Lusztig cycles, is often referred to as weave calculus, cf.~e.g.~\cref{ssec:primer_weavecalculus}.

Weave calculus has been successfully used to tackle different problems, some of which include the construction of cluster algebras on braid varieties \cite{CGGLSS25,CGGJ24,casals2025comparingclusteralgebrasbraid}, a proof that the Muller-Speyer twist is Donaldson-Thomas and that source and target positroids seeds quasi-coincide \cite{casals2023demazureweavesreducedplabic}, the Floer theoretic characterization of Gaiotto-Moore-Neitzke spectral networks \cite{casals2025spectralnetworksbettilagrangians,GNMSN,GMN14_Snakes}, and the microlocal study of Legendrian links \cite{CasalsWeng22,casals2022conjugate,CasalsLi24}. In interesting cases, including all finite and affine Dynkin types, all cluster seeds and the cluster automorphism groups can be realized by weaves \cite{AnBaeLee21ADE,Hughes21A,Hughes21D}, and see also \cite{asplund2025decompositionsaugmentationvarietiesweaves,CapovillaHughesWeng25,hughes2025legendriandoublestwistspuns} for more applications of weave calculus. We emphasize that in all these results, the study of the Lusztig cycles for the corresponding weaves is a key ingredient, as it provides all the necessary data for the cluster seeds, including the quivers and both types of cluster variables.\\

The main goal of this article is to develop the foundations of categorical weave calculus, which is centered around the categorical Lusztig cycles. The applications above hint towards the existence of such a diagrammatic categorical calculus, but at the same time indicate that its development requires new ideas and techniques. Indeed, on the one hand, such calculus would provide a diagrammatic presentation of many interesting Fukaya-Seidel categories and parts of their spaces of stability conditions, functorially interact with the additive cluster categorifications, and generalize the categorification of higher-rank Fock-Goncharov theory of local systems, cf.~\cite{Chr25c,KL25}, to irregular Betti surfaces. On the other hand, at a technical level, it should provide a categorical lift of the tropical Lusztig moves for weighted braid words, the intersection theory of weave cycles in Demazure weaves from \cite{CasalsZaslow,CasalsWeng22,CGGLSS25}, and feature a triangulated counterpart of the theory of highest weight abelian categories.\\

\noindent The manuscript addresses many of these aspects by developing new ideas and techniques. These include the construction of weave perverse schobers, the study of rigid filtered objects associated to weighted braid words, and new homological insights in the context of $k$-linear stable $\infty$-categories, such as the categorical lift of the tropical Lusztig rules in terms of the propagation of sections of weaves schobers, and the relation between mutations of exceptional collections and mutations for their associated simple-minded and silting collections.\\

Regarding the technical implementation of such ideas, the categorical weave calculus we present uses the formalism of perverse schobers, introduced in \cite{KS14}. Intuitively, the perverse schobers we use are a certain class of constructible sheaves of $k$-linear stable $\infty$-categories on a surface. Specifically, we use the theory of perverse schobers parametrized by ribbon graphs embedded in surfaces, introduced in \cite[Sections 3\&4]{Chr22}, \cite[Section 3]{CHQ23}, see also \cref{ssec:perverse_schober_preliminaries} below, which provides the key ingredients on perverse schobers needed to develop our results.\\

\subsection{The main results} 

A central point illustrated by our results is that

\begin{center}
{\it There exists a categorical enhancement of weave calculus, constructed using perverse schobers. In it, Lusztig cycles and their calculus have canonical enhancements to a relative CY3 category.}
\end{center}

\noindent As evidence, we now highlight three of our main results, as follows.\\

{\bf (A)} Let us first start with the categorical enhancement of the Lusztig cycles associated to a Demazure weave $\w$. By definition, a weighted braid word $(\beta,{\bf a})$ is a pair consisting of a positive braid word $\beta=s_{i_0}\cdots s_{i_l}$ and a collection of non-negative integers ${\bf a}\coloneqq (a_0,\ldots,a_l)$, where $s_j$ denotes the $j$th Artin generator of the braid group in Type A, i.e.~corresponding to the $j$th vertex of the Dynkin diagram. First, as summarized in \cref{sssec:Lusztigcycles}, each non-categorical Lusztig cycle of a weave $\w$ can be understood as a certain sequence of weighted braid words $(\beta_1,{\bf a}_1)\to\cdots\to(\beta_k,{\bf a}_k)$ whose weights propagate downwards through $\w$ using the following tropical Lusztig rules:
\begin{multline}
\label{eq:intro_lusztig_trop}
\qquad \qquad \qquad \qquad \quad \Phi_1:(a_1,a_2)\mapsto \min(a_1,a_2), \quad \Phi_2:(a_1,a_2)\mapsto (a_2,a_1),\\ \qquad\qquad \Phi_3: (a_1,a_2,a_3)\mapsto  \left(a_2+a_3-\min(a_1,a_3),\min(a_1,a_3),a_1+a_2-\min(a_1,a_3)\right),
\end{multline}
where $\Phi_1$ is used at trivalent vertices, the braid word changing according to $s_i^2\lr s_i$, $\Phi_2$ at tetravalent vertices, the braid word changing via $s_is_k=s_ks_i$ with $i,k$ non-adjacent in the Dynkin diagram, and $\Phi_3$ at hexavalent vertices, where the braid word changes via $s_is_js_i=s_js_is_j$ with $i,j$ adjacent.\\

\noindent Therefore, there are two tasks to obtain the right notion of a categorical Lusztig cycle:

\begin{itemize}
    \item[(i)] The first task is to construct the categorical version of a weighted braid word $(\beta,{\bf a})$. This is achieved in \cref{sec:filtrations_from_braid_words}, where we construct a canonical rigid filtered object $L_{\beta,{\bf a}}$ associated to any weighted braid word $(\beta,{\bf a})$ arising from Lusztig cycles:
    \begin{equation}\label{eq:categorical_weighted_braid_intro}
        \fbox{$\exists!$ rigid filtered object $L_{\beta,{\bf a}}\in\D$ whose graded pieces recover $\bw$}
    \end{equation}
    This is established in \cref{thm:unique_rigid_resolutions}, which refines the above statement. Specifically, the object $\Lbw$ belongs to the derived $\infty$-category $\D=\D^{\on{perf}}(\Pi_2(A_{n}))$ of the additive derived preprojective algebra of type $A_n$.     The existence and uniqueness of the object $\Lbw$ is a novel result within the representation theory of $\Pi_2(A_{n})$. For its construction, we will rely on a partial geometric model for $\D$ in terms of matching $2$-spheres and matching paths in the disc, which allows us to express $\Lbw$ as a sequence of resolutions for the intersections of matching paths.

    \item[(ii)]     The second task is to use the sequence $(\beta_1,{\bf a}_1)\to\cdots\to(\beta_k,{\bf a}_k)$ associated to a Lusztig cycle to glue these filtered objects $L_{\beta_i,{\bf a}_i}$ together into an object of a certain category. This should be done in a manner that the tropical Lusztig rules displayed in \cref{eq:intro_lusztig_trop} are recovered when considering graded pieces. This step requires constructing such a category from a Demazure weave $\w$ first: we achieve this in \cref{sec:weave_schobers} by constructing a perverse schober $\F_\w$ associated to a weave $\w$. In this context, the required gluing of the filtered objects is achieved by establishing certain fiber and cofiber sequences between them, as discussed in \cref{sec:lusztigcycles}. For instance, if $(\beta_i,{\bf a}_i)\to (\beta_{i+1},{\bf a}_{i+1})$ is a trivalent vertex, \Cref{prop:semi_twist_Lusztig_cycle_at_trivalent_vertex} will establish a cofiber sequence
    \begin{equation}\label{eq:cofibersequence_trivalent_intro}
        \fbox{At a trivalent: $\tau_{\geq 0}\on{Mor}(S_{p},L_{\beta_i,{\bf a}_i})\otimes S_{p} \longrightarrow L_{\beta_i,{\bf a}_i}\longrightarrow L_{\beta_{i+1},{\bf a}_{i+1}}$}
    \end{equation}
    allowing us to glue the associated filtered objects at that trivalent vertex, as local sections of the weave schober. Similarly, if $(\beta_i,{\bf a}_i)\to (\beta_{i+1},{\bf a}_{i+1})$ is a hexavalent or tetravalent vertex then gluing occurs via the following equivalences:
    \begin{equation}\label{eq:filteredobjects_hexavalent_intro}
    \fbox{At a non-trivalent vertex: $L_{\beta_i,{\bf a}_i}\simeq L_{\beta_{i+1},{\bf a}_{i+1}}$}
    \end{equation}
    The equivalences in (\ref{eq:filteredobjects_hexavalent_intro}) are of the underlying objects, which suffices to glue, but the filtrations themselves are typically different. Such equivalences are established in \cref{prop:4_and_6_moves}. The result of such gluings, using (\ref{eq:categorical_weighted_braid_intro}), (\ref{eq:cofibersequence_trivalent_intro}) and (\ref{eq:filteredobjects_hexavalent_intro}), is so that the categorical Lusztig cycles shall be objects in the category of global sections $\glsecF$, cf.~\cref{def:categorical_Lusztig_cycles}, which will be a smooth and proper $k$-linear stable $\infty$-category with a relative $3$-Calabi--Yau structure.
\end{itemize}

In a nutshell, these constructions in Sections \ref{sec:weave_schobers}, \ref{sec:lusztigcycles} and \ref{sec:filtrations_from_braid_words} show that categorical Lusztig cycles exist and abide by a categorical enhancement of the tropical Lusztig propagation rules for weighted braid words. Though establishing these results requires developing new results that might be of interest on their own, as illustrated by the content of these sections, the first main result can be summarized as follows:\\

\begin{fthm}[Categorical Lusztig cycles]\label{thm:main1} Let $\w:\beta\to\delta(\beta)$ be a Demazure weave with $m$ trivalent vertices and $\F_\w$ its weave schober. Then:

\begin{enumerate}[$(1)$]

\item Let $(\beta,{\bf a})$ be a weighted braid word arising from a horizontal slice of a Lusztig cycle of $\w$. Then, up to equivalence, there exists a unique $(\beta,{\bf a})$-filtered object whose colimit $\Lbw \in \D$ is rigid. In addition, its graded pieces uniquely recover the weighted braid word $(\beta,{\bf a})$.

\item Let $(\beta_i,{\bf a}_i)\to (\beta_{i+1},{\bf a}_{i+1})$ be the two weighted braid words arising from a horizontal slice of a Lusztig cycle of $\w$ above and under a vertex $p\in\Wgraph$. Then

\begin{equation}\label{eq:main_equivalences}
L_{\beta_{i+1},{\bf a}_{i+1}}\simeq\begin{cases}
\cof(\tau_{\geq 0}\on{Mor}(S_{p},L_{\beta_i,{\bf a}_i})\otimes S_{p} \longrightarrow L_{\beta_i,{\bf a}_i}) & \mbox{ if }p\mbox{ is trivalent},\\
L_{\beta_i,{\bf a}_i} & \mbox{ if }p\mbox{ is non-trivalent}.
\end{cases}
\end{equation}

\noindent In addition, the categorical equivalences \eqref{eq:main_equivalences} recover the tropical Lusztig propagation rules at the vertex $p$ by considering the associated graded pieces of $\mathcal{L}_{\beta_i,{\bf a}_i}$ and $\mathcal{L}_{\beta_{i+1},{\bf a}_{i+1}}$.

\item For each $i\in[1,m]$, there exists a global section $\bm L_i^\w\in\glsecF$ whose value at the stalk of height $h$ is given by $L_{\beta,{\bf a}}$, where $(\beta,{\bf a})$ is the weighted braid word associated to the combinatorial Lusztig weave cycle $\g^\w_i$ at height $h$.
\end{enumerate}
\noindent By definition, the objects $\{\LC_i\}_{i\in[1,m]}$ are referred to as the categorical Lusztig cycles.
\end{fthm}

\noindent Theorem \ref{thm:main1} is proven in Sections \ref{sec:lusztigcycles}, \ref{sec:filtrations_from_braid_words} and \ref{sec:weave_propagation_LC}, once the weave schober $\F_\w$ is introduced and studied in \cref{sec:weave_schobers}. Three comments on Theorem \ref{thm:main1} for now:
\begin{itemize}
    \item[$(i)$] We also show that the rigid geometric realization $\Lbw$ in Theorem \ref{thm:main1}.(1) belongs to the subcategory $\mathcal{M}\sse\D$ of matching objects, thus giving a geometric model for $\Lbw$, cf.~\cref{sec:filtrations_from_braid_words}.
    
    \item[$(ii)$] The morphism in \eqref{eq:main_equivalences}, whose cofiber measures how the objects $\Lbw$ propagate through a trivalent vertex of $\w$, is the non-negative truncation $\tau_{\geq0}(\on{ev})$ of the evaluation map of the derived Hom. Its cofiber categorically represents a truncated spherical twist, referred to as a semi-twist. Different versions of semi-twists are ubiquitous in categorical weave calculus, as they categorically capture the tropical propagation rules for weave cycles. 
    
    \item[$(iii)$] Intuitively, Theorem \ref{thm:main1}.(3) states that the canonical rigid objects $L_{\beta_1,{\bf a}_1},\ldots,L_{\beta_k,{\bf a}_k}$ arising via Theorem \ref{thm:main1}.(1) from scanning a given Lusztig cycle of $\w$, can be glued together to a global section $\bm L_i^\w\in\glsecF$ of the weave schober $\F_\w$. The equivalences in Theorem \ref{thm:main1}.(2) are a key part of this gluing and the well-definedness of the resulting categorical Lusztig cycles.
\end{itemize}


{\bf (B)} Let us now discuss the categorical enhancement of weave calculus. Theorem \ref{thm:main1} provides the categorical Lusztig cycles, with the basic properties that indeed make them categorical avatars of the classical Lusztig weave cycles. A central feature of weave calculus are the notions of weave equivalences and weave mutations. These are summarized in \cref{sssec:weave_equivalences_mutations}, see particularly Figures \ref{fig:weave_equivalences} and \ref{fig:weave_mutation}. In this regard, the first task to start developing a categorical weave calculus is to ensure that:

\begin{itemize}
    \item[(i)] The relevant categories, objects and properties associated to a Demazure weave $\w$ should be {\it invariant under weave equivalences}. In particular, the category $\glsecF$ and its categorical Lusztig cycles $\{\LC_i\}_{i\in[1,m]}$ featuring in Theorem \ref{thm:main1} should be invariant under weave equivalences.

    \item[(ii)] There must be a categorical enhancement for weave mutation, so that we have precise formulas describing how the relevant categories, objects and homological properties associated to a weave $\w$ change under weave mutations. In particular, we should understand how the category $\glsecF$ and its categorical Lusztig cycles $\{\LC_i\}_{i\in[1,m]}$ change under weave mutations.
\end{itemize}

We will establish these two items throughout the manuscript, laying the foundations for a categorical weave calculus. Specifically, though the weave schober $\F_\w$ itself is not directly a weave invariant, we show in \cref{sec:weave_schobers} that its $\infty$-category of global sections $\glsecF$ is invariant under both weave equivalences and weave mutations, cf.~\cref{thm:weaveequivglsec}. The behavior of categorical Lusztig cycles is harder to establish, cf.~\cref{sec:proof_Lusztig_cycles}, but the outcome is neat:

\begin{center}
{\it The categorical Lusztig cycles form a simple-minded collection in $\glsecF$, and are invariant under weave equivalences. Plus, weave mutations induce mutations of such simple-minded collections.}
\end{center}

\noindent The proof that the categorical Lusztig cycles form a simple-minded collection, in the sense of \Cref{def:recollection_exceptional_SMC_silting}.(4), is contentful. It requires us to construct certain exceptional collections in $\glsecF$ from which $\{\LC_i\}_{i\in[1,m]}$ can be constructed, which we do in \cref{sec:weave_thimbles} and \cref{ssec:LC_from_thimbles}. The relation between the exceptional collection and the categorical Lusztig cycles is similar to the relation between the standard objects and the simple objects in a highest weight abelian category. (The precise relation to the theory of highest weight categories is explored in Appendix \ref{ssec:AppendixA}.) Nevertheless these exceptional collections are {\it not} invariant under weave equivalence and thus, in a sense, are not well-defined from the viewpoint of weave calculus. We overcome this fact in Sections \ref{ssec:Lusztigcycles_weaveequivalence} and \ref{ssec:Lusztigcycles_weavemutation} by studying how mutations of exceptional collections affect the associated simple-minded collections, which can be seen as a homological result of interest on its own.\\

At core, the second result establishing a categorical weave calculus can be summarized as follows:\\

\begin{fthm}[Categorical weave calculus]\label{thm:main2} Let $\w:\beta\to\delta(\beta)$ be a Demazure weave and $\F_\w$ its weave schober. Then the following holds:

\begin{enumerate}[$(1)$]
    \item The stable $\infty$-category of global sections $\glsecF$ is invariant under weave equivalence and weave mutation. In particular, it is an invariant of the boundary braid $\beta$.

    \item The categorical Lusztig cycles $\{\LC_i\}_{i\in[1,m]}$ form a simple-minded collection in $\glsecF$.

    \item Let $\w\sim\w'$ be a weave equivalence. Then the collections of categorical Lusztig cycles $\{\LC_i\}$ and $\{\bm L_i^{\w'}\}$ coincide, up to a permutation of the indices $i\in[1,m]$.

    \item Let $\w\to\w'$ be a forward weave mutation. Then $\{\bm L_i^{\w'}\}$ is obtained from $\{\LC_i\}$ by a forward mutation of simple-minded collections.
\end{enumerate}
\end{fthm}

\noindent Theorem \ref{thm:main2} is established in \cref{sec:proof_Lusztig_cycles}. Note that Theorem \ref{thm:main2}.(4) implies that $\{\LC_i\}$ undergoes a backward mutation of simple-minded collections under a backward weave mutation of $\w$. In a nutshell, Theorem \ref{thm:main2} establishes that $\{\LC_i\}$, built via Theorem \ref{thm:main1}, form a simple-minded collection in $\glsecF$, well-defined for the weave equivalence class of $\w$, and weave mutations lead to a form of tilting for $\{\LC_i\}$.\\

{\bf (C)} Finally, we also introduce and develop the dual notion of categorical Lusztig cycles in \cref{sec:silting_dualLusztig}, which we refer to as categorical dual Lusztig cycles. Note that the notion of dual Lusztig cycles is new even within non-categorical weave calculus, as no weave cycles for the cluster $\mathcal{A}$-variables were introduced in \cite{CGGLSS25}.\footnote{Intuitively, Lusztig cycles yield cluster $\mathcal{X}$-variables and dual Lusztig cycles correspond to cluster $\mathcal{A}$-variables. Categorically, the categorical Lusztig cycles additively categorify the $\mathcal{X}$-variables, and their categorical duals additively categorify the $\mathcal{A}$-variables, via the cluster tilting object.} Since the Lusztig cycles propagate downwards through the weave and are, coarsely said, as small as possible, the dual Lusztig cycles must propagate the weave upwards and be as large as possible, e.g.~they never have compact support.\\

Categorically, the canonical dual of the categorical Lusztig cycles $\{\LC_i\}$, considered as a simple-minded collection, is the unique Koszul dual silting collection $\{\lsilt^\w_i\}$ determined by
\begin{equation}\label{eq:dual_collection}
\on{Ext}^l(\lsilt^\w_i,\LC_j)\simeq \begin{cases} k & l=0\text{ and }i=j \\ 0 & \text{else}\,.\end{cases}
\end{equation}
The Koszul duality referred to above consists of a bijection between silting collections and simple-minded collections which interchanges the respective notions of mutations, see also \Cref{thm:silting_vs_SMC} and \cite{KN_unpublished}.\\

Nevertheless, the description via \eqref{eq:dual_collection} does not allow us to compute or manipulate categorical dual Lusztig cycles via weave calculus, as it is significantly unclear how $\{\lsilt^\w_i\}$ would be described from the given weave $\w$. In fact, it can be shown that the objects $\{\lsilt^\w_i\}$ are typically \textit{not} represented by weave cycles in $\w$, since global sections of $\F_\w$ propagate downwards. We resolve this mismatch, by introducing a duality on weave schobers. Specifically, for each weave $\w$, we introduce a different weave schober $\mathcal{F}_\w^{\uparrow}$ together with two equivalences of global sections
\begin{equation}\label{eq:intro_ringel_duality}
(\iota^{\uparrow})^R\circ \iota^{\downarrow},\ (\iota^{\uparrow})^L\circ \iota^{\downarrow}:\glsecF\lr \Glsec{\rgraph_\w^{\uparrow}}{\mathcal{F}_\w^{\uparrow}}\,,
\end{equation}
which intuitively capture the categorical action of rotating counterclockwise and clockwise by $\pi$. Diagram \eqref{eq:diagram_twosided_intro} provides the categories and functors involved to define the equivalences in \eqref{eq:intro_ringel_duality}. We will show that the equivalences \eqref{eq:intro_ringel_duality} behave like a version of the Ringel duality functor from highest weight theory in the context of stable $\infty$-categories of global sections of perverse schobers, see \cref{lem:Ringel_duality_and_silting_objects,lem:standard_vs_costandard_upwards_silt} or the discussion in \Cref{ssec:tiltingtheory_exceptionalcollections}. 

\begin{equation}\label{eq:diagram_twosided_intro}
\begin{tikzcd}
\fib(\on{ev}_{e_\infty})\simeq\Glsec{\rgraph_\w^\downarrow}{\mathcal{F}_\w^\downarrow} \arrow[r, "\iota^{\downarrow}", hook] \arrow[rr, "(\iota^{\uparrow})^L\circ \iota^{\downarrow}"', bend right=49]\arrow[rr, "\simeq", bend right=49] \arrow[rr, "(\iota^{\uparrow})^R\circ \iota^{\downarrow}", bend left=49] \arrow[rr, "\simeq"', bend left=49] & \Glsec{\rgraph_\w^{\updownarrow}}{\mathcal{F}_\w^{\updownarrow}} \arrow[l, "(\iota^{\downarrow})^L"', two heads, bend right] \arrow[l, "(\iota^{\downarrow})^R", two heads, bend left] \arrow[r, "(\iota^{\uparrow})^L"', bend right] \arrow[r, "(\iota^{\uparrow})^R", two heads, bend left] & \Glsec{\rgraph_\w^{\uparrow}}{\mathcal{F}_\w^{\uparrow}}\simeq\fib(\on{ev}_{e_1}) \arrow[l, "\iota^{\uparrow}"', hook']
\end{tikzcd}
\end{equation}

\noindent The remarkable fact is then that the image $\rsilt^\uparrow_i$ of $\lsilt^\w_i$ under the duality $(\iota^\uparrow)^R\circ \iota^\downarrow$ is representable by an upwards facing weave cycle on $\w$, even if $\lsilt^\w_i$ is not. In fact, the proofs of our results provide an explicit formula for the weights of this weave cycle, cf.~\cref{rem:description_weave_cycle_of_silt}. The main statement regarding the silting collection of categorical dual Lusztig cycles is summarized as:

\begin{fthm}[Weave realization of the silting collection]\label{thm:main3} Let $\w:\beta\to\delta(\beta)$ be a Demazure weave and $\F_\w$ its weave schober. Let $\{\lsilt^\w_i\}$ be the silting collection Koszul dual to the simple-minded collection of categorical Lusztig cycles $\{\LC_i\}$. Then, for any $i\in[1,m]$, the following holds:

\begin{enumerate}[$(1)$]
    \item The Ringel-type dual $\rsilt^\uparrow_i$ is representable by a positive weave cycle $\rho^\w_i:E(\w)\lr\Z_{\geq0}$.
    
     \item In addition, $\rho^\w_i$ can be obtained by starting at the trivalent vertex $p_i\in\w$ and propagating upwards with the tropical Lusztig rules, where at a trivalent vertex $p$ with southern edge $\on{e}_{s}\in E(\w)$, the upwards right edge $\on{e}_{ne}$ has weight
\begin{equation}\label{def:y_intro}
\framebox{$\rho^\w_i(\on{e}_{ne})=\rho^\w_i(\on{e}_{s})+\on{dim}\on{Ext}^1_\D(S_{p},\on{tot}(\C_{\ast>b(p)}))$}\,
\end{equation}
while the upward left edge has weight $\rho^\w_i(\on{e}_{s})$.
    
    \item There exists an equivalence $\on{ev}_{e_1}(\lsilt_i^\w)\simeq \on{ev}_{e_\infty}(\rsilt_i^\uparrow)$ in $\D$. In particular, the evaluation of $\lsilt_i^\w$ at the bottom edge $e_1$ is a rigid object which admits a $(\beta,{\bf a})$-filtration, where the weights ${\bf a}={\bf a}(\rho^\w_i)$ for $\beta$ are determined by the weave cycle $\rho^\w_i$ at the top horizontal slice.
\end{enumerate}
\end{fthm}

 \noindent In Theorem \ref{thm:main3}.(2), $S_{p}\in\D$ is the 2-spherical object associated to the trivalent vertex $p$, cf.~\cref{def:sphericalobj}, and $\on{tot}(\C_{\ast>b(p)})\in \D$ is an explicit object obtained from the weighted braid word right below $p$ and the corresponding weights of $\rho_i$, cf.~\cref{def:y}. We emphasize that the dual Lusztig cycles $\{\rho_i^\w\}$ in Theorem \ref{thm:main3}, which are (non-categorical) weave cycles, are already new from the perspective of non-categorical weave calculus. The correction term \eqref{def:y_intro} is a key new contribution that allows for such dual Lusztig cycles to be described explicitly as weave cycles. This is an instance where the categorical weave calculus, which leads naturally to \eqref{def:y_intro}, answers an unresolved question within the non-categorical weave calculus.\\

\noindent Theorem \ref{thm:main3}.(3) allows us to compute certain evaluations of the silting collection $\{\lsilt_i^\w\}$ using their Ringel-type duals and their corresponding weave cycles. It is desirable to be able to compute such evaluations explicitly, as Theorem \ref{thm:main3}.(3) does via the evaluation $\on{ev}_{e_\infty}(\rsilt_i^\uparrow)$ of the Ringel-type dual. Specifically, the evaluations $\{\on{ev}_{e_1}(\lsilt_i^\w)\}$ determine the cluster tilting object used to additively categorify the cluster algebra associated with the braid variety, as will be explained in \cite{CasalsChrist2}. 

\subsection{A comment on applications}

Theorems \ref{thm:main1}, \ref{thm:main2}, \ref{thm:main3}, and the new techniques developed to prove them are a core contribution of this manuscript, establishing the foundations of categorical weave calculus. At the same time, these contributions lead to a number of new results and applications.  Indeed, in the same manner that the non-categorical weave calculus leads to new results in a variety of contexts, cf.~\cref{ssec:scientific_context} above, the uses of this categorical enhancement are manifold. In particular:\\

\noindent {\bf $(i)$ Symplectic topology}. A symplectic geometric realization of the cluster algebras constructed in \cite{CGGLSS25} stems from the study of Lagrangian fillings of Legendrian links, see e.g.~\cite{CasalsWeng22,CasalsZaslow}. This is a geometric problem about certain surfaces in 4-dimensions, and a Demazure weave $\w$ yields such a Lagrangian filling. Intuitively, the cluster algebra is the coordinate ring of functions for a moduli of Lagrangian fillings, with different Lagrangian fillings (e.g.~coming from Demazure weaves) providing the cluster charts. It is much desirable to have a symplectic realization for categorical weave calculus, e.g.~for the weave schober $\F_\w$, its saturated relative 3-Calabi--Yau category $\glsecF$, and the categorical Lusztig cycles and their duals.
    
    For that, we shall first construct a symplectic 6-manifold $(Y_\w,\la_\w)$ with a Lefschetz fibration $\pi:Y_\w\lr\cC$ such that $\glsecF$ will be the partially wrapped Fukaya category associated to $\pi$. Then we will explain how $(Y_\w,\la_\w)$ relates to the Lagrangian fillings in 4-dimensions by using matching 2-disks and the $\mathbb{L}$-compressing systems from \cite{CasalsGao24}. In particular, this yields a rigorous symplectic geometric model for these cluster algebras, from the 4-dimensional study of Lagrangian fillings, and for the cluster categories, here obtained via the Fukaya-Seidel category of a 6-dimensional manifold, with a geometric map connecting these two models. In this symplectic context, we shall explain the meaning of the categorical Lusztig cycles in Theorem \ref{thm:main1}, in terms of Lagrangian skeleta, of the dual ones in Theorem \ref{thm:main3}, in terms of Lagrangian cocores, and of the cluster algebra from \cite{CGGLSS25} as it relates to the Fukaya-Seidel category.\\
    
\noindent {\bf $(ii)$ Quivers with potential}. By combining \cite[Section 2]{CasalsGao24} and \cite[Section 7.4]{CGGLSS25}, there is a canonical ice quiver with potential $(Q_\w,W_\w)$ associated to a Demazure weave $\w$. We shall establish a relation between the $\infty$-category of global sections $\glsecF$ of the weave schober $\F_\w$, and the perfectly valued derived $\infty$-category of the relative 3-Calabi--Yau Ginzburg dg-algebra of such $(Q_\w,W_\w)$. The relative Ginzburg algebra is obtained from the categorical dual Lusztig cycles $\{ \lsilt_i^\w \}$, and its simple modules identify with the categorical Lusztig cycles $\{ \LC_i\}$. In particular, the description in terms of the quiver potential shows that the weave cycle intersection pairing of the Lusztig cycles, see \cite[Sections 4.5\&4.6]{CGGLSS25}, computes the Euler characteristics of the $\on{Ext}$-groups between the categorical Lusztig cycles. We will also prove that the Donaldson-Thomas transformation of the quiver can be realized via a sequence of weave mutations and that, in categorical weave calculus, the associated transformation coincides with the Ringel-type duality in \eqref{eq:intro_ringel_duality}.\\

\noindent {\bf $(iii)$ Additive categorifications of cluster algebras}. The results of this manuscript lead to additive categorifications of the cluster algebras associated to braid varieties. In the reduced case $\beta=w_0w$, $w\in W$ a reduced word, such additive categorifications have been previously studied via the abelian categories of modules over preprojective algebras and their exact subcategories, see for instance~\cite{GLS06,GLS08,BIRS09}. Our contributions relate to these previous works as follows.\\

First, we tackle the general case where $\beta$ is any braid word with $\delta(\beta)=w_0$, not necessarily the reduced case. The silting collections we construct yield cluster tilting objects in the Higgs category or, equivalently by \cite[Ex.~8.19]{Wu21}, in the module category of the preprojective algebra. These cluster tilting objects generalize those associated to the maximal rigid modules from \cite{GLS06,GLS08,BIRS09}. The so-called layer modules in these works are a special case of the $2$-spherical objects associated with the trivalent vertices of the weave, see \Cref{def:sphericalobj} and  \Cref{rmk:reduced_case_sphericaltwists}. In this sense, the filtrations of the maximal rigid modules by layer modules, cf.~\cite[Cor.~2.9]{IR11}, are generalized to arbitrary weaves via the filtrations in Theorem \ref{thm:main1} and Theorem \ref{thm:main3}.(3). In addition, \cite{GLS07,GLS11,IR11} show that the endomorphism algebra of the cluster tilting object in the reduced case is quasi-hereditary, i.e.~its module category is a highest weight abelian category. The triangulated categories constructed in our paper thus generalize the derived categories of these quasi-hereditary algebras. We allow the associated category to be triangulated, not just abelian, and equip them with a full exceptional collection, in place of a highest weight structure, cf.~\cref{ssec:tiltingtheory_exceptionalcollections} and Appendix \ref{ssec:AppendixA}.\\

\noindent Second, even in the reduced case and for a given reduced expression, we construct cluster tilting objects corresponding to many cluster seeds, not just one. Specifically, we recover the seed associated to a maximal rigid module  by considering the right-inductive weave. The above discussions concern Type A, and details will appear in \cite{CasalsChrist2}. We expect this relation to hold in the case of a general finite Dynkin diagram, cf.~\cref{rmk:Dynkincase}.

\begin{remark}
Relatedly, in the context of local systems with irregular singularities, non-categorical weave calculus allows for the construction of cluster algebras for the coordinate rings of wild character varieties, see e.g.~\cite[Section 1.15]{FG06} and \cite[Section 3]{CasalsZaslow}, or \cite{casals2025spectralnetworksbettilagrangians}. Similar to (2) above, there is an important dissonance between the real 4-dimensional geometry that leads to the non-abelianization maps of \cite{GNMSN,GMN14_Snakes} and the cluster algebras in \cite{CasalsWeng22,CGGLSS25}, versus the real 6-dimensional geometry that leads to the 3-Calabi--Yau categorical techniques used to categorify such cluster algebras, see e.g.~\cite[Section 6]{Gon17} or \cite{Chr22,Chr25c}. The categorical weave calculus developed here, in conjunction with application (i) above, resolves this picture by providing a geometric factorization of the appropriate symplectic 3-Calabi--Yau manifold as a high-dimensional branched cover over the symplectic push-off of the Lagrangian filling associated to $\w$, explaining geometrically how these two approaches are related.
\qed
\end{remark}

Finally, putting (i)-(ii) above together yields a form of homological mirror symmetry, established through the global sections $\glsecF$ of the weave schober $\F_\w$: in the A-side there is a Fukaya-Seidel category and its associates, such as the relative Rabinowitz category, and in the B-side there is the relative Ginzburg dg-algebra and its associates, such as the cluster category or Higgs category, cf.~\cite[Corollary 5.26]{Wu21}. Here the A-side is readily symplectic, whereas these categories in the B-side are understood as a form of non-commutative crepant resolutions. These further results and applications of categorical weave calculus will be established in \cite{CasalsChrist2}.\\


\noindent {\bf Organization of the paper}: The manuscript proves the above results as follows. \cref{sec:preliminaries} provides the necessary ingredients on weaves and perverse schobers, in \cref{ssec:primer_weavecalculus} and \Cref{ssec:perverse_schober_preliminaries} respectively. The weave schober $\F_\w$ in Theorems \ref{thm:main1} and \ref{thm:main2} is constructed in Section \ref{ssec:defining_weave_schobers} and Theorem \ref{thm:main2}.(1) is proven in \cref{ssec:weaveschober_equivalences_mutations}, cf.~\cref{thm:weaveequivglsec}. The categorical Lusztig cycles are constructed in Section \ref{sec:lusztigcycles}, using results latter established in \cref{sec:weave_propagation_LC}. Theorem \ref{thm:main1} is also proven in \cref{sec:weave_propagation_LC}. The proofs in \cref{sec:weave_propagation_LC} base on results and techniques on categorical weighted braid words developed in \cref{sec:filtrations_from_braid_words}. Theorem \ref{thm:main2}.(2), (3) and (4) are then proven in \cref{sec:proof_Lusztig_cycles}, using the constructions from \cref{sec:weave_thimbles}, and Theorem \ref{thm:main3} in \cref{sec:silting_dualLusztig}. Finally, discussions on the relation with the theory of highest weight abelian categories can be found in \cref{ssec:tiltingtheory_exceptionalcollections} and Appendix \ref{ssec:AppendixA}.\qed

\vspace{0.5cm}

\noindent {\bf Acknowledgements}: We are grateful to Daping Weng, Jos\'e Simental, Yoon Jae Nho, Melissa Sherman-Bennet, and James Hughes for discussions on weave calculus. We thank Ricardo Canesin, Johannes Flake, Jonathan Gruber, Bernhard Keller, Henning Krause and Catharina Stroppel for helpful discussions on highest weight categories. R.C.~is supported by the National Science Foundation under the grant DMS-2505760, and a UC Davis College of L\&S Dean's Fellowship. M.C.~is a member of the Hausdorff Center for Mathematics at the University of Bonn (DFG GZ 2047/1, project ID 390685813).\qed

\vspace{0.5cm}

\noindent {\bf Notation on $\infty$-categories}: We use the theory of stable $\infty$-categories as developed in \cite{HTT,HA}. The $\infty$-category of small $\infty$-categories is denoted by $\on{Cat}_\infty$, as in \cite[Chapter 3]{HTT}. The $\infty$-category of small stable $\infty$-categories is denoted $\on{St}$, cf.~\cite[Section 1.1.4]{HA}. All stable $\infty$-categories appearing in this paper will be small. We denote by $\on{St}^{\on{idem}}\subset \on{St}$ the full subcategory of idempotent complete stable $\infty$-categories. $\on{Ind}$-completion defines an equivalence between $\on{St}^{\on{idem}}$ and the $\infty$-category of compactly generated, presentable $\infty$-categories and compact objects preserving, left adjoint functors $\mathcal{P}r^{L,\omega}_{\on{St}}$. The symmetric monoidal $\infty$-structure on $\mathcal{P}r^{L,\omega}_{\on{St}}$ thus induces a symmetric monoidal structure on $\on{St}^{\on{idem}}$. Given a field $k$, $\D^{\on{perf}}(k)\in \on{St}^{\on{idem}}$ defines an algebra object and we define $\on{LinCat}_k^{\on{sm}}\coloneqq \on{Mod}_{\D^{\on{perf}}(k)}(\on{St})$ to be the $\infty$-category of small, idempotent complete, $k$-linear stable $\infty$-categories. Note that the forgetful functor $\on{LinCat}_k^{\on{sm}}\to \on{St}^{\on{idem}}$ preserves limits by \cite[Cor.~3.4.3.2]{HA}. All limits computed in $\on{St}^{\on{idem}}$ in this paper will be preserved by the forgetful functor $\on{St}^{\on{idem}}\to \on{St}$, as a $k$-linear stable $\infty$-category with a full exceptional collection is automatically idempotent complete. The forgetful functor $\on{St}\to\on{Cat}_\infty$ also preserves limits.

Given a $k$-linear stable $\infty$-category $\C\in \on{LinCat}_k^{\on{sm}}$ and two objects $X,Y\in \C$, we write $\on{Mor}_\C(X,Y)\in \D(k)$ for the morphism object in the $\D(k)$-linear $\infty$-category $\on{Ind}\C$, see \cite[Def.~4.2.1.28]{HA}, i.e.~the $k$-linear $\infty$-categorical version of the derived Hom. There is also an isomorphism of $k$-vector spaces
\[
\on{Ext}^{i}_\C(X,Y)\simeq H_{0}\on{Mor}(X,Y[i])\,.
\]
Finally, all chain complexes are graded homologically and $k[i]$ denotes a chain complex concentrated in degree $i$, as in \cite{HA}. In particular, for the non-negative truncation functor $\tau_{\geq 0}$, which we use often for semi-twists, we have $\tau_{\geq0}(k[-1]\oplus k[2])\simeq k[2]$.
\qed

\section{Preliminary Ingredients}\label{sec:preliminaries}

This section introduces the necessary ingredients to establish the results in this manuscript. Specifically, \cref{ssec:primer_weavecalculus} provides a summarizing account of non-categorical weave calculus, including the definition of weave equivalences and mutations, and the non-categorical Lusztig cycles. 
\cref{ssec:perverse_schober_preliminaries} recalls the notion of a perverse schober parametrized by a ribbon graph and provides a few facts about perverse schobers that we use throughout the article.



\begin{center}
	\begin{figure}[h!]
		\centering
		\includegraphics[scale=1]{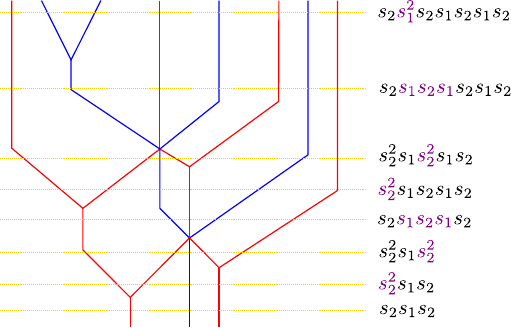}
		\caption{A Demazure weave $\w:\beta\lr\delta(\beta)$ from the 3-stranded braid word $\beta=s_2s_1(s_1s_2)^3$ to $\beta'=s_2s_1s_2$, where blue weave edges are labeled by $s_1\in S_3$ and red weave edges by $s_2\in S_3$. Generic horizontal slices are depicted in yellow dashed lines, and the positive braid word associated to each such slice is written to its right. We have highlighted in purple the letters of each braid word to which the next move will be applied when scanning the Demazure weave downward.}
		\label{fig:weave_intro_example}
	\end{figure}
\end{center}

\subsection{A primer on weave calculus}\label{ssec:primer_weavecalculus}

\begin{center}
	\begin{figure}[h!]
		\centering
		\includegraphics[scale=2.1]{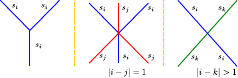}
		\caption{The three types of vertices allowed in the diagrammatic version of a weave $\w$. (Left) Trivalent vertex, where all edges are decorated with the same permutation $s_i\in S_n$. It represents the move $s_i^2\to s_i$. (Center) Hexavalent vertex, with the edge decorations alternating between $s_i,s_j$ with $|i-j|=1$. It represents the braid move $s_is_js_i\to s_js_is_j$. (Right) Tetravalent vertex, with edge decorations as drawn with $|i-k|>1$. It represents the braid move $s_is_k=s_ks_i$.}
		\label{fig:weaves_types_vertices}
	\end{figure}
\end{center}

The focus of the article will be on Demazure weaves, whose non-categorical weave calculus we now briefly summarize. For specificity, we focus on Type $A_{n}$, though similar considerations apply to all types, cf.~\cite[Section 4]{CGGLSS25}. Let $\beta,\beta'$ be two positive braid words on $n$-strands, i.e.~words on the positive Artin generators $s_1,\ldots,s_{n-1}$ of the braid group on $n$-strands. By definition, a Demazure weave $\w:\beta\lr\beta'$ is a sequence of positive braid words $\beta\to\beta_1\to\beta_2\to\cdots\to\beta'$ such that each move $\beta_j\to\beta_{j+1}$ is either a braid move in the positive braid monoid or the move $s_i^2\to s_i$. \cref{fig:weave_intro_example} depicts an example of a Demazure weave $\w$ diagrammatically. Specifically, by drawing a positive braid word $s_{i_1}\cdots s_{i_l}$ as a set of parallel vertical strands on the plane, with the $j$th strand labeled by $s_{i_j}$, these moves can be depicted diagrammatically as in \cref{fig:weaves_types_vertices}. A core aspect of weave calculus are the Lusztig cycles of a Demazure weave, their intersection theory, and their behavior under weave equivalences and mutations, which we will momentarily summarize.

\noindent These diagrams for Demazure weaves are always read top-to-bottom, the diagram itself is often referred to as a weave. Considered as a planar graph, a Demazure weave admits the three types of vertices as in \cref{fig:weaves_types_vertices}, and the connected components of the complement of the vertices are said to be weave edges, cf.~\cite[Definition 2.1]{CasalsWeng22} or \cite[Section 2]{CasalsZaslow} for the general definition. The trivalent vertices of a weave $\w$ will be denoted $p_1,\dots,p_m$, ordered bottom-to-top, e.g.~with $p_i$ at height $i$. Each horizontal slice not intersecting a vertex determines a positive braid word by reading the weave edges left-to-right, as illustrated in \cref{fig:weave_intro_example}.\\

\begin{center}
	\begin{figure}[h!]
		\centering
		\includegraphics[scale=1.2]{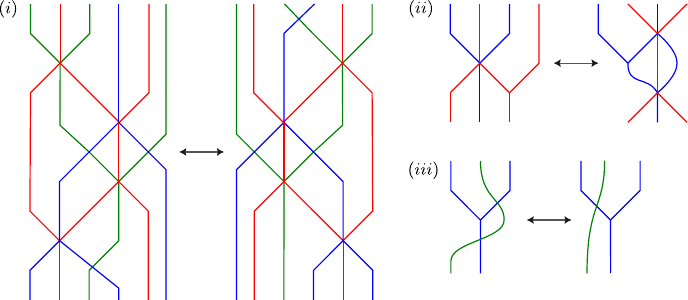}
		\caption{A few weave equivalences, in all of them blue indicates $s_1\in S_3$, red $s_2$, and green $s_3$. $(i)$ This weave equivalence encodes a syzygy between braid moves, it is also known as the Zamolodchikov relation. $(ii)$ A weave equivalence known as the push-through move, as it allows for pushing trivalent vertices through hexavalent vertices. $(iii)$ A weave equivalence allowing for distant colors to pass through vertices of each other.}
		\label{fig:weave_equivalences}
	\end{figure}
\end{center}


\subsubsection{Weave equivalence and mutations}\label{sssec:weave_equivalences_mutations} The notion of weave mutation was introduced in \cite[Section 4.8]{CasalsZaslow}. Equivalences between weaves, also known as moves, were discussed in \cite[Theorem 1.1]{CasalsZaslow}. See also \cite[Section 4]{CGGJ24}. We consider the following equivalence relation on weaves:
\begin{itemize}
    \item[(i)] Two weaves $\w,\w':\beta\to\beta'$ consisting only of braid moves are declared to be equivalent, see e.g.~\cref{fig:weave_equivalences}.(i).
    \item[(ii)] If $|i-j|=1$, the weaves $\w:s_is_js_is_j\to s_js_is_js_j\to s_js_is_j$ and $\w':s_is_js_is_j\to s_is_is_js_i\to s_is_js_i\to s_js_is_j$  are declared to be equivalent, cf.~\cref{fig:weave_equivalences}.(ii).
    \item[(iii)] If $|i-j|\geq2$, the weaves $\w:s_is_js_i\to s_is_is_j\to s_is_j\to s_js_i$ and $s_is_js_i\to s_js_is_i\to s_js_i$ are declared to be equivalent, cf.~\cref{fig:weave_equivalences}.(iii).
\end{itemize}
\noindent By definition, two weaves $\w,\w':\beta\to\beta'$ are said to be equivalent if there is a sequence of local equivalences as above starting at $\w$ and ending at $\w'$. Note that the equivalences in (ii) and (iii) are parameterized by rank $2$ subdiagrams of types $A_2$ and $A_1\times A_1$ respectively.\\

\begin{center}
	\begin{figure}[h!]
		\centering
		\includegraphics[scale=1.5]{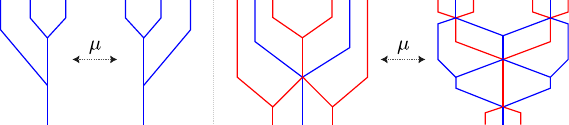}
		\caption{(Left) An elementary weave mutation, between the only two Demazure weaves from $\w,\w':s_1^3\to s_1$. (Right) A more elaborate form of weave mutation, obtained by composing weave equivalences with an elementary weave mutation.}
		\label{fig:weave_mutation}
	\end{figure}
\end{center}

\begin{center}
	\begin{figure}[h!]
		\centering
        \includegraphics[scale=1.2]{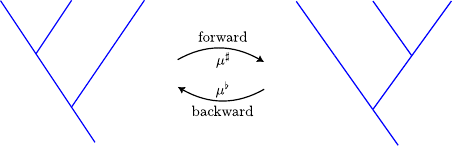}
		\caption{A redrawing of \cref{fig:weave_mutation} (left), with new notation to specify direction. Recall that a weave mutation is given by a sequence of weave equivalences, one local exchange of the weave according to the two diagrams depicted (with any color $s_i$), and another sequence of weave equivalences. In a forward weave mutation, the piece of weave on the left is exchanged by the piece of weave on the right. A backward weave mutation is the inverse process.}\label{fig:WeaveMutation_ForwardBackward}
	\end{figure}
\end{center}

\noindent By definition, the two weaves $\w,\w':s_i^3\to s_i$ are said to be are related by an elementary weave mutation, cf.~\cref{fig:weave_mutation} (left) or \cref{fig:WeaveMutation_ForwardBackward}. They are not considered to be weave equivalent. In general, two weaves are said to be related by a weave mutation if they can be related by a sequence of weave equivalences, an elementary weave mutation, and another sequence of weave mutations. See for instance cf.~\cref{fig:weave_mutation} (right). It is proven in \cite[Theorem 4.6]{CGGJ24} that any two Demazure weaves $\w,\w':\beta\to\delta(\beta)$ are related by a sequence of equivalence moves and elementary mutations, where we have fixed a braid word for the Demazure product $\delta(\beta)$ of $\beta$.

\noindent In order to match the notation used in the tilting theory of finite length hearts, we specifically distinguish the direction of a weave mutation by using the adjectives forward and backward according to \cref{fig:WeaveMutation_ForwardBackward}.


\subsubsection{Weighted braid words and Lusztig propagation rules}\label{sssec:Lusztigcycles} Let $\w:\beta\to\beta'$ be a Demazure weave and $E(\w)$ its set of weave edges. By definition, a weave cycle on $\w$ is a function $\g: E(\w) \to \Z_{\geq 0}$ that assigns a non-negative integer to each edge of the weave. The values of $\g$ are said to be the weights of the edges for the cycle $\w$, cf.~\cite[Section 4.4]{CGGLSS25}.

In the same manner that each horizontal slice of a weave yields a braid word $s_{i_1}\cdots s_{i_k}$, cf.~\cref{fig:weave_intro_example}, given a weave cycle $\g$ and a horizontal slice we obtain a braid word $s_{i_1}\cdots s_{i_k}$ and a collection of non-negative integers ${\bf a}\coloneqq (a_1,\ldots,a_k)$, where $a_j\coloneqq \g(e_{i_j})$ is the weight of the $j$th weave edge $e_{i_j}\in E(\w)$ being intersected in that horizontal slice, corresponding to the letter $s_{i_j}$ in the braid word $s_{i_1}\cdots s_{i_k}$. By definition, we refer to such a pair $(s_{i_1}\cdots s_{i_k},{\bf a})$ as a weighted braid word. Of particular importance are the following weave cycles associated to the trivalent vertices of the weave $\w$:

\begin{definition}[Lusztig cycles]\label{def:Lusztig_cycles}
Let $\w$ be a Demazure weave and $p\in\w$ a trivalent vertex. By definition, the Lusztig cycle $\g_p$ associated to $p$ is the unique weave cycle $\g_p: E(\mathfrak{W}) \to \Z_{\ge 0}$ such that:
\begin{itemize}
    \item[$(i)$] $\gamma_p(e_\mathsf{s})=1$ for the unique south edge $e_\mathsf{s}$ at the trivalent vertex $p$, and $\gamma_p(e)=0$ if $e$ is any weave edge that contains any points above $p$.
    \item[$(ii)$] At any other trivalent vertex with top edges $e_1,e_2$ and south edge $e$, $\g_p$ satisfies
    \begin{equation}\label{eq:3_move}\g_p(e)=\min(\gamma_p(e_1),\g_p(e_2))\,.\end{equation}
    \item[$(iii)$] For a 4-valent vertex  with top edges $e_1,e_2$ and  bottom edges $e_1',e_2'$, $\g_p$ satisfies
    \[(\g_p(e'_1),\g_p(e'_2))=(\g_p(e_2),\g_p(e_1))\,.\]
    \item[$(iv)$] For a 6-valent vertex  with top edges $e_1,e_2,e_3$ and bottom edges $e_1',e_2',e'_3$, $\g_p$ satisfies
    \begin{align*}
        \g_p(e'_1) & = \g_p(e_2)+\g_p(e_3)-\min(\gamma_p(e_1),\g_p(e_3))\,,\\
        \g_p(e'_2) & = \min(\gamma_p(e_1),\g_p(e_3))\,,\\
        \g_p(e'_3) & = \g_p(e_2)+\g_p(e_1)-\min(\gamma_p(e_1),\g_p(e_3))\,.\qed
    \end{align*}
\end{itemize}
\end{definition}

\begin{center}
	\begin{figure}[h!]
		\centering
        \includegraphics[scale=1.6]{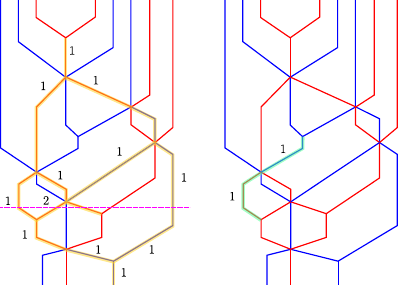}
		\caption{Two instances of Lusztig cycles on a Demazure weave $\w:(s_1^2s_2^2)^2\to s_1s_2s_1$. (Left) The Lusztig cycle $\g_1$ associated to the first (top) trivalent vertex. The non-zero weave edges for $\g_1$ are highlighted in yellow, and the specific weight is written next to each weave edge. For instance, at the horizontal pink slice, the weighted braid word is $((s_2^2s_1)^2,(1,2,0,1,0,1))$. (Right) The Lusztig cycle $\g_2$ associated to the second trivalent vertex, with the non-zero weave edges for $\g_2$ highlighted in cyan.}\label{fig:weave_Lusztigcycle_example}
	\end{figure}
\end{center}

See \cref{fig:weave_intro_example2} and \cref{fig:weave_Lusztigcycle_example} for instances of Lusztig cycles; in the former a highlighted edge denotes weight 1 and a non-highlighted edge weight 0, in the latter the weights of the weave edges are specified. \cref{def:Lusztig_cycles} can be also described in terms of weighted braid words as follows. To construct $\g_p$, first consider a horizontal weave slice right below the height of $p$, which yields a certain braid word $\beta(p)\coloneqq s_{i_1}\cdots s_{i_{k(p)}}$ and assign the weights ${\bf a}(p)\coloneqq (0,0,\ldots,0,1,0,\ldots,0)$, where the only non-zero weight $1$ is assigned to the weave edge south of $p$. Then, propagate this weighted braid word $(\beta(p),{\bf a}(p))$ downwards by using the following rules for the weights:
\begin{multline}
\label{eq: lusztig trop}
\qquad \qquad \qquad \qquad \quad \Phi_1:(a_1,a_2)\mapsto \min(a_1,a_2), \quad \Phi_2:(a_1,a_2)\mapsto (a_2,a_1),\\ \qquad\qquad \Phi_3: (a_1,a_2,a_3)\mapsto  \left(a_2+a_3-\min(a_1,a_3),\min(a_1,a_3),a_1+a_2-\min(a_1,a_3)\right),
\end{multline}
where $\Phi_1$ is used at trivalent vertices, where the braid word changes according to $s_i^2\lr s_i$, $\Phi_2$ at tetravalent vertices, where the braid word changes as $s_is_k=s_ks_i$ with $i,k$ non-adjacent, and $\Phi_3$ at hexavalent vertices, where the braid word changes via $s_is_js_i=s_js_is_j$ with $i,j$ adjacent.

\begin{remark} The downward propagation rules in \cref{eq: lusztig trop} above are referred to as tropical Lusztig rules, as they are obtained by tropicalizing the relations
\begin{multline*}
\label{eq: lusztig trop2}
\qquad \qquad \qquad \qquad \quad x_i(t_1)x_i(t_2)=x_i(t_1+t_2), \quad x_i(t_1)x_k(t_2)=x_k(t_2)x_i(t_1),\\ x_i(t_1)x_j(t_2)x_i(t_3)=
x_j\left(\frac{t_2t_3}{t_1+t_3}\right)x_i(t_1+t_3)x_j\left(\frac{t_1t_2}{t_1+t_3}\right), \hskip 2.8cm
\end{multline*}
\noindent where $x_i(t)\coloneqq \exp(E_{i,i+1}t)$ is the 1-parameter subgroup in $\on{SL}_n$ corresponding to the positive simple root $\alpha_i$. Here the indices $i,j$ are adjacent, and $i,k$ non-adjacent, in the Dynkin diagram.\qed
\end{remark}

Finally, of key importance are the intersection theory of weave cycles and their associated variables, especially for the Lusztig cycles in \cref{def:Lusztig_cycles}. The precise formulas for their intersections are detailed in \cite[Sections 4.5 \& 4.6]{CGGLSS25}, and the variables are described in \cite[Section 5.2]{CGGLSS25}, cf.~also \cite[Section 4]{CasalsWeng22}. In particular, the collection of Lusztig cycles determines an intersection quiver, as illustrated in \cref{fig:weave_intro_example2}.


\subsection{Ingredients on perverse schobers}\label{ssec:perverse_schober_preliminaries}

The weave schober $\F_\w$ associated to a Demazure weave $\w$, introduced in \cref{sec:weave_schobers}, will be a perverse schober parametrized by a ribbon graph $\rgraph$. Let us recall the key definitions and properties of such perverse schobers. We refer to \cite{KS14}, \cite[Sections 3.1 \& 3.2]{Chr22b}, or \cite[Section 3]{CHQ23} for further details.\\


\subsubsection{\texorpdfstring{$\rgraph$}{G}-parametrized perverse schobers}\label{sssec:G_parametrized_perverse_schobers} Let us describe $\rgraph$-parametrized perverse schobers, and the ingredients needed to define them.

\begin{definition}\label{def:exit_path_category}
Let $\rgraph$ be a ribbon graph. By definition, its exit path category $\on{Exit}(\rgraph)\in \on{Cat}_\infty$ is the nerve of the $1$-category whose objects are the vertices and edges of $\rgraph$, and whose non-identity morphisms are of the form $v\to e$, where $e$ is an edge of $\rgraph$ incident to a vertex $v$.\qed
\end{definition}

By the exodromy equivalence, cf.~e.g.~\cite[Theorem 5.17]{PT22}, constructible sheaves on $\rgraph$ (understood as a stratified space) valued in the $\infty$-category of small stable $\infty$-categories $\on{St}$ can equivalently be described as functors $\on{Exit}(\rgraph)\to \on{St}$. Thus we can use $\on{Exit}(\rgraph)$ in \cref{def:exit_path_category} and have $\on{Fun}(\on{Exit}(\rgraph),\on{St})$ as a model for constructible sheaves on $\rgraph$ with values in the $\infty$-category of stable $\infty$-categories. The notion of perversity in classical sheaf theory, cf.~\cite{BBD82}, has the following incarnation in the context of such sheaves of categories on ribbon graphs.\\ 

\noindent For each $n\in\mathbb{N}$, let $\rgraph_{n}$ be the ribbon graph with a single vertex $v$ and $n$ incident external edges $e_1,\ldots,e_n$. The ribbon graph $\rgraph_n$ is said to be the $n$-spider $\rgraph_{n}$. A perverse schober on $\rgraph_{n}$ is defined as follows:

\begin{definition}\label{def:schobernspider}
A perverse schober on the $n$-spider $\rgraph_{n}$ is a functor $\mathcal{F}\colon \on{Exit}(\rgraph_n)\to \on{St}$ such that:
\begin{enumerate}
\item[(1)] If $n=1$, the functor
\[ F\colon \mathcal{V}\longrightarrow \mathcal{N}\,,\qquad\mbox{where } \mathcal{V}\coloneqq F(v),\quad \mathcal{N}\coloneqq F(e),\quad F\coloneqq \F(v
\to e),\]
given by restricting to the unique edge $e$ is required to be spherical.\footnote{$F$ is spherical if it admits a right adjoint $G$, and both the cotwist functor $T_{\mathcal{V}}\coloneqq \on{fib}(\on{id}_{\mathcal{V}}\xrightarrow{\on{unit}}GF)\in \on{Fun}(\mathcal{V},\mathcal{V})$ and twist functor $T_{\mathcal{N}}\coloneqq \on{cof}(FG\xrightarrow{\on{counit}}\on{id}_{\mathcal{N}})\in \on{Fun}(\mathcal{N},\mathcal{N})$ are auto-equivalences. \\
Previous papers of the second author used different conventions for the twist and cotwists functors. The convention employed here is the standard in the literature on spherical objects. It is employed to simplify the notation in the computations involving spherical objects.}

\item[(2)] If $n\geq 2$, the collection of functors 
\[ (F_i\colon \mathcal{V}^n\longrightarrow \mathcal{N}_i)_{i\in \mathbb{Z}/n\mathbb{Z}},\qquad\mbox{where }\mathcal{V}^n\coloneqq F(v),\quad \mathcal{N}_i\coloneqq F(e_i),\quad F_i\coloneqq \F(v
\to e_i),\]
satisfies that
\begin{enumerate}
	\item there exist adjunctions $F_i\dashv G_i\dashv H_i$ for some $G_i,H_i$,
    \item $G_i$ is fully faithful, which is equivalent to $F_iG_i\simeq \on{id}_{\mathcal{N}_i}$ via the counit,
    \item $F_{i}\circ G_{i-1}$ is an equivalence of $\infty$-categories,\footnote{The convention used in \cite{CDW23} is that $F_iG_{i+1}$ is an equivalence of $\infty$-categories instead. This difference corresponds to reversing the orientation of the plane containing the $n$-spider. In the convention used here, the wrapping in the resulting partially wrapped Fukaya categories will be counterclockwise.}
    \item $F_i\circ G_j\simeq 0$ if $j\neq i,i-1$,
    \item $\on{fib}(H_{i-1})=\on{fib}(F_{i})$ as full subcategories of $\mathcal{V}^n$.\qed
\end{enumerate}
\end{enumerate}
\end{definition}

In \cref{def:schobernspider}.(1), the functor $F$ uniquely determines $\F$ as $v,e$ are the only objects and $v\to e$ is the only 1-morphism in the 1-category defining $\on{Exit}(\rgraph_n)$. Similarly, in \cref{def:schobernspider}.(2), the functors $(F_1,\ldots,F_n)$ uniquely determine $\F$ for $n\geq2$. \cref{def:schobernspider} provides local perversity conditions at $\rgraph_n$. Passing to Grothendieck groups, these conditions imply that the arising diagram encodes the sections of a perverse sheaf on the disc with support on $\rgraph_n$, see \cite{KS16}. This notion is globalized to a ribbon graph as follows:

\begin{definition}[Perverse schober on $\rgraph$]\label{def:schober}
Let $\rgraph$ be a ribbon graph. By definition, $\rgraph$-parametrized perverse schober is a functor
\[\mathcal{F}\colon \on{Exit}(\rgraph)\rightarrow \on{St}\]
such that, for each $n$-valent vertex $v$ of $\rgraph$, the restriction of $\mathcal{F}$ to the coslice category $\on{Exit}(\rgraph)_{v/}$ determines a perverse schober on the $n$-spider, as in \Cref{def:schobernspider}.\qed
\end{definition}

We will mostly be interested in $k$-linear perverse schobers, meaning that the functor $\mathcal{F}\colon \on{Exit}(\rgraph)\rightarrow \on{St}$ factors through $\on{LinCat}_k^{\on{sm}}\to \on{St}$. 

The $\infty$-category of global sections of a $\rgraph$-parametrized perverse schober is defined as follows:

\begin{definition}\label{def:sections}
Let $\mathcal{F}$ be a $\rgraph$-parametrized perverse schober.
\begin{itemize}
    \item[(i)] By definition, the $\infty$-category of global sections of $\mathcal{F}$ is
\[\Glsec{\rgraph}{\F}\coloneqq \on{lim}(\mathcal{F})\in \on{St}\,.\]
Note that $\Glsec{\rgraph}{\F}$ can be identified with the $\infty$-category of coCartesian sections of the Grothendieck construction $p\colon \Gamma(\mathcal{F})\to \on{Exit}(\rgraph)$ of $\mathcal{F}$, see \cite[\href{https://kerodon.net/tag/05RX}{Prop.~05RX}]{Ker}. Since $\mathcal{F}$ arises from a strictly commuting diagram $\on{Exit}(\rgraph)\to \on{Set}_\Delta$, there is an explicit simplicial set version of the Grothendieck construction \cite[\href{https://kerodon.net/tag/025X}{Def.~025X}]{Ker}.

\item[(ii)] By definition, the $\infty$-category $\Losec{\rgraph}{\mathcal{F}}$ of lax sections of $\mathcal{F}$ is the $\infty$-category of all sections of the Grothendieck construction $p\colon \Gamma(\mathcal{F})\rightarrow \on{Exit}(\rgraph)$ of $\mathcal{F}$. Note that $\Losec{\rgraph}{\mathcal{F}}$ describes the lax limit of $\mathcal{F}$ in the $(\infty,2)$-category of stable $\infty$-categories.\qed
\end{itemize}
\end{definition}

\noindent In \cref{def:sections}, $\Glsec{\rgraph}{\F}\subset \Losec{\rgraph}{\F}$ is a full subcategory. The focus in this manuscript is on $\Glsec{\rgraph}{\F}$, but $\Losec{\rgraph}{\F}$ is used for gluing constructions of global sections, see for instance \cref{construction:local_sections}.

\begin{remark}
A vertex $v$ of a $\rgraph$-parametrized perverse schober $\mathcal{F}$ is said to be a singularity if the fiber of the functor
$$\prod_{i=2}^{n} \mathcal{F}(v\to e_i):\mathcal{F}(v)\lr \prod_{i=2}^n\mathcal{F}(e)$$
is not equivalent to zero. For $n=1$, the fiber above is taken to be $\F(v)$. This fiber is also called the $\infty$-category of vanishing cycles of $\F$ at $v$, and it is equivalent to the fiber of any other collection of $n-1$ of the $n$ functors $\{\mathcal{F}(v\to e_i)\}_{1\leq i\leq n}$.\qed
\end{remark}

Given a ribbon graph $\rgraph$ together with a chosen collection of internal edges, we can contract these internal edges to produce a new ribbon graph $\rgraph'$. This is denoted by $c\colon \rgraph \to\rgraph'$, indicating that $\rgraph'$ arises via contraction from $\rgraph$. In order to understand the behavior of weave schobers under weave equivalences and weave mutations, we need the following result on how $\rgraph$-parametrized perverse schober behave under such contractions:

\begin{proposition}[{\!\!\cite[Prop.~4.28]{Chr22b} \& \cite[Lemma 3.16]{CHQ23}}]\label{prop:contraction}
Let $c\colon \rgraph \to \rgraph'$ be a contraction of ribbon graphs.
\begin{enumerate}[(1)]
    \item Let $\mathcal{F}$ be a $\rgraph$-parametrized perverse schober and assume that $c$ does not contract two singularities of $\mathcal{F}$ into a single vertex. Then there exists a canonical $\rgraph'$-parametrized perverse schober $c_*(\mathcal{F})$ together with a canonical equivalence of $\infty$-categories
    \[ \Glsec{\rgraph}{\mathcal{F}}\simeq \Glsec{\rgraph'}{c_*(\mathcal{F})}\,.\]
     \item Let $\mathcal{F}'$ be a $\rgraph'$-parametrized perverse schober and choose for every singularity $v\in \rgraph_0'$ of $\mathcal{F}$ a preimage under $c$ in $\rgraph_0$. Then there exists a canonical $\rgraph$-parametrized perverse schober $c^*(\mathcal{F})$ together with a canonical equivalence of perverse schobers $c_*c^*(\mathcal{F}')\simeq \mathcal{F}'$.\qed
\end{enumerate}
\end{proposition}

\noindent \cref {prop:contraction} will be used shortly when stating and proving \cref{lem:exchangesingularities} below.

\subsubsection{A model for \texorpdfstring{$\rgraph_2$}{G2}-parametrized perverse schobers}\label{sssec:model_perverseschobers}

For weave schobers, the ribbon graphs we consider have mostly 2-valent vertices. There is a useful construction of perverse schobers on  the $2$-spider $\rgraph_2$ from spherical functors, as follows. Given a spherical functor $F\colon \V\to \N$, let $\V\overset{\rightarrow}{\times}_F \N$ be defined by the following pullback of $\infty$-categories: 
\[
\begin{tikzcd}
\V\overset{\rightarrow}{\times}_F \N \arrow[d] \arrow[r] \arrow[rd, "\lrcorner", phantom] & {\on{Fun}(\Delta^1,\N)} \arrow[d, "\on{ev}_0"] \\
\V \arrow[r, "F"]                                                                         & \N    ,                                        
\end{tikzcd}
\]
where $\on{ev_0}$ is evaluation of a functor $\Delta^1\to\N$ at the vertex $\{0\}$ of $\Delta^1$. Note that $\V\overset{\rightarrow}{\times}_F \N$ describes a model for the $(\infty,2)$-categorical lax limit of the diagram given by $F$. Concretely, an object in $\V\overset{\rightarrow}{\times}_F \N$ can be identified with a tuple
\begin{equation}\label{eq:description_VFN}
(X,Y,\eta),\qquad \mbox{where }X\in \V,\ Y\in \N,\text{ and } \eta\colon F(X)\to Y\in \V \mbox{ is a morphism}.
\end{equation}

The following construction of $\rgraph_2$-parametrized perverse schobers from a spherical functors is proven in \cite[Prop.~3.7]{Chr22b}:

\begin{lemma}[Perverse schobers on $\rgraph_2$ from spherical functors]\label{lem:schober2gon}
Let $F\colon \V\to \N$ be a spherical functor. Then the functors
\[
\varrho_1\colon \V\overset{\rightarrow}{\times}_F \N \to \on{Fun}(\Delta^1,\N)\xrightarrow{(\{1\}\subset \Delta^1)^*} \N=\N_1
\]
and
\[
\varrho_2\colon \V\overset{\rightarrow}{\times}_F \N \to \on{Fun}(\Delta^1,\N)\xrightarrow{\on{cof}} \N=\N_2
\]
together define a perverse schober $\F_F$ on $\rgraph_2$, with $\F_F(v)\coloneqq \V\overset{\rightarrow}{\times}_F \N$ and $\F_F(e_i)\coloneqq \N$ for $i=1,2$.\qed
\end{lemma}

In the description of \cref{eq:description_VFN}, the functors $\varrho_1,\varrho_2$ from \cref{lem:schober2gon} act on objects as
\begin{equation}\label{eq:rho_maps}
\varrho_1(X,Y,\eta)=Y\,,\quad \varrho_2(X,Y,\eta)=\on{cof}(\eta)\,.
\end{equation}
Later, we will use the left and right adjoints $\varrho_1^L,\varrho_1^R:\N_1\to \V\overset{\rightarrow}{\times}_F \N$ of $\varrho_1$, which are given on objects as
\[\varrho_1^L(Y)=(0,Y,0),\qquad \varrho_1^R(Y)=(G(Y),Y,\on{counit}\colon FG(Y)\to Y)\,,\]
where $G$ is the right adjoint of $F$.\\

\begin{notation}\label{not:schobers}
We use the following notation for perverse schobers on ribbon graphs with $2$-valent and $1$-valent vertices, based on the local model from \cref{lem:schober2gon}. Specifically, we can describe a perverse schober in such a ribbon graph by labeling each vertex by a spherical functor $F\colon \V\to \N$ and each adjacent edge either by $\varrho_1$ or $\varrho_2$, in the $2$-valent case, or $\varrho_1=F$ in the $1$-valent case, possibly composed with an autoequivalence $E:\N\to\N$. For instance, we can draw the following diagram:
\begin{equation}\label{eq:notation_schober1}
\begin{tikzcd}[column sep=large]
{} & F \arrow[r, "{(E\circ \varrho_2,\varrho_1)}"', no head] \arrow[l, "\varrho_1", no head] & F' \arrow[r, "\varrho_2"', no head] & {}
\end{tikzcd}
\end{equation}
This notation in \cref{eq:notation_schober1} thus describes a perverse schober on a $2$-valent graph with $2$ vertices. In the exit path description, this schober corresponds to the following diagram
\[
\begin{tikzcd}
\N & \V\overset{\rightarrow}{\times}_F \N \arrow[l, "\varrho_1"'] \arrow[r, "E\circ \varrho_2"] & \N & \V\overset{\rightarrow}{\times}_{F'} \N \arrow[r, "\varrho_2"] \arrow[l, "\varrho_1"'] & \N\,.
\end{tikzcd}
\]
\qed
\end{notation}

\begin{remark}\label{rem:lift_amounts_to_relative_triangles}
Let $F\colon \V\to \N$ be a spherical functor and $Y_1,Y_2\in \N$ two objects. We can ask what are the lifts of $Y_1,Y_2$ to an object $W\in \V\overset{\rightarrow}{\times}_F \N$ satisfying that $\varrho_1(W)\simeq Y_1$, $\varrho_2(W)\simeq Y_2$. To specify such an object $W=(X,Y_1,\eta)$, we must supply an object $X\in \V$, together with a fiber and cofiber sequence in $\N$ 
\begin{equation}\label{eq:Frelative}
F(X)\xlongrightarrow{\eta} Y_1\longrightarrow Y_2\,.
\end{equation}
Thus, the $\infty$-category $\V\overset{\rightarrow}{\times}_F \N$ can be understood as parametrizing $F$-relative triangles in $\N$, as in \cref{eq:Frelative}. Lifting a collection of local sections on the edges of $\rgraph$ of a $\rgraph$-parameterized perverse schober, which is locally modeled as above, to a global section thus amounts to specifying such a relative triangle for each vertex of $\rgraph$.\qed
\end{remark}


\begin{center}
	\begin{figure}[h!]
		\centering
		\includegraphics[scale=0.8]{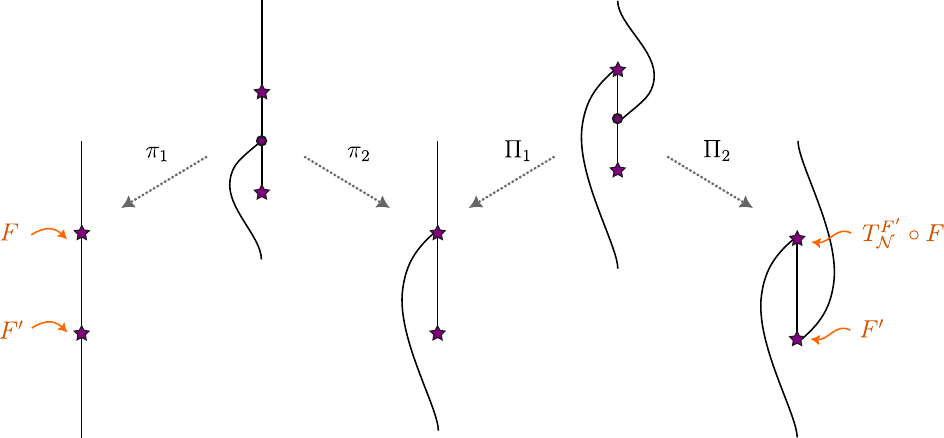}
		\caption{A sequence of spans of contractions realizing the clockwise $\pi$-rotation for the leftmost ribbon graph $\rgraph$. The key functors for the schober $\F$ on the left are indicated in orange, and those for the resulting schober $(\Pi_2)_*\Pi_1^*(\pi_2)_*\pi_1^*(\mathcal{F})$, depicted on the right.
        }\label{fig:ZigZagCorrespondence_Clockwise}
	\end{figure}
\end{center}

\subsubsection{Hurwitz moves for perverse schobers}

\begin{center}
	\begin{figure}[h!]
		\centering
		\includegraphics[scale=0.8]{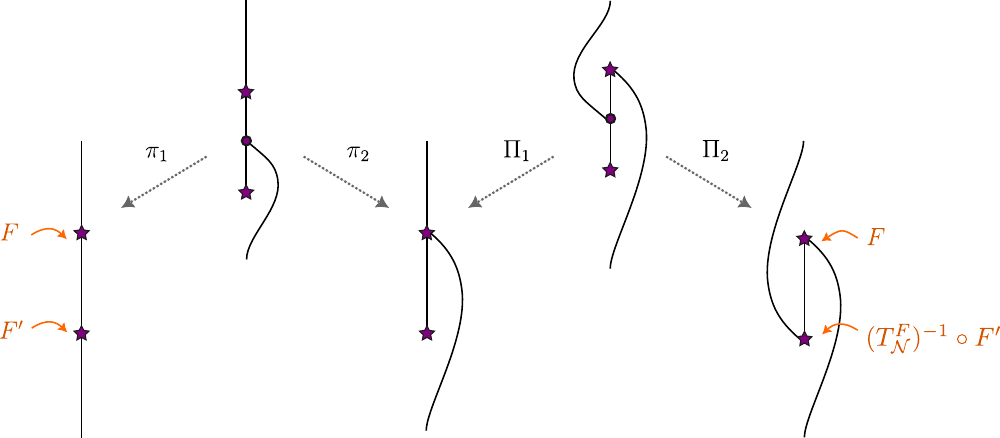}
		\caption{Similar to \Cref{fig:ZigZagCorrespondence_Clockwise}, this depicts a sequence of spans of contractions realizing the counterclockwise $\pi$-rotation for the leftmost ribbon graph $\rgraph$. The schober $\F$ on the left is written in orange, and the resulting schober $(\Pi_2)_*\Pi_1^*(\pi_2)_*\pi_1^*(\mathcal{F})$ is depicted on the right.
        }\label{fig:ZigZagCorrespondence_CounterClockwise}
	\end{figure}
\end{center}

Consider a ribbon graph $\rgraph$ with two adjacent $2$-valent vertices $v,v'$, cf.~\cref{fig:ZigZagCorrespondence_CounterClockwise} (left). We can perform two rotations by $\pi$ exchanging $v$ and $v'$, either counterclockwise, as in \cref{fig:ZigZagCorrespondence_CounterClockwise}, or clockwise, as illustrated in \cref{fig:ZigZagCorrespondence_Clockwise}. The following lemma describes how a $\rgraph$-parametrized perverse schober changes when performing such rotations, which is a schober analogue of the classical Hurwitz moves, cf.~\cite[Section II.6]{Hurwitz1891} or \cite[Section (16d)]{Sei08}. In its statement and proof we use \cref{not:schobers} above.

\begin{lemma}\label{lem:exchangesingularities} 
\begin{enumerate}[$(1)$]
\item Consider the two spans of contractions of ribbon graphs in \Cref{fig:ZigZagCorrespondence_Clockwise}. Denote the leftmost ribbon graph by $\rgraph$ and the rightmost ribbon graph by $\rgraph'$. Let $F,F'\colon \V\lr\N$ be two spherical functors and $\mathcal{F}$ the $\rgraph$-parametrized perverse schober
\[
\mathcal{F}\quad\coloneqq \begin{tikzcd}
{} \arrow[r, "\rho_1", no head] & F \arrow[r, "{(\rho_2,\rho_1)}", no head] & F' & {} \arrow[l, "\rho_2"', no head].
\end{tikzcd}
\]
Then the pull-push $(\Pi_2)_*\Pi_1^*(\pi_2)_*\pi_1^*$ along these contractions via \Cref{prop:contraction} yields, up to equivalence, the following $\rgraph'$-parametrized perverse schober
\[
(\Pi_2)_*\Pi_1^*(\pi_2)_*\pi_1^*(\mathcal{F})\quad=\begin{tikzcd}
{} & F' \arrow[l, "{\rho_1}"', no head] & T_{\mathcal{N}}^{F'}\circ F \arrow[l, "{(\rho_2,\rho_1)}"', no head] \arrow[r, "\rho_2", no head] & {}
\end{tikzcd}
\]
where $T_{\mathcal{N}}^{F'}=\on{cof}((F')^R\circ F'\xrightarrow{\on{cu}}\on{id}_\N)$ is the twist functor. 

\item Consider the two spans of contractions of ribbon graphs in \Cref{fig:ZigZagCorrespondence_CounterClockwise}. Denote the leftmost ribbon graph by $\rgraph$ and the rightmost ribbon graph by $\rgraph'$. Let $F,F':\V\lr\N$ be two spherical functors and $\mathcal{F}$ the $\rgraph$-parametrized perverse schober
\[
\mathcal{F}\quad\coloneqq \begin{tikzcd}
{} \arrow[r, "\rho_1", no head] & F \arrow[r, "{(\rho_2,\rho_1)}", no head] & F' & {} \arrow[l, "\rho_2"', no head]
\end{tikzcd}
\]
Then the pull-push $(\Pi_2)_*\Pi_1^*(\pi_2)_*\pi_1^*$ along these contractions via \Cref{prop:contraction} yields, up to equivalence, the following $\rgraph'$-parametrized perverse schober
\[
(\Pi_2)_*\Pi_1^*(\pi_2)_*\pi_1^*(\mathcal{F})\quad=\begin{tikzcd}
{} & (T_{\mathcal{N}}^{F})^{-1}\circ F' \arrow[l, "{\rho_1}"', no head] & F \arrow[l, "{(\rho_2,\rho_1)}"', no head] \arrow[r, "\rho_2", no head] & {}
\end{tikzcd}
\]
\end{enumerate}
\qed
\end{lemma}

\begin{proof}[Proof of \Cref{lem:exchangesingularities}.]
Let us focus on proving (1), as (2) follows from (1). Throughout the proof we will denote the twist functor $T_\N^{F'}$ by $T^{F'}$ and the twist functor of $F\dashv F^r$ by $T^F$. First, the perverse schober $\F$ is equivalent to

\begin{equation}\label{eq:proof_schober1}
\F=
\begin{tikzcd}
{} \arrow[r, "\rho_1", no head] & F \arrow[r, "{(\rho_2,\rho_1)}", no head] & F' & {} \arrow[l, "\rho_2"', no head]
\end{tikzcd}
\simeq 
\begin{tikzcd}
{} \arrow[r, "\rho_1", no head] & F \arrow[r, "{(\rho_2,\rho_2)}", no head] & F' & {} \arrow[l, "{T^{F'}\rho_1}"', no head]
\end{tikzcd}
\end{equation}

\noindent Indeed, \cref{eq:proof_schober1} follows from the local rotation formula for a perverse schober, cf.~\cite[Prop.~3.11]{Chr22}. Let us expand the right hand side of \cref{eq:proof_schober1} by adding a non-singular vertex between the vertices associated to $F$ and $F'$, according to \cref{fig:ZigZagCorrespondence_Clockwise} so as to model $\pi_1^*(\F)$. We denote the resulting schober by using \cref{not:schobers} analogously for schobers on a 3-spider, see also \cite[Lemma 3.3]{Chr22b} for the definition of the $\varrho_i$ functors, so that $\pi_1^*(\F)$ can be diagrammatically expressed as

\begin{equation}\label{eq:proof_schober2}
\pi_1^*(\F)\coloneqq 
\begin{tikzcd}
{} \arrow[r, "\varrho_1", no head] & F \arrow[r, "{(\rho_2,\rho_2)}"', no head] & 0 \arrow[rr, "T^{F'}\rho_1", no head, bend right=49] & F' \arrow[l, "{(\rho_3,\rho_1)}"', no head] & {}
\end{tikzcd}
\simeq 
\begin{tikzcd}
{} \arrow[r, "{(T^{F})^{-1}\varrho_2}", no head] & F \arrow[r, "{(\rho_1,\rho_3)}"', no head] & 0 \arrow[rr, "T^{F'}\rho_2", no head, bend right=49] & F' \arrow[l, "{(\rho_1[1],\rho_1)}"', no head] & {}
\end{tikzcd}
\end{equation}

\noindent where, as before, \cref{eq:proof_schober2} is the equivalence given by the local rotation formula, cf.~\cite[Prop.~3.11]{Chr22}. By collapsing the non-singular vertex, cf.~\cite[Lemma 4.26]{Chr22}, the pushforward schober $(\pi_2)_*\pi_1^*\F$ thus becomes

\begin{equation}\label{eq:proof_schober3}
(\pi_2)_*\pi_1^*\F=
\begin{tikzcd}
{} \arrow[r, "{(T^{F})^{-1}\varrho_3}"', no head] & F \arrow[rr, "T^{F'}\rho_2", no head, bend right=49] & F' \arrow[l, "{(\rho_1[1],\rho_1)}"', no head] & {}
\end{tikzcd}
\simeq 
\begin{tikzcd}
{} \arrow[r, "{\varrho_1[1]}"', no head] & F \arrow[rr, "T^{F'}\rho_3", no head, bend right=49] & F' \arrow[l, "{(\rho_2[1],\rho_1)}"', no head] & {}
\end{tikzcd}
\end{equation}
\noindent Let us iterate this expansion-contraction with $\Pi_1,\Pi_2$. The pull-back $\Pi_1^*(\pi_2)_*\pi_1^*(\F)$ is obtained by adding a non-singular vertex according to \cref{fig:ZigZagCorrespondence_Clockwise}, which yields
\begin{equation}\label{eq:proof_schober4}
\Pi_1^*(\pi_2)_*\pi_1^*(\F)=
\begin{tikzcd}[column sep=50]
{} \arrow[rr, "{\varrho_1[1]}", no head, bend left=30] & F \arrow[r, "{(\varrho_1,\varrho_3)}", no head] & 0 & F' \arrow[l, "{(\rho_2[1],\rho_1)}"', no head] & {} \arrow[lll, "T^{F'}\rho_2", no head, bend left=20]
\end{tikzcd}
\end{equation}
\noindent By rotating along the non-singular vertex, we have the equivalence
\begin{equation}\label{eq:proof_schober4}
\Pi_1^*(\pi_2)_*\pi_1^*(\F)\simeq
\begin{tikzcd}[column sep=50]
{} \arrow[rr, "{\varrho_2}", no head, bend left=30] & F \arrow[r, "{(\varrho_1,\varrho_1)}", no head] & 0 & F' \arrow[l, "{(\rho_3,\rho_1)}"', no head] & {} \arrow[lll, "T^{F'}\rho_2", no head, bend left=20]
\end{tikzcd}
\end{equation}\label{eq:proof_schober5}
Therefore, the contraction along $\Pi_2$ yields
\begin{equation}\label{eq:proof_schober5}
\left(\Pi_2\right)_*\Pi_1^*(\pi_2)_*\pi_1^*(\F)=
\begin{tikzcd}[column sep=large]
{} \arrow[rr, "\varrho_2", no head, bend left=30] & F \arrow[r, "{(\varrho_1,\varrho_1)}", no head] & F' &  {} \arrow[ll, "T^{F'}\rho_2", no head, bend left=20]
\end{tikzcd}
=
\begin{tikzcd}
{} \arrow[r, "\varrho_2", no head] & F' \arrow[r, "{(\varrho_1,\varrho_1)}", no head] & F \arrow[r, "T^{F'}\varrho_2", no head] & {}
\end{tikzcd}
\end{equation}
where the right hand side of \cref{eq:proof_schober5} is notational re-writing, and thus an actual identity, as opposed to a more general equivalence. By first rotating along the leftmost vertex $F'$ of \cref{eq:proof_schober5} and passing to the inverse we have the equivalences

\begin{equation*}
\left(\Pi_2\right)_*\Pi_1^*(\pi_2)_*\pi_1^*(\F)
\simeq \begin{tikzcd}[column sep=50]
{} \arrow[r, "\varrho_1", no head] & F' \arrow[r, "{((T^{F'})^{-1}\varrho_2,\varrho_1)}", no head] & F \arrow[r, "T^{F'}\varrho_2", no head] & {}
\end{tikzcd} \simeq \begin{tikzcd}[column sep=large]
{} \arrow[r, "\varrho_1", no head] & F' \arrow[r, "{(\varrho_2,T^{F'}\varrho_1)}", no head] & F \arrow[r, "T^{F'}\varrho_2", no head] & {}
\end{tikzcd}
\end{equation*}
\noindent Therefore we conclude that
\[\left(\Pi_2\right)_*\Pi_1^*(\pi_2)_*\pi_1^*(\F)\simeq \begin{tikzcd}
{} \arrow[r, "\varrho_1", no head] & F' \arrow[r, "{(\varrho_2,\varrho_1)}", no head] & T^{F'}\circ F \arrow[r, "\varrho_2", no head] & {}
\end{tikzcd}\,.
\]
\end{proof}

\section{The perverse schober of a Demazure weave}\label{sec:weave_schobers}

This section introduces the perverse schober $\F_\w$ associated to a Demazure weave $\w$ in \cref{ssec:defining_weave_schobers}, and establishes the behavior of its global sections under weave equivalence and weave mutation, in \cref{ssec:weaveschober_equivalences_mutations}.


\begin{center}
	\begin{figure}[h!]
		\centering
		\includegraphics[scale=1.2]{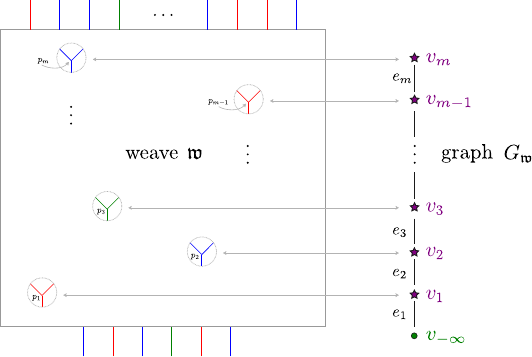}
		\caption{A weave $\w$, depicted left, and its associated graph $\Wgraph$, drawn to its right. The trivalent vertices $p_1,\ldots,p_m$ of $\w$ and the associated vertices $v_1,\ldots,v_m$ are also depicted, with $v_i$ at the same height as $p_i$, $i\in[1,m]$. The edges $e_1,\ldots,e_m$ of $\Wgraph$ are also drawn, with $e_i$ denoting the edge between $v_{i-1}$ and $v_i$. The weave schober $\wsch$ is a perverse schober on $\Wgraph$.}\label{fig:Weave_GraphForWeaveSchober}
	\end{figure}
\end{center}


\subsection{Definition of weave schobers}\label{ssec:defining_weave_schobers}

Let $\w:\beta\lr\delta(\beta)$ be an $n$-stranded Demazure weave. The height function endows the trivalent vertices of $\w$ with a total order, increasing as we read bottom-to-top. We assume that $\w$ is generic, in that no two vertices of $\w$ are at the same height. For instance, if we consider $\w$ as a raked weave as in \cite[Section 5.1.1]{CGGLSS25}, the horizontal right semirays $\{r_v\}_{v\in V(\w)}$ do not intersect each other. These semirays $r_v$ are depicted in an example as the dashed yellow segments in \Cref{fig:Example_WeaveSchober1}. As before, we let $p_1,\dots,p_m$ be the trivalent vertices of $\w$, ordered from bottom-to-top by height, cf.~e.g.~\Cref{fig:Weave_GraphForWeaveSchober}. Let us now start describing the perverse schober $\F_\w$ associated to $\w$.\\

\noindent {\bf Domain of definition}. Let $\Wgraph$ be the linear spanning graph of the $1$-gon with $m$ vertices, as depicted in \Cref{fig:Weave_GraphForWeaveSchober}. In practice, we draw the graph $\Wgraph$ such that:

\begin{itemize}
    \item[(i)] $\Wgraph$ is drawn vertically to the right side of $\w$,

    \item[(ii)] The $m$ vertices of $\Wgraph$ are labeled by $v_1,\dots,v_m$ from bottom-to-top, and the edges by $e_{1},\dots,e_m$, with $e_1$ being the bottom external (non-compact) edge.

    \item[(iii)] The $i$th vertex $v_i$ of $\Wgraph$ is drawn at the same height as the trivalent vertex $p_i\in \w$, so that $v_m$ lies on top, $v_1$ on the bottom, and the external edge lies below $v_1$.
\end{itemize}

\begin{center}
	\begin{figure}[h!]
		\centering
		\includegraphics[scale=1.6]{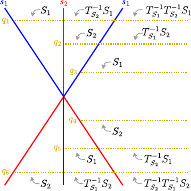}
		\caption{The local model for a 6-valent vertex at a weave $\w$, analyzed in \cref{ex:Objects_Weave_HexavalentVertex}. The objects $S_{q_i}$ associated to the points $q_i\in \w$ are depicted. If $q_i$ has color $i_j\in S_n$, these objects $S_{q_i}$ are obtained by following the yellow dashed lines to the right from $q_i$, transforming the initial $S_{i_j}$ at $q_i$ by twisting by $T_{S_{i_k}}$ every time a weave $i_k$-line is crossed.}\label{fig:Weave_HexavalentVertex}
	\end{figure}
\end{center}

\begin{center}
	\begin{figure}[h!]
		\centering
		\includegraphics[scale=1.2]{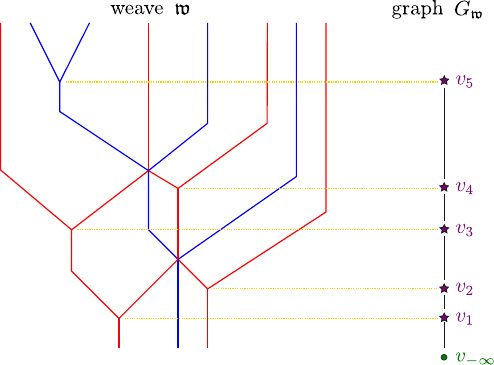}
		\caption{The 3-stranded weave $\w$ for \cref{ex:running_example1}. The positive braid word at the top is $\beta=s_2s_1(s_1s_2)^3$, the braid word at the bottom is its Demazure product $\delta(\beta)=\s_2\s_1\s_2$. The graph $G_\w$ over which the weave schober $\F_\w$ is depicted to the right of $\w$. The vertices of $G_\w$ are labeled $v_1,v_2,v_3,v_4,v_5$ bottom to top. The south point $\vinf$ is depicted in a round green circle.}\label{fig:Example_WeaveSchober1}
	\end{figure}
\end{center}

\noindent See Figures \ref{fig:Weave_HexavalentVertex} and \ref{fig:Example_WeaveSchober1} for illustrations of such $\Wgraph$. We also consider a generic point $\vinf$ in the non-compact edge $e_1$. By construction, the perverse schober $\F_\w$ associated to $\w$ will be a perverse schober parametrized by $\gw$.\\

\noindent {\bf Generic stalks}. We choose the generic stalks of the perverse schober $\F_w$, i.e.~the stalks at the points of any edge, to be the perfect derived $\infty$-category $\D\coloneqq\D^{\on{perf}}(\Pi_2(A_{n}))$ of the (additive) derived preprojective algebra $\Pi_2(A_{n})$ of type $A_n$, also known as the $2$-Calabi--Yau completion of type $A_n$ \cite{Kel11}. This choice is in part motivated by the equivalence between $\D$ and the wrapped Fukaya category of the 4-dimensional $A_n$-Milnor fibre, cf.~e.g.~\cite[Section 3]{LU21}.\footnote{The additive derived preprojective algebra $\Pi_2(A_{n})$ is equivalent to the multiplicative derived preprojective algebra of type $A_n$, and this is true irrespective of the characteristic of the field $k$.} We denote by $S_1,\dots,S_n\in \D$ the 1-dimensional modules associated with the $n$ vertices of the $A_n$-quiver, i.e.~the simple modules of the vertices in the standard heart of $\D$. We note that $S_i\in \D$ is $2$-spherical for all $1\leq i\leq n$, that $\on{Ext}^1(S_i,S_j)\simeq k$ if $|i-j|=1$, and $\on{Ext}^1(S_i,S_j)\simeq 0$ if $|i-j|\not =1$.\\

\noindent {\bf Spherical functors}. The spherical functors describing $\F_\w$ will be defined by spherical objects, which are themselves constructed as follows. Let $p\in \w$ be a point on the Demazure weave, consider the horizontal right semiray $r_p$ based at $p$, and suppose that $r_p$ does not intersect any vertex of $\w$. The intersections of $r_p$ with $\w$ determine a positive braid word $\beta_p=s_{i_0}s_{i_1}\dots s_{i_{l(p)}}$ by reading the $s_i$-labels in the associated weave edges left-to-right, cf.~\cref{fig:Weave_HexavalentVertex}. In these hypotheses and notation, we define:

\begin{definition}[Spherical objects from weave points]\label{def:sphericalobj}
The spherical object associated to $p\in\w$ is
\[ S_p\coloneqq  S_{\beta_p}=T^{-1}_{S_{i_{l(p)}}}\cdots T^{-1}_{S_{i_{2}}}T^{-1}_{S_{i_{1}}}(S_{i_0})\in \D \]
where $T_{S_j}=\on{cofib}(\on{Mor}_\D(S_j,\mhyphen)\otimes S_j\to \on{id}_\D)$ is the twist functor.\qed
\end{definition}

\noindent In the case that $p\in\w$ is a trivalent vertex, $S_p$ in \cref{def:sphericalobj} coincides with the $\w$-transport of $S_{i_0}$ along the right thimble $\wp_r(p)$, see \cref{ssec:weave_transport} below. More generally, \cref{def:sphericalobj} should be understood as the transport of $S_{i_0}$ along the horizontal right semiray $r_p$. In \cref{def:sphericalobj}, the inverse $T_{S_{i_j}}^{-1}$ of the twist functor $T_{S_{i_j}}$ exists because $S_{i_j}$ is a 2-spherical object, and thus $T_{S_{i_j}}$ is an auto-equivalence by \cite[Prop.~2.10]{ST01}. In this case we have the formula $T_{S_{i_j}}^{-1}\simeq \on{fib}(\on{id}_\D\to\on{Mor}(\mhyphen,S_{i_j})^\ast\otimes S_{i_j})$, cf.~loc.~cit.

\begin{example}\label{ex:Objects_Weave_HexavalentVertex}
Consider the 3-stranded weave $\w$ in \Cref{fig:Weave_HexavalentVertex}, which is the local model for a 6-valent vertex in a weave $\w$. The positive braid word at the top is $s_1s_2s_1$, and the one at the bottom is $s_2s_1s_2$. Given the points $q_i\in\w$ and their associated dashed yellow lines $r_{q_i}$ as depicted, $i\in[1,6]$, the associated spherical objects are

\begin{align*} 
S_{q_1}&\simeq \itw1\itw2S_1\simeq S_2 & S_{q_2}&\simeq \itw1S_2 & S_{q_3}&\simeq S_1\\ 
S_{q_4}&\simeq S_2 & S_{q_5}&\simeq \itw2S_1 & S_{q_6}&\simeq \itw2\itw1S_2\simeq S_1,
\end{align*}
where we have used the equivalence $\tw1 S_2\simeq T_{S_2}^{-1}S_1$. These are the objects at the rightmost end of the dashed yellow lines in \Cref{fig:Weave_HexavalentVertex}.\qed
\end{example}

\noindent {\bf Weave schobers}. Given a spherical object $S\in \D$, we denote by $\mhyphen\otimes S\colon \D^{\on{perf}}(k)\to \D$ the corresponding spherical functor. Let us use \Cref{not:schobers} to define the perverse schober $\F_\w$ of $\w$:

\begin{definition}[Perverse schober of a weave]\label{def:weave_schobers}
Let $\w$ be a Demazure weave with trivalent vertices $p_1,\dots,p_m$ and associated graph $\gw$. By definition, the perverse schober $\F_\w$ associated to $\w$ is the $\gw$-parametrized perverse schober given by
\[
\begin{tikzcd}
\mhyphen\otimes S_{p_m} \arrow[d, "{(\varrho_1,\varrho_1)}", no head] \\
\mhyphen\otimes S_{p_{m-1}} \arrow[d, "{(\varrho_2,\varrho_1)}", no head] \\
\vdots \arrow[d, "{(\varrho_2,\varrho_1)}", no head]                   \\
\mhyphen\otimes S_{p_1} \arrow[d, "\varrho_2", no head]               \\
{}                                                                   
\end{tikzcd}
\]
We refer to any such perverse schober as a weave schober.\qed
\end{definition}

\begin{center}
	\begin{figure}[h!]
		\centering
		\includegraphics[scale=1.2]{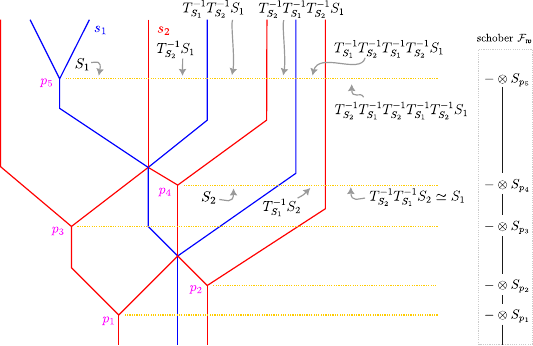}
		\caption{The weave $\w$ from \cref{fig:Example_WeaveSchober1} featuring in \cref{ex:running_example1}. The left-to-right transport of the simple $S_2\in\D$, resp.~$S_1$, for the weave trivalent vertex $p_4$, resp.~$p_5$, is illustrated.\\}
        \label{fig:Example_WeaveSchober2}
	\end{figure}
\end{center}

\begin{example}\label{ex:running_example1} Consider the 3-stranded weave $\w$ in \cref{fig:Example_WeaveSchober1}, so $n=3$ and $m=5$. The graph $G_\w$ is depicted to its right, with the horizontal slicings depicted in dashed yellow lines. These lines match each trivalent vertex $p_i$ of the weave $\w$ to a vertex $v_i$ of the graph $G_\w$, $i\in[1,m]$. The weave schober $\F_\w$ is determined by the spherical objects $S_{p_i}$ at the vertex $v_i$. The spherical object $S_{p_i}$ is determined by $\w$:  starting with the simple $S_{i_j}$ at the trivalent $p_i$, where $i_j$ is the weave color of $p_i$, we transport $S_{i_j}$ to the right along the dashed yellow lines. By \cref{def:sphericalobj}, the rule for such transport is that we apply the inverse twist $T_{S_{i_k}}^{-1}$ as we transport across a weave $i_k$-line, cf.~\cref{ssec:weave_transport}. Such transport for $i=4,5$ is depicted in \cref{fig:Example_WeaveSchober2}. Specifically, the spherical objects for this weave schober $\F_\w$ are

\[S_{p_1}\simeq T^{-1}_{S_2}T^{-1}_{S_1}(S_2)\simeq S_1,\qquad S_{p_2}\simeq S_2,\qquad S_{p_3}\simeq (T_{S_1}T_{S_2})^{-2}(S_2)\simeq T^{-1}_{S_2}(S_1[1])\,,\]
\[S_{p_4}\simeq T^{-1}_{S_2}T^{-1}_{S_1}(S_2)\simeq S_1, \qquad S_{p_5}\simeq T^{-1}_{S_2}(T_{S_2}T_{S_1})^{-2}(S_1)\simeq T^{-1}_{S_2}T^{-1}_{S_1}(S_2[1])\simeq S_1[1]\,\]

\noindent where we have used the three equivalences $T_{S_2}^{-1}T_{S_1}^{-1}(S_2)\simeq S_1$, or equivalently $T_{S_2}(S_1)\simeq T_{S_1}^{-1}(S_2)$, $T_{S_1}^{-1}T_{S_2}^{-1}(S_1)\simeq S_2$, and $T_{S_i}^{-1}(S_i)\simeq S_i[1]$, for $i=1,2$.\qed
\end{example}

\begin{remark}
    The stable $\infty$-category of global sections $\glsecF$ can be more concretely understood as the perfect derived $\infty$-category of the directed dg category on the object $S_{p_1},\dots,S_{p_m}$ in the dg avatar of $\D$ in the sense of \cite{Sei08}, see also the notes \cite{Christ_excep_coll} for more details. 
\end{remark}

We denote the left adjoint of the (spherical) evaluation functor $\on{ev}_{e_1}\colon \glsecF\to \mathcal{F}(e_1)$ of global section to the bottom edge $e_1$ by $\on{ind}^L_{e_1}$.

\begin{proposition}\label{prop:rel_3_CY_str}
Let $\w$ be a Demazure weave. The functor 
\[ \on{ind}^L_{e_1}\colon \D=\mathcal{F}_\w(e_1)\to \glsecF\] admits a weak relative left $3$-Calabi--Yau structure and the functor 
\[
\on{ev}_{e_1}\colon \glsecF\to \mathcal{F}_\w(e_1)=\D
\]
admits a weak relative right $3$-Calabi--Yau structure. 
\end{proposition}

\begin{proof}
For each $2$-spherical object $S\in \D$, by \cite[Lem.~6.4]{Chr23}, the functor $\on{Mor}_\D(S,\mhyphen)\colon \D\to \D^{\on{perf}}(k)$ admits a weak relative left $3$-Calabi--Yau structure, and the functor $(\mhyphen)\otimes S\colon \D^{\on{perf}}(k)\to \D$ admits a weak relative right $3$-Calabi--Yau structure. The assertion then follows from \cite[Prop.~5.2]{Chr23} and \cite[Thm.~5.7]{Chr23}. 
\end{proof}

\begin{remark}[On generalizing to finite Dynkin type]\label{rmk:Dynkincase} The calculus for Demazure weaves is defined for any finite Dynkin type $\mathfrak{g}$, cf.~\cite[Sections 4 \& 6]{CGGLSS25}. Specifically, the simply-laced cases are constructed directly, and the non-simply-laced cases via folding, cf.~\cite[Section 6.2]{CGGLSS25}. In line with this, the weave schober $\F_\w$ in \cref{def:weave_schobers} can be constructed analogously in any of the simply-laced cases, where the stalk category $\D$ is taken to be the derived $\infty$-category of the derived preprojective algebra of the corresponding type $\mathfrak{g}$.\footnote{If $\mbox{char}(k)=0$, the additive and multiplicative preprojective algebras are isomorphic for simply-laced Dynkin cases.} The folding technique applies to extend the construction of the weave schober in the non-simply-laced cases. In this more general framework, it is our conceptual understanding that the main results in this manuscript generalize to any finite Dynkin type in characteristic 0. A number of technical details would need to be written and verified, and thus we leave the development of categorical weave calculus for any finite Dynkin type as a conjectural generalization.\qed
\end{remark}


\subsection{The \texorpdfstring{schober $\F_\w$}{weave schober} under weave equivalences and mutations}\label{ssec:weaveschober_equivalences_mutations}

Consider a weave equivalence $\w\simeq\w'$ between Demazure weaves $\w,\w'$, cf.~\cref{sssec:weave_equivalences_mutations} or \cite[Section 4.2]{CGGLSS25}. By definition, such an equivalence is said to be height-preserving if the height-induced order of the trivalent vertices of $\w$ coincides with that of the trivalent vertices of $\w'$. Any weave equivalence can be expressed as a (possibly alternating) composition of height-preserving weave equivalences and the local weave equivalence depicted in \Cref{fig:Weave_HeightChangeEquivalence}. We refer to these particular weave equivalences as pulley moves: specifically, a forward pulley move exchanges \Cref{fig:Weave_HeightChangeEquivalence} (left) for \Cref{fig:Weave_HeightChangeEquivalence} (right), and a backwards pulley move does the reverse.\\

\noindent The distinction between height-preserving weave equivalences and pulley moves is crucial in the study of weave schobers. Indeed, the exchange of heights for trivalent vertices of a weave, whether caused by a pulley move or a weave mutation, is what modifies the weave schober in an interesting manner. The precise statement reads as follows:

\begin{center}
	\begin{figure}[h!]
		\centering
		\includegraphics[scale=1.2]{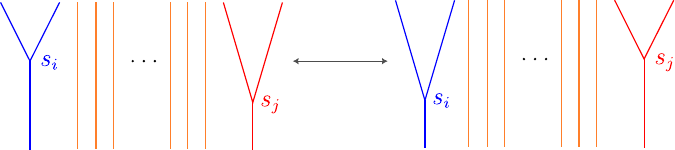}
		\caption{The local weave equivalence that exchanges the height of trivalent vertices, referred to as a pulley move. In this equivalence, there are {\it no} trivalent vertices whose height is between those of the two trivalent vertices depicted in this picture. The weave lines in color orange indicate arbitrary and possibly different $s_k\in W$ labels, and the two permutations $s_i,s_j\in W$ are also arbitrary.\\}
        \label{fig:Weave_HeightChangeEquivalence}
	\end{figure}
\end{center}

\begin{theorem}\label{thm:weaveequivglsec} In the notation above, the following holds:

\begin{enumerate}[(1)]
\item Let $\w\simeq \w'$ be a height-preserving weave equivalence. Then, there is an equivalence of weave schobers $\mathcal{F}_{\w}\simeq \mathcal{F}_{\w'}$.\\

\item Let $\w\simeq \w'$ be a forward pulley move or $\w,\w'$ be related by a forward weave mutation, exchanging the heights of the vertices $p_{i-1}$ and $p_i$. Then, the weave schobers $\mathcal{F}_{\w},\mathcal{F}_{\w'}$ are related as follows:
\begin{equation}\label{eq:forward_schober}
\mathcal{F}_{\w}= \begin{tikzcd}
\mhyphen\otimes S_{p_m} \arrow[d, "{(\varrho_1,\varrho_1)}", no head]     \\
\vdots \arrow[d, "{(\varrho_2,\varrho_1)}", no head]                       \\
\mhyphen\otimes S_{p_{i}} \arrow[d, "{(\varrho_2,\varrho_1)}", no head]   \\
\mhyphen\otimes S_{p_{i-1}} \arrow[d, "{(\varrho_2,\varrho_1)}", no head] \\
\vdots          
\end{tikzcd}\qquad\longrightarrow\qquad
\mathcal{F}_{\w'}\simeq \mu_i^\sharp(\mathcal{F}_\w)\coloneqq  \begin{tikzcd}
\mhyphen\otimes S_{p_m} \arrow[d, "{(\varrho_1,\varrho_1)}", no head]                         \\
\vdots \arrow[d, "{(\varrho_2,\varrho_1)}", no head]                                           \\
\mhyphen\otimes S_{p_{i-1}} \arrow[d, "{(\varrho_2,\varrho_1)}", no head]                     \\
\qquad \mhyphen\otimes T_{S_{p_{i-1}}}(S_{p_{i}}) \arrow[d, "{(\varrho_2,\varrho_1)}", no head] \\
\vdots\end{tikzcd}
\end{equation}

\item Let $\w,\w':\beta\to\delta(\beta)$ be two Demazure weaves and $\F_\w,\F_{\w'}$ their associated weave perverse schobers. Then there is an equivalence of $\infty$-categories of global sections
\begin{equation}\label{eq:equiv_globalsections}
\glsecF\simeq R\Gamma(\rgraph_{\w'},\mathcal{F}_{\w'})\,.
\end{equation}
\end{enumerate}
In particular, the $\infty$-category of global sections of a weave schober $\F_\w$ is an invariant of $\beta$ for any Demazure weave $\w:\beta\to\delta(\beta)$, independent of the choice of $\w$, up to equivalence.
\qed
\end{theorem}

\noindent The weave schober in the right of \Cref{eq:forward_schober} will also be denoted by $\mu_i^\sharp(\mathcal{F}_\w)$ and referred to as the schober obtained from $\mathcal{F}_\w$ via a forward or left (Hurwitz) move. The equivalence of global sections between the two perverse schobers $\mathcal{F}_\w$, $\mathcal{F}_{\w'}$, with $\w$ and $\w'$ related by a pulley move, is a schober-theoretic version of \cite[Lem.~5.22]{Sei08}.

\begin{remark}[Backward move for schobers]\label{rmk:backward_pulley_schober}
The analogue of \cref{eq:forward_schober} if $\w$ and $\w'$ are related by a backward pulley move or a backward mutation, exchanging the heights of the vertices $p_{i-1}$ and $p_i$, is as follows:
\begin{equation}\label{eq:backward_schober}
\mathcal{F}_{\w}= \begin{tikzcd}
\mhyphen\otimes S_{p_m} \arrow[d, "{(\varrho_1,\varrho_1)}", no head]     \\
\vdots \arrow[d, "{(\varrho_2,\varrho_1)}", no head]                       \\
\mhyphen\otimes S_{p_{i}} \arrow[d, "{(\varrho_2,\varrho_1)}", no head]   \\
\mhyphen\otimes S_{p_{i-1}} \arrow[d, "{(\varrho_2,\varrho_1)}", no head] \\
\vdots          
\end{tikzcd}\qquad\longrightarrow\qquad
\mathcal{F}_{\w'}\simeq \mu_i^\flat(\mathcal{F}_\w)\coloneqq \begin{tikzcd}
\mhyphen\otimes S_{p_m} \arrow[d, "{(\varrho_1,\varrho_1)}", no head]                         \\
\vdots \arrow[d, "{(\varrho_2,\varrho_1)}", no head]                                           \\
\qquad \mhyphen\otimes T^{-1}_{S_{p_{i}}}(S_{p_{i-1}}) \arrow[d, "{(\varrho_2,\varrho_1)}", no head]                     \\
\mhyphen\otimes S_{p_{i}} \arrow[d, "{(\varrho_2,\varrho_1)}", no head] \\
\vdots\end{tikzcd}
\end{equation}
This is proven similarly to \cref{thm:weaveequivglsec}.(2), or it can also be logically deduced from inverting it. The weave schober on the right of \Cref{eq:backward_schober} is denoted by $\mu_i^\flat(\mathcal{F}_\w)$ and referred to as the schober obtained from $\mathcal{F}_\w$ via a backward or right (Hurwitz) move.\qed
\end{remark}


\begin{center}
	\begin{figure}[h!]
		\centering
		\includegraphics[scale=0.7]{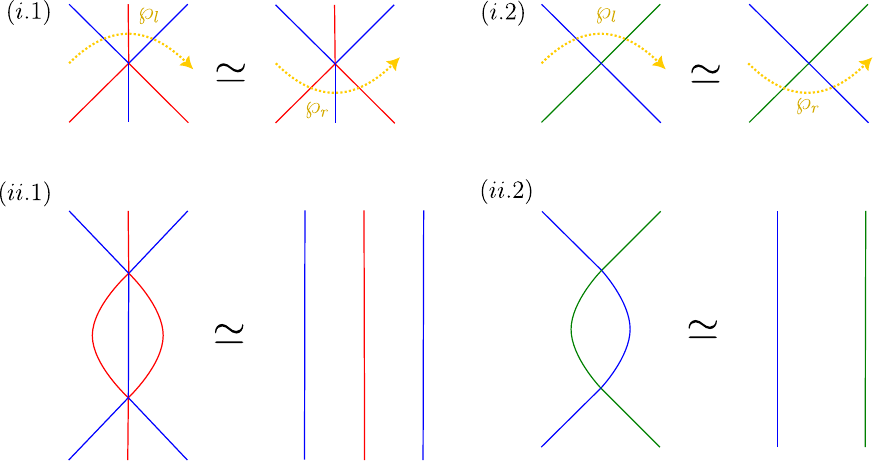}
		\caption{(i) A height exchange for a 6-valent vertex, in $(i.1)$, and for a 4-valent vertex, in $(i.2)$. (ii) Canceling pairs, for 6-valent vertices in $(ii.1)$, and for 4-valent vertices in $(ii.1)$.}
        \label{fig:WeaveEquivalences_Schober}
	\end{figure}
\end{center}

\subsubsection{Proof of \cref{thm:weaveequivglsec} and \ref{thm:main2}.(1)}  

\noindent The proofs of Parts (1) and (2) and \Cref{thm:weaveequivglsec} are mutually independent, and Part (3) follows from combining them. We refer to \cref{sssec:weave_equivalences_mutations} for the necessary ingredients on weave equivalences and weave mutations.\\

\begin{center}
	\begin{figure}[h!]
		\centering
		\includegraphics[scale=0.8]{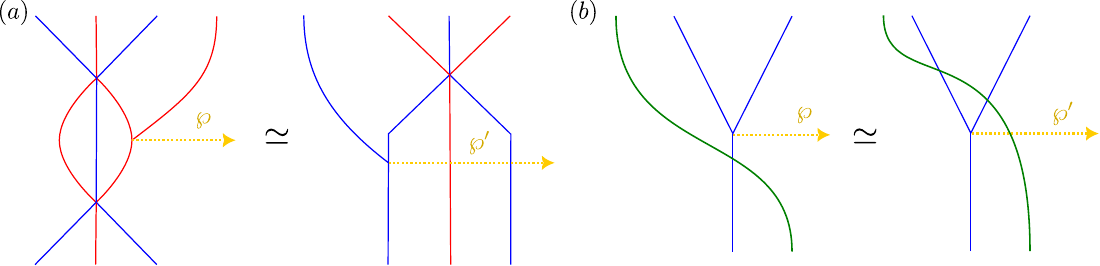}
		\caption{(a) The pushthrough move, where blue color for a weave line indicates $s_i$-label and red color indicates $s_{i+1}$-label. (b) A case of commutation of distant colors, with a weave line labeled by $s_j$, in green, passing through a trivalent $s_i$-vertex, $|i-j|\geq2$.\\}
        \label{fig:WeaveEquivalences2_Schober}
	\end{figure}
\end{center}

\noindent {\bf Proof of Part (1)}. It suffices to establish the equivalence of weave schobers by verifying it for the weave equivalences in \cite[Section 4.2]{CGGJ24}, cf.~also \cite[Theorem 4.2]{CasalsZaslow}. Specifically, these are:
\begin{itemize}
    \item[(i)] {\it Height exchange of 4- and 6-vertices}, as depicted in \Cref{fig:WeaveEquivalences_Schober}.(i).

    \item[(ii)] {\it Canceling pairs of 4- and 6-vertices}, as depicted in \Cref{fig:WeaveEquivalences_Schober}.(ii).

    \item[(iii)] {\it Commutation of distant colors}, which involve an $s_i$-weave line passing above a weave vertex involving $s_j$-labels with $|i-j|\geq2$, see e.g.~\Cref{fig:WeaveEquivalences2_Schober}.(b). Confer also Moves VI and VI' in \cite[Figures 24\&25]{CasalsZaslow}, or \cite[Subsection 4.2.3]{CGGJ24}.

    \item[(iv)] {\it Pushthrough move}, as depicted in \Cref{fig:WeaveEquivalences2_Schober}.(a).

    \item[(v)] {\it Zamolodchikov relation}, as depicted in \cite[Figure 106]{CasalsZaslow}, cf.~also Move IV in Figure 22 in ibid, or \cite[Subsection 4.2.6]{CGGJ24}.
\end{itemize}

\noindent For $(i)$ in the case of a 6-valent vertex, the invariance of the weave schober follows from the fact that the spherical objects $S_{p_i}$ do not change up to equivalence due to the braid relation
\begin{equation}\label{eq:braidrelation_twistfunctors}
T_{S_i}^{-1}T_{S_{i+1}}^{-1}T_{S_i}^{-1}\simeq T_{S_{i+1}}^{-1}T_{S_i}^{-1}T_{S_{i+1}}^{-1}.
\end{equation}
Indeed, if one considers the path $\wp_l$, resp.~$\wp_r$ in the left, resp.~right, of \Cref{fig:WeaveEquivalences_Schober}.(i.1), then for any object $L\in\D$ at the left endpoint of such paths the resulting object when traversing these paths left-to-right are equivalent by \Cref{eq:braidrelation_twistfunctors}.\footnote{In the notation of $\w$-transports, introduced in \Cref{def:combinatorial_transport}, this equivalence reads $\mbox{Hol}_\w(L,\wp_l)\simeq \mbox{Hol}_\w(L,\wp_r).$} It thus follows that any spherical object $S_{p_i}$ does not change when the height of a 6-valent vertex of $\w$ is changed. For $(i)$ in the case of a 4-valent vertex, the same argument applies but using the commutation relation
\begin{equation}\label{eq:commutationrelation_twistfunctors}
T_{S_i}^{-1}T_{S_j}^{-1}\simeq T_{S_j}^{-1}T_{S_i}^{-1},\quad |i-j|\geq 2\,.
\end{equation}

\noindent The invariance of $\F_\w$ under the weave equivalences $(ii)$, $(iii)$ for 4 and 6-valent vertices, and $(v)$ is implied similarly by Equations \ref{eq:braidrelation_twistfunctors} and \ref{eq:commutationrelation_twistfunctors}, as they do not involve 3-valent vertices. The remaining cases are $(iii)$, in the case of a 3-valent vertex, and the pushthrough move $(iv)$, as depicted in \Cref{fig:WeaveEquivalences2_Schober}. In the former case, invariance of $\F_\w$ follows from $T_{S_j}(S_i)\simeq S_i$ for $|i-j|\geq2$. In the latter case, it follows from the equivalence $T_{S_i}^{-1}T_{S_{i+1}}^{-1}S_{i}\simeq S_{i+1}$.

\noindent {\bf Proof of Part (2)}. In either a forward pulley move or a forward weave mutation, only the heights of the two involved trivalent vertices are exchanged. It thus follows that the spherical objects associated to any other trivalent vertices (other than these two) remain identical. By direct computation, the effect of such height exchange on the two spherical objects associated to these two trivalent vertices being exchanged yields \cref{eq:forward_schober}.\\

\noindent {\bf Proof of Part (3)}. By \cite[Lemma 4.4]{CGGLSS25} or \cite[Theorem 4.6]{CGGJ24} any two such weaves $\w,\w'$ are related by a sequence of weave equivalences and weave mutations. Therefore the statement follows from Part (1), and Part (2), combined with \cref{lem:exchangesingularities}.(1). Indeed, \cref{lem:exchangesingularities}.(1) supplies the required equivalence of global sections by noting that $\mu_i^\sharp(\mathcal{F}_\w)$ is the result of pull-push along the zig-zag realizing the clockwise $\pi$-rotation of the ribbon graph in \cref{fig:ZigZagCorrespondence_Clockwise}.\qed\\

\noindent Finally, we emphasize that \cref{thm:weaveequivglsec}.(3) is equivalent to Theorem \ref{thm:main2}.(1). Thus Theorem \ref{thm:main2}.(1) is now proven. The proofs of Parts (2), (3) and (4) of Theorem \ref{thm:main2} are more challenging, and are the content of \cref{sec:proof_Lusztig_cycles}.



\section{Categorical Lusztig cycles}\label{sec:lusztigcycles}

The object of this section is to construct the categorical Lusztig cycles associated to a weave $\w$, which will be global sections of the weave schober $\F_\w$, i.e.~objects in $\glsecF$. Their construction requires as input a number of new results that will be proven in \cref{sec:weave_propagation_LC}. Assuming these results for now, the construction is presented in \Cref{sssec:categoricalLusztigcycles_definition1}, while \Cref{ssec:results_on_categorical_lusztig_cycles} collects the main statements about them, which are proven in \cref{sec:proof_Lusztig_cycles}. We however first summarize the construction of the categorical Lusztig cycles as follows: \\

{\bf Outline of the construction}. Let $\w$ be a Demazure weave and $\{\g_i^\w\}$ its collection of Lusztig cycles, $i\in[1,m]$. For a fixed Lusztig cycle $\g_i^\w$, any horizontal slice $h\sse\mathbb{D}$ that does not intersect a weave vertex specifies a weighted braid word $(\b_h,{\bf a}_h)$. In a nutshell, the categorical Lusztig cycle ${\bm L}_i^\w$ is constructed by associating to each slice $h$ a local object $L_{\b_h,{\bf a}_h}$ with a filtration adapted to the weighted braid word $(\b_h,{\bf a}_h)$, and assembling such $L_{\b_h,{\bf a}_h}$ together as the horizontal slice $h$ scans the weave from the vertex $v_i$ downwards. Specifically, these two steps are as follows:\\

\noindent {\bf Step 1}. Given a weighted braid word $(\beta,{\bf a})$ arising as a slice of a Lusztig cycle, we will show:

\begin{equation}\label{eq:categorical_weighted_braid}
\fbox{$\exists!$ rigid filtered object $L_{\beta,{\bf a}}\in\D$ whose graded pieces recover $\bw$}
\end{equation}

\noindent This is established in \cref{thm:unique_Lusztig_filtered_object}, which refines the above statement.\\

\noindent {\bf Step 2}. Consider a vertex $p$ of the weave $\w$ and a Lusztig cycle $\g_i^\w$. Let us denote the rigid filtered objects associated to the slice above and below a trivalent vertex by $L_{\beta,{\bf a}}$ and $L_{\beta',{\bf a'}}$. These rigid objects are provided by Statement (\ref{eq:categorical_weighted_braid}) above. The relation between $L_{\beta,{\bf a}}$ and $L_{\beta',{\bf a'}}$ is as follows:

\begin{enumerate}
    
    \item If $p$ is trivalent, we will show in that there exists a cofiber sequence
    \begin{equation}\label{eq:cofiberseq_trivalent0}
    \tau_{\geq 0}\on{Mor}(S_{p},L_{\beta,{\bf a}})\otimes S_{p} \stackrel{\eta}{\longrightarrow} L_{\beta,{\bf a}}\longrightarrow L_{\beta',{\bf a'}}\,.
    \end{equation}
    Such sequence is proven in \Cref{prop:semi_twist_Lusztig_cycle_at_trivalent_vertex} and captures the data for a propagation of $L_{\beta,{\bf a}}$ past the trivalent, see also \Cref{rem:lift_amounts_to_relative_triangles} for more explanation.
\begin{equation}
\fbox{At a trivalent, $L_{\beta,{\bf a}}$ and $L_{\beta',{\bf a'}}$ are glued via the cofiber sequence (\ref{eq:cofiberseq_trivalent0}).}\\
\end{equation}

    \item If $p$ is hexavalent or tetravalent, we will show that $L_{\beta,{\bf a}}\simeq L_{\beta',{\bf a'}}$, cf.~\cref{prop:4_and_6_moves}. Then:
\begin{equation}\label{eq:filteredobjects_hexavalent}
\fbox{At a non-trivalent vertex, $L_{\beta,{\bf a}}$ and $L_{\beta',{\bf a'}}$ are identified via $L_{\beta,{\bf a}}\simeq L_{\beta',{\bf a'}}$.}\\
\end{equation}
\end{enumerate}

\noindent Let us for now assume Statement (\ref{eq:categorical_weighted_braid}), the cofiber sequence (\ref{eq:cofiberseq_trivalent0}), and the equivalence in Statement (\ref{eq:filteredobjects_hexavalent}), and proceed to construct the categorical Lusztig cycles.


\subsection{Definition of the categorical Lusztig cycles}\label{sssec:categoricalLusztigcycles_definition1}

Let $\w$ be a Demazure weave with $m$ trivalent vertices $p_1,\dots,p_m$, ordered from bottom to top. The vertices $v_1,\ldots,v_m$ and edges $e_1,\ldots,e_m$ of the graph $\rgraph_w$ are also ordered bottom to top, cf.~\cref{fig:Weave_GraphForWeaveSchober}. Given the weave schober $\F_\w$, we describe a global section $\mathfrak{s}\in\glsecF$ by specifying:

\begin{enumerate}
    \item Its value $\mathfrak{s}(v_j)$ at the vertices $v_j$, i.e.~an object in $\F_\w(v_j)\simeq \V\times_{F_j}\D\simeq \D^{\on{perf}}(k)\times_{F_j}\D$, if $i<m$, or in $\D^{\on{perf}}(k)$, if $i=m$, according to \cref{sssec:model_perverseschobers} and \cref{sec:weave_schobers}, where $F_j\coloneqq \mhyphen\otimes S_{p_j}\colon \D^{\on{perf}}(k)\to \D$.

    \item Its value $\mathfrak{s}(e_j)$ at the edges $e_j$, i.e.~an object in $\F_\w(e_j)\simeq\D$.

    \item Its value $\mathfrak{s}(v_j\to e_j)$ on the generization morphisms $v_j\to e_{j}$, which is an object in $\F_w(v_j\to e_j)$, i.e.~a functor $\F_w(v_j)\to\F_w(e_j)$, and similarly for $v_j\to e_{j+1}$, $j\in[1,m-1]$, and $v_m\to e_m$.\\
\end{enumerate}

\noindent These values for $\mathfrak{s}$ need to be compatible with each other and $\F_\w$ in order for them to define an actual global section in $\glsecF$, e.g.~so it is coCartesian and not just lax, cf.~\cref{def:sections}.

\begin{definition}[Categorical Lusztig cycles]\label{def:categorical_Lusztig_cycles}
Let $\w$ be a Demazure weave with trivalent vertices $p_1,\ldots,p_m$ and fix $i\in[1,m]$. By definition, the categorical Lusztig cycle $\LC_i\in \glsecF$ is the global section of $\F_w$ characterized as follows:

\begin{enumerate}
    \item Its value at the vertex $v_j$ is
$$\LC_i(v_j)\coloneqq\begin{cases}
(v_j;\tau_{\geq 0}\on{Mor}(S_{p_j},L_{\beta_{j+1},{\bf a}_{j+1}}),L_{\beta_{j+1},{\bf a}_{j+1}},\on{ev})\, & \text{if } j\not= i, \\
(v_i;k[-1],0,0), & \text{if }i=j\text{ and }i<m,\\
(v_m;k), & \text{if }i=j=m,\\
\end{cases}$$
where $(\b_j,{\bf a}_j)$ is the weighted braid word associated to the Lusztig weave cycle $\g_i^\w$ and any choice of generic horizontal weave slice lying below $p_j$ (but above $p_{j-1}$).

\item Its value at the edge $e_j$ is $\LC_i(e_j)\coloneqq (e_j;L_{(\b_j,{\bf a}_j)})$.

\item The coCartesian edges $\LC_i(v_j\to e_{j+1})$ 
are determined by the trivial identity 
\[
\F_\w(v_j\to e_{j+1})(\LC_i(v_j))=\varrho_1((v_j;\tau_{\geq 0}\on{Mor}(S_{p_j},L_{\beta_{j+1},{\bf a}_{j+1}}),L_{\beta_{j+1},{\bf a}_{j+1}},\on{ev}))= L_{\beta_{j+1},{\bf a_{j+1}}}= \LC_i(e_{j+1}).
\]

\noindent The coCartesian edges $\LC_i(v_j\to e_{j})$ are determined by the equivalence
\[
\F_\w(v_j\to e_{j})(\LC_i(v_j))=\varrho_2((v_j;\tau_{\geq 0}\on{Mor}(S_{p_j},L_{\beta_{j+1},{\bf a}_{j+1}}),L_{\beta_{j+1},{\bf a}_{j+1}},\on{ev}))=\on{cof}(\on{ev})\simeq L_{\beta_{j},{\bf a}_{j}}=\LC_i(e_j)\,,
\]
arising from the fiber and cofiber sequence \eqref{eq:cofiberseq_trivalent0} and the equivalences \eqref{eq:filteredobjects_hexavalent}.
\end{enumerate}
\qed
\end{definition}

In \cref{def:categorical_Lusztig_cycles}.(2), we use the equivalence in Statement (\ref{eq:filteredobjects_hexavalent}) so as to have a well-defined value at the edges. Indeed, two generic horizontal weave slices for the edge $e_j$ might yield different weighted braid words and different filtrations, and so we must ensure the equivalence of the colimit of these filtered objects. Note that, by construction, each ${\bm L}_i^\w$ in \cref{def:categorical_Lusztig_cycles} is indeed a global section of $\F_\w$, since the generization maps connecting vertices and edges are coCartesian.

\begin{remark}\label{rmk:support_categorical_Lusztig_cycles}
A given Lusztig weave cycle $\g_i^\w$ is supported at and below the vertex $v_i$, and thus the associated weighted braid words $(\beta_{j},{\bf a}_{j})$ have all weights zero if $j>i$. Since, in general, $L_{\beta,{\bf a}}\simeq 0$ if all the weights ${\bf a}$ are zero, it follows that the categorical Lusztig cycle ${\bm L}_i^\w$ is also supported at and below $v_i$, i.e.~${\bm L}_i^\w(v_j)\simeq 0$ and ${\bm L}_i^\w(e_j)\simeq 0$ if $i<j$.\qed
\end{remark}

\noindent For any weave trivalent vertex $p_i$ with $i<m$, we can depict part of the categorical Lusztig cycle $\LC_i$ in \cref{def:categorical_Lusztig_cycles} near the associated vertex $v_i$ of $\Wgraph$ as follows:
\[
\begin{tikzcd}[column sep=-10]
      & {(v_{i+1};0,0,0)} \arrow[rd] \arrow[ld] &             & {(v_i;k[-1],0,0)} \arrow[ld] \arrow[rd] &               & {(v_{i-1};\tau_{\geq 0}\on{Mor}(S_{p_{i-1}},S_{p_i}),S_{p_i},\eta)} \arrow[ld] \arrow[rd] &                                   & \dots \arrow[ld] \\
\dots &                                                    & (e_{i+1};0) &                                                    & (e_i;S_{p_i}) &                                                                                                         & {(e_{i-1};L_{\beta_{i-1},{\bf a}_{i-1}})} &                      
\end{tikzcd}
\]
where, by \cref{rmk:support_categorical_Lusztig_cycles}, all the values to the left of $(v_i;k[-1],0,0)$ must vanish, i.e.~they are of the form $(v_j;0,0,0)$ and $(e_j;0)$.

\begin{remark}\label{rmk:weightedbraidword_from_filtration}
For a Demazure $\w$, the categorical Lusztig cycle $L_i^w$ allows to recover 
the weave Lusztig cycle $\g_i^\w$ for all $i\in[1,m]$. Indeed, let $\on{ev}_{e_j}\colon \mathcal{H}(\rgraph,\mathcal{F}_{\w})\to \mathcal{F}(e_j)\simeq\D$ be the evaluation functor of global sections at the edge $e_j$, $j\in[1,m]$. It follows from \cref{def:categorical_Lusztig_cycles} that $\on{ev}_{e_j}({\bm L}_i^\w)\simeq L_{\beta_j,{\bf a}_j}$. By \Cref{lem:uniqueness_of_weight}, the object $L_{\beta_j,{\bf a}_j}\in \D$ uniquely determines the weight ${\bf a}_j$ for $\beta_j$. Concretely, the grading pieces of the filtered object $L_{\beta_j,{\bf a}_j}$ recover the weighted braid word $(\beta_j,{\bf a}_j)$, and the collection of such $(\beta_j,{\bf a}_j)$ uniquely determine $\g_i^\w$.\qed
\end{remark}

\begin{remark}
The weave standard thimbles define a full exceptional collection in $\glsecF$, see \Cref{prop:weave_full_exceptional_collection} below. The stable $\infty$-category $\glsecF$ thus admits a conical iterated recollement into $m$-many $\D(k)$'s, with the object $k$ in the $i$-th component identifying with $\Delta_i^\w$. The restriction of a global section to the $j$-th component of the iterated recollement can be expressed as applying the functor $\on{Mor}(\Delta_j^\w,\mhyphen)$. The $i$-th restriction of $\LC_i$ is given by $k$. The $j$-th restriction of $\LC_i$ to $\D(k)$, with $j\not=i$, is given by the $k$-module $\tau_{\geq 0}\on{Mor}(S_{p_j},L_{\beta_{j+1},{\bf a}_{j+1}})\in \D(k)$ appearing in \Cref{def:categorical_Lusztig_cycles} in the definition of $\LC_i(v_j)$.
\end{remark}


\subsection{The main results on categorical Lusztig cycles}\label{ssec:results_on_categorical_lusztig_cycles}

Let us collect the main results that we will establish on the categorical Lusztig cycles from \cref{def:categorical_Lusztig_cycles}. To ease notation, given a Demazure weave $\w$ with $m$ trivalent vertices, we write $\{\LC_i\}$ implicitly understanding that the index always ranges in $i\in[1,m]$.\\

{\bf (1)} First, we show that the collection of categorical Lusztig cycles is invariant under weave equivalence. Specifically, let $\w,\w'$ be two Demazure weaves and $\F_\w,\F_{\w'}$ their associated weave schobers. Recall that given a weave equivalence $\w\simeq\w'$, \Cref{thm:weaveequivglsec}.(3) provides an equivalence of $\infty$-categories between the global sections $\F_\w$ and $\F_{\w'}$, cf.~\cref{eq:equiv_globalsections}. The precise invariance reads as follows:

\begin{theorem}[Categorical Lusztig cycles are invariant under weave equivalence]\label{thm:Lusztigcycles_weaveequivalence}
Let $\w,\w'$ be two Demazure weaves, $\w\simeq\w'$ a weave equivalence and $\e:\glsec(\rgraph_{\w},\mathcal{F}_\w)\stackrel{\simeq}{\lr} \glsec(\rgraph_{\w'},{\mathcal{F}_{\w'}})$ the associated equivalence of global sections of their weave schobers $\F_\w$ and $\F_{\w'}$. Consider the trivalent vertices $p_1,\dots,p_m$ of $\w$ and the permutation $\sigma\colon [m]\to [m]$ such that $p_{\sigma(1)}',\dots,p_{\sigma(m)}'$ are the ordered trivalent vertices of $\w'$ after applying the weave equivalence $\w\simeq\w'$. Then
\begin{equation*}\label{eq:equivalence_LC}
\e(\LC_i)\simeq {\bm{L}}_{\sigma(i)}^{\w'},\quad \forall i\in[1,m].
\end{equation*}
\qed
\end{theorem}

\noindent In a nutshell, \cref{thm:Lusztigcycles_weaveequivalence} states that weave equivalences yield the same collections of categorical Lusztig cycles, up to reordering. Conceptually, it makes the collection of categorical Lusztig cycles well-defined as a collection of objects assigned to a weave equivalence class, in the same manner that the isomorphism type of the $\infty$-category $\glsecF$ is well-defined on the equivalence class of $\w$, even if the weave schober $\F_\w$ itself is not. \cref{thm:Lusztigcycles_weaveequivalence} is proven in \cref{ssec:Lusztigcycles_weaveequivalence}.\\

{\bf (2)} Second, we show that the collection of categorical Lusztig cycles $\{\LC_i\}$, $i\in[1,m]$, associated to a weave $\w$ forms a simple-minded collection. Easing the notation by $\mathcal{H}_\w\coloneqq \glsecF$, this means that each categorical Lusztig cycle $\LC_i$ is indecomposable and that they satisfy the following: 
\begin{itemize}
        \item The derived Hom $\on{Mor}_{\mathcal{H}_\w}(\LC_i,\LC_i)$ is coconnective for all $i\in[1,m]$,
        \item The derived Hom $\on{Mor}_{\mathcal{H}_\w}(\LC_i,\LC_j)[1]$ is coconnective, for all $i<j$, $i,j\in[1,m]$,
        \item $\bigoplus_{i=1}^m \LC_i$ generates $\mathcal{H}_\w$ under finite limits, colimits, and direct summands. 
    \end{itemize}
\noindent See \cref{def:recollection_exceptional_SMC_silting}
, or \cite[Section 3.2]{KY14} and references therein, for more details on simple-minded collections. In \cref{ssec:LC_from_thimbles} we prove the following result:

\begin{theorem}[Categorical Lusztig cycles are simple-minded]\label{thm:categorical_Lusztig_cycles_are_SMC}
Let $\w$ be a Demazure weave. Then the categorical Lusztig cycles $\{\LC_i\}$ form a simple-minded collection in $\glsecF$.\qed
\end{theorem}

\noindent A consequence of \cref{thm:categorical_Lusztig_cycles_are_SMC} is that, intuitively, $\{\LC_i\}$ behave as the simple objects in an abelian category. Specifically, \cref{thm:categorical_Lusztig_cycles_are_SMC} is equivalent to the following result:

\begin{corollary}\label{cor:categorical_Lusztig_cycles_are_SMC}
Let $\w$ be a Demazure weave. Then there exists a bounded $t$-structure on $\glsecF$ such that the simple objects in its heart are given by the categorical Lusztig cycles $\{\LC_i\}$.\qed
\end{corollary}

\noindent \cref{cor:categorical_Lusztig_cycles_are_SMC} follows from \cref{thm:categorical_Lusztig_cycles_are_SMC} by \cite[Theorem 6.1]{KY14}, see also \cite[Corollary 9.5]{KN13} or \cite[Section 2]{AlNofayee09_simpletstructures}.\\

{\bf (3)} Third, we show that a mutation of Demazure weaves, from $\w$ to $\w'$, induces a mutation of simple-minded collections, from $\{\LC_i\}$ to $\bm{L}_i^{\w'}$. See \cref{def:forward_mutation_of_simple_minded_collections} and \cref{prop:backwards_mutation_along_rigid_simple} for the definition and explicit description of mutations for simple-minded collections. In particular, weave mutation induces a tilting of the associated hearts. The specific result we prove in \cref{ssec:Lusztigcycles_weavemutation} reads as follows:

\begin{theorem}\label{thm:weavemutation_SMCmutation1}
Let $\w$ be a Demazure weave. Then the following holds:
\begin{enumerate}
\item[$(i)$] Let $\w\to \w'$ be a backward weave mutation, and suppose it exchanges the height of the $i$-th and $(i-1)$-th trivalent vertices of $\w$. Then the simple-minded collection $\{\bm{L}_{j}^{\w'}\}$ is the backward simple-minded mutation of $\{\LC_j\}$ at $\LC_{i}$, i.e.~

\[\bm{L}_{j}^{\w'}\simeq\begin{cases}
\LC_i[-1] & \mbox{ if }i=j,\\
\on{cof}(\on{Ext}^1(\LC_i,\LC_j)\otimes \LC_i[-1] \to \LC_j) & \mbox{ else }.
\end{cases}
\]

\item[$(ii)$]  Let $\w\to \w'$ be a forward weave mutation, and suppose it exchanges the height of the $i$-th and $(i-1)$-th trivalent vertices. Then the simple-minded collection $\{\bm{L}_{j}^{\w'}\}$ is the forward simple-minded mutation of $\{\LC_j\}$ at $\LC_i$, i.e.~

\[\bm{L}_{j}^{\w'}\simeq\begin{cases}
\LC_i[1] & \mbox{ if }i=j,\\
\on{fib}(\LC_j\to \on{Ext}^1(\LC_j,\LC_i)^*\otimes \LC_i[1]) & \mbox{ else }.
\end{cases}
\]

\qed
\end{enumerate}
\end{theorem}

\noindent Theorem \ref{thm:main2} in the introduction is a consequence of \cref{thm:Lusztigcycles_weaveequivalence}, \cref{thm:categorical_Lusztig_cycles_are_SMC} and \cref{thm:weavemutation_SMCmutation1} above. The proof of \cref{thm:Lusztigcycles_weaveequivalence}, \cref{thm:categorical_Lusztig_cycles_are_SMC} and \cref{thm:weavemutation_SMCmutation1} is in itself interesting, as it involves introducing a full exceptional collection $\mathcal{E}_\w$ in the $\infty$-category of global sections $\glsecF$, and developing a recipe to construct the categorical Lusztig cycles $\{\LC_i\}$ from such a full exceptional collection $\mathcal{E}_\w$. The construction of $\mathcal{E}_\w$ is achieved in \Cref{sec:weave_thimbles}, and see also \Cref{prop:simplesfromstandard} and \Cref{thm:categorical_Lusztig_cycles_from_exceptionals} for its relation to $\{\LC_i\}$. Similarly to $\F_\w$, though the full exceptional collection $\mathcal{E}_\w$ is useful for computations, it should not be considered as intrinsically attached to the Demazure weave $\w$ since it is not invariant under weave equivalences.\\

\begin{remark}
Finally, as stated in Theorem \ref{thm:main3}, we also show in \cref{thm:description_of_silting_dual_Lusztig_cycles} that the Koszul dual collection to the categorical Lusztig cycles $\{\LC_i\}$ can be realized by weave cycles once a Ringel-type duality has been applied, cf.~\cref{prop:upwards_downwards_thimbles_duality} and \cref{lem:Ringel_duality_and_silting_objects}. These facts are proven in detail in \cref{sec:silting_dualLusztig} and consequently establish the analogous results for such dual silting collection as we have for the categorical Lusztig cycles $\{\LC_i\}$.\qed
\end{remark}

\section{Filtered objects from weighted braid words}\label{sec:filtrations_from_braid_words}

The goal of this section is to construct filtered objects $\lbw\in\D$ from weighted braid words $\bw$.\\

\noindent The geometric realizations of these filtered objects are key for the definition of the categorical Lusztig cycles in \cref{def:categorical_Lusztig_cycles}, see e.g.~Statement (\ref{eq:categorical_weighted_braid}). That said, the theory we develop in this section is more general and thus we present it in a self-contained manner, entirely within the $\infty$-category $\D$, and detached from the theory of weaves and their Lusztig cycles. We shall only specialize to weighted braid words $\bw$ arising from Lusztig cycles in \cref{sec:weave_propagation_LC}, where Statement (\ref{eq:categorical_weighted_braid}), in the more precise form of \Cref{thm:unique_Lusztig_filtered_object}, is proven using the general techniques developed in the present section.\\

\noindent The input ingredients for this section are:
\begin{enumerate}
    \item The perfect derived $\infty$-category $\D\coloneqq\D^{\on{perf}}(\Pi_2(A_{n}))$ of the additive derived preprojective algebra $\Pi_2(A_{n})$ of type $A_n$. Its simples are denoted by $S_1,\dots,S_n\in \D$, as in \cref{ssec:defining_weave_schobers}.

    \item Given a positive braid word $\beta\coloneqq s_{i_0}s_{i_1}\dots s_{i_{l}}$, we define its associated spherical object $S_\beta$ to be
    \begin{equation}\label{eq:sphericalobj_braid}
        S_{\beta}\coloneqq T^{-1}_{S_{i_{l}}}\cdots T^{-1}_{S_{i_{2}}}T^{-1}_{S_{i_{1}}}(S_{i_0})\in \D
    \end{equation}
    where $T_{S_j}=\on{cofib}(\on{Mor}_\D(S_j,\mhyphen)\otimes S_j\to \on{id}_\D)$ is the twist functor along the 2-spherical object $S_j$. This is in line with \cref{def:sphericalobj}.
\end{enumerate}
\noindent The section is organized as follows:
\begin{itemize}
    \item[$(a)$] \cref{ssec:spherical_objects_braids} establishes homological bounds on $\Mor_\D(S_\a,S_\b)$ for certain types of positive braid words $\a,\b$, and shows that the 2-spherical semi-twists preserve rigidity. Both these technical ingredients are crucial for the existence and uniqueness of the filtered objects $\lbw\in\D$ associated to weighted braid words $\bw$.

    \item[$(b)$] \cref{ssec:matching_2spheres} provides a diagrammatic model, based on the study of matching paths, to describe a subcategory $\sM\sse\D$. \cref{lem:matching_curves_stable_under_semitwist} will establish that $\M$ is stable under a class of semi-twists. These results in \cref{ssec:matching_2spheres} and the subcategory $\sM$ are used to prove the existence of the filtered objects $\lbw\in\D$ associated to certain weighted braid words $\bw$.

    \item[$(c)$] \cref{ssec:filteredobject_from_weightedbraid} proves the uniqueness of the rigid filtered objects $\lbw$ in $\D$, should they exist, and establishes a criterion, based on the diagrammatic model of \cref{ssec:matching_2spheres}, for the existence of such $\lbw$ given a weighted braid word $\bw$.
\end{itemize}

\noindent The main result that we prove in this section, and use subsequently, is the following:

\begin{theorem}[Uniqueness of filtered objects $\lbw$]\label{thm:unique_rigid_resolutions}
Let $\bw$ be a weighted braid word. Suppose that there exists a $(\beta,{\bf a})$-filtered object $\lbw\colon N(\mathbb{Z}_{\geq 0})\to \D$ whose geometric realization is rigid and lies in $\mathcal{M}$. Then any other rigid $(\beta,{\bf a})$-filtered object in $\D$ is equivalent to $\lbw$.
\end{theorem}

\noindent The necessary definitions for all the notions featured in \cref{thm:unique_rigid_resolutions}, including that of a $(\beta,{\bf a})$-filtered object, the subwords $\b_i,\b_j$, the degree 1 crossings, and the subcategory $\sM\sse\D$, will be provided in \cref{ssec:matching_2spheres} and \cref{ssec:filteredobject_from_weightedbraid} below. \cref{thm:unique_rigid_resolutions} is proven in \cref{ssec:filteredobject_from_weightedbraid}.\\

\subsection{Bounds on \texorpdfstring{$\on{Ext}$}{Ext}-groups between braid \texorpdfstring{$2$}{2}-spherical objects \texorpdfstring{$S_\b$}{}}\label{ssec:spherical_objects_braids}

This subsection proves \cref{lem:connectiveHom} and \cref{lem:semi_twists_preserve_rigidity}, establishing the vanishing of certain Ext-groups between certain objects of $\D$.\\

First, the silting object $\Pi_2(A_n)\in \D$ determines a $t$-structure, cf.~\Cref{lem:t_structure_from_silting_object}, which we refer to as the standard $t$-structure. The objects $S_\b$ typically do not belong to its associated heart $\D^\heartsuit \sse\D$. For instance, if $\beta=s_1^k$ then $S_{s_1^k}=S_1[k-1]\not\in\D^\heartsuit$ for any $k\in[2,\infty)$. That said, it is well known that

\begin{equation}\label{eq:reduced_in_heart}
\beta\mbox{ reduced }\Longrightarrow S_\b\in\D^\heartsuit.
\end{equation}

\begin{remark}\label{rmk:2sphericalobjects_tstructures}
\begin{enumerate}[(1)]
    \item Suppose that $\beta=s_{i_l}\dots s_{i_0}$ is reduced. Since the Grothendieck group $K_0(\D^\heartsuit)$ is freely generated by the $n$ simples $S_1,\dots,S_n$, it canonically identifies with the root lattice in type $A_n$. Then the module $S_\b\in \D^\heartsuit$ projects to the positive root $s_{i_l}\dots s_{i_1}(\a_{i_0})$ in $K_0(\D^\heartsuit)$. 
    
    \item By \cite[Lemma 3.3]{HW22}, if an object lies in the homologically non-positive aisle of the standard $t$-structure of $\D$, then so do its images under the action of the negative spherical twists. Therefore the spherical object $S_\beta\in \D$ lies in the homologically negative aisle of the $t$-structure of $\D$.\qed 
\end{enumerate}
\end{remark}

\begin{remark}\label{rmk:reduced_case_sphericaltwists}
For a reduced braid word $\beta$, the repeated (positive) spherical twists $\{S_{\alpha}\}_{\alpha<\beta}\subset \D^\heartsuit$ are also known as the layer modules, and arise as the graded pieces of a filtration of the corresponding preprojective algebra, see \cite[Prop.~2.2]{AIRT12} and \cite[Rem.~4.4]{Can25}. We also note that negative spherical twists correspond to positive spherical twists in the opposite category, and that since $\D\simeq \D^{\on{op}}$, the usage of negative or positive spherical twists can be considered a convention.
\end{remark}

\noindent Given two positive braid words $\alpha,\beta$, let us denote by $\alpha<\beta$ if $\alpha$ is a strict suffix of $\beta$, i.e.~if there exist $i_0,\dots,i_m$ such that $\beta=s_{i_0}s_{i_1}\dots s_{i_m}\alpha$ and $\alpha\neq\beta$. Note that for reduced permutation braids, the corresponding partial order is the Bruhat order. The following vanishing result for the morphism objects (i.e.~derived Homs) between such two $2$-spherical objects will be crucial.

\begin{lemma}\label{lem:connectiveHom}
Let $\alpha,\beta$ be two braid words such that $\alpha<\beta$. Then:
\begin{itemize}
    \item[$(i)$] The morphism object $\on{Mor}_\D(S_{\alpha}[-1],S_{\beta})\in \D(k)$ is connective, i.e.~
\[\on{Ext}^i_\D(S_{\alpha},S_\beta)\simeq 0,\quad i\in[2,\infty)\,.\]

    \item[$(ii)$] The morphism object $\on{Mor}_\D(S_{\beta}[-1],S_{\alpha})\in \D(k)$ is coconnective, i.e.~
    \[ \on{Ext}_\D^i(S_{\beta},S_{\alpha})\simeq 0,\quad i\in(-\infty,0]\,. \] 
\end{itemize}
\end{lemma}

\begin{proof} \noindent For Part (1), after possibly acting with a self-equivalence of $\D$, we can and do assume that $\alpha=s_{i_l}$ and $\beta=s_{i_0}s_{i_1}\dots s_{i_l}$. In the case that $\beta$ is reduced, the statement follows from Statement (\ref{eq:reduced_in_heart}). Indeed, there are no extensions higher than degree $2$ between the spherical objects in the standard heart $\D^\heartsuit$ and a simple $S_{i_l}$. Thus $\on{Ext}^i_\D(S_{i_l},S_\beta)\simeq 0$ for $i\in[2,\infty)$ since $S_\beta\in\D^\heartsuit$ in this reduced case. Let us now prove the general non-reduced case by an induction on the length $l$.\\

\noindent Suppose $\beta$ is not reduced, and thus we can choose an expression for $\beta$ such that $s_{i_k}=s_{i_{k+1}}$ for some index $k\in [0,l-1]$. Let us denote
\[ \beta^k\coloneqq s_{i_0}\cdots s_{i_k}\,,\quad \beta^{k+1}\coloneqq s_{i_0}\cdots s_{i_k}s_{i_{k+1}}\,,\quad \beta^\dagger\coloneqq s_{i_{k+2}}\cdots s_{i_l}\,,
\]
so that $\beta=\beta^{k+1}\beta^\dagger$. The induction hypothesis implies that the statement holds for $\beta^k\beta^\dagger$, which has length strictly less than $\beta$. That is, the induction hypothesis implies that $\on{Mor}_\D(S_{\alpha}[-1],S_{\beta^k\beta^\dagger})$ is connective, since $\alpha=s_{i_l}$ is still a suffix for $\beta^k\beta^\dagger$. The goal is to show that $\on{Mor}_\D(S_{i_l}[-1],S_\beta)$ is connective.\\

\noindent For the induction step, we use that $S_{\beta^{k+1}}\simeq T_{S_{i_{k+1}}}^{-1}(S_{\b^k})$ and the corresponding cofiber sequence 
\begin{equation}\label{eq:cofibseq_negtwist}
S_{\beta^{k+1}}\longrightarrow S_{\beta^k}\longrightarrow 
\on{Mor}_\D(S_{\beta^k},S_{i_{k+1}})^*\otimes S_{i_{k+1}}
\end{equation}
induced by the negative twist of $S_{\beta^k}$ along $S_{i_{k+1}}$. By mapping from 
\[S^{\beta_\dagger}[-2]\coloneqq T_{S_{i_{k+2}}}\cdots T_{S_{i_{l-1}}}(S_{i_{l}})[-2]\simeq T_{S_{i_{k+2}}}\cdots T_{S_{i_{l-1}}}T_{S_{i_{l}}}(S_{i_{l}})[-1]\] 
and considering $\on{Ext}$-groups, 
the sequence (\ref{eq:cofibseq_negtwist}) gives rise to a long exact sequence

\begin{equation*}\label{eq:cofibseq_negtwist2}
\begin{tikzcd}[row sep=2em, column sep=2.5em]
    & \cdots \arrow[r] 
    & \textcolor{purple}{\on{Ext}^{j-1}(S^{\beta_\dagger}[-2],\on{Mor}_\D(S_{\beta^k},S_{i_{k+1}})^*\otimes S_{i_{k+1}})} \arrow[dll] \\
    \on{Ext}^{j}(S^{\beta_\dagger}[-2],S_{\beta^{k+1}}) \arrow[r] 
    & \on{Ext}^{j}(S^{\beta_\dagger}[-2],S_{\beta^{k}}) \arrow[r] 
    & \textcolor{purple}{\on{Ext}^{j}(S^{\beta_\dagger}[-2],\on{Mor}_\D(S_{\beta^k},S_{i_{k+1}})^*\otimes S_{i_{k+1}})}
\end{tikzcd}
\end{equation*}

\noindent for $j\in(0,\infty)$. Since we have the equivalences
\[\on{Mor}_\D(S^{\beta_\dagger}[-2],S_{\beta^{k+1}})\simeq \on{Mor}_\D(S_{i_l}[-1],S_\beta)\quad\mbox{and}\quad \on{Mor}_\D(S^{\beta_\dagger}[-2],S_{\beta^{k}})\simeq \on{Mor}_\D(S_{i_l}[-1],S_{\beta^k\beta^\dagger})\,,\]
and $\on{Mor}_\D(S_{i_l}[-1],S_{\beta^k\beta^\dagger})$ is connective by the inductive hypothesis, it suffices to show that the two terms highlighted in purple vanish. We claim this is the case, and prove it as follows. On the one hand, since $S_{i_k}\simeq S_{i_{k+1}}$, as $i_k=i_{k+1}$, we obtain that
\[
\on{Mor}_\D(S_{\beta^k},S_{i_{k+1}})^*\simeq \on{Mor}_\D(S_{i_{k+1}},S_{\beta^k})[2]\simeq \on{Mor}_\D(S_{i_{k}}[-1],S_{\beta^{k}})[1]
\]
must be concentrated in homological degrees in $[1,\infty)$, by applying the induction hypothesis to $\beta^k$. On the other hand, we also have that
\[ \on{Mor}_\D(S^{\beta_\dagger}[-2],S_{i_{k+1}})\simeq \on{Mor}_\D(S_{i_l}[-1],S_{s_{i_{k+1}}\beta^\dagger})\]
is connective, by the inductive hypothesis applied to $s_{i_{k+1}}\beta^\dagger$. Therefore, the two purple outer entries of the above long exact sequence vanish, as claimed. Therefore, the long exact sequence above and the inductive hypothesis imply that $\on{Mor}_\D(S_{i_l}[-1],S_\beta)$ is connective, as required. This establishes Part (1). Part (2) follows from Part (1) because $\D$ is a 2-Calabi--Yau category.
\end{proof}

\begin{example} (1) Let us consider $\alpha=s_1$ and $\beta=s_is_1$, for some $i\in[1,n]$. Then
\[
S_{s_is_1}\simeq\begin{cases}
S_1[1], & \mbox{for }i=1,\\
\fib(S_2\to S_1[1]), & \mbox{for }i=2,\\
S_i, & \mbox{for }i\in[3,n],
\end{cases}
\mbox{ and thus }
\on{Mor}_\D(S_1[-1],S_{s_is_1})\simeq\begin{cases}
k\oplus k[2], & \mbox{for }i=1,\\
k[1], & \mbox{for }i=2,\\
0, & \mbox{for }i\in[3,n].
\end{cases}
\]
Therefore $\on{Mor}_\D(S_1[-1],S_{s_is_1})$ are indeed connective, in accordance with \cref{lem:connectiveHom}.(1).\\

\noindent (2) Let $\alpha=s_1$ and $\beta=(s_1s_2s_1s_1s_2s_2)s_1$. Then $S_\beta\simeq S_2$ and $\on{Mor}_\D(S_1[-1],S_{\beta})\simeq k$. More generally, there exist braid words $\beta$ of arbitrarily large length such that $\on{Mor}_\D(S_1[-1],S_{\beta})\simeq k$.\qed
\end{example}


\subsubsection{$2$-spherical semi-twists preserve rigidity}

Let us now show that semi-twists in $\D$ preserve rigidity. Semi-twists are defined in \Cref{def:semitwists}. For instance, the negative inverse semi-twist of $X$ by $Y$ is
    \begin{equation}\label{def:negative_semitwist}
    T^{-,\leq 0}_{Y}(X)\coloneqq \on{fib}(X\to \tau_{\leq 0}(\on{Mor}_\C(X,Y)^*)\otimes Y)\,,
    \end{equation}
    where $\tau_{\leq 0}(\on{Mor}_\C(X,Y)^*)\subset \on{Mor}_\C(X,Y)^*$ denotes the maximal summand concentrated in homologically negative degrees. An important homological property of semi-twists is as follows:

\begin{lemma}[semi-twists preserve rigidity]\label{lem:semi_twists_preserve_rigidity}
Let $\D$ be a $k$-linear $2$-Calabi--Yau\footnote{Either weak left or triangulated $2$-Calabi--Yau would be sufficient.} stable $\infty$-category. Let $S\in \D$ be a $2$-spherical object. Consider further two compact objects $A,B\in \D$ with $\on{Ext}^1_\D(A,B)\simeq 0$. Then
\begin{align}
\on{Ext}^1_\D(T_S^{-,\leq 0}(A),T_S^{-,\leq 0}(B))&\simeq 0 \label{eq:Ext_1_twist_1}\\
\on{Ext}^1_\D(T_S^{+,\geq 0}(A),T_S^{+,\geq 0}(B))&\simeq 0\label{eq:Ext_1_twist_2}\\
\on{Ext}^1_\D(T_S^{-,> 0}(A),T_S^{-,>0}(B))&\simeq 0\label{eq:Ext_1_twist_3}\\
\on{Ext}^1_\D(T_S^{+,< 0}(A),T_S^{+,<0}(B))&\simeq 0\label{eq:Ext_1_twist_4}\,.
\end{align}
In particular, the functors $T_S^{-,\leq 0}(\mhyphen),T_S^{+,\geq 0}(\mhyphen),T_S^{-,>0}(\mhyphen),T_S^{+,< 0}(\mhyphen)$ preserve rigid objects.
\end{lemma}

\begin{proof}
The equivalence \eqref{eq:Ext_1_twist_2} follows from the equivalence \eqref{eq:Ext_1_twist_1} in $\D^{\on{op}}$ and, similarly, \eqref{eq:Ext_1_twist_4} follows from \eqref{eq:Ext_1_twist_3}. It thus suffices to prove \eqref{eq:Ext_1_twist_1} and \eqref{eq:Ext_1_twist_3}.\\

Let us start to prove the equivalence \eqref{eq:Ext_1_twist_1}. To ease notation, we denote morphism objects in $\D$ by $\on{Mor}(\mhyphen,\mhyphen)$ instead of $\on{Mor}_\D(\mhyphen,\mhyphen)$. By the definition of $T^{-,\leq 0}$ in (\ref{def:negative_semitwist}), the morphism object $\on{Mor}(T_S^{-,\leq 0}(A),T_S^{-,\leq 0}(B))$ is obtained as the totalization  of the following square:
\[
\begin{tikzcd}
{\on{Mor}(A,\tau_{\leq 0}\on{Mor}(B,S)^*\otimes S)}                                               & {\on{Mor}(A,B)} \arrow[l]                                               \\
{\on{Mor}(\tau_{\leq 0}\on{Mor}(A,S)^*\otimes S,\tau_{\leq 0}\on{Mor}(B,S)^*\otimes S)} \arrow[u] & {\on{Mor}(\tau_{\leq 0}\on{Mor}(A,S)^*\otimes S,B)} \arrow[u] \arrow[l]
\end{tikzcd}
\]
The grading in the totalization is chosen, so that $\on{Mor}(T_S^{-,\leq 0}(A),T_S^{-,\leq 0}(B))$ is obtained as the cofiber in the vertical direction and fiber in the horizontal direction. The above square is equivalent to the following square: 
\begin{equation}\label{eq:diagram_semitwist1}
\begin{tikzcd}
{\on{Mor}(A,S)\otimes \tau_{\leq 0}\on{Mor}(B,S)^*}                                               & {\on{Mor}(A,B)} \arrow[l]                                               \\
{(k\oplus k[-2])\otimes \tau_{\geq 0}\on{Mor}(A,S) \otimes \tau_{\leq 0}\on{Mor}(B,S)^*} \arrow[u] & {\tau_{\geq 0}\on{Mor}(A,S)\otimes \on{Mor}(S,B)} \arrow[u] \arrow[l,"f"']
\end{tikzcd}
\end{equation}
where we have used that $S\in\D$ is a 2-spherical object, and thus $\on{Mor}(S,S)\simeq k\oplus k[-2]$, for the bottom left entry of (\ref{eq:diagram_semitwist1}). The $2$-Calabi--Yau structure of $\D$ provides an equivalence $\varepsilon:\on{Mor}(B,S)^*[-2]\to\on{Mor}(S,B)$, which we can use to change the lower right corner $\on{Mor}(S,B)$ of (\ref{eq:diagram_semitwist1}) to $\on{Mor}(B,S)^*[-2]$. In addition, the functoriality of the Calabi--Yau duality implies that the composition of the inclusion
$$\tau_{\geq 0}\on{Mor}(A,S) \otimes \tau_{\leq 0}\on{Mor}(B,S)^*[-2]\hookrightarrow \tau_{\geq 0}\on{Mor}(A,S)\otimes \on{Mor}(B,S)^*[-2]$$
with the morphism
$$\tau_{\geq 0}\on{Mor}(A,S)\otimes \on{Mor}(B,S)^*[-2]\stackrel{f\circ(\on{id}\otimes\varepsilon)}{\lr} (k\oplus k[-2])\otimes \tau_{\geq 0}\on{Mor}(A,S) \otimes \tau_{\leq 0}\on{Mor}(B,S)^*$$

\noindent restricts to the second summand to an equivalence. Now, by moving the upper right term of (\ref{eq:diagram_semitwist1}) to the lower left, the totalization of (\ref{eq:diagram_semitwist1}) is equivalent to the totalization of the following square:
\small
\begin{equation}\label{eq:Mor_square}
\begin{tikzcd}
{\on{Mor}(A,S)\otimes \tau_{\leq 0}\on{Mor}(B,S)^*}                                               & {0} \arrow[l]                                               \\
\on{Mor}(A,B)\oplus {(k\oplus k[-2])\otimes \tau_{\geq 0}\on{Mor}(A,S) \otimes \tau_{\leq 0}\on{Mor}(B,S)^*} \arrow[u] & {\tau_{\geq 0}\on{Mor}(A,S)\otimes \on{Mor}(B,S)^*[-2]} \arrow[u] \arrow[l]
\end{tikzcd}
\end{equation}
\normalsize
In order to compute the totalization of (\ref{eq:Mor_square}) we consider the following two trivial squares
\begin{equation}\label{eq:trivial_diagram1}
\begin{tikzcd}
{\tau_{\geq 0}\on{Mor}(A,S)\otimes \tau_{\leq 0}\on{Mor}(B,S)^*}           & 0 \arrow[l]           \\
{\tau_{\geq 0}\on{Mor}(A,S)\otimes \tau_{\leq 0}\on{Mor}(B,S)^*} \arrow[u,"\on{id}"] & 0 \arrow[l] \arrow[u]
\end{tikzcd}
\end{equation}
and
\begin{equation}\label{eq:trivial_diagram2}
\begin{tikzcd}
0                                                                                & 0 \arrow[l]                                                                                           \\
{ \tau_{\geq 0}\on{Mor}(A,S) \otimes \tau_{\leq 0}\on{Mor}(B,S)^*[-2]} \arrow[u] & { \tau_{\geq 0}\on{Mor}(A,S) \otimes \tau_{\leq 0}\on{Mor}(B,S)^*[-2]} \arrow[u] \arrow[l, "\on{id}"']
\end{tikzcd}
\end{equation}
\noindent The totalizations of (\ref{eq:trivial_diagram1}) and (\ref{eq:trivial_diagram2}) vanish, as the identity map is an isomorphism. Since every object in $\D(k)$ is equivalent to its homology, we can decompose $\on{Mor}(A,S)\simeq \tau_{\geq 0}\on{Mor}(A,S)\oplus \tau_{\leq 0}\on{Mor}(A,S)$. Therefore we have a morphism of squares from (\ref{eq:trivial_diagram1}) to \eqref{eq:Mor_square}. The cofiber of such a morphism is the square
\begin{equation}\label{eq:cofiber_square_morphism}
\begin{tikzcd}
{\tau_{<0}\on{Mor}(A,S)\otimes \tau_{\leq 0}\on{Mor}(B,S)^*}                                         & 0 \arrow[l]                                                                 \\
{\on{Mor}(A,B)\oplus  \tau_{\geq 0}\on{Mor}(A,S) \otimes \tau_{\leq 0}\on{Mor}(B,S)^*[-2]} \arrow[u] & {\tau_{\geq 0}\on{Mor}(A,S)\otimes \on{Mor}(B,S)^*[-2]} \arrow[u] \arrow[l]
\end{tikzcd}
\end{equation}
Since the totalization of \eqref{eq:trivial_diagram1} vanishes, the totalization of \eqref{eq:cofiber_square_morphism} is equivalent to the totalization of \eqref{eq:Mor_square}. Similarly, there is a morphism from \eqref{eq:cofiber_square_morphism} to \eqref{eq:trivial_diagram2}, whose fiber is

\begin{equation}\label{eq:fiber_square_morphism}
\begin{tikzcd}
{\tau_{<0}\on{Mor}(A,S)\otimes \tau_{\leq 0}\on{Mor}(B,S)^*}                                         & 0 \arrow[l]                                                                 \\
{\on{Mor}(A,B)} \arrow[u] & {\tau_{\geq 0}\on{Mor}(A,S)\otimes \tau_{>0}\on{Mor}(B,S)^*[-2]} \arrow[u] \arrow[l]
\end{tikzcd}
\end{equation}

\noindent As before, the totalization of \eqref{eq:cofiber_square_morphism} is equivalent to the totalization of \eqref{eq:fiber_square_morphism}, thus we focus on the latter. Passing to the totalization of \eqref{eq:fiber_square_morphism}, the term $\tau_{<0}\on{Mor}(A,S)\otimes \tau_{\leq 0}\on{Mor}(B,S)^*$ only contributes in degrees $i\leq -2$, and t$\tau_{\geq 0}\on{Mor}(A,S)\otimes \tau_{>0}\on{Mor}(B,S)^*[-2]$ only contributes in degrees $i\geq 0$. Hence, $\on{Ext}^1(T_S^{-,\leq 0}(A),T_S^{-,\leq 0}(B))$ vanishes, as required. This establishes the vanishing \eqref{eq:Ext_1_twist_1}.\\

\noindent For the equivalence \eqref{eq:Ext_1_twist_3}, the argument is similar to the above. Namely, the morphism object $\on{Mor}(T_S^{-,>0}(A),T_S^{-,>0}(B))$ is equivalent to the totalization of the following square in $\D(k)$:
\[
\begin{tikzcd}
{\on{Mor}(A,S)\otimes \tau_{>0}\on{Mor}(B,S)^*}                                               & {0} \arrow[l]                                               \\
\on{Mor}(A,B)\oplus {(k\oplus k[-2])\otimes \tau_{<0}\on{Mor}(A,S) \otimes \tau_{>0}\on{Mor}(B,S)^*} \arrow[u] & {\tau_{<0}\on{Mor}(A,S)\otimes \on{Mor}(B,S)^*[-2]} \arrow[u] \arrow[l]
\end{tikzcd}
\]
Its totalization is in turn equivalent to the totalization of the square
\[
\begin{tikzcd}
{\tau_{\geq 0} \on{Mor}(A,S)\otimes \tau_{>0}\on{Mor}(B,S)^*}                                               & {0} \arrow[l]                                               \\
\on{Mor}(A,B) \arrow[u] & {\tau_{<0}\on{Mor}(A,S)\otimes \tau_{\leq 0}\on{Mor}(B,S)^*[-2]} \arrow[u] \arrow[l]
\end{tikzcd}
\]
which vanishes in degree $-1$, thus establishing the equivalence \eqref{eq:Ext_1_twist_3}. Finally, the fact that each of the four functors $T_S^{-,\leq 0}(\mhyphen),T_S^{+,\geq 0}(\mhyphen),T_S^{-,>0}(\mhyphen),T_S^{+,< 0}(\mhyphen)$ preserves rigid objects follows from the equivalences in the statement by taking $A=B$.
\end{proof}

\begin{remark}
A result in line with \Cref{lem:semi_twists_preserve_rigidity}, in the case of objects in the module $1$-category of the preprojective algebra, appears in \cite[Prop.~5.2]{GLS08}. A further similar result appears in \cite[Lem.~4.5]{PYK23}.\qed
\end{remark}


\subsection{A geometric model for $\D$}\label{ssec:matching_2spheres}

In this subsection we explain how to describe certain objects in $\D$ and the morphisms between them in terms of immersed paths in the plane. Such objects and their geometric description will be used in \cref{ssec:filteredobject_from_weightedbraid} to construct the filtered objects $\lbw$. In brief, this subsection develops the classical theory of matching cycles in the context of a perverse schober whose $\infty$-category of global sections is equivalent to $\D$.\\

\begin{center}
	\begin{figure}[h!]
		\centering
		\includegraphics[scale=1.2]{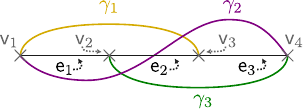}
		\caption{The graph $\rgraph_\D$ in the case of $n=3$, with its four vertices $\vd_1,\vd_2,\vd_3,\vd_4$ (gray) and its three edges $\ed_1,\ed_2,\ed_3$ (black). There are three $\pi_\D$-matching paths $\g_1,\g_2,\g_3$, depicted in yellow, purple and gray respectively.}
        \label{fig:MatchingPaths_GraphD}
	\end{figure}
\end{center}

Let $\rgraph_\D$ be the linear ribbon graph with $n+1$ vertices, i.e.~an $A_{n+1}$-Dynkin diagram. See \cref{fig:MatchingPaths_GraphD} for the case $n=3$. Consider the spherical functor $\phi \colon \D^{\on{perf}}(k)\lr\D^{\on{perf}}(k[t^{\pm1}])$, with $k[t^{\pm1}]$ the ungraded $k$-algebra of Laurent polynomials, determined by mapping $k$ to the simple representation with $k$ as its underlying vector space and $t$ acting by $1$. Given such a functor $\phi$ we define, using \Cref{not:schobers}, the following $k$-linear $\rgraph_\D$-parametrized perverse schober:

\begin{equation}\label{eq:schoberD}
\begin{tikzcd}
\F_\D\coloneqq & \phi \arrow[r, "{(\varrho_1,\varrho_1)}", no head] & \phi \arrow[r, "{(\varrho_2,\varrho_1)}", no head] & \dots \arrow[r, "{(\varrho_2,\varrho_1)}", no head] & \phi
\end{tikzcd}
\end{equation}

\noindent It follows from \cite[Thm.~4.19]{Chr21b} that there exists an equivalence of $k$-linear $\infty$-categories $\D\simeq \Glsec{\rgraph_\D}{\mathcal{F}_\D}$, i.e.~$\D$ is equivalent to the category of global sections of the schober $\F_\D$ in \eqref{eq:schoberD}.

\begin{remark}
Specifically, this partial geometric model that we will describe captures, in the context of perverse schobers, aspects of \cite[Section (16g)]{Sei08}, and some of the algebraic structures associated to the symplectic geometry of Lagrangian 2-spheres. See also \cite[Sections 3\&4]{KS02}, \cite[Section 5]{Chr21b} and \cite[Section 4]{CHQ23}. \qed
\end{remark}


\subsubsection{Matching objects in $\D$}\label{sssec:2spherical_matching_objects}

Let $\mathbb{D}\subset \mathbb{R}^2$ the open disk with center $(\frac{n+3}{2},0)$ and diameter $n+2$. Consider the points $\vd_i=(i,0)\in \mathbb{D}$, $i\in[1,n+1]$, and the punctured disk $\tilde{\mathbb{D}}\coloneqq\mathbb{D}\backslash \{\vd_1,\dots,\vd_{n+1}\}$. By definition, a grading structure $\nu$ on $\mathbb{D}$ is a homotopy class of sections of the projectivised tangent bundle $\mathbb{P}(T\tilde{\mathbb{D}})$; a section of $\mathbb{P}(T\tilde{\mathbb{D}})$ is also known as a line field. We equip the disk $\mathbb{D}$ above with the unique grading structure satisfying that the winding number around each puncture is given by $2$. The corresponding line field $\nu$ can then be chosen to be constant, pointing everywhere in the horizontal direction. In particular, $\nu$ extends to all of $\mathbb{D}$.\\

Let us fix an embedding of the geometric realization of $\rgraph_\D$ into $\mathbb{D}$ such that the $i$-th vertex $\vd_i$ of $\rgraph_\D$ is mapped to $\vd_i\in \mathbb{D}$. To ease notation, we do not distinguish the vertices of $\rgraph_\D$ from their images in $\mathbb{D}$ and equally refer to them as $\vd_i$. Note that in this model we are depicting $\rgraph_\D$ and the schober $\F_\D$ horizontally, cf.~\cref{fig:MatchingPaths_GraphD} and \cref{eq:schoberD}: this is done intentionally in order to contrast with $\Wgraph$ and $\F_\w$, which we always draw vertically, cf.~\cref{ssec:defining_weave_schobers}.\\

\begin{definition}\label{def:pi_D_matching_path}
By definition, an immersed curve $\gamma\colon [-1,1]\lr \mathbb{D}$ is a $\pi_\D$-matching path if 
\begin{itemize}
    \item The image $\gamma((-1,1))$ is disjoint from the vertices of $\rgraph_\D$, and its endpoints $\gamma(-1)$ and $\gamma(1)$ must be vertices of $\rgraph_\D$. 
    \item The complement $\mathbb{D}\setminus \gamma(-1,1)$ does not contain any once-punctured disks without marked points on the boundary, and whose boundary entirely consists of points in the image of $\gamma$.\footnote{In particular, if $\gamma$ wraps around a vertex, it does so at most by one full turn and then necessarily returns from the same direction that it approached the vertex.} 
    \item $\gamma$ is not homotopic relative the vertices of $\rgraph_\D$ to a constant curve.\qed
\end{itemize}
\end{definition}
\noindent The $\pi_\D$-matching paths in \cref{def:pi_D_matching_path} will be always considered up to homotopy relative to the vertices of $\rgraph_\D$. \cref{fig:MatchingPaths_GraphD} depicts three embedded $\pi_\D$-matchings paths in the case $n=3$.
\begin{definition}\label{def:grading_matchingpath}
    A grading of a $\pi_\D$-matching path $\gamma\colon [-1,1]\to \mathbb{D}$ consists of a homotopy class of paths in $\Gamma([-1,1],\gamma^*\mathbb{P}(T\tilde{\mathbb{D}}))$ from the line field $\gamma^*\nu$ to the line field of tangents $\dot{\gamma}$.\qed
\end{definition}

\begin{remark}\label{rem:canonical_grading}
Gradings of a $\pi_\D$-matching path $\gamma$, as in \cref{def:grading_matchingpath}, form a $\mathbb{Z}$-torsor. We can and typically do choose the line field $\nu$ that is parallel to the edges of $\rgraph_\D$. Note that choosing $\gamma(0)$ to lie on an edge $\ed$ and $\dot{\gamma}(0)$ to be parallel to $\ed$ determines a grading, i.e.~a base point for the $\mathbb{Z}$-torsor.\qed
\end{remark}

\begin{center}
	\begin{figure}[h!]
		\centering
		\includegraphics[scale=1.3]{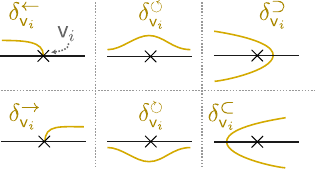}
		\caption{The six different types of local pieces for a $\pi_\D$-matching path, as defined in \cref{def:pi_D_matching_path}.}
        \label{fig:WeaveVertex_Segments2}
	\end{figure}
\end{center}

Each $\pi_\D$-matching path from \cref{def:pi_D_matching_path} can be uniquely written as the composite of six different homotopy classes of segments, up to reversal of orientation and smooth isotopy relative to its endpoints, with each segment lying in a contractible neighborhood of a vertex $\vd_i\in\rgraph_\D$. These six types of segments are depicted in \Cref{fig:WeaveVertex_Segments2}, paired in three different columns. The segments depicted in the $k$th column of \Cref{fig:WeaveVertex_Segments2} are said to be segments of the $k$th type. For instance, the three $\pi_\D$-matchings paths depicted in \cref{fig:MatchingPaths_GraphD} only use segments of the first and second type.

\noindent In our context, the purpose of the $\pi_\D$-matching paths from \cref{def:pi_D_matching_path} is that we can associate a global section $X_\g$ of $\mathcal{F}_\D$, i.e.~an object $X_\g\in\D$, to each such graded $\pi_\D$-matching path $\g$. This can be seen as an instance of the constructions given in \cite[Constr.~5.13]{Chr21b} and \cite[Section 4.3]{CHQ23}.\footnote{Note that $\mathcal{F}_\D$ admits a canonical arc system kit in the sense of \cite{CHQ23}.} The construction can be described as follows.\\

\begin{construction}[Matching objects]\label{construction:matching_objects}

Let us label the vertices of $\rgraph_\D$ by $\vd_1,\dots,\vd_{n+1}$ left-to-right and the incident edges by $\ed_1,\dots,\ed_{n}$, also left-to-right. Fix a $\pi_\D$-matching path $\g$: the object $X_\gamma\in\D$ will be obtained by gluing local sections associated to the different types of segments for $\gamma$, as in \cref{fig:WeaveVertex_Segments2}.\\

{\bf (\ref{construction:matching_objects}.1).} First, let us choose the grading $\nu$ for $\gamma$ as in \Cref{rem:canonical_grading}: such grading of $\gamma$ determines an integer $d$. From this, we associate a degree $d(\delta)\in\Z_{\leq0}$ to every type of segment $\delta$ for $\gamma$, as follows. Denote the segments of $\gamma$ by $\delta_{-a},\dots,\delta_0,\dots,\delta_b$, with $a,b>0$, where $\delta_0$ the segment for the starting point of $\gamma$ and such that they are ordered compatibly with the way they are horizontally concatenated to form $\gamma$. In this case, we set $d(\delta_0)\coloneqq d$ and for the next segment in $\delta_0,\delta_{-1},\dots,\delta_{-a}$ (resp.~$\delta_0,\delta_1,\dots,\delta_b$), we recursively increase (resp.~decrease) the degree by $1$ if the previous segment was of the third type going clockwise (resp.~counterclockwise), and otherwise the degree remains the same.\\

{\bf (\ref{construction:matching_objects}.2).} Second, we associate to each type of a segment in \cref{fig:WeaveVertex_Segments2} a local section. Specifically, with $\rgraph_\D$ and $\F_\D$ as above, we denote by $p\colon \Gamma(\mathcal{F}_D)\to \on{Exit}(\rgraph_\D)$ the Grothendieck construction of $\F_\D$. Objects in $\Gamma(\mathcal{F}_\D)$ can be described as pairs $(\vd;O_\vd)\in\Gamma(\mathcal{F}_\D)$, where $\vd\in \on{Exit}(\rgraph_{\D})$ and $O_\vd\in\F_\D(\vd)$ are objects. By \eqref{eq:schoberD} and \cref{sssec:model_perverseschobers} we have

\begin{equation*}
\F_\D(\vd_i)\simeq
\begin{cases}
\D^{\on{perf}}(k)\overset{\rightarrow}{\times}_\phi \D^{\on{perf}}(k[t^{\pm1}]) & \mbox{ if }i\in[2,n],\\
\D^{\on{perf}}(k) \mbox{ if }i\in\{1,n+1\}.
\end{cases}
\end{equation*}
\noindent Let us consider the following partial sections $M_{\delta_{\vd_i}^{\dagger}}:\Delta^0=\{\vd_i\}\lr\Gamma(\F_\D)$ of $p$, where the superscript $\dagger$ denotes $\dagger\in\{\leftarrow,\rightarrow,\circlearrowleft,\circlearrowright,\supset,\subset\}$, defined as follows:
\begin{align*}
&M_{\delta_{\vd_i}^{\circlearrowleft}}({\vd_i})\coloneqq (0,\phi(k),0),\qquad\qquad\qquad\qquad\quad
M_{\delta_{\vd_i}^{\circlearrowright}}({\vd_i})\coloneqq (k\oplus k[-1],\phi(k),\phi\phi^R\phi(k)\xrightarrow{\on{counit}}\phi(k))\,,\\
& \\
&M_{\delta_{\vd_i}^{\leftarrow}}({\vd_i})\coloneqq \begin{cases} 
    k & \mbox{ if }i=1 \\ 
    (k;\phi(k),\on{id}_{\phi(k)}) & \mbox{ else}\,
    \end{cases},\quad
M_{\delta_{\vd_i}^{\rightarrow}}({\vd_i})\coloneqq \begin{cases} 
    k & \mbox{ if } i=n+1 \\ 
    (k[-1];0,0) & \mbox{ else}
    \end{cases}\\
& \\
&M_{\delta_{\vd_i}^{\supset}}(\vd_i)\coloneqq \begin{cases}
    k\oplus k[1] & i=1\\
    (k\oplus k[1],0,0)& \mbox{ else}\,
    \end{cases},\qquad
M_{\delta_{\vd_i}}^{\subset}(\vd_i)\coloneqq \begin{cases}
    k\oplus k[1] & i=n+1\\
    (k\oplus k[1],\phi(k)\oplus \phi(k)[1],\eta)& \mbox{ else}\,,
    \end{cases}\\
\end{align*}
where the morphism $\eta$ for $M_{\delta_{\vd_i}}^{\subset}(\vd_i)$ is
$$\eta\coloneqq \begin{pmatrix}\on{id}_{\phi(k)} & 0\\  \alpha & \on{id}_{\phi(k)[1]} \end{pmatrix}$$
with $\alpha\colon \phi(k)\to \phi(k)[1]$ a non-zero morphism. In these descriptions above, for $M_{\delta_{\vd_i}}^{\supset}(\vd_i)$ and $M_{\delta_{\vd_i}}^{\subset}(\vd_i)$ we assume the segment is oriented clockwise. In these cases, if $\delta$ is a segment of the third type with the other orientation, then we set $M_{\delta}\coloneqq M_{\delta^{\on{rev}}}[-1]$ with $\delta^{\on{rev}}$ the segment with the reversed orientation.\\

Let us extend these partial sections $M_{\delta_{\vd_i}^{\dagger}}:\Delta^0\to\Gamma(\F_\D)$ above to local lax sections of $\mbox{Exit}(\rgraph_\D)$ via $p$-relative left Kan extensions, i.e.~fitting the diagram

\begin{equation*}\label{eq:relative_leftKanExtensions}
\begin{tikzcd}
 & \Gamma(\mathcal{F}_\D) \arrow[d,"p"'] \\
\Delta^0 \arrow[ru,"M_{\delta_{\vd_i}^{\dagger}}"] \arrow[r,hook,"\iota"] & \on{Exit}(\rgraph_\D) \arrow[u,bend right,blue,"\mbox{Lan}^p_\iota(M_{\delta_{\vd_i}^{\dagger}})"']
\end{tikzcd}
\end{equation*}

\noindent This leads to the following lax sections associated to segments:

\begin{definition}\label{def:localsection_from_segment}
Consider an ungraded segment $\delta$ of type $\dagger$ near the vertex $\vd_i\in\rgraph_\D$. By definition, the local section $M_{\delta}\in\mathcal{L}(\rgraph_\D,\mathcal{F}_\D)$ associated to $\delta$ is the $p$-relative Kan extension of $M_{\delta^\dagger_{\vd_i}}$.\qed
\end{definition}

{\bf (\ref{construction:matching_objects}.3).} Finally, the object $X_\g$ is obtained by gluing the lax sections of its segments constructed in \cref{def:localsection_from_segment}. Specifically, the gluing of these local sections amounts to repeated pushouts in the $\infty$-category of lax sections $\mathcal{L}(\rgraph_\D,\mathcal{F}_\D)$, as follows. Let $M_{\ed_i}\in \mathcal{L}(\rgraph_\D,\mathcal{F}_\D)$ be the lax section determined by
$$M_{\ed_i}(\ed_i)=\phi(k)\in\F_\D(\ed_i),\qquad M_{\ed_i}(x)\simeq 0\mbox{ for all }x\not=\ed_i,\quad x\in \on{Exit}(\rgraph_\D).$$
Let $\ed_{j_i}$ be the edge where the segment $\delta_i$ of $\gamma$ starts, and consider the associated conical inclusion of lax sections $M_{\ed_{j_i}}[d(\delta_i)]\hookrightarrow M_{\delta_i}[d(\delta_i)]$. By construction, we always glue in a segment $\delta_i$ of $\gamma$ at a time by passing to the pushout along the morphism $M_{\ed_{j_i}}[d(\delta_i)]\hookrightarrow M_{\delta_i}[d(\delta_i)]$. Note that, if $\delta_i$ is a segment of the first or second type, then this inclusion is determined by evaluating to an equivalence at $\ed_{j_i}$. The result of this gluing process is an object $X_\gamma\in\mathcal{L}(\rgraph_\D,\mathcal{F}_\D)$ lying in the full subcategory $\Glsec{\rgraph_\D}{\mathcal{F}_\D}\subset \Losec{\rgraph_\D}{\mathcal{F}_\D}$ of global sections, which we summarize as follows:

\begin{definition}[Matching objects]\label{def:object_from_matching}
Let $\g$ be a graded $\pi_\D$-matching path. By definition, the object $X_\g\in \Glsec{\rgraph_\D}{\mathcal{F}_\D}$ is the result of gluing the local lax sections $M_\delta$ of the segments $\delta$ associated to $\g$ along the conical inclusions $M_{\ed_{\delta}}[d(\delta)]\hookrightarrow M_{\delta}[d(\delta)]$ of lax sections, where $\ed_{\delta}$ is the edge corresponding to $\delta$, as in (\ref{construction:matching_objects}.3) above.\qed
\end{definition}

{\bf \noindent End of Construction \ref{construction:matching_objects}.}
\end{construction}

\noindent We will use the category additively spanned by the matching objects from \cref{def:object_from_matching}:

\begin{definition}\label{def:M_geometric_objects} By definition, the $\infty$-category $\mathcal{M}\coloneqq \on{Add}(\{X_{\gamma}\})\subset \D$ is the additive hull of objects $X_\g\in\D$ associated with graded $\pi_\D$-matching paths. By definition, an object $X\in\D$ is said to be matching if $X\in\mathcal{M}$.\qed
\end{definition}

\begin{example}[Simples $S_i$ are matching]\label{ex:simple_as_matching_sphere}
Each of the 2-spherical simples $S_i\in\D$ is matching, i.e.~$S_i\in \mathcal{M}$. Specifically, $S_i$ arises as the global section $\F_\D$ associated with the $\pi_\D$-matching path $\ell_i$ given by a (straight) line $\ed_i$ connecting the $i$th and $(i+1)$-the vertices of $\rgraph_\D$. In this case, $\gamma$ consists of two segments $\delta_1\coloneqq \delta^{\leftarrow}_{\vd_i}$ and $\delta_2\coloneqq \delta^{\rightarrow}_{\vd_{i+1}}$ of the first type, composed at the edge $\ed_i$. Thus, $X_{\ell_i}$ is the pushout of a diagram of the following form:
\begin{equation}\label{eq:Si_geometric}
\begin{tikzcd}
{M_\ed} \arrow[d, hook] \arrow[r, hook] & {M_{\delta_2}} \\
{M_{\delta_1}}                        &                            
\end{tikzcd}
\end{equation}
where the local lax sections $M_e,M_{\delta_1},M_{\delta_2}$ are as in Construction \ref{construction:matching_objects} above. We can depict these lax sections as follows:
\[
M_{\delta_1}\simeq\begin{tikzcd}[column sep=tiny]
\dots \arrow[rd] &               & {(\vd_i;k,\phi(k),\on{id})} \arrow[rd] \arrow[ld] &                   & {(\vd_{i+1};0,0,0)} \arrow[rd] \arrow[ld] &               & \dots \arrow[ld] \\
                 & (\ed_{i-1};0) &                                                     & (\ed_i;\phi(k)) &                                           & (\ed_{i+1};0) &                 
\end{tikzcd}
\]
\[
M_{\delta_2}\simeq\begin{tikzcd}[column sep=tiny]
\dots \arrow[rd] &               & {(\vd_i;0,0,0)} \arrow[rd] \arrow[ld] &                   & {(\vd_{i+1};k[-1],0,0)} \arrow[rd] \arrow[ld] &               & \dots \arrow[ld] \\
                 & (\ed_{i-1};0) &                                                     & (\ed_i;\phi(k)) &                                           & (\ed_{i+1};0) &                 
\end{tikzcd}
\]
\[
M_{\ed}\simeq\begin{tikzcd}[column sep=tiny]
\dots \arrow[rd] &               & {(\vd_i;0,0,0)} \arrow[rd] \arrow[ld] &                   & {(\vd_{i+1};0,0,0)} \arrow[rd] \arrow[ld] &               & \dots \arrow[ld] \\
                 & (\ed_{i-1};0) &                                                     & (\ed_i;\phi(k)) &                                           & (\ed_{i+1};0) &                 
\end{tikzcd}
\]
The pushout $X_{\ell_i}$ of \eqref{eq:Si_geometric} is the following global, i.e.~coCartesian, section of $p\colon \Gamma(\mathcal{F}_\D)\to \on{Exit}(\rgraph_\D)$:
\[
X_{\gamma}=\begin{tikzcd}[column sep=tiny]
\dots \arrow[rd] &               & {(\vd_i;k,\phi(k),\on{id})} \arrow[rd] \arrow[ld] &                   & {(\vd_{i+1};k[-1],0,0)} \arrow[rd] \arrow[ld] &               & \dots \arrow[ld] \\
                 & (\ed_{i-1};0) &                                                     & (\ed_i;\phi(k)) &                                           & (\ed_{i+1};0) &                 
\end{tikzcd} 
\]
It indeed follows that $X_{\ell_i}\simeq S_i$ and thus $S_i\in\mathcal{M}$.\qed
\end{example}

\begin{remark}\label{rmk:canonical_grading2}
Let $\gamma$ be a $\pi_\D$-matching path without segments of the third type. Then there exists a canonical grading for $\gamma$ such that $X_\gamma\in \D^\heartsuit$ lies in the standard heart $\D^\heartsuit\subset \D$.\qed
\end{remark}

Applying a spherical twist to a matching object $X_\g\in\mathcal{M}$ along any of the simples $S_1,\dots,S_n$ can be described in terms of $\pi_\D$-matching paths, as follows. For that, we consider the half-twists $f_i\colon \mathbb{D}\to \mathbb{D}$, considered up to compactly supported isotopies of $\mathbb{D}$, cf.~\cite[Figure 14.2]{Donaldson_RiemannSurfaces}, \cite[Section 3b]{KS02} or \cite[Section (16d)]{Sei08}. The mapping class of such $f_i\colon \mathbb{D}\to \mathbb{D}$, $i\in[1,n]$, is represented by any compactly supported diffeomorphism that acts on the edge $\ed_i\in\rgraph_\D$ connecting the $i$th and $(i+1)$th vertices $\vd_i,\vd_{i+1}\in\rgraph_\D$ by a counter-clockwise half-rotation (thus exchanging $\vd_i$ and $\vd_{i+1}$) and fixes the other vertices $\vd_j$ of $\rgraph_\D$ for $j\in[1,n+1]$ and $j\not\in\{i,i+1\}$. Given a graded $\pi_\D$-matching path, the grading for the $\pi_\D$-matching path $f_i\circ \gamma$ is obtained by composition of $f_i^*(\dot{\gamma})\simeq (f_i\circ \gamma)^*(\nu)$ with the canonical homotopy between line fields $f_i^*(\nu)\simeq \nu$. Spherical twist along the simples $S_1,\ldots,S_n$ are then described as follows:

\begin{lemma}\label{lem:twist_via_homeomorphism}
Let $\gamma$ be a graded $\pi_\D$-matching path, $X_\g\in\D$ its matching object, and $f_i\in\mbox{Diff}_c(\mathbb{D})$ the $i$th half-twist, $i\in[1,n]$. Then 
$$T_{S_i}(X_{\gamma})\simeq X_{f_i\circ \gamma}.$$
\end{lemma}

\noindent \cref{lem:twist_via_homeomorphism} follows from the same argument used to prove \cite[Thm.~7.5]{Chr21b}, see also \cite[Remark 16.14]{Sei08}. \Cref{lem:twist_via_homeomorphism} implies that the objects associated to embedded $\pi_\D$-matching paths can be constructed out of the $2$-spherical objects $S_1,\dots,S_n$ via iterated spherical twists. In addition, the association of embedded graded $\pi_\D$-matching paths to equivalence classes of global sections of $\mathcal{F}_\D$ is injective, see e.g.~the proof of \cite[Cor.~5.1]{CHQ23}.

\begin{example}\label{ex:matching_paths1}
Consider the three $\pi_\D$-matching paths $\g_1,\g_2,\g_3$ in \cref{fig:MatchingPaths_GraphD}, with $\rgraph_\D$ corresponding to the $A_3$-Dynkin diagram. By \cref{lem:twist_via_homeomorphism} we can describe their associated matching objects as
$$X_{\g_1}\simeq T_{S_1}(S_2),\qquad X_{\g_2}\simeq T^{-1}_{S_3}T_{S_2}(S_1), \qquad X_{\g_3}\simeq T_{S_3}(S_2).$$
\noindent In this instance, we have graded the $\pi_\D$-matching paths via the canonical grading, cf.~\cref{rem:canonical_grading}, and since there are no segments of the third type the resulting matching objects are indeed in $\D^\heartsuit$, as per \cref{rmk:canonical_grading2}.\qed
\end{example}
 
\begin{remark}\label{rem:curves_for_self_extensions}
There can be distinct {\it immersed} graded $\pi_\D$-matching paths which give rise to equivalent global sections. 
For instance, let $\tau_1$ be the immersed $\pi_\D$-matching path starting at $\vd_i$, encircling $\vd_{i+1}$ exactly once, and ending at $\vd_i$ again. Similarly, let $\tau_2$ be the immersed $\pi_\D$-matching path starting at $\vd_{i+1}$, encircling $\vd_{i}$ exactly once, and ending at $\vd_{i+1}$ again. Then both $\pi_\D$-matching paths $\tau_1,\tau_2$, when suitably graded, give rise to equivalent objects $X_{\tau_1}\simeq X_{\tau_2}$ in $\D$, namely the extension $X_{\tau_1}\simeq \on{fib}(S_i\to S_i[2])$. For the repeated fibers $\on{fib}(S_i\to \on{fib}(S_i[2]\to \dots \on{fib}(S_i[2(m-1)]\to S_i[2m])))$, there are similarly two corresponding $\pi_\D$-matching paths.\qed
\end{remark}


\subsubsection{Morphisms between matching objects in $\D$}\label{sssec:morphisms_2spherical_matching_objects}

Let us now provide a description of the morphism objects between the matching objects from \cref{sssec:2spherical_matching_objects} in terms of counts of graded intersections between their corresponding $\pi_\D$-matching paths. We distinguish two types of intersections: endpoint intersections and crossings. Crossings are transverse intersections of two $\pi_\D$-matching paths away from their endpoints. Both types of intersections are considered as directed: the order of the curves at the intersection point is specified.\\

{\bf Index of intersection points}. First, we must associate an index $d(x)\in\mathbb{Z}$ to each type of directed intersection $x$, as follows. For the case of crossings, let $\gamma_1,\gamma_2$ be two $\pi_\D$-matching paths and $c_i\coloneqq  \dot{\gamma_1}\simeq (\gamma_i)^*\nu$ their canonical gradings, $i\in[1,2]$. Suppose that $\g_1,\g_2$ intersect transversely at an interior point $p=\gamma_1(t_1)=\gamma_2(t_2)$, for some $t_1,t_2\in(-1,1)$. Then we define the index $d(x)\in \pi_1(\mathbb P(T_p{\bf S}))\simeq \mathbb Z$ to be the class $c_2(t_2)\cdot\kappa\cdot c_1(t_1)^{-1}$, where $\kappa$ is a path from $\dot{\gamma}_1(t_1)$ to $\dot{\gamma}_2(t_2)$ given by clockwise rotation in $T_p\mathbb{D}$ by an angle less than $\pi$.\footnote{Note that we indeed have $\pi_1(\mathbb P(T_p{\bf S}))\simeq \mathbb Z$, and we choose the isomorphism to be determined by mapping $1$ to a clockwise rotation.} In this case, the crossing is said to be in degree $d(x),d(x)+1$, since these are the degrees of the non-trivial extensions to which it will contribute.\\

\noindent For the case of endpoint intersections, let $\gamma_1,\gamma_2$ be two $\pi_\D$-matching paths with an endpoint intersection $x=\vd_i$ for some $i\in[1,n+1]$. We follow \cite[Eq.~(2.5)]{IQZ20} to define the intersection index of $x$. For that, we choose a small circle around $\vd_i$ and let $\alpha\colon [0,1]\to \mathbb{D}$ be a curve going fully around the circle in the clockwise direction. Let $y_1$ be the crossing from $\gamma_1$ to $\alpha$, and $y_2$ the crossing from $\gamma_2$ to $\alpha$, so that $y_1,y_2$ are crossing intersections. The index of the intersection point $x$ is then defined as $d(x)\coloneqq d(y_1)-d(y_2)$, where we have chosen an arbitrary grading for the circle $\alpha$.

\begin{example}\label{ex:endpointintersection_index}
If $\gamma_1=\ed_1$ and $\gamma_2=\ed_2$ are the first two edges of the geometric realization of $\rgraph_\D$ equipped with their canonical gradings, their endpoint intersection has index $1$.\qed
\end{example}

\begin{center}
	\begin{figure}[h!]
		\centering
		\includegraphics[scale=1.5]{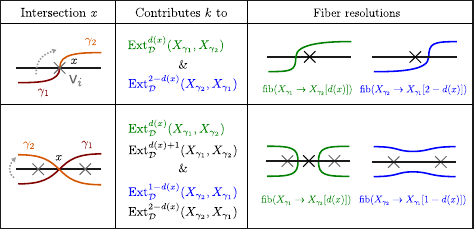}
		\caption{A summary of the key cases in \cref{prop:morphisms_via_intersections}. The first column depicts the type intersection from $\g_1$ to $\g_2$, endpoint above and crossing below. The second column records its contributions to the morphisms between the matching objects $X_{\g_1}$ and $X_{\g_2}$. The third depicts, when possible, a description of the corresponding extensions in terms of matching objects.} 
        \label{fig:WeaveVertex_Ext0}
	\end{figure}
\end{center}

Let $\gamma_1,\gamma_2$ be two graded $\pi_\D$-matching paths, our goal is to describe the morphism object $\on{Mor}_\D(X_{\gamma_1},X_{\gamma_2})$ in terms of the count of graded intersections of $\gamma_1,\gamma_2$. We will only need the case where at least one of $\gamma_1$ or $\gamma_2$ is embedded, and thus we assume that. The following result provides such a description of $\on{Mor}_\D(X_{\gamma_1},X_{\gamma_2})$, where its two key cases are summarized in \Cref{fig:WeaveVertex_Ext0}:

\begin{proposition}\label{prop:morphisms_via_intersections}
Let $\gamma_1,\gamma_2$ be two graded $\pi_\D$-matching paths whose underlying curves are distinct and $X_{\g_1},X_{\g_2}\in\D$ their associated matching objects. Suppose that either $\gamma_1$ or $\gamma_2$ is embedded and choose the order so that $\g_1$ preceeds $\g_2$. Then there exists an equivalence in $\D^{\on{perf}}(k)$
\begin{equation}\label{eq:morphisms_from_intersections}
\on{Mor}_\D(X_{\gamma_1},X_{\gamma_2})\simeq \bigoplus_{x\text{ end.~int.}}k[-d(x)]\oplus \bigoplus_{x\text{ crossing}}k[-d(x)]\oplus k[-d(x)-1]\,.
\end{equation}
\end{proposition}

\begin{proof}
In the case of embedded $\pi_\D$-matching paths, \eqref{eq:morphisms_from_intersections} is proved in different levels of generality and contexts in \cite{KS02,IQZ20,Chr21b}. Specifically, in the context of perverse schobers, \cite[Thm.~6.4]{Chr21b}, implies \eqref{eq:morphisms_from_intersections} for graded $\pi_\D$-matching paths without segments of the third type. In order to reduce the general case to the case of $\pi_\D$-matching paths without segments of the third type, we proceed as follows.\\

Swapping $\gamma_1$ and $\gamma_2$ if necessary, we assume that $\gamma_1$ is embedded. By acting with an autoequivalence of $\D$, we can and do assume that $X_{\gamma_1}\simeq S_i$ for some $i\in[1,n]$. Let us choose a sequence of graded $\pi_\D$-matching paths $\gamma_1',\dots,\gamma_m'$ without segments of the third type together with iterated fiber and cofiber sequences $Y_{i+1}\to X_{\gamma_{i+1}'}\to X_{\gamma_i'}$ such that $Y_m\simeq X_{\gamma_2}$. The number $m$ is the number of segments of the third type of $\gamma_2$. We do choose this sequence such that the morphism $X_{\gamma_{i+1}'}\to X_{\gamma_i'}$ arises from an endpoint intersection away from the two endpoints of $\gamma_1$ and geometrically $\gamma_{i}',\gamma_{i+1}',\gamma_1$ do not form a triangle involving two crossings with $\gamma_1$ and the endpoint intersection. Then the arising fiber and cofiber sequences, which read
\[
\on{Mor}_\D(X_{\gamma_1},X_{\gamma_{i}'}[-1])\to \on{Mor}_\D(X_{\gamma_1},Y_{i+1})\to \on{Mor}_\D(X_{\gamma_1},X_{\gamma_{i+1}'}),
\]
all split. Furthermore the number of intersections of $\gamma_1$ with $\gamma_2$ is equal to the number of intersections of $\gamma_1$ with $\gamma_1',\dots,\gamma_m'$. This reduces the general case to cases with no segments of the third type.
\end{proof}

\begin{center}
	\begin{figure}[h!]
		\centering
		\includegraphics[scale=1.5]{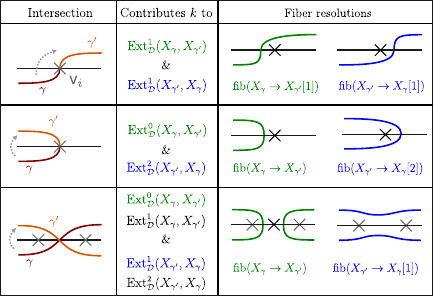}
		\caption{The intersections of $\pi_\D$-matching paths and their contribution to morphisms between the associated matching objects in the case of the canonical grading of \cref{rem:canonical_grading}. For instance, by \cref{ex:endpointintersection_index} the first row here must coincide with the first row of \cref{fig:WeaveVertex_Ext1} with the endpoint intersection $\vd_i=x$ indexed by $d(x)=1$.}
        \label{fig:WeaveVertex_Ext1}
	\end{figure}
\end{center}

\begin{example}\label{ex:matching_paths2} Let us continue with \cref{ex:matching_paths1}, cf.~\cref{fig:MatchingPaths_GraphD}. By \cref{ex:simple_as_matching_sphere}, $S_i\simeq X_{\ed_i}$, where the edge $\ed_i$ is understood as a straight $\pi_\D$-matching path. By \cref{ex:endpointintersection_index}, the endpoint intersection index of $\ed_1$ and $\ed_2$ is 1. By \cref{prop:morphisms_via_intersections},
$$\on{Mor}(X_{\ed_1},X_{\ed_2})\simeq k[-1],$$
which indeed recovers the fact that $\mbox{Ext}^1(S_1,S_2)\simeq k$ and $\mbox{Ext}^i(S_1,S_2)\simeq0$ if $i\neq1$. Similarly, for $\g_2$ as in \cref{fig:MatchingPaths_GraphD}, the morphisms between the matching object $X_{\g_2}$ and $S_2$ are given by
$$\on{Mor}(X_{\g_2},X_{\ed_2})\simeq k\oplus k[-1],$$
since their only intersection point is a crossing $x$ of degree $d(x)=0$, cf.~third row of \cref{fig:WeaveVertex_Ext1}.\qed
\end{example}

\begin{remark}\label{rmk:immersed_fiber_resolution}
(1) The resolution of an intersection can fail to be union of $\pi_\D$-matching paths because of the second assumption in \Cref{def:pi_D_matching_path} stating that $\pi_\D$-matching paths should not cut out any once-punctured discs. In this case, the $\pi_\D$-matching paths describing the fiber are obtained by applying the local move in \cref{fig:MatchingPath_resolution}.

\begin{center}
	\begin{figure}[h!]
		\centering
		\includegraphics[scale=1.5]{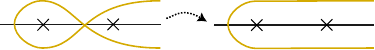}
		\caption{The resolution in \cref{rmk:immersed_fiber_resolution}, describing fibers as matching objects.}
        \label{fig:MatchingPath_resolution}
	\end{figure}
\end{center}

(2) For crossings with immersed $\pi_\D$-matching paths arising from self-extensions as in \Cref{rem:curves_for_self_extensions}, the fibers of both morphisms arising from the crossing can be matching (and be computed by choosing a different representative). However, if the fiber of a morphism between two matching objects arising from an intersection is not matching, it must arise from a crossings and lie in the larger of the two degrees of the crossing.\qed
\end{remark}

Finally, we will need to use that semi-twists preserve objects in the subcategory $\mathcal{M}\sse\D$. This is the content of the following:

\begin{lemma}[semi-twists preserve $\mathcal{M}$]\label{lem:matching_curves_stable_under_semitwist}
Let $\gamma,\gamma'$ be graded $\pi_\D$-matching paths and $X_\g,X_{\g'}\in\mathcal{M}$ their associated matching objects. Suppose that $\gamma'$ is embedded. Then both semi-twists $T_{X_{\gamma'}}^{-,\leq 0}(X_{\gamma})$ and $T_{X_{\gamma'}}^{+,\geq 0}(X_{\gamma})$ are matching, i.e.~they are objects in $\mathcal{M}$.
\end{lemma}

\begin{proof}
By definition, $T_{X_{\gamma'}}^{+,\geq 0}(X_{\gamma})=\on{cof}(\tau_{\geq 0}(\on{Mor}_\D(X_{\gamma'},X_{\gamma}))\otimes X_{\gamma'}\to X_{\gamma})$. By \cref{prop:morphisms_via_intersections}, the dimension of the homology of the morphism object $\on{Mor}_\D(X_{\gamma'},X_{\gamma})$ counts endpoint intersections and crossings of $\gamma,\gamma'$. If an endpoint intersection $x\in \g\cap\g'$ gives rise to a negative extension, then the passage to the cofiber of the corresponding morphism  $X_{\gamma'}[i]\stackrel{x}{\to} X_{\gamma}$, with $i\geq 0$, has the effect of resolving that endpoint intersection of $\gamma,\gamma'$, cf.~first row of \cref{fig:WeaveVertex_Ext0}.\\

\noindent If the intersection $x\in \g\cap\g'$ is a crossing contributing negative extensions, there are two options. First, if both of the corresponding extensions are negative, then the passage to the cofiber $X_{\gamma'}[i]\oplus X_{\gamma'}[i+1]\to X_{\gamma}$ has the effect of applying a Dehn-twist around $\gamma'$ to $\gamma$ at the crossing, for $i\in[0,\infty)$. Second, if only one of the corresponding extensions is negative, then the cofiber of the morphism $X_{\gamma}\to X_{\gamma'}$ is given by the direct sum of the two $\pi_\D$-matching paths arising from resolving the crossing. See the second row of \cref{fig:WeaveVertex_Ext0} for a depiction of these cases. Therefore, for either types of intersections, the passage to the cofiber defining the semi-twist remains in $\mathcal{M}$, which proves the statement for $T_{X_{\gamma'}}^{+,\geq 0}(X_{\gamma})$. The proof that $T_{X_{\gamma'}}^{-,\leq 0}(X_{\gamma})\in\mathcal{M}$ is analogous. 
\end{proof}


\subsection{The filtered object $\lbw$ associated with a weighted braid word $\bw$}
\label{ssec:filteredobject_from_weightedbraid} 
The object of this subsection is to prove \cref{thm:unique_rigid_resolutions}. Let us formally define the notions featured in its statement, including that of a weighted braid word $\bw$ and a $\bw$-filtration.

\begin{definition}[Weighted braid words]\label{def:weighted_braidword}
Let $\beta=\sigma_{i_0}\dots \sigma_{i_l}$ be a positive braid word of length $l+1$. By definition, a weight ${\bf a}=(a_0,\dots,a_l)$ for $\beta$ consists of a sequence of non-negative integers $a_i\geq 0$. A weighted braid word is a pair $(\beta,{\bf a})$ consisting of a positive braid word $\b$ and a weight for $\b$.\qed
\end{definition}

\noindent In the case of permutations, \cref{def:weighted_braidword} appeared in \cite[Definition 6.13]{casals2025comparingclusteralgebrasbraid} in the form of weighted expressions and their associated Lusztig data. Recall from \cite[Definition 1.2.2.9]{HA} that a $\Z_{\geq0}$-filtered object in $\D$ is a functor $N(\Z_{\geq0})\to\D$, where $N(\Z_{\geq0})$ is the nerve of $\Z_{\geq0}$, with $\Z_{\geq0}$ being understood as a poset category.
Given a positive braid word $\beta=\sigma_{i_0}\dots \sigma_{i_l}$ we denote
\begin{equation}\label{eq:betaj}
\beta_j\coloneqq \sigma_{i_{j}}\sigma_{i_{j+1}}\dots \sigma_{i_l},\quad j\in[0,l]
\end{equation}
\noindent Note that $\beta_0=\beta$ and $\beta_l=\sigma_{i_l}$.

\begin{definition}[$\bw$-filtered objects]\label{def:filtered_by_weighted_braid_word}
Let $(\beta,{\bf a})$ be a weighted braid word of length $l+1$ and $X\in \D$ an object. By definition, a $(\beta,{\bf a})$-filtration of $X$ is a $\mathbb{Z}_{\geq 0}$-filtered object $X(\mhyphen)\colon N(\mathbb{Z}_{\geq 0})\to \D$ such that:
\begin{enumerate}
    \item the $i$th graded piece $\on{gr}_i(X(\mhyphen))\coloneqq \on{cof}(X(i-1)\to X(i))$ is

    \begin{equation}
    \on{gr}_i(X(\mhyphen))\simeq
    \begin{cases}
        S_{\beta_i}^{\oplus a_i} & \mbox{ if }i\in[0,l],\\
        0 & \mbox{ if }i\in[l+1,\infty).
    \end{cases}
    \end{equation}

    \item As objects in $\D$, there exists an equivalence $\on{colim} X(\mhyphen)\simeq X$.
\end{enumerate}
By definition, a $\bw$-filtered object $X\in\D$ is an object $X\in\D$ equipped with a $\bw$-filtration.\qed
\end{definition}

\begin{remark}\label{rmk:geometric_realization_filtered_object}
Given a $\mathbb{Z}_{\geq 0}$-filtered object $X(\mhyphen)$ in $\D$, the colimit $\on{colim} X(\mhyphen)$ is also called the geometric realization of $X(\mhyphen)$, as it corresponds to the (simplicial) geometric realization under the Lurie--Dold--Kan correspondence, cf.~\cite[Rem.~1.2.4.2]{HA}. In the setting of \Cref{def:filtered_by_weighted_braid_word}, $\on{colim} X(\mhyphen)$ is equivalent to $X(l)$.\qed 
\end{remark}

\noindent In general, there might be multiple $\bw$-filtered objects in $\D$. To obtain uniqueness, as in \cref{thm:unique_rigid_resolutions}, we focus on rigid objects. Recall that an object $X\in\D$ is called rigid if $\on{Ext}^1_\D(X,X)\simeq 0$.

\begin{example}
    Let $\beta$ be a braid word of length $l+1$ and ${\bf a}=(0,\dots,1,\dots,0)$ the weight with a unique non-zero entry $1$ in position $i\in[0,l]$. Then there exists a unique $(\beta,{\bf a})$-filtered object in $\D$ up to equivalence: it is given by $S_{\beta_i}$ equipped with the filtration 
    \[
    0\to \dots \to 0 \to S_{\beta_i}\xrightarrow{\on{id}} S_{\beta_i} \xrightarrow{\on{id}}S_{\beta_i} \xrightarrow{\on{id}}\dots\,,
    \]
    where the first appearance of $S_{\beta_i}$ is at the $i$th entry.\qed
\end{example}

Given an object $X$ with a $(\beta,{\bf a})$-filtration, we can recover the weight from $X$ and $\beta$:

\begin{lemma}\label{lem:uniqueness_of_weight}
    Let $\beta$ be a positive braid word. For any object $X\in \D$ there exists at most one weight ${\bf a}$ for $\beta$, such that $X$ admits a $(\beta,{\bf a})$-filtration.
\end{lemma}    

\begin{proof}
Suppose that $X$ admits a $(\beta,{\bf a})$ and a $(\beta,{\bf b})$-filtration with ${\bf a}\not ={\bf b}$. Let $l+1$ be the length of $\beta$. There exists a maximal $0\leq i\leq l$ such that $b_i\not=a_i$. Switching ${\bf a},{\bf b}$ if necessary, we may assume that $b_i>a_i$. Passing to repeated fibers of morphisms $X\to S_{\beta_j}^{\oplus a_j}$, $j\geq i$, we may assume that $a_j=b_j=0$ for $j>i$ and $a_i=0$ but $b_i>0$. Let $\C_\ast$ be the coherent chain complex corresponding to the $(\beta,{\bf b})$-filtration of $X$. We choose any non-zero morphism $S_{\beta_i}[-2]\to S_{\beta_i}^{\oplus b_i}$ and extend it using \Cref{lem:connectiveHom} to a morphism between coherent chain complexes $(S_{\beta_i}[-2-i])_\ast\to \C_\ast$, where $(S_{\beta_i}[-2-i])_\ast$ is considered in degree $i$ with value $S_{\beta_i}[-2-i]$. Passing to the cofiber totalization defines a non-zero morphism $S_{\beta_i}[-2]\to X$. Using that $X$ admits a $(\beta,{\bf a})$-filtration with $a_j=0$ for $j\geq i$ and \Cref{lem:connectiveHom}, we find that $\on{Ext}^2_\D(S_{\beta_i},X)\simeq 0$, which is a contradiction. We must thus have ${\bf a}={\bf b}$.
\end{proof}


\subsubsection{Preparatory lemmas}\label{sssec:preparatory_lemmas} We use the following two lemmas, \cref{lem:non_rigidity_preserved_under_resolution} and \cref{lem:resolving_crossings_non_geometric}, in the proof of \Cref{thm:unique_rigid_resolutions}.

\begin{lemma}\label{lem:non_rigidity_preserved_under_resolution}
    Let $S\in \D$ be a $2$-spherical object and $A\in \D$ be such that $\on{Ext}^i_\D(S,A)\simeq 0$ for all $i\in[2,\infty)$. Then, for any non-zero morphism $f:S[-1]\to A$,
    \[ \on{dim}\on{Ext}^1_\D(\cof(f),\cof(f))\geq \on{dim}\on{Ext}^1_\D(A,A)\,.\] 
    In particular, if $A$ is not rigid, then the cofiber $\cof(f)$ is not rigid.
\end{lemma}

\begin{proof}
    Let us denote $A'\coloneqq \cof(f)$. Since the morphism object functor is exact in both entries, the morphism object $\on{Mor}_\D(A',A')$ is equivalent to the horizontal cofiber and vertical fiber of the following square in $\D(k)$:
    \[
    \begin{tikzcd}
{\on{Mor}_\D(A,S[-1])} \arrow[r] \arrow[d] & {\on{Mor}_\D(A,A)} \arrow[d] \\
{\underbrace{\on{Mor}_\D(S[-1],S[-1])}_{\simeq k\oplus k[-2]}} \arrow[r]       & {\on{Mor}(S[-1],A)}         
\end{tikzcd}
    \]
where $\on{Mor}_\D(S[-1],S[-1])\simeq k\oplus k[-2]$ follows from the assumption that $S$ is 2-spherical.
By hypothesis, the bottom right entry $\on{Mor}(S[-1],A)$ is connective, and by the $2$-Calabi--Yau property on $\D$ the top left entry $\on{Mor}(A,S[-1])$ is concentrated in homological degrees $i\in(-\infty,-2]$. It thus follows that the above square contains a direct summand in the stable $\infty$-category of squares in $\D(k)$ of the form
\[
\begin{tikzcd}
0 \arrow[r] \arrow[d] & {\on{Ext}^1_\D(A,A)[-1]} \arrow[d] \\
0 \arrow[r]           & 0                             
\end{tikzcd}
\]
This implies that there is a direct summand $\on{Ext}^1_D(A,A)\subset \on{Ext}^1_\D(A',A')$, and the statement follows.
\end{proof}

\noindent For the next lemma and its proof, cf.~\cref{def:M_geometric_objects} for the subcategory $\mathcal{M}\sse\D$ additively spanned by matching objects, and \cref{sssec:morphisms_2spherical_matching_objects} for the degrees between crossings of the associated $\pi_\D$-matching paths.

\begin{lemma}\label{lem:resolving_crossings_non_geometric}
Let $S\in \mathcal{M}$ be a $2$-spherical object and $X\in \D$ an object that does not lie in $\mathcal{M}$. Suppose that $\on{Mor}_\D(S,X)$ is concentrated in degrees $i\in[-2,\infty)$. Then, for any morphism $f:S[-1]\to X$, the cofiber $\cof(f)$ does not lie in $\mathcal{M}$.
\end{lemma}

\begin{proof}
By contradiction, suppose that the cofiber $Y\coloneqq \cof(f)$ of the morphism $S[-1]\to X$ lies in $\mathcal{M}$. Then the cofiber $X$ of the morphism $Y[-1]\to S[-1]$ is not matching and it thus has a contribution from a crossing from $Y$ to $S$ of degrees $-1$ and $0$, cf.~\cref{fig:WeaveVertex_Ext0}. By the $2$-Calabi--Yau property of $\D$, there is also a crossing from $S$ to $Y$ with degrees $2$ and $3$. This is in contradiction with the fact that, by hypothesis, the morphism object
\[
\on{Mor}_\D(S,Y)\simeq \on{cof}(\underbrace{\on{Mor}_\D(S,S[-1])}_{\simeq k[-1]\oplus k[-3]}\to \on{Mor}_\D(S,X))
\]
is concentrated in degrees $i\in[-2,\infty)$.
\end{proof}

\begin{corollary}\label{lem:pieces_of_filtration_are_rigid_and_geometric}
Let $\bw$ be a weighted braid word, $\mathcal{L}(\mhyphen)\colon N(\mathbb{Z}_{\geq 0})\to \D$ a $(\beta,{\bf a})$-filtered object with geometric realization $L$.
\begin{itemize}
    \item[(1)] If $L$ is rigid, then $\mathcal{L}(i)$ is also rigid for all $i\in[0,\infty)$
    \item[(2)] If $L\in\mathcal{M}$, then also $\mathcal{L}(i)\in \mathcal{M}$ for all $i\in [0,\infty)$.
\end{itemize}
\end{corollary}

\begin{proof}
Using \Cref{lem:connectiveHom}, part (1) follows from \Cref{lem:non_rigidity_preserved_under_resolution}, and part (2) follows from \Cref{lem:resolving_crossings_non_geometric}.
\end{proof}


\subsubsection{Resolving endpoint intersections}\label{sssec:endpoint_intersections} 

The remaining ingredient before proving \cref{thm:unique_rigid_resolutions} is used to construct the rigid filtered objects $\lbw$: it is the construction of certain morphisms whose cofibers are rigid, as established in \cref{prop:unique_morphism_with_rigid_cofiber} below. For that, we start with the cyclic order on endpoint intersections at a vertex of $\rgraph_\D$, defined as follows:

\begin{definition}[Cyclic order at vertices]\label{def:cyclic_order_endpoint_intersections}
Let $\vd\in\rgraph_\D$ be a vertex and $C\coloneqq \{\gamma_1,\ldots,\g_m\}$ be a collection of distinct $\pi_\D$-matching paths all of which have an endpoint at $\vd\in\rgraph_\D$, arriving at $\vd$ with different slopes. By definition, the intersection order on $C$ is the cyclic order on $C$ obtained by choosing representatives of the $\pi_\D$-matching paths in $C$ which have minimal numbers of intersections with each other, and declaring the cyclic order to be induced by the clockwise orientation. That is, in the intersection order, $\gamma_i<\gamma_j<\gamma_k$ 
if and only if the angle with which $\gamma_j$ approaches $\vd$ lies between the angles with which $\gamma_k$ and $\gamma_i$ approach $\vd$.\qed
\end{definition}

\noindent In \cref{def:cyclic_order_endpoint_intersections}, angles are understood to lie in $[0,2\pi)$, and cf.~\cite[Section 3a]{KS02} for a discussion on representatives of paths with minimal number of intersections.

\begin{remark}\label{rem:cyclic_order_homological_description}
An equivalent and more algebraic version of \Cref{def:cyclic_order_endpoint_intersections} is that $\gamma_i<\gamma_j<\gamma_k$ if and only if the composite of the non-zero morphisms $X_{\gamma_i}\to X_{\gamma_j}[a]$ and the $[a]$ shift of $X_{\gamma_j}\to X_{\gamma_k}[b]$ arising from the endpoint intersections of $\g_i$ and $\g_j$, and $\g_j$ and $\g_k$, yields a non-zero morphism $X_{\g_i}\to X_{\gamma_k}[a+b]$.\qed
\end{remark}

\begin{proposition}\label{prop:unique_morphism_with_rigid_cofiber}
    Let $X\in \M$ be rigid object and $S\in \M$ a 2-spherical object. Suppose that $\on{Ext}^i_\D(S,X)\simeq 0$, $i\in[2,\infty)$, and that there are no degree $1$ crossings from $S$ to $X$. Then there exists a unique morphism $f:S[-1]\to X$, unique up to postcomposition with automorphisms of $X$, whose cofiber $\cof(f)$ is rigid. Furthermore, this cofiber is matching, i.e.~$\cof(f)\in \M$. 
\end{proposition}

\begin{proof}
Let us write $Y\coloneqq \cof(f)$, and express $S=X_{\gamma}$ for some $\pi_\D$-matching path $\gamma$. Let $\vd,\vd'\in\rgraph_\D$ be the two vertices of where $\gamma$ ends. Consider the induced intersection order, cf.~\cref{def:cyclic_order_endpoint_intersections}, on the set $C$, resp.~$C'$, consisting of those $\pi_\D$-matching paths describing $X\in\M$ with an endpoint at $\vd$, resp.~$\vd'$, and $\g$ itself. We turn this cyclic orders on $C$ and $C'$ into total orders by declaring $\gamma$ to be the largest element in this total order. In each of these two total orders, consider the smallest element such that the corresponding endpoint intersection with $\gamma$ is of degree $1$. These yield two morphisms $a,a'\colon S[-1]\to X$, one from such smallest element in $C$ and the corresponding one from $C'$.\\

\noindent Now, we claim that the required morphism is their sum $f\coloneqq a+a'$, $f: S[-1]\to X$. Let us argue that $\cof(f)$ must be rigid. First, note that the cofiber of any morphism $S[-1]\to X$ not containing contributions from both $a$ and $a'$ cannot itself be rigid. Indeed, passing to the cofiber amounts to resolving the endpoint intersections, cf.~\cref{fig:WeaveVertex_Ext1}, and resolving degree $1$ endpoint intersections different than the minimal degree $1$-ones leads to a self-crossing of the cofiber near $\vd$, or $\vd'$, of degrees $0$ and $1$, which contradicts rigidity. Now, to show that $\cof(f)$ is rigid, note that if $a$ or $a'$ corresponds to an endpoint intersection with non-minimal elements (though minimal of degree $1$), the assumption $\on{Ext}^j_\D(S,X)\simeq 0$, $j\in[2,\infty)$, and the rigidity $\on{Ext}^1(X,X)\simeq 0$ imply that either $a$ or $a'$ must factor as $S[-1]\to X[-i]\to X$ for some $i\in(-\infty,-2]$. That said, the degree $(-i+1)$-endpoint intersection leads to self-crossings of $\on{cof}(a+a')$ in degrees $(-i+1),(-i+2)$ and $i,(i+1)$, and since these are not degree $1$, we conclude $\on{cof}(f)$ is rigid.\\

Finally, for uniqueness, note that the degree $0$ endpoint intersections of $X$ with itself at $\vd$ and $\vd'$ give rise to automorphisms of $X$ which allow to remove all contributions from degree $1$ endpoint intersections that are larger in the total order, as follows from \Cref{rem:cyclic_order_homological_description}. Thus the morphism is uniquely determined up to such automorphisms, as required.
\end{proof}

\begin{center}
	\begin{figure}[h!]
		\centering
		\includegraphics[scale=1.5]{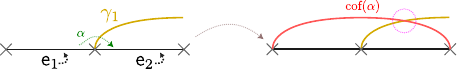}
		\caption{(Left) Pieces of $\pi_\D$-matching paths for \cref{ex:endpoint_resolutions}, with $\g_1$ indicating only a piece that continues to other parts of the graph $\rgraph_\D$, and $\ed_1,\ed_2$ as $\pi_\D$-matching paths inside a neighborhood containing $\vd_1,\vd_2$ and $\vd_3$. The morphism $\alpha:X_{\ed_1}[-1]\lr X_{\ed_2}$ of degree 1 is indicated in green. (Right) The endpoint resolution of the morphism $\alpha$, which yields a crossing intersection with the matching object $X_{\g_1}$. Note that the $\pi_\D$-matching path for $\cof(\alpha)$ is given by the image of the inverse half-twist $f_2^{-1}(\ed_1)$, cf.~\cref{lem:twist_via_homeomorphism}.}
        \label{fig:crossing_resolutions}
	\end{figure}
\end{center}

\begin{example}\label{ex:endpoint_resolutions} (1) Consider the case $n=3$, the two $\pi_D$-matching paths $\ed_1$ and $\ed_2$, and choose $S\coloneqq X_{\ed_1}\simeq S_1$ and $X\coloneqq X_{\ed_2}\simeq S_2$, cf.~\cref{ex:simple_as_matching_sphere}. The hypothesis of \cref{prop:unique_morphism_with_rigid_cofiber} are satisfied, and the conclusion is straightforward in this example. Indeed, there is a unique degree 1 morphism $f:S[-1]\lr X$, see e.g.~\cref{ex:endpointintersection_index}, and its cofiber $\cof(f)\simeq T_X^{-1}(S)$ is indeed rigid and, by \cref{lem:twist_via_homeomorphism}, matching, see e.g.~\cref{fig:crossing_resolutions} (right).\\

\noindent (2) \cref{fig:crossing_resolutions} depicts an instance illustrating that resolving an endpoint intersection differently than in the proof of \cref{prop:unique_morphism_with_rigid_cofiber} yields a non-rigid object. Consider $S\coloneqq X_{\ed_1}$ and $X=X_{\g_1}\oplus X_{\ed_2}$, such that the two $\pi_\D$-matching paths for $X$ are $\g_1$ and $\ed_2$, as in \cref{fig:crossing_resolutions} (left). Suppose that $X$ is rigid, which is for instance the case if the $\pi_\D$-matching path $\g_1$ only intersects $\ed_2$ at $\vd_2$, by \cref{fig:WeaveVertex_Ext1}. In this situation, reading clockwise, we have a degree 1 morphism from $S$ to $X_{\g_1}$, a degree 0 morphism from $X_{\g_1}$ to $X_{\ed_2}$, and a degree 1 morphism from $X_{\ed_2}$ to $S$. Their intersection order at $\vd_2$ is $X_{\g_1}<X_{\ed_2}<S$ and according to the proof of \cref{prop:unique_morphism_with_rigid_cofiber} we should resolve the degree 1 morphism from $S$ to $X_{\g_1}$ if we want to obtain a rigid object.\\

If, for illustrative purposes, we instead resolved the degree 1 morphism $\alpha$ from $S$ to $X_{\ed_2}$, depicted in green in \cref{fig:crossing_resolutions} (left), then its resolution would lead to an object $X_{\g_1}\oplus\cof(\alpha)$. That said, the direct sum $X_{\g_1}\oplus\cof(\alpha)$ would not be rigid due to the presence of the crossing between $\cof(\alpha)$ and $X_{\g_1}$, which is circled in pink in \cref{fig:crossing_resolutions} (right). Indeed, this crossing contributes to $\on{Ext}^1(\cof(\alpha),X_{\g_1})$ by \cref{fig:WeaveVertex_Ext1}, and thus the direct sum $X_{\g_1}\oplus\cof(\alpha)$ is not rigid.\qed
\end{example}

\subsubsection{Proof of \Cref{thm:unique_rigid_resolutions}}\label{sssec:proof_unique_rigid_resolutions} 
\begin{proof}[Proof of \Cref{thm:unique_rigid_resolutions}]
Let us suppose that the $\bw$-filtered object $\lbw$ exists and prove its uniqueness, up to equivalence. Let $\lbw'$ be another $(\beta,{\bf a})$-filtered object with rigid geometric realization. We now prove that the restrictions of $\lbw$ and $\lbw'$ to $N(\mathbb{Z}_{[0,i]})$ are equivalent, by induction on $i\in [0,\infty)$. The base case $i=0$ holds by construction. For the inductive step, we suppose that the assertion holds for $i-1$, and show it for $i$, as follows. By definition, the $i$th object $\lbw(i)$ of the filtered object $\lbw$ arises as the cofiber of a morphism $S_{\beta_i}^{\oplus a_i}[-1]\to \lbw(i-1)$ or, equivalently, as the $a_i$-many iterated cofibers along morphisms with domain $S_{\beta_i}[-1]$. By \Cref{lem:pieces_of_filtration_are_rigid_and_geometric}, $\lbw(i)$ is rigid and lies in $\mathcal{M}$, and thus $S_{\beta_i}$ and $\lbw(i-1)\simeq \lbw'(i-1)$ cannot have a crossing contributing to degrees $0$ and $1$. \Cref{prop:unique_morphism_with_rigid_cofiber} thus shows that there exists a unique morphism $f\colon S_{\beta_i}^{\oplus a_i}[-1]\to \lbw(i-1)$, up to equivalence, whose cofiber is rigid. Therefore, the morphisms $\lbw(i-1)\to \lbw(i)$ and $\lbw(i-1)'\to \lbw(i)'$ must both be equivalent to the cofiber morphism $\lbw(i-1)\to \on{cof}(f)$ of $f$. We thus obtain that the restrictions of $\lbw$ and $\lbw'$ to $N(\mathbb{Z}_{[0,i]})$ are also equivalent, concluding the induction step. This establishes the uniqueness of $\lbw$, up to equivalence.\end{proof}

\noindent For context, thought it will not be used in this article, we record our expectation that the following variation of \cref{thm:unique_rigid_resolutions} might hold:

\begin{conjecture}\label{conj:filtered_by_weighted_braid_words}
Let $\bw$ be a weighted braid word. Then there exists a unique $\bw$-filtered object $\lbw\in\D$ whose geometric realization is rigid.\qed
\end{conjecture}

\noindent \cref{conj:filtered_by_weighted_braid_words} makes no assertion that $\lbw\in\mathcal{M}$, cf.~\cref{thm:unique_rigid_resolutions}. An approach to such conjecture via our argument for \cref{thm:unique_rigid_resolutions} would require understanding how to generalize the endpoint resolutions, cf.~the proof of \cref{prop:unique_morphism_with_rigid_cofiber}, to interact with degree 1 crossings while preserving rigidity, cf.~\cref{ex:endpoint_resolutions}.


\subsubsection{A few comments on \texorpdfstring{$\lbw$}{the filtered object Lba}}\label{sssec:comments_filteredobject}

{\bf (A) The recursive form for $\lbw$.} The proof of \cref{thm:unique_rigid_resolutions} provides an explicit construction of the filtered object $\lbw$, if it exists, as follows. It is given by the $\D$-valued $\Z_{\geq0}$-filtration
\begin{equation}\label{eq:bw_filtration1}
\lbw(0)\stackrel{f_1}{\lr} \lbw(1)\stackrel{f_2}{\lr}\cdots \lbw(l-1)\stackrel{f_l}{\lr} \lbw(l)
\end{equation}
recursively defined by
\begin{equation}\label{eq:bw_filtration2}
\lbw(i)\coloneqq \cof(S_{\b_i}^{\oplus a_i}[-1]\stackrel{\vartheta_i}{\lr}\lbw(i-1)),\quad i\in[1,l],\qquad 
\lbw(0)\coloneqq S_{\b_0}^{\oplus a_0},
\end{equation}
and constant equal to $\lbw(i)\coloneqq \lbw(l)$ for all $i\in[l,\infty)$. By \cref{prop:unique_morphism_with_rigid_cofiber}, the morphism
\[\vartheta_i\colon S_{\b_i}^{\oplus a_i}[-1]\lr\lbw(i-1)\]
corresponds to the (up to) $a_i$ endpoint intersections of degree $1$, at each of the endpoints of $S_{\beta_i}$, which are minimal in the order described in \Cref{def:cyclic_order_endpoint_intersections}. The morphisms $f_i$ of the filtration (\ref{eq:bw_filtration1}) are given by the canonical inclusion of $\lbw(i-1)$ into the cofiber of $\vartheta_i$.

\begin{remark}
By the proof of \cref{prop:unique_morphism_with_rigid_cofiber}, the cofiber
\[\cof(S_{\b_i}^{\oplus a_i}[-1]\stackrel{\vartheta_i}{\lr}\lbw(i-1))\]
can be geometrically understood as a resolution of endpoint intersections of the matching paths describing $S_{\b_i}^{\oplus a_i}$, which are readily obtained by using \cref{lem:twist_via_homeomorphism}, and the matching paths defining $\lbw(i-1)$. In this sense, $\lbw$ can be understood as an iterated resolution of endpoint intersections of matching paths.\qed
\end{remark}


{\bf (B) The coherent chain complex $\C_{\beta,{\bf a}}$ corresponding to $\lbw$.} The $\bw$-filtered object $\lbw$ can be alternatively described in terms of a coherent chain complex, as follows from the $\infty$-categorical Dold--Kan correspondence, see also \Cref{rem:Dold_Kan}. Appendix \ref{sec:CoherentChainComplexes} provides a summary on coherent chain complexes. Intuitively, the role of coherent chain complexes in a stable $\infty$-category is akin to that of twisted complexes of $A_\infty$-modules in a $A_\infty$-category. In a nutshell, the graded pieces of any $\mathbb{Z}_{\geq 0}$-filtered object $X(\mhyphen)$ in $\D$ determine a $\D$-valued coherent chain complex. Specifically, the $i$-th term of such a coherent chain complex is the shift of the corresponding piece of the associated graded, i.e.~the cofiber $\cof(f_i)[-i]$, where $f_i\colon X(i-1)\to X(i)$ is the morphism defining the filtered object.
In our case this coherent chain complex is described as follows.\\

\noindent Let $(\beta,{\bf a})$ be a weighted braid word arising from a horizontal slice of a Lusztig cycle. We denote by $\cbw\in\on{Ch}_{\geq 0}(\D)$ the coherent chain complex determined by the filtered object $\lbw$. By construction, cf.~\cref{eq:bw_filtration2}, the $i$th graded piece $\cof(f_i)$ of the filtration (\ref{eq:bw_filtration1}) is equivalent to $S_{\beta_i}^{\oplus a_i}$, and thus $\cbw$ is concentrated in degrees $0$ to $l$ and has the form
    \begin{equation}\label{eq:bw_coherentchaincomplex}
\C_{\beta,{\bf a}}=\left(S_{\beta_l}^{\oplus a_l}[-l]\lr S_{\beta_{l-1}}^{\oplus a_{l-1}}[-l+1]\lr \dots \lr S_{\beta_{2}}^{\oplus a_{2}}[-2]\lr S_{\beta_{1}}^{\oplus a_{1}}[-1]\lr S_{\beta_0}^{\oplus a_0}\right)
\end{equation}

\noindent where we have drawn arrows for the differentials of degree $(-1)$, but we did not draw those corresponding to the higher coherence data. The geometric realization $L_{\beta,{\bf a}}=\on{colim}\lbw(i)$ is equivalent to the total cofiber of $\C_{\ba}$, which can be computed as an iterated cofiber, see \cref{eq:bw_filtration2} above and \Cref{rem:iterated_cofiber}. In particular, $L_{\beta,{\bf a}}$ is given by an iterated extension between the $S_{\beta_i}^{\oplus a_i}$'s.

\begin{remark}
Even though the coherent chain complexes $\cbw$ contains the same data as the filtered objects $\lbw$, it is sometimes convenient to use them when establishing results related to $\lbw$. In particular, we use this presentation in terms of coherent chain complexes in the proofs of some of the results on categorical Lusztig cycles, cf.~the proofs of \cref{prop:4_and_6_moves} and \cref{prop:semi_twist_Lusztig_cycle_at_trivalent_vertex}.\qed
\end{remark}


{\bf (C) The role of rigidity.} We conclude this section with an example illustrating the behavior of non-rigid $\bw$-filtered objects under the tropical propagation rules for weave cycles. These weighted braid words do arise as slices of weave cycles, though not as slices of Lusztig cycles. We hope this example highlights the importance of rigidity in our constructions.

\begin{example}\label{ex:6move_consistency}
To ease notation in this example, a weighted braid word $\bw$ where ${\bf a}$ takes only the values $0$ and $1$ will be denoted by writing the braid word $\beta$ with dots above the letters where ${\bf a}$ evaluates to $1$. For instance, $s_1\dot{s}_2s_1$ denotes $\bw=(s_1s_2s_1,(0,1,0))$.\\

Let us consider the follow three weighted braid words:
\[ (\beta,{\bf a})=\sigma_1\dot{\sigma_2}\sigma_1\sigma_3\dot{\sigma_2}\sigma_1\sim (\beta',{\bf a'})=\dot{\sigma_2}\sigma_1\dot{\sigma_2}\sigma_3\dot{\sigma_2}\sigma_1\sim  (\beta'',{\bf a''})=\dot{\sigma_2}\sigma_1\sigma_3\dot{\sigma_2}\sigma_3\sigma_1.\]
Note that the corresponding braid words are braid equivalent, differing only via braid moves, and the weights have been propagated according to the tropical propagation rules from \cref{sssec:Lusztigcycles}, cf.~\cref{eq: lusztig trop}. For the first weighted braid word $\bw=\sigma_1\dot{\sigma_2}\sigma_1\sigma_3\dot{\sigma_2}\sigma_1$, the two spherical objects corresponding to the first and second $\dot{\sigma_2}$ are $T_{S_1}^{-1}T_{S_2}^{-1}T_{S_3}^{-1}T_{S_1}^{-1}(S_2)$ and $T_{S_1}^{-1}(S_2))$, respectively. Note that the former is equivalent to $T_{S_3}(S_2)$, and thus the $\pi_\D$-matching paths associated to $T_{S_3}(S_2)$ and $ T_{S_1}^{-1}(S_2))$ have a degree $1$ crossing, but no endpoint intersections. On the one hand, resolving all endpoint intersections yields the direct sum $T_{S_3}(S_2)\oplus  T_{S_1}^{-1}(S_2)$, which is a non-rigid object with a $(\beta,{\bf a})$-filtration. On the other hand, the cofiber of the crossing morphism $T_{S_3}(S_2)[-1]\to T_{S_1}^{-1}(S_2)$ is a rigid but non-matching object with a $(\beta,{\bf a})$-filtration.\\

For the second weighted braid word $(\beta',{\bf a'})$, \cref{eq: lusztig trop} and the construction of $\lbw$ imply that applying the (weighted) braid move $\sigma_1\dot{\sigma_2}\sigma_1\to\dot{\sigma_2}\sigma_1\dot{\sigma_2}$ from $\bw$ has the effect of expressing $T_{S_3}(S_2)$ as $\on{cof}(S_3[-1]\to S_2)$, see e.g.~\cref{ex:R3move_filteredobjects} below. Similarly, for the third weighted braid word $(\beta'',{\bf a''})$, the corresponding weighted braid move $\dot{\sigma_2}\sigma_3\dot{\sigma_2}\to \sigma_3\dot{\sigma_2}\sigma_3$ from $(\beta',{\bf a'})$ has the effect of expressing $T_{S_1}^{-1}T_{S_3}^{-1}(S_2)$ as $\on{cof}(T_{S_1}^{-1}(S_2)[-1]\to S_3)$. The spherical objects associated to the first and second $\dot{\sigma_2}$ in $(\beta'',{\bf a''})$ are
$$T_{S_1}^{-1}T_{S_3}^{-1}T_{S_2}^{-1}T_{S_3}^{-1}T_{S_1}^{-1}(S_2)\simeq S_2,\quad\mbox{ and }\quad T_{S_1}^{-1}T_{S_3}^{-1}(S_2),$$
respectively. The former again has a degree $1$ crossing with the latter, but no endpoint intersection. Therefore $S_2\oplus T_{S_1}^{-1}T_{S_3}^{-1}(S_2)$ is a non-rigid matching $(\beta'',{\bf a''})$-filtered object, whereas $\on{cof}(S_2[-1]\to T_{S_1}^{-1}T_{S_3}^{-1}(S_2))$ is a rigid non-matching $(\beta'',{\bf a''})$-filtered object. Finally, note that the two rigid objects arising from $\bw$ and $(\beta'',{\bf a''})$ are
\[
\on{cof}(T_{S_3}(S_2)[-1]\to T_{S_1}^{-1}(S_2)) \simeq \on{cof}(S_2[-1]\to T_{S_1}^{-1}T_{S_3}^{-1}(S_2)),
\]
which are equivalent. In contrast, the non-rigid objects we obtained from $\bw$ and $(\beta'',{\bf a''})$ are
\[
T_{S_3}(S_2)\oplus  T_{S_1}^{-1}(S_2) \not\simeq S_2\oplus  T_{S_1}^{-1}T_{S_3}^{-1}(S_2)\,,
\]
which are not equivalent. Thus in the non-rigid case, the objects themselves change, and not just their filtrations, if one uses the tropical propagation rules from \cref{eq: lusztig trop}.\qed
\end{example}


\section{Weave propagation of rigid filtered objects}\label{sec:weave_propagation_LC}

The goal of this section is to establish how a rigid $\bw$-filtered object $\lbw$, in the sense of \Cref{def:filtered_by_weighted_braid_word}, changes when either a weighted braid move, or a weighted $s_i^2\to s_i$ transformation, are applied to $\bw$, cf.~\cref{sssec:Lusztigcycles}. Specifically, we prove the equivalence in Statement (\ref{eq:filteredobjects_hexavalent}) and prove the cofiber sequence from \cref{eq:cofiberseq_trivalent0}, in \cref{prop:4_and_6_moves} and \cref{prop:semi_twist_Lusztig_cycle_at_trivalent_vertex} respectively. The former implies that the geometric realization of $\lbw$ is invariant under a weighted braid move, whereas the latter exhibits the homological relation between $\lbw$ and the filtered object associated to the weighted braid word obtained from $\bw$ via a weighted $s_i^2\to s_i$ transformation.\\

Given a weighted braid word $\bw$ such that a rigid $\bw$-filtered object $\lbw$ exists, we denote by
\[L_{\beta,{\bf a}}\coloneqq \mbox{colim }\lbw\]
the geometric realization of $\lbw$, cf.~\cref{rmk:geometric_realization_filtered_object}.


\subsection{Propagation via weighted braid moves}\label{sssec:invariance_of_L_beta_a_under_braid_moves}

The weighted analogue of the braid moves is as follows:

\begin{definition}[Weighted braid moves]\label{def:weighted_braid_moves} Let $\bw$ be a weighted braid word. By definition, a weighted braid move $\bw\to(\b',{\bf a'})$ is a braid move $\beta\to\beta'$ such that:

\begin{enumerate}
    \item If $\beta\to\beta'$ is of the form $s_is_{i+1}s_i\leftrightarrow s_{i+1}s_is_{i+1}$, applied to the $j,(j+1)$ and $(j+2)$ crossings of $\b$, then the weights satisfy
    $$(a'_j,a'_{j+1},a'_{j+2})=\left(a_{j+1}+a_{j+2}-\min(a_{j},a_{j+2}),\min(a_{j},a_{j+2}),a_{j}+a_{j+1}-\min(a_{j},a_{j+2})\right),$$
    and else $a'_k=a_k$ if $k\not\in\{j,j+1,j+2\}$.\\

    \item If $\beta\to\beta'$ is of the form $s_is_{q}\leftrightarrow s_{q}s_i$, $|q-i|\geq2$, applied to the $j,(j+1)$st crossings of $\b$, then the weights satisfy
        $$(a'_j,a'_{j+1})=\left(a_{j+1},a_{j}\right),$$
    and else $a'_k=a_k$ if $k\not\in\{j,j+1\}$.
\end{enumerate}
In short, a weighted braid move is a braid move along with a change of the associated weights according to the tropical propagation rules in \cref{eq: lusztig trop}.\qed
\end{definition}

\noindent By the construction of Lusztig cycles, cf.~\cref{sssec:Lusztigcycles}, the weighted braid words associated to a Lusztig cycle before and after a 4-valent or 6-valent vertex of a weave $\w$ are related by a single weighted braid move as in \cref{def:weighted_braid_moves}. Consequently, the move in \cref{def:weighted_braid_moves}.(1) is referred to as a 6-valent move, or a weighted Reidemeister 3 move, and the move in \cref{def:weighted_braid_moves}.(2) is said to be 4-valent move, or a weighted commutation move.\\

Let us now show that, under the assumption of rigidity, being filtered by a weighted braid word is preserved under weighted braid moves. This provides a categorical interpretation of the tropical propagation rules for Lusztig cycles, as described in \cref{sssec:Lusztigcycles}.

\begin{proposition}\label{prop:4_and_6_moves}
Let $(\beta,{\bf a})$ be weighted braid word and $X\in \D$ a rigid $(\beta,{\bf a})$-filtered object. Suppose that $(\beta,{\bf a})\to(\beta',{\bf a}')$ is a weighted braid move. Then $X$ also admits a $(\beta',{\bf a'})$-filtration.
\end{proposition}

\begin{proof}
Let $\mathcal{L}$ be the $(\beta,{\bf a})$-filtered object with geometric realization $X$. Let $\mathcal{C}$ be the coherent chain complex corresponding to $\mathcal{L}$, as in \cref{sssec:comments_filteredobject}.(B), so that the totalization $\on{tot}(\mathcal{C})\simeq X$ is equivalent to the geometric realization $X$ of $\mathcal{L}$. Let us prove the statement for the two weighted braid moves $\bw\to(\beta',{\bf a'})$ from \cref{def:weighted_braid_moves}.\\

{\bf 6-valent move.} Let us suppose that the braid $\b$ is of the form $\beta=\beta_A (\sigma_{i+1}\sigma_{i}\sigma_{i+1}) \beta_B$, where $\beta_A,\beta_B$ are arbitrary positive braid words, and denote by $\kappa\coloneqq \ell(\beta_A)$ the length of $\beta_A$. Thus the resulting braid $\beta'$ reads $\beta'=\beta_A \sigma_{i}\sigma_{i+1}\sigma_{i} \beta_B$, with its weights ${\bf a'}$ depending on the weights ${\bf a}$ according to \cref{def:weighted_braid_moves}.(1).\\

First, we can reduce to the case where the suffix $\beta_B=1$ is the trivial braid word as follows. For $\beta_B=\sigma_{i_0}\dots\sigma_{i_k}$, we write $T_{\beta_B}=T_{S_{i_k}}\dots T_{S_{i_0}}$. By applying the iterated inverse twist $T_{\beta_B}^{-1}[-\kappa]$ to the coherent chain complex $\C$, we obtain a coherent chain complex of the form
\[
T_{\beta_B}^{-1}(\C)[-\kappa] \simeq  ( \dots \rightarrow S_{i+1}^{\oplus a_{\kappa+2}}[-2] \xrightarrow{\dd_{\kappa+2}} \on{fib}(S_{i}\to S_{i+1}[1])^{\oplus  a_{\kappa+1}}[-1]\xrightarrow{\dd_{\kappa+1}} S_{i}^{\oplus a_{\kappa}} \rightarrow \dots )\,.
\] 
In fact, both depicted differentials vanish, i.e.~$\dd_{\kappa+1}=0$ and $\dd_{\kappa+2}=0$, because of the vanishing of the corresponding $\on{Ext}^1$-groups. It thus suffices to work directly with the arising $3$-term coherent chain complex
\begin{equation}\label{eq:3-term_chain_cplx} S_{i+1}^{\oplus a_{\kappa+2}}[-2] \xrightarrow{0} \on{fib}(S_{i}\to S_{i+1}[1])^{\oplus  a_{\kappa+1}}[-1]\xrightarrow{0} S_{i}^{\oplus a_{\kappa}},\end{equation}
which is effectively setting $\beta_B=1$ and focusing on the crossings where the braid move occurs. Note that the presence of the prefix $\beta_A$ only shifts the subindices by its length $\kappa$. \Cref{lem:connectiveHom,lem:non_rigidity_preserved_under_resolution} together with the rigidity of $X$ imply that the totalization of \eqref{eq:3-term_chain_cplx} is rigid.\\

\noindent Second, our task is to compare the totalization of \eqref{eq:3-term_chain_cplx} and the totalization of the following coherent chain complex
\begin{equation}\label{eq:3-term_chain_cplx2}
S_{i}^{\oplus a_{\kappa+2}'}[-2] \xrightarrow{0}  \on{fib}(S_{i+1}\to S_{i}[1])^{\oplus  a_{\kappa+1}'}[-1]\xrightarrow{0} S_{i+1}^{\oplus a_{\kappa}'}
\end{equation}
given by specifying the morphism
\[ S_{i}^{\oplus a_{\kappa+2}'}[-1]\subset \on{cof}(S_{i}^{\oplus a_{\kappa+2}'}[-2] \xrightarrow{0}  \on{fib}(S_{i+1}\to S_{i}[1])^{\oplus  a_{\kappa+1}'}[-1])\to  S_{i+1}^{\oplus a_{\kappa}'}\]
to be any morphism of maximal rank, when considered as a $a_{\kappa+2}'\times a_{\kappa}'$-matrix via the equivalence $\on{Ext}^1_\D(S_i^{a_{\kappa+2}'},S_{i+1}^{a_\kappa})\simeq k^{a_{\kappa}'\cdot a_{\kappa+2}'}$.\\

\noindent Third, we compute the totalization of \eqref{eq:3-term_chain_cplx} as follows. Since the cofiber of the morphism $S_{i+1}^{\oplus a_{\kappa+2}}[-2] \xrightarrow{0} \on{fib}(S_{i}\to S_{i+1}[1])^{\oplus  a_{\kappa+1}}[-1]$ is equivalent to $S_{i+1}^{\oplus a_{\kappa+2}}[-1] \oplus \on{fib}(S_{i}\to S_{i+1}[1])^{\oplus  a_{\kappa+1}}[-1]$, the totalization of \eqref{eq:3-term_chain_cplx} is thus given by the cofiber of the arising morphism 
\[
S_{i+1}^{\oplus a_{\kappa+2}}[-1] \oplus \on{fib}(S_{i}\to S_{i+1}[1])^{\oplus  a_{\kappa+1}}[-1] \xlongrightarrow{(\alpha,0)} S_{i}^{\oplus a_{\kappa}}\,.
\]\\
We present $\alpha$ as a $(a_{\kappa+2} \times a_{\kappa})$-matrix and observe that the totalization of \eqref{eq:3-term_chain_cplx} can only be rigid only if this matrix is of maximal rank, and so it is. Therefore, the totalization of \eqref{eq:3-term_chain_cplx} is given by the direct sum:
\begin{equation}\label{eq:totalization1}
\textcolor{blue}{\on{fib}(S_{i}\to S_{i+1}[1])^{\oplus  a_{\kappa+1}}}\oplus\textcolor{red}{\on{cof}(S_{i+1}[-1]\to S_{i})^{\oplus \on{min}(a_{\kappa+2},a_{\kappa})}}\oplus\textcolor{cyan}{S_i^{\oplus [a_{\kappa}-a_{\kappa+2}]_+}}\oplus\textcolor{orange}{S_{i+1}^{\oplus[-a_{\kappa}+a_{\kappa+2}]_+}}
\end{equation}

\noindent where we use the notation $[a]_+\coloneqq \max(a,0)$. Note that either $[a_{\kappa}-a_{\kappa+2}]_+=0$ or $[-a_{\kappa}+a_{\kappa+2}]_+=0$, and thus either the third term (in cyan) or the fourth one (in orange) vanish, which is essential for rigidity to hold, as $S_i\oplus S_{i+1}$ is not rigid.\\

Similarly, the totalization of the coherent complex \eqref{eq:3-term_chain_cplx2} reads
\begin{equation}\label{eq:totalization2}
\textcolor{blue}{\on{cof}(S_{i}[-1]\to  S_{i+1})^{\oplus \on{min}(a_{\kappa+2}',a_{\kappa}')}}\oplus\textcolor{red}{\on{fib}(S_{i+1}\to S_{i}[1])^{\oplus  a_{\kappa+1}'}}\oplus \textcolor{cyan}{S_{i}^{\oplus [-a_{\kappa}'+a_{\kappa+2}']_+}}\oplus\textcolor{orange}{S_{i+1}^{\oplus [a_{\kappa}'-a_{\kappa+2}']_+}}
\end{equation}

\noindent It now suffices to compare the summands of \eqref{eq:totalization1} and \eqref{eq:totalization2}. For that, let us use that for this weighted braid move the weights ${\bf a'}$, as in \cref{def:weighted_braid_moves}.(1), indeed satisfy:
\[ \on{min}(a_{\kappa+2}',a_{\kappa}')= a_{\kappa+1}\,,\quad a_{\kappa+1}'=\on{min}(a_{\kappa+2},a_{\kappa})\,,\]
\[ \on{max}(0,a_{\kappa}'-a_{\kappa+2}')=\on{max}(0,-a_{\kappa}+a_{\kappa+2})\,, \]
\[ \on{max}(0,-a_{\kappa}'+a_{\kappa+2}')=\on{max}(0,a_{\kappa}-a_{\kappa+2})\,.\]

\noindent Therefore the third terms of \eqref{eq:totalization1} and \eqref{eq:totalization2}, in cyan, match, as do their fourth terms, in orange. The first $\textcolor{blue}{\on{fib}(S_{i}\to S_{i+1}[1])^{\oplus  a_{\kappa+1}}}$ of \eqref{eq:totalization1}, in blue, is equivalent to
$$\on{fib}(S_{i}\to S_{i+1}[1])^{\oplus  a_{\kappa+1}}\simeq \on{cof}(S_{i}[-1]\to S_{i+1})^{\oplus  a_{\kappa+1}}\simeq \on{cof}(S_{i}[-1]\to S_{i+1})^{\oplus  \on{min}(a_{\kappa+2}',a_{\kappa}')},$$
and thus it is indeed equivalent to the first term of \eqref{eq:totalization2}. Similarly, the second term $\textcolor{red}{\on{fib}(S_{i+1}\to S_{i}[1])^{\oplus  a_{\kappa+1}'}}$ of \eqref{eq:totalization2}, in red, is equivalent to
$$\on{fib}(S_{i+1}\to S_{i}[1])^{\oplus  a_{\kappa+1}'}\simeq \on{cof}(S_{i+1}[-1]\to S_{i})^{\oplus  a_{\kappa+1}'}\simeq \on{cof}(S_{i+1}[-1]\to S_{i})^{\oplus  \on{min}(a_{\kappa+2},a_\kappa)}$$
and hence it is equivalent to the second term of \eqref{eq:totalization1}. In conclusion, the objects \eqref{eq:totalization1} and \eqref{eq:totalization2} are equivalent. By grafting twice, using \Cref{lem:graftcplx}, we conclude that the rigid object $X$ is also $(\beta',{\bf a}')$-filtered.\\

{\bf 4-valent move.} Let us now suppose that $\beta=\beta_A \sigma_i\sigma_j\beta_B$, $|i-j|\geq 2$. The resulting word is $\beta'=\beta_A \sigma_j \sigma_i \beta_B$, with the weights ${\bf a}$ to ${\bf a'}$ changing according to \cref{def:weighted_braid_moves}.(2), i.e.~$a_{\kappa}'=a_{\kappa+1}$ and $a_{\kappa+1}'=a_{\kappa}$ where, as before, $\kappa\coloneqq \ell(\beta_A)$ is the length of $\beta_A$. In this case we must compare the totalizations of the original coherent complex
\begin{equation}\label{eq:totalization3}
\C\coloneqq \qquad \dots \to \C_{\kappa+2}\to \C_{\kappa+1}\xrightarrow{\dd_{\kappa+1}} \C_{\kappa}\to \C_{\kappa-1}\to \dots
\end{equation}
for $\beta$, and the totalization of the, yet to be defined, complex $\C'$ for $\beta'$, which reads
\begin{equation}\label{eq:totalization4}
\C'\coloneqq \qquad \dots \to \C_{\kappa+2}\to \C_{\kappa}[-1]\xrightarrow{\dd'_{\kappa+1}} \C_{\kappa+1}[1]\to \C_{\kappa-1}\to \dots
\end{equation}

The key fact is that $\on{Ext}^1(S_{\sigma_j\beta_B},S_{\sigma_i\sigma_j\beta_B})\simeq \on{Ext}^1(S_{\sigma_j},S_{\sigma_i})\simeq 0$ vanishes, since $|i-j|\geq2$ , and therefore the differential $\dd_{\kappa+1}:\C_{\kappa+1}\lr\C_\kappa$ in the coherent complex $\C$ vanishes, so that we have
\[ S_{\sigma_j\beta_B}^{\oplus a_{\kappa+1}}[-\kappa-1]\simeq \C_{\kappa+1}\xrightarrow{0}  \C_{\kappa}\simeq S_{\sigma_i\sigma_j\beta_B}^{\oplus a_{\kappa}}[-\kappa]\] 
in degrees $\kappa$ and $\kappa+1$ of $\C$. Similarly $\dd'_{\kappa+1}\simeq0$ for $\C'$ in \eqref{eq:totalization4}. The cofibers of $\dd_{\kappa+1}$ and ${\dd'_{\kappa+1}}$ are thus equivalent to the direct sum $\C_{\kappa+1}[1]\oplus \C_{\kappa}$. Therefore, by grafting via \Cref{lem:graftcplx}, the totalization $X$ of $\C$ in \eqref{eq:totalization3} is equivalent to the totalization of the complex
\begin{equation}\label{eq:totalization5}
\dots \to \C_{\kappa+3}[1] \to \C_{\kappa+2}[1]\to \textcolor{blue}{\C_{\kappa+1}[1]\oplus \C_{\kappa}}\to \C_{\kappa-1}\to \dots
\end{equation}
where the grafting part is indicated in blue in \eqref{eq:totalization5}. The chain complex $\C'$ is obtained from \eqref{eq:totalization5} using grafting as well, and hence their totalizations coincide. Thus $\on{tot}(\C')\simeq\on{tot}(\C)\simeq X$. In consequence, $X\simeq \on{tot}(\C')\simeq\mbox{colim}\mathcal{L}_{\b',{\bf a'}}$, which shows that $X$ is also $(\beta',{\bf a})$-filtered, as required.
\end{proof}

\begin{corollary}[Invariance of $\Lbw$ under weighted braid moves]\label{cor:4_and_6_moves} Let $\w$ be a Demazure weave and $\g_i$ a Lusztig cycle. Suppose that $(\beta,{\bf a})$ and $(\beta',{\bf a}')$ are two weighted braid words obtained from two different horizontal slices of $\g_i$ such that there are no trivalent vertices of $\w$ between them. Then there exists an equivalence $L_{\beta,{\bf a}}\simeq L_{\beta',{\bf a}'}$ in $\D$.
\end{corollary}

\begin{proof}
By \Cref{thm:unique_rigid_resolutions}, $(\beta,{\bf a})$-filtered objects are unique. By hypothesis, $(\beta,{\bf a})$ and $(\beta',{\bf a}')$ are related by weighted braid moves, and thus \Cref{thm:unique_rigid_resolutions} and \cref{prop:4_and_6_moves} imply that $L_{\beta,{\bf a}}\simeq L_{\beta',{\bf a}'}$.
\end{proof}


\begin{example}\label{ex:R3move_filteredobjects} Let us use the dot notation for weighted braid words with binary weights $0$ or $1$, as in \cref{ex:6move_consistency}, so $s_1\dot{s}_2s_1$ denotes $\bw=(s_1s_2s_1,(0,1,0))$ and $\dot{s}_2s_1\dot{s}_2$ denotes $\bw=(s_2s_1s_2,(1,0,1))$.\\

\noindent $(i)$ Consider the weighted braid word $\bw=s_1\dot{s}_2s_1$, so that $\beta_0=s_1s_2s_1$, $\beta_1=s_2s_1$ and $\beta_2=s_1$. Then the associated filtered object has $\lbw(0)=0$, as $a_0=0$, $\lbw(1)=S_{\b_1}=T_{S_1}^{-1}(S_2)$, since $a_1=1$, and $\lbw(i)\simeq\lbw(1)$ for all $i\in[1,\infty)$, as $a_2=0$. Thus the filtration in \cref{thm:unique_rigid_resolutions}, cf.~\cref{eq:bw_filtration1}, reads

\begin{equation}\label{eq:example1_filtration}
\{0\}\lr T_{S_1}^{-1}(S_2)\stackrel{\mbox{id}}{\lr} T_{S_1}^{-1}(S_2).    
\end{equation}

\noindent From the viewpoint of \cref{sssec:comments_filteredobject}.(B), the corresponding coherent chain complex $\cbw$ is given by
\begin{equation}\label{eq:example1_complex}
\C_{s_1\dot{s}_2s_1}=\left(0\stackrel{d_2}{\lr} T_{S_1}^{-1}(S_2)[-1]\stackrel{d_1}{\lr} 0\right),
\end{equation}
where the only non-zero entry is in homological degree 1. Note that the cofiber totalization of $\C_{s_1\dot{s}_2s_1}$ in \cref{eq:example1_complex} is equivalent to the cofiber of the zero map $d_1:T_{S_1}^{-1}(S_2)[-1]\lr 0$, which indeed coincides with $L_{s_1\dot{s}_2s_1}=T_{S_1}^{-1}(S_2)$ from \cref{eq:example1_filtration}.\\

\noindent $(ii)$ Consider the weighted braid word $\bw=\dot{s}_2s_1\dot{s}_2$, so that $\beta_0=s_2s_1s_2$, $\beta_1=s_1s_2$ and $\beta_2=s_2$. Then $S_{\beta_0}= T_{S_2}^{-1}T_{S_1}^{-1}S_2\simeq S_1$ and $S_{\beta_2}=S_2$ and the associated filtered object reads
\begin{equation}\label{eq:example2_filtration}
S_1\stackrel{\on{id}}{\lr} S_1\stackrel{\theta}{\lr} T_{S_1}^{-1}(S_2),
\end{equation}
where $\theta:S_1\lr T_{S_1}^{-1}(S_2)$ is identified via $T_{S_1}^{-1}(S_2)\simeq T_{S_2}S_1$ with the canonical map $S_1\lr T_{S_2}S_1$ to the cofiber of the evaluation map from $S_2$ to $S_1$. The morphism also $\theta$ corresponds to an endpoint intersection.\\

\noindent In this case, the coherent chain complex $\cbw$ is given by

\begin{equation}\label{eq:example2_complex}
\C_{\dot{s}_2s_1\dot{s}_2}  =
\begin{tikzcd}
{S_2}[-2] \arrow[r, "0"] \arrow[rr, "\iota\neq 0", bend left] & 0 \arrow[r, "0"] & {S_1}\,,
\end{tikzcd}
\end{equation}
\noindent where $\iota:S_2[-1]\lr S_1$ is a non-zero morphism in $\on{Mor}_\D(S_2,S_1)\simeq k[-1]$. The cofiber totalization of $\C_{\dot{s}_2s_1\dot{s}_2}$ in \cref{eq:example2_complex} is
$$\cof(S_2[-1]\lr\cof(0\lr S_1))\simeq \cof(S_2[-1]\lr S_1)\simeq \cof(\Mor_\D(S_2,S_1)\otimes S_2\lr S_1)\simeq T_{S_2}S_1.$$
Since $T_{S_2}S_1\simeq T_{S_1}^{-1}(S_2)$, this coincides with $L_{\dot{s}_2s_1\dot{s}_2}=T_{S_1}^{-1}(S_2)$ from \cref{eq:example2_filtration}.\\

\noindent $(iii)$ From $(i)$ and $(ii)$ above, we directly deduce that $L_{s_1\dot{s}_2s_1}\simeq L_{\dot{s}_2s_1\dot{s}_2}$, as they are both equivalent to $T_{S_1}^{-1}(S_2)$. This is indeed in accordance with \cref{prop:4_and_6_moves} and \cref{cor:4_and_6_moves}. Note that the induced filtrations are nevertheless different, cf.~filtrations (\ref{eq:example1_filtration}) and (\ref{eq:example2_filtration}).\qed
\end{example}


\subsection{Propagation across trivalent vertices}\label{sssec:categoricalLusztigcycles_definition}

The goal of this subsection is to prove the cofiber sequence \eqref{eq:cofiberseq_trivalent0}, which we establish in \cref{prop:semi_twist_Lusztig_cycle_at_trivalent_vertex}. Specifically, given a weighted braid word $\bw$ with $\Lbw$ rigid, we want to relate $\Lbw$ and $L_{\b',{\bf a'}}$ where $(\b',{\bf a'})$ is obtained from $\bw$ via the $3$-valent move, meaning the weighted braid operation $s_{i}^2\to s_i$ applied to two consecutive and equal crossings of $\b$, where the weights change according to \cref{eq: lusztig trop}. That is, all weights remain equal except the weight of the new crossing, which must be the minimum of the two weights associated to the two crossings to which we applied the operation.\\

For context, we recall from \cref{sec:lusztigcycles} that the cofiber sequence \eqref{eq:cofiberseq_trivalent0} is a central ingredient in the construction of the categorical Lusztig cycles. In particular, by \Cref{rem:lift_amounts_to_relative_triangles}, it supplies a collection of relative triangles which are used in \cref{def:categorical_Lusztig_cycles} to specify the required global sections of the weave schober. The statement we need is as follows:

\begin{proposition}[Filtered objects across trivalents]\label{prop:semi_twist_Lusztig_cycle_at_trivalent_vertex}
Let $\w$ be a Demazure weave, $p\in\w$ a trivalent vertex and $\g$ the Lusztig cycle of a trivalent vertex above $p$. Consider the weighted braid words $(\beta,{\bf a})$ and $(\beta',{\bf a'})$ associated to $\g$ from horizontal slices right above and right below $p$, respectively, and suppose that a rigd $\bw$-filtered object $\mathcal{L}_{\beta,{\bf a}}$ exists. Then a rigid $(\beta',{\bf a'})$-filtered object $\mathcal{L}_{\beta',{\bf a'}}$ exists and there exists a cofiber sequence in $\D$
\begin{equation}\label{eq:cofiberseq_trivalent}
\tau_{\geq 0}\on{Mor}(S_{p},L_{\beta,{\bf a}})\otimes S_{p} \stackrel{ev}{\longrightarrow} L_{\beta,{\bf a}}\longrightarrow L_{\beta',{\bf a'}}\,
\end{equation}
where $\tau_{\geq 0}\on{Mor}(S_{p},L_{\beta,{\bf a}})$ denotes the maximal summand of the morphism object contained in homologically positive degrees and the left morphism $ev$ indicates the corresponding evaluation morphism.
\end{proposition}

\begin{proof} Let us use the coherent chain complex $\C_{\beta,{\bf a}}$ in $\D$ associated to the filtered object $\lbw$, cf.~\cref{sssec:comments_filteredobject}.(B). To ease notation we set $\C\coloneqq \C_{\beta,{\bf a}}$, which reads
\begin{equation}\label{eq:trivalent_totalization1}
\begin{tikzcd}[column sep=tiny]
\C_\ast\coloneqq \qquad \dots \arrow[r] & {\C_{\k+1}} \arrow[r] & { S_p^{\oplus a_\k}[-\k]} \arrow[r]                            & {S_p^{\oplus a_{\k-1}}[2-\k]} \arrow[r] & {\C_{\k-2}} \arrow[r]                                  & \dots
\end{tikzcd}
\end{equation}
where $\kappa-1$ and $\kappa$ denote the indices of the two equal crossings $s_{i_{\k-1}}=s_{i_\k}$ of $\beta$ to which the operation $s_{i_\k}^2\to s_{i_\k}$ is being applied. Note that in \eqref{eq:trivalent_totalization1} we have used that the $\kappa-1$ and $\kappa$ crossings are equal to express the degree $(\k-1)$ term of $\C_\ast$ as
$$S_{\beta_{\k-1}}[-(\k-1)]\simeq (T_{s_{i_\k}}^{-1} S_{\b_\k})([-(\k-1)])\simeq S_{\b_\k}[2-\k] \simeq S_p[2-\k],$$
since we have $S_p\simeq S_{\b_\k}$ by construction.\\

First, we consider a coherent chain complex $\mathcal{E}_\ast$ endowed with a coherent chain map $\mathcal{E}_\ast\to \C_\ast$, as follows. In degrees $0$ to $(\k-1)$, we define
$$\mathcal{E}_{\leq \k-1}\coloneqq \tau_{\geq 0}\on{Mor}(S_p[2-\k],\C_{\leq \k-1})\otimes S_p[2-\k].$$

\noindent In degree $\k$, we set $\mathcal{E}_{\k}\coloneqq S_p^{\oplus a_\k-\on{min}(a_\k,a_{\k-1})}[-\k]$ and use the zero morphism $\mathcal{E}_{\k}\to \on{tot}(\mathcal{E}_{\leq \k-1})$ to extend $\mathcal{E}_{\leq \k-1}$ to $\mathcal{E}_{\leq \k}$. In degrees strictly larger than $\k$ we declare $\mathcal{E}_{j}\simeq 0$, for all $j\in[\k+1,\infty)$. By using evaluation in degrees $0$ to $\k-1$, and a monomorphism in degree $\k$, there is a canonical chain map
\begin{equation}\label{eq:mapEtoC}
f:\mathcal{E}_\ast \to \C_\ast.
\end{equation}
By construction, the chain map $f_{\leq\k-1}:\mathcal{E}_{\leq \k-1}\to \C_{\leq \k-1}$ induces after totalization the evaluation morphism
\begin{equation}\label{eq:fk_totalization_evaluation}
\on{tot}(f_{\leq\k-1})\simeq ev:\tau_{\geq 0}\on{Mor}(S_p,\on{tot}(\C_{\leq \k-1}))\otimes S_p\to \on{tot}(\C_{\leq \k-1}).
\end{equation}
\noindent Now, we claim that $\SE_\ast$ is equivalent to the coherent chain complex
\begin{equation}\label{eq:E_equivalent}
\E_\ast\simeq\quad\dots \lr 0 \lr  S_p^{\oplus a_b-\on{min}(a_\k,a_{\k-1})}[-\k] \stackrel{0}{\lr} {S_p^{\oplus a_{\k-1}}[2-\k]} \lr \on{Mor}(S_p,\C_{\k-2})\otimes S_p \lr \dots
\end{equation}
\noindent Indeed, \Cref{lem:connectiveHom} implies that the morphism object $\on{Mor}(S_p[2-\k],\C_j)$ is connective for all $j<\k-1$ from which \eqref{eq:E_equivalent} follows.\\

Second, let us prove that $\on{tot}(\mathcal{E}_\ast)\simeq \tau_{\geq 0}\on{Mor}(S_p,L_{\beta,{\bf a}})\otimes S_p$, as follows. \Cref{lem:connectiveHom} implies that $\on{Mor}(S_p,S_{\beta_j})$ is concentrated in strictly negative degree for all $j\in(\k,\infty)$, and thus we have an equivalence
\begin{equation}\label{eq:equivalence_truncation}
\tau_{\geq 0}\on{Mor}(S_p,L_{\beta,{\bf a}})\simeq \tau_{\geq 0} \on{Mor}(S_p,\on{tot}(\C_{\leq \k})).
\end{equation}
By the construction of $\C$, cf.~\eqref{eq:trivalent_totalization1}, we have a cofiber sequence
\begin{equation}\label{eq:cofiberseq_totD1}
\on{tot}(\C_{\leq k-1}) \to \on{tot}(\C_{\leq \k})\to S_p^{\oplus a_{\k}}.
\end{equation}
Applying the functor $\on{Mor}(S_p,\mhyphen)$ to \eqref{eq:cofiberseq_totD1} leads to the cofiber sequence
\begin{equation}\label{eq:cofiberseq_totD2}
\on{Mor}(S_p,\on{tot}(\C_{\leq \k-1}))\to \on{Mor}(S_p,\on{tot}(\C_{\leq \k})) \to \underbrace{\on{Mor}(S_p,S_p)^{\oplus a_\k}}_{\simeq (k\oplus k[-2])^{\oplus a_\k}}\,.
\end{equation}
\noindent By \Cref{lem:connectiveHom}, and the fact that the differential $\dd_{\k+1}:S_p^{\oplus a_\k}[-\k]\to S_p^{\oplus a_{\k-1}}[2-\k]$ in \eqref{eq:trivalent_totalization1} has maximal rank, since $L_{\beta,{\bf a}}$ is rigid, we obtain that
\[ \on{Ext}^1(S_p,\on{tot}(\C_{\leq \k-1}))\simeq \on{Ext}^1(S_p,S_p^{\oplus a_{\k-1}}[1])\simeq k^{\oplus a_{\k-1}}\]
and that the connecting homomorphism
$$k^{\oplus a_\k}\to \on{Ext}^1(S_p,\on{tot}(\C_{\leq \k-1}))\simeq k^{\oplus a_{\k-1}}$$
for \eqref{eq:cofiberseq_totD2} is of maximal rank, i.e.~of rank $\on{min}(a_\k,a_{\k-1})$. In consequence, we obtain from \eqref{eq:cofiberseq_totD2} the (split) cofiber sequence
\begin{equation}\label{eq:totD_cofiberseq}
\tau_{\geq 0}\on{Mor}(S_p,\on{tot}(\C_{\leq \k-1}))\to \tau_{\geq 0}\on{Mor}(S_p,\on{tot}(\C_{\leq \k})) \to k^{\oplus a_\k-\on{min}(a_\k,a_{\k_1})}\,.
\end{equation}
Together with the equivalence \eqref{eq:equivalence_truncation}, \eqref{eq:totD_cofiberseq} implies that
\begin{equation}\label{eq:totD_truncatedMor}
\on{tot}(\mathcal{E}_\ast) \simeq \tau_{\geq 0}\on{Mor}(S_p,L_{\beta,{\bf a}})\otimes S_p\,.
\end{equation}

Third, by considering the cofiber of the map $f:\E_\ast\lr\C_\ast$ in \eqref{eq:mapEtoC}, we obtain the following diagram:
\[
\begin{tikzcd}[column sep=tiny]
\color{blue}\E_\ast\simeq \arrow[d,"f"] & \color{blue}\dots \arrow[r] & \color{blue}0 \arrow[d] \arrow[r]                                        & \color{blue}{ S_p^{\oplus a_\k-\on{min}(a_\k,a_{\k-1})}[-\k]} \arrow[d] \arrow[r, "0"] & \color{blue}{S_p^{\oplus a_{\k-1}}[2-\k]} \arrow[d] \arrow[r] & \color{blue}{\on{Mor}(S_p,\C_{\k-2})\otimes S_p} \arrow[d] \arrow[r] & \color{blue}\dots \\
\color{orange}\C_\ast\simeq \arrow[d] & \color{orange}\dots \arrow[r] &\color{orange} {\C_{\k+1}} \arrow[d, "\on{id}"'] \arrow[r] &\color{orange} { S_p^{\oplus a_\k}[-\k]} \arrow[d] \arrow[r]                            &\color{orange} {S_p^{\oplus a_{\k-1}}[2-\k]} \arrow[d] \arrow[r] &\color{orange} {\C_{\k-2}} \arrow[d] \arrow[r]                                  & \color{orange}\dots \\
\color{purple}\on{cof}(f)\simeq & \color{purple}\dots \arrow[r] & \color{purple}{\C_{\k+1}} \arrow[r]                       & \color{purple}{ S_p^{\oplus \on{min}(a_\k,a_{\k-1})}[-\k]} \arrow[r]         & \color{purple}0 \arrow[r]                                           & \color{purple}{T_{S_p}\C_{\k-2}} \arrow[r]                                 & \color{purple}\dots
\end{tikzcd}
\]
By \eqref{eq:fk_totalization_evaluation} and the equivalence \eqref{eq:totD_truncatedMor}, the totalization $\on{tot}(f):\on{tot}(\E_\ast)\lr\on{tot}(\C_\ast)$ is equivalent to the evaluation map
$$ev:\tau_{\geq 0}\on{Mor}(S_p,L_{\beta,{\bf a}})\otimes S_p\to L_{\beta,{\bf a}}.$$
In order to conclude the statement of the proposition, it therefore suffices to show that the totalization $\on{tot}{\on{cof}(f)}$ of the bottom row $\cof(f)$ of the diagram above, in purple, is indeed equivalent to $L_{\b',{\bf a'}}$. Indeed, it follows from grafting, see \Cref{lem:graftcplx}, that there is a coherent chain complex
\[
\C_{\beta',{\bf a}'}\simeq (\dots \to \C_{\k+1}[1] \to S_{p}^{\on{min}(a_\k,a_{\k-1})}[-\k+1]\to T_{S_p}\C_{\k-2}\to \dots),
\]
whose totalization equivalent to $\on{tot}(\cof(f))$. Thus
\[\on{cof}(ev)\simeq\on{tot}(\cof(f))\simeq \on{tot}(\C_{\beta',{\bf a}'})\,.\]
The desired filtered object $\mathcal{L}_{\beta',{\bf a'}}$ is the one corresponding to the coherent chain complex $\C_{\beta',{\bf a}'}$.
\end{proof}

\begin{remark}
In terms of semi-twists, cf.~\Cref{def:semitwists}, \Cref{prop:semi_twist_Lusztig_cycle_at_trivalent_vertex} shows that
\begin{equation}\label{eq:semitwist_trivalent}
L_{\beta',{\bf a'}}\simeq T^{+,\geq 0}_{S_{p}}(L_{\beta,{\bf a}})
\end{equation}
under the given hypotheses. By \cref{lem:semi_twists_preserve_rigidity}, semi-twists preserve rigidity and thus $L_{\beta,{\bf a}}$ being rigid implies that $L_{\beta',{\bf a'}}$ must be rigid as well.\qed
\end{remark}

\begin{example}\label{ex:cofiberseq_trivalent} Let us provide a key instance of the cofiber sequence (\ref{eq:cofiberseq_trivalent}). For that, we consider the Demazure weave $\w:s_1^2\lr s_1$ consisting of a trivalent vertex $p
\in\w$. Consider an arbitrary weighted braid word $\bw$ above the trivalent vertex, so that $\beta=s_1^2$ and $(a_0,a_1)=(b,a)$ for some $a,b\in\Z_{\geq0}$. Since $S_p\simeq S$, $S_{\beta_0}=S[1]$ and $S_{\beta_1}=S$, and the associated filtered object $\lbw$ is
\begin{equation}\label{eq:cofiberseq_trivalent_above}
\lbw(0)=S[1]^a,\quad L_{\beta,{\bf a}}\simeq \lbw(1)=\cof(S[-1]^b\stackrel{f}{\lr} S[1]^a)
\end{equation}
\noindent where $f:S[-1]^b\lr S[1]^a$ can be described by a linear map $M_f\in \on{Mat}_{a\times b}(k)$ tensored with $S$, since $\mor(S,S)\simeq k\oplus k[-2]$ as $S$ is 2-spherical. Here we used the notation $S^a\coloneqq S^{\oplus a}$ for ease. Due to its
construction in terms of endpoint resolutions, $M_f$ is of maximal rank $\on{rk}(M_f)=\min(a,b)$. By using the cofiber sequence describing $L_{\beta,{\bf a}}$ in \cref{eq:cofiberseq_trivalent_above}, we obtain
$$\mor(S,L_{\ba})\simeq k[-2]^b\oplus k[-1]^{a-\min(a,b)}\oplus k^{b-\min(a,b)}\oplus k[1]^a$$
and thus
\begin{equation}
\tau_{\geq0}\mor(S,L_{\beta,{\bf a}})\otimes S\simeq S^{b-\min(a,b)}\oplus S[1]^a.
\end{equation}

\noindent Consider the weighted braid word $(\b',a')=(s_1,\min(a,b))$ below the trivalent vertex, as specified by \cref{eq:3_move}. The associated filtered object is $L_{\b',a'}\simeq S^{\min(a,b)}$. The cofiber sequence in \cref{eq:cofiberseq_trivalent} therefore reads
\begin{equation}\label{eq:cofiberseq_trivalent_example}
S^{b-\min(a,b)}\oplus S[1]^a \stackrel{g}{\lr} \cof(S[-1]^b\stackrel{f}{\lr} S[1]^a)\stackrel{h}{\lr} S^{\min(a,b)}\,.
\end{equation}
Specifically, if $b=\min(a,b)$, then the cofiber sequence (\ref{eq:cofiberseq_trivalent}) coincides with the cofiber sequence
$$S[1]^a \stackrel{g}{\lr} \cof(S[-1]^b\stackrel{f}{\lr} S[1]^a)\stackrel{h}{\lr} S^{b}\,
$$
of the map $f$, with $h=\delta$ being the natural connecting morphism $\delta:\cof(S[-1]^b\stackrel{f}{\lr} S[1]^a)\lr S^{b}$. Since $\delta$ is a connecting morphism, we have $\fib(h)=S[1]^a$, and this is indeed a cofiber sequence. If $a=\min(a,b)$, then the cofiber sequence (\ref{eq:cofiberseq_trivalent}) reads 

$$S^{b-\min(a,b)}\oplus S[1]^a \stackrel{g}{\lr} \cof(S[-1]^b\stackrel{f}{\lr} S[1]^a)\stackrel{h}{\lr} S^{\min(a,b)}\,,$$
with $h=(p\otimes S)\circ\delta$ being the composition of the connecting morphism $\delta$ and any surjective map $p:k^b\lr k^a$ tensored with $S$. In this case, $\fib(h)$ is obtained as the trivial extension of $\fib(p\otimes S)\simeq \ker(p)\otimes S\simeq S^{b-\min(a,b)}$ by $\fib(\delta)\simeq S[1]^a$, and the above is indeed a cofiber sequence.\qed
\end{example}

\begin{corollary}\label{thm:unique_Lusztig_filtered_object}
Let $(\beta,{\bf a})$ be a weighted braid word arising from a horizontal slice of a Lusztig cycle of a Demazure weave. Then, up to equivalence, there exists a unique $(\beta,{\bf a})$-filtered object $\lbw \in \D$ whose geometric realization $\Lbw$ is rigid. In addition, $\Lbw\in\mathcal{M}$.
\end{corollary}

\begin{proof}
Let $p\in\w$ be the vertex associated to the Lusztig cycle $\g=\g_p$, so that $\g$ is born at $p$ and propagates downwards, cf.~\cref{def:Lusztig_cycles}. By construction, the weights of the Lusztig cycle above $p$ vanish, and therefore the statement trivially holds for horizontal slices above $p$. For the horizontal slice right below $p$ and the Lusztig cycle $\g_p$, the associated weighted braid word $(\beta_p,{\bf a}_p)$ has a unique non-zero weight $a_j$, for some $j\in[0,\ell(\beta)]$, and this weight has value $a_j=1$. Thus $L_{\beta_p,{\bf a}_p}\simeq S^{\oplus a_j}_{\beta_j}\simeq S_{\beta_j}$, which is indeed rigid and belongs to $\mathcal{M}$, since $S_{\beta_j}\in\mathcal{M}$ by \cref{lem:twist_via_homeomorphism}. The object $L_{\beta_p,{\bf a}_p}$ is trivially $(\beta,{\bf a})$-filtered and \cref{thm:unique_rigid_resolutions} implies it is unique. The statement thus holds for the horizontal slice right below $p$ by taking $\mathcal{L}_{\beta_p,{\bf a}_p}$ with the trivial filtration and colimit $L_{\beta_p,{\bf a}_p}$.\\

Let us now argue that for any horizontal slice $h$ of the Demazure weave $\w$ below $p$. By construction, the weighted braid word $(\b_h,{\bf a}_h)$ associated to $\g_p$ and such horizontal slice $h$ is obtained from $(\b_p,{\bf a}_p)$ by a combination of weighted braid moves and (weighted) trivalent moves $s_i^2\to s_i$, typically intertwined. For the weighted braid moves, \cref{prop:4_and_6_moves} and \cref{cor:4_and_6_moves} imply the required statement. For the weighted trivalent moves, i.e.~propagating downwards through a 3-valent vertex, \cref{prop:semi_twist_Lusztig_cycle_at_trivalent_vertex} shows that the resulting filtered object below exists and it is obtained by applying a positive semi-twist to the filtered object above, cf.~equivalence \eqref{eq:semitwist_trivalent}. \cref{lem:semi_twists_preserve_rigidity} shows that rigidity is therefore preserved when scanning downward past a trivalent vertex, and \Cref{lem:matching_curves_stable_under_semitwist} proves that being matching is also preserved, i.e.~the result filtered object has colimit in $\mathcal{M}$. This establishes the existence of the required $(\beta,{\bf a})$-filtered rigid object $\lbw \in \D$, with geometric realization in $\mathcal{M}$. Its uniqueness follows from \cref{thm:unique_rigid_resolutions}.
\end{proof}

This concludes the proofs for all the ingredients needed in the construction of the categorical Lusztig cycles in \cref{sec:lusztigcycles}. In particular, the results thus far prove Theorem \ref{thm:main1}. For reference, we record the proof as follows.

\begin{proof}[Proof of Theorem \ref{thm:main1}] Part (1) is stated and proven in \cref{thm:unique_Lusztig_filtered_object}. The proof for Part (2) follows from combining \cref{cor:4_and_6_moves} and \cref{prop:semi_twist_Lusztig_cycle_at_trivalent_vertex}, cf.~also \cref{rmk:weightedbraidword_from_filtration}. Part (3) is obtained from the construction of the categorical Lusztig cycles in \cref{sssec:categoricalLusztigcycles_definition1}, which is detailed and shown to be well-defined in the results of \cref{sec:filtrations_from_braid_words} and this section.
\end{proof}


\section{Matching sections for weave schobers and weave thimbles}\label{sec:weave_thimbles}

The goal of this section is to construct a full exceptional collection $\stds^\w\coloneqq \{\st1^\w,\ldots,\Delta_m^\w\}$ in the $k$-linear stable $\infty$-category of global sections $\glsecF$, with $\F_\w$ the weave schober of a Demazure weave $\w$. The exceptional collection $\stds^\w$ will be used to prove homological properties of the categorical Lusztig cycles $\{\LC_i\}$. In particular, we use them to prove Parts (2), (3) and (4) of Theorem \ref{thm:main2} and Theorem \ref{thm:main3}.\\

It is important to emphasize that the exceptional collection $\stds^\w$ is {\it not} invariant under weave equivalences: there exist weave equivalences $\w\simeq\w'$ such that $\stds^\w\neq\stds^{\w'}$. Thus these exceptional collections $\stds^\w$ are not intrinsically well-defined within the context of weave calculus. That said, the categorial Lusztig cycles $\{\LC_i\}$ can be expressed in terms of the collection $\stds^\w$ in a manner that $\{\LC_i\}$ are invariant under weave equivalences applied to $\w$.\\

\noindent This section is organized as follows:

\begin{enumerate}
    \item \cref{subsec:matching} provides a general construction, in line with \cref{sssec:2spherical_matching_objects}, to obtain global sections in $\glsecF$ from a certain type of matching paths on the plane of the weave $\w$. Note that, in contrast with \cref{ssec:matching_2spheres}, the matching paths in \cref{subsec:matching} are drawn on the plane where the weave $\w$ and the graph $\Wgraph$ are themselves drawn, and not in a different plane and with a different graph $\rgraph_\D$ as it occurred in \cref{ssec:matching_2spheres}. The resulting objects in this section belong to $\glsecF$ and are, in a sense, 3-dimensional, cf.~\cref{lem:machting_spheres_and_thimbles}; the objects built in \cref{sssec:2spherical_matching_objects} belonged to $\D$ and were, in a sense, 2-dimensional, cf.~\cref{sssec:morphisms_2spherical_matching_objects}.

    \item \cref{ssec:weave_transport} provides a specific type of paths associated to a weave and discusses the notion and properties of weave parallel transport. The key geometric idea is to think of the weave $\w$ as a set of branch cuts. Even if it is a highly non-minimal branch cut, its rigid algebraic structure leads to remarkable relations between the different objects associated to the different types of weave paths.

    \item \cref{ssec:weave_thimbles} uses the general construction from \cref{subsec:matching} applied to the weave paths in \cref{ssec:weave_transport} to define the weave thimbles $\stds^\w$ of a Demazure weave $\w$, and their duals. A significant observation will be that these weave thimbles identify with the inductions of the vanishing cycles in the sense of \cite{Chr25b}. \cref{ssec:weavethimbles_exceptional} proves that the weave thimbles $\stds^\w$ form an exceptional collection, and establishes their behavior under weave equivalences and weave mutations.
\end{enumerate}


\subsection{Global sections of \texorpdfstring{$\F_\w$}{weave schobers} from matching paths}\label{subsec:matching}

Let $\w$ be a Demazure weave and $p_1,\dots,p_m$ be its ordered trivalent vertices, ordered from bottom to top. Following \cref{ssec:defining_weave_schobers}, we have the associated ribbon graph $\Wgraph$ with vertices $v_1,\dots,v_m$ and edges $e_,\ldots,e_m$, see \cref{fig:Weave_GraphForWeaveSchober}. As before, we use $\vinf\in e_1$ to denote a generic point in the unique non-compact edge $e_1$. To ease notation, we will often denote the trivalent vertices $p_1,\dots,p_m$ of the weave $\w$ by $v_1,\dots,v_m$, identifying $p_i=v_i$ if the context allows. Note that the graph $\Wgraph$ is always drawn vertically, in contrast to $\rgraph_\D$ in \cref{ssec:matching_2spheres}, which was drawn horizontally, cf.~\cref{sssec:2spherical_matching_objects}. Both $\w$ and $\Wgraph$ are drawn in a plane, equivalently a disk, that we denote by $\mathbb{D}$.

\begin{definition}[$\w$-matching paths]\label{def:matchingpath}
Let $\w$ be a Demazure weave with trivalent vertices $p_1,\dots,p_m$. By definition, a $\w$-matching path $\gamma$ in $\mathbb{D}$ is an embedded curve $\gamma\colon [0,1]\to \mathbb{D}$, such that:
\begin{enumerate}[$(1)$]
\item $\gamma(0)=v_i$ and $\gamma(1)=v_k$, for some $k\in \{-\infty,1,\dots,m-1\}$ and $i\in(k,m]$.
\item Near $\gamma(0)$ the image of $\g$ must coincide with either \cref{fig:Weave_ParallelTransportModels}.(1) or \cref{fig:Weave_ParallelTransportModels}.(2). That is, either point straight to the right, or point immediately to the left after having started with a velocity vector pointing upwards.
\item $\gamma$ goes only downwards after its start, i.e.~with the exception of the local model near $\gamma(0)$ in the previous item, the path $\gamma$ goes downwards.
\item The interior of its image $\gamma(0,1)$ is disjoint from $\{v_i\}_{1\leq i\leq m}$,
\item If $j\not =-\infty$, there must exist an equivalence
\[S_{b}\simeq \mathcal{F}^\rightarrow(\gamma|_{(0,1)})(S_{a})[l]\]
in $\D$ for some $l\in \mathbb{Z}$, where $a,b$ are the indices such that the weave edges incident to $p_i$ are labeled by $a$, i.e.~by the permutation $s_a$, and the weave edges incident to $p_k$ are labeled by $b$. Here we also denoted $\mathcal{F}^\rightarrow(\gamma|_{(0,1)})\colon \mathcal{F}_\w(e_i)\simeq \mathcal{F}_\w(e_{k+1})$ as the transport of $\mathcal{F}$ along $\gamma$, cf.~\cite[Section 4.3]{Chr25b}. 
\end{enumerate}
Such $\w$-matching paths are considered up to ambient compactly supported smooth isotopy of $\mathbb{D}$ relative to $\{v_{-\infty},v_1,\dots,v_n,\}$.\qed
\end{definition}

\begin{center}
	\begin{figure}[h!]
		\centering
		\includegraphics[scale=1.1]{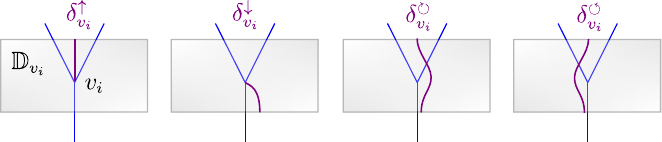}
		\caption{Four types of $\w$-matching segments, cf.~\cref{def:mathcing_segments} and \cref{rem:segments}.}
        \label{fig:WeaveVertex_Segments}
	\end{figure}
\end{center}

\noindent Note that \cref{def:matchingpath} is purposefully ambiguous given the abuse of notation $p_i=v_i$ identifying trivalents of $\w$ with vertices of $\Wgraph$. Indeed, it could either suggest that we are drawing a path $\g$ with endpoints on the trivalent vertices of $\w$, or with endpoints on the vertices of $\Wgraph$. Both viewpoints are valuable: the former is more intuitively geometric and matches the classical theory of matching paths, where the endpoints are branch points, whereas the latter is technically more accurate in the context of building sections for the weave schober $\F_\w$, which is actually defined on $\Wgraph$. \cref{fig:WeaveVertex_Segments} is an instance of the former viewpoint, drawing endpoints at the trivalents of weave $\w$ itself. 

\begin{remark} The condition \cref{def:matchingpath}.(3) exists to fix the correct shift: e.g.~the different angles in which a path leaves a trivalent of $\w$ only affect the resulting object by a shift. Also, in contrast to \cref{sssec:morphisms_2spherical_matching_objects}, we do not need to consider gradings for the $\w$-matching paths in \cref{def:matchingpath}.  The condition \cref{def:matchingpath}.(5) is an algebraic version of the requirement that parallel transport along $\gamma$ connects the two vanishing spheres at the endpoints of $\gamma$, see \cref{subsubsec:glsecfrommatchingpath}.(iii).\qed
\end{remark}

\begin{definition}\label{def:mathcing_segments}
Given $i\in[1,m]$, let $\mathbb{D}_{v_i}\subset \mathbb{D}$ be a closed horizontal band whose upper boundary, resp.~lower boundary, coincides with a horizontal slices immediately above, resp.~below, the vertex $v_i$. 
\begin{enumerate}[(1)]
\item For $i\in[1,m-1]$, a $\w$-matching segment in $\mathbb{D}_{v_i}$ is defined to be an embedded curve $\gamma\colon [0,1]\to \mathbb{D}_{v_i}$ such that 
\begin{itemize}
\item $\gamma$ is smoothly isotopic to a path that only goes downwards, relative to its endpoints.
\item $\gamma$ does not intersect $\{v_i\}$ except for at most at one endpoint.
\item All endpoints of $\gamma$ lie in $(\partial\mathbb{D}_{v_i}\cap \mathbb{D}^\circ)\cup \{v_i\}$, i.e.~on the top or bottom bounding slice or at $v_i$. 
\end{itemize}
Such $\w$-matching segments are considered up to compactly supported smooth isotopy of $\mathbb{D}_{v_i}$ relative $\mathbb{D}_{v_i}\cup \{v_i\}$.

\item If $i=m$, a matching segment $\gamma$ in $\mathbb{D}_{v_m}$ is defined similarly, now imposing that $\gamma(0)=v_m$ and that $\gamma(1)$ lies on the lower boundary slice of $\mathbb{D}_{v_m}$.\qed
\end{enumerate}
\end{definition}

\noindent \cref{fig:WeaveVertex_Segments} depicts four $\w$-matching segments illustrating \cref{def:mathcing_segments}.

\begin{remark}\label{rem:segments}
There is a unique matching segment in $\mathbb{D}_{v_m}$, denoted $\delta_{v_m}^\downarrow$. For $i\in[2,n]$, there are exactly four matching segments $\delta_{v_i}^\dagger$ in $\mathbb{D}_{v_i}$ up to compactly supported isotopy, where $\dagger\in\{\uparrow,\downarrow,\circlearrowright,\circlearrowleft\}$, cf.~\cref{fig:WeaveVertex_Segments}. Indeed, there are two such segments with an endpoint at $v_i$, starting either at the top boundary, denoted $\delta_{v_i}^\uparrow$, or ending at the bottom boundary, denoted $\delta_{v_i}^\downarrow$. There are two $\w$-matching segments starting at the top boundary and ending at the bottom boundary: going around $v_i$ either clockwise\footnote{Note that clockwise segments turn right at a $2$-valent vertex, and are thus isotopic to local clockwise trajectories considered in \cite{Chr25b}, since our orientations of the plane here are reversed from \cite{Chr25b}.}, denoted $\delta_{v_i}^\circlearrowright$, or counterclockwise, denoted $\delta_{v_i}^\circlearrowleft$.\qed
\end{remark}

\noindent Note also that every $\w$-matching path in \cref{def:matchingpath} arises uniquely as the composite of matching segments.


\subsubsection{Local lax sections from matching segments}\label{construction:local_sections} 

Let $v_i$ be a vertex of $\rgraph_\w$ and $\delta$ a $\w$-matching segment in $\mathbb{D}_{v_i}$, as in \cref{def:mathcing_segments}. We suppose further that we are given the following data, denoted $L$ (also called a local value):
\begin{itemize}
    \item If $\delta=\delta_{v_i}^\downarrow,\delta_{v_i}^\uparrow$, an element $L= k[l]\in \D^{\on{perf}}(k)$, i.e.~an integer $l\in \mathbb{Z}$. 
    \item If $\delta=\delta_{v_i}^\circlearrowleft,\delta_{v_i}^\circlearrowright$, an element $L\in \D=\mathcal{F}(e_{i+1})$.
\end{itemize}
We now proceed with defining an associated local lax section $M^L_{\delta}\in \mathcal{L}(\rgraph_\w,\mathcal{F}_{\w})$, in the sense of \Cref{def:sections}, supported at $v_i$ and its incident edges $e_{i+1},e_i$, cf.~\cref{def:local_lax_sections} below. For that, let $p\colon \Gamma(\mathcal{F}_\w)\to \on{Exit}(\rgraph_{\w})$ be the Grothendieck construction of $\mathcal{F}_\w$. As in \cref{sssec:2spherical_matching_objects}, we describe objects in $\Gamma(\mathcal{F}_\w)$ as pairs $(v;O_v)\in\Gamma(\mathcal{F}_\w)$, where $v\in \on{Exit}(\rgraph_{\w})$ and $O_v\in\F_\w(v)$ are objects. Note that $\F_\w(v_i)=\D^{\on{perf}}(k)\overset{\rightarrow}{\times}_{(\mhyphen)\otimes S_{p_i}} \D$ if $i\in[1,m-1]$, whose elements we write as tuples as in \eqref{eq:description_VFN}, and $\F_\w(v_m)=\D^{\on{perf}}(k)$. Consider the following five partial sections of $p$:
\begin{align*}
&\mathcal{M}^\downarrow_{L,i}:\Delta^0=\{v_i\}\lr\Gamma(\F_w),\quad   \mathcal{M}^\downarrow_{L,i}(v_i)\coloneqq(v_i;k[l-1],0,0)\,,\stepcounter{equation}\tag{\theequation}\label{eq:lax_down_value_at_vertex}\\ 
&\mathcal{M}^\uparrow_{L,i}:\Delta^0=\{v_i\}\lr\Gamma(\F_w),\quad \mathcal{M}^\uparrow_{L,i}(v_i)\coloneqq(v_i;k[l],S_{p_i}[l],\on{id}_{S_{p_i}[l]})\,,\\
&\mathcal{M}^\circlearrowright_{L,i}:\Delta^0=\{v_i\}\lr\Gamma(\F_w),\quad \mathcal{M}^\circlearrowright_{L,i}(v_i)\coloneqq(v_i;0,L,0)\,,\stepcounter{equation}\tag{\theequation}\label{eq:standard_value_at_vertex}\\
&\mathcal{M}^\circlearrowleft_{L,i}:\Delta^0=\{v_i\}\lr\Gamma(\F_w),\quad \mathcal{M}^\circlearrowleft_{L,i}(v_i)\coloneqq(v_i;\on{Mor}(S_{p_i},L),L, \on{Mor}(S_{p_i},L)\otimes S_{p_i}\stackrel{\on{ev}}{\shortrightarrow} L)\,,\stepcounter{equation}\tag{\theequation}\label{eq:costandard_value_at_vertex}\\
&\mathcal{M}^\downarrow_{L,m}:\Delta^0=\{v_m\}\lr\Gamma(\F_w),\quad  \mathcal{M}^\downarrow_{L,m}(v_m)\coloneqq(v_m;L)=(v_m;k[l]).\\
\end{align*}
The index $i$ for the first four functors ranges in $i\in[1,m-1]$. These functors effectively assign the value of a yet to be constructed lax section at the corresponding vertex $v_i$ or $v_m$, which depends on the nature of the matching segment. Intuitively, each of these can be made into a lax section of the weave schober by extending with zeros. However, for gluing an alternative extension to a lax section will be more convenient, given universally via the left Kan extension.

\begin{equation}\label{eq:relative_leftKanExtensions}
\begin{tikzcd}
 & \Gamma(\mathcal{F}_\w) \arrow[d,"p"'] \\
\Delta^0 \arrow[ru,"\mathcal{M}^\dagger_{L,i}"] \arrow[r,hook,"\iota"] & \on{Exit}(\rgraph_\w) \arrow[u,bend right,blue,"\mbox{Lan}^p_\iota(\mathcal{M}^\dagger_{L,i})"']
\end{tikzcd}
\end{equation}

\begin{definition}[Local lax sections]\label{def:local_lax_sections}
Let $
\iota:\Delta^0\lr\on{Exit}(\rgraph_\w)$ denote the inclusion of a vertex $v_i\in \rgraph_\w$. By definition, the local lax section  
\[ M^L_{\delta^\dagger_{v_i}}\coloneqq\on{Lan}^p_\iota(\mathcal{M}^\dagger_{L,i})\in\losec(\rgraph_\w,\mathcal{F}_\w)\]
for the segment $\delta^\dagger_{v_i}$ with value $L$ is the $p$-relative left Kan extension $\on{Lan}^p_\iota(\mathcal{M}^\dagger_{L,i})$ of $\mathcal{M}^\dagger_{L,i}$ along $\iota$.\qed
\end{definition}

\noindent In \cref{def:local_lax_sections}, the superscript $\dagger$ denotes $\dagger\in\{\uparrow,\downarrow,\circlearrowright,\circlearrowleft\}$ and the lax sections are local in that they are supported at $v_i$ and their adjacent edges $e_i,e_{i+1}$, meaning that $M^L_{\delta^\dagger_{v_i}}(x)\simeq (x;0)$ for all $x\not =v_i,e_i,e_{i+1}$. The diagram in \cref{def:local_lax_sections} depicts such lax sections $M^L_{\delta^\dagger_{v_i}}$ from \cref{def:local_lax_sections} as relative left Kan extensions, emphasized in blue. For instance, in the case of $v_m$, this implies that 
\[ M_{\delta_{v_m}^\downarrow}^L(v_m\to e_m)\colon (v_m;k[l])\to (e_m;S_{p_m}[l])\] 
is the coCartesian edge in $\Gamma(\mathcal{F}_\w)$ encoding the application of the functor $(\mhyphen)\otimes S_{p_m}$. In general, we compute the values of $M^L_{\delta}$ at the adjacent edges $e_{i+1}$ and $e_i$ by applying the functors $\varrho_1$, respectively $\varrho_2$, from \Cref{lem:schober2gon}. This leads to the following values at the adjacent edges:
\begin{equation*}
 M^L_{\delta_{v_i}^\uparrow}(e_{i+1})\simeq S_{p_i}[l]\,, \quad\quad \quad M^L_{\delta_{v_i}^\uparrow}(e_i)\simeq 0
\end{equation*}
\begin{equation}\label{eq:matchingsection_down}
 M^L_{\delta_{v_i}^\downarrow}(e_{i+1})\simeq 0\,, \quad\quad \quad M^L_{\delta_{v_i}^\downarrow}(e_i)\simeq S_{p_i}[l]
\end{equation}
\begin{equation}\label{eq:cwtransport}
 M^L_{\delta_{v_i}^\circlearrowright}(e_{i+1})\simeq L\,, \quad\quad \quad M^L_{\delta_{v_i}^\circlearrowright}(e_i)\simeq L
\end{equation}
\begin{equation}\label{eq:ccwtransport}
 M^L_{\delta_{v_i}^\circlearrowleft}(e_{i+1})\simeq L\,, \quad\quad \quad M^L_{\delta_{v_i}^\circlearrowleft}(e_i)\simeq T_{S_{p_i}}(L)\,.
\end{equation}

\begin{remark}
These are the schober local models categorifying the appearances of Dehn twists for parallel transports in Picard-Lefschetz theory, which in the case of weave schobers are combinatorialized in \cref{ssec:weave_thimbles}. To wit, Equations (\ref{eq:cwtransport}) and (\ref{eq:ccwtransport}) encode the fact that, for a regular fiber fixed east of the critical points, counterclockwise encircling a critical value leads to a positive twist along the associated vanishing cycle. Specifically, the parametrized schober refines this information by declaring that passing down to the right of a critical value is trivialized to be the identity, as in \Cref{eq:cwtransport}, and passing down to left of it results in the application of the positive twist associated to its vanishing cycle, as in \Cref{eq:ccwtransport}.\qed
\end{remark}


\subsubsection{Global sections from $\w$-matching paths}\label{subsubsec:glsecfrommatchingpath}
\label{constr:glsecfrommatchingpath}
Let us now glue the local lax sections from \cref{def:local_lax_sections} to obtain a global section of $\F_\w$. For that, let $\gamma$ be a $\w$-matching path in $\mathbb{D}$, as in \cref{def:matchingpath}, starting at a vertex $v_i$ and ending at a vertex $v_{i'}$, for some $i,i'\in[1,m]\cup\{-\infty\}$ with $i'<i$. Let us express $\g$ as a composition of $\w$-matching segments $\delta_i,\dots,\delta_{i'}$, cf.~\cref{def:mathcing_segments}, or until $\delta_{1}$ if $i'=-\infty$. We decorate each such a segment $\delta_i$ with a local object $L_i$ as follows: 
\begin{itemize}
    \item[$(i)$] For the initial segment $\delta_i=\delta_{v_i}^\downarrow$, we set $L_i\coloneqq k\in \D^{\on{perf}}(k)$.
    \item[$(ii)$] For $j\in(i',i)$, we recursively declare $L_{j}\coloneqq M_{\delta_{j+1}}^{L_{j+1}}(e_{j})$, using \cref{def:local_lax_sections}.
    \item[$(iii)$] If $i'\not = -\infty$, the condition \Cref{def:matchingpath}.(5) is equivalent to the assertion that $M_{\delta_{j+1}}^{L_{j+1}}(e_{j})\simeq S_p[l]$ for some $l\in \mathbb{Z}$. We then declare $L_{i'}\coloneqq k[l] \in \D^{\on{perf}}(k)$.
\end{itemize}

\noindent In order to obtain the global section $M_\gamma\in \mathcal{H}(\rgraph_\w,\mathcal{F}_{\w})$ associated to the $\w$-matching path $\g$, we glue together the compatible family of local sections $\{M_{\delta_{j}}^{L_{j}}\}_{j\in[i',i]}$, which is rigorously done as follows.\\

For an edge $e_{j}$ and object $Y\in \D= \mathcal{F}_\w(e_{j})$, let $\iota_{e_{j}}(Y)\in \mathcal{L}(\rgraph_\w,\mathcal{F}_\w)$ be the lax section supported at $e_{j}$ with value $Y$. For all $j\in(i',i]$, we choose morphisms
\begin{equation}\label{eq:gluing_morphisms}
\iota_{j}(L_{j})\longrightarrow M_{\delta_{j}}^{L_{j}},M_{\delta_{j-1}}^{L_{j-1}},\quad j\in(i',i]
\end{equation}
in $\mathcal{L}(\rgraph_\w,\mathcal{F}_\w)$ which restrict at the edge $e_{j}$ to an equivalence. Let $Z_{\gamma}\subset \on{Exit}(\rgraph_{\w})$ be the full subcategory spanned by the vertices $v_i,v_{i+1},\dots,v_{i'}$, omitting $v_{i'}$ if $i'=-\infty$, and their connecting edges $e_i,\dots, e_{i'+1}$.

\begin{definition}[Global sections from $\w$-matching paths]\label{def:globalsec_from_matchingpath} Let $\w$ be a Demazure weave, $\F_\w$ its weave schober and $\g$ a $\w$-matching path. By definition, the global section $M_\g\in\glsecF$ is defined as the colimit
$$M_\g\coloneqq\mbox{colim}(D_{\gamma}\colon Z^{\on{op}}_{\gamma}\longrightarrow \mathcal{L}(\rgraph_\w,\mathcal{F}_\w))$$
where the diagram $D_\g$ is given by $D_\g(v_{j})\coloneqq M_{\delta_{j}}^{L_{j}}$ and $D_\g(e_{j})\coloneqq \iota_{e_{j}}(L_{j})$, where we implicitly use the morphisms from \eqref{eq:gluing_morphisms}.\qed
\end{definition}

\noindent By inspecting its construction, we deduce that $M_\g\in\glsecF$ in \cref{def:globalsec_from_matchingpath} is coCartesian and, though the diagram $D_\g$ depends on the morphisms \eqref{eq:gluing_morphisms}, its colimit $M_\g$ is independent of these choices.

\begin{definition}\label{eq:matching_spheres_thimbles} Let $\w$ be a Demazure weave with $m$ trivalent vertices, and $\gamma\colon [0,1]\to \mathbb{D}$ be a $\w$-matching path. By definition, the object $M_\g\in\glsecF$ is said to be
\begin{itemize}
    \item a $\w$-matching sphere if the endpoint $\g(1)=v_i$ is a vertex $i\in[1,m]$,
    \item a $\w$-thimble if instead $\g(1)=\vinf$.\qed
\end{itemize}
\end{definition}

\noindent The notation in \cref{eq:matching_spheres_thimbles} is motivated by the following computation:

\begin{lemma}\label{lem:machting_spheres_and_thimbles}
    Let $\gamma\colon [0,1]\to \mathbb{D}$ be a matching path. 
    \begin{enumerate}[(1)]
    \item If $\gamma(1)=v_i$ for some $i\geq 1$, then
        \[ \on{Mor}_{\glsecF}(M_\gamma,M_\gamma)\simeq k\oplus k[-3]\,,\]
    meaning that $M_{\gamma}$ is a $3$-spherical object.
    \item If $\gamma(1)=v_{-\infty}$, then 
    \[ \on{Mor}_{\glsecF}(M_\gamma,M_\gamma)\simeq k\,,\]
    meaning that $M_{\gamma}$ is exceptional.
    \end{enumerate}
\end{lemma}

\begin{proof} By \cref{def:globalsec_from_matchingpath}, $M_{\gamma}\simeq \displaystyle\on{colim}_{Z^{\on{op}}_\gamma} D_{\gamma}$ and thus
\[
\on{Mor}_{\glsecF}(M_\gamma,M_\gamma)\simeq \displaystyle\lim_{x\in Z_\gamma} \on{Mor}_{\mathcal{L}(\rgraph_\w,\mathcal{F}_\w)}(D_\gamma(x),M_\gamma)\,.
\]
Since $D_{\gamma}(x)$ is the $p$-relative left Kan extension of its restriction to $x\in Z_\gamma\subset \on{Exit}(\rgraph_\w)$, it follows that 
\[
\on{Mor}_{\mathcal{L}(\rgraph_\w,\mathcal{F}_\w)}(D_\gamma(x),M_\gamma)\simeq \on{Mor}_{\mathcal{F}_\w(x)}(D_\gamma(x)({\bf x}),M_\gamma({\bf x}))\,
\]
where here we have used ${\bf x}$ to indicate the vertex or edge ${\bf x}\in\Wgraph$ tautologically corresponding to the object $x\in Z_\g$.
By unraveling the diagrams involved, we obtain that $\on{Mor}_{\glsecF}(M_\gamma,M_\gamma)$ is the limit of the diagram
\[
\begin{tikzcd}
k \arrow[rd, hook] &                 & {k\oplus k[-2]} \arrow[rd, "\simeq"] \arrow[ld, "\simeq "] &        & {k\oplus k[-2]} \arrow[rd, "\simeq "] \arrow[ld, "\simeq "] &                 & k \arrow[ld, hook] \\
                   & {k\oplus k[-2]} &                                                            & \dots  &                                                             & {k\oplus k[-2]} &                   
\end{tikzcd}
\]
if $\gamma(1)=v_k$ with $k\geq 1$, and otherwise, if $\gamma(1)=v_{-\infty}$, it is the limit of the diagram
\[
\begin{tikzcd}
k \arrow[rd, hook] &                 & {k\oplus k[-2]} \arrow[rd, "\simeq"] \arrow[ld, "\simeq "] &        & {k\oplus k[-2]} \arrow[rd, "\simeq "] \arrow[ld, "\simeq "] &                 \\
                   & {k\oplus k[-2]} &                                                            & \dots  &                                                             & {k\oplus k[-2]}
\end{tikzcd}
\]
In the former case, we obtain
$$\on{Mor}_{\glsecF}(M_\gamma,M_\gamma)\simeq\mbox{lim}\left(
\begin{tikzcd}
 & k \arrow[d, hook]\\
k \arrow[r, hook] & {k\oplus k[-2]}
\end{tikzcd}\right)\simeq k\oplus k[-3],$$
\noindent since both inclusions are to the first factor of $k\oplus k[-2]$. In the latter case, 
\[ \on{Mor}_{\glsecF}(M_\gamma,M_\gamma)\simeq \mbox{lim}\left(
\begin{tikzcd}
 & k\oplus k[-2] \arrow[d, "\simeq"]\\
k \arrow[r, hook] & {k\oplus k[-2]}
\end{tikzcd}\right) \simeq k\,,\]
as required.
\end{proof}


\begin{center}
	\begin{figure}[h!]
		\centering
        \includegraphics[scale=1.1]{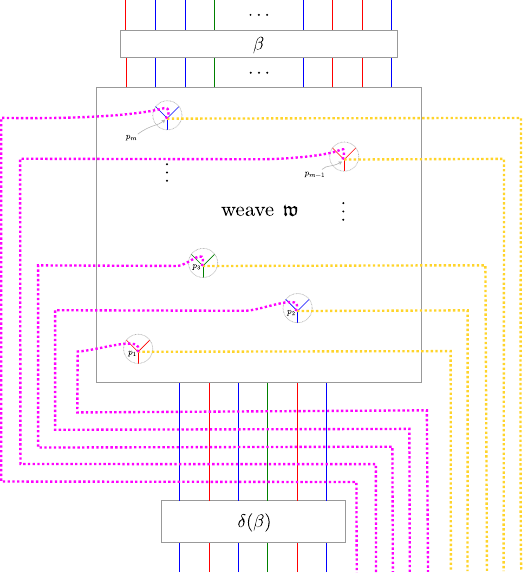}
		\caption{A general Demazure weave with its weave vanishing paths, as introduced in \cref{def:weave_thimbles}. The left vanishing paths $\wp_l$ are depicted in pink, and the right vanishing paths $\wp_r$ in yellow.}\label{fig:WeaveThimbles}
	\end{figure}
\end{center}

\subsection{Weave vanishing paths and transport}\label{ssec:weave_transport} Let $\w:\beta\to\delta(\beta)$ be a Demazure weave. We now introduce a certain type of $\w$-matching path, as in \cref{def:matchingpath}, which will be used to construct the exceptional collection $\stds^\w$ and its dual.

\begin{center}
	\begin{figure}[h!]
		\centering
        \includegraphics[scale=1.2]{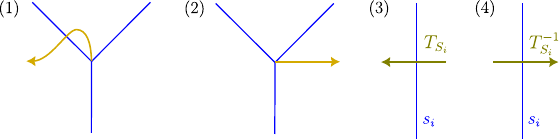}
		\caption{(1) The start of the path for a left vanishing path at the trivalent vertex of a weave. (2) The start of the path for a right vanishing path. (3) and (4) depict the rules for parallel transport in the presence of a weave, cf.~\Cref{def:combinatorial_transport} .}\label{fig:Weave_ParallelTransportModels}
	\end{figure}
\end{center}

\begin{definition}[Weave vanishing paths]\label{def:weave_thimbles} Let $p\in\w$ be a trivalent vertex. By definition:
\begin{enumerate}
    \item The right weave vanishing path $\wp_r(p)\sse\mathbb{D}$ is the embedded path given by concatenating the straight horizontal path from $p$ to the right, until $\w$ is entirely to the left of the path, with the straight vertical path heading downwards, as in \Cref{fig:WeaveThimbles} (yellow). The path $\wp_r$ starts as in \Cref{fig:Weave_ParallelTransportModels}.(2).

    \item The left weave vanishing path $\wp_l(p)\sse\mathbb{D}$ is the embedded path given by concatenating the following four paths. First, the horizontal path from $p$ to the left starting with the upwards hook, as in \Cref{fig:Weave_ParallelTransportModels}.(1), until $\w$ is entirely to the right of the path. Second, the straight vertical path heading downwards, until we descend to a height where the horizontal weave slice spells $\delta(\beta)$. Third, the straight horizontal path to the right from that point, and fourth the straight vertical path heading downwards. Left weave vanishing paths are depicted as in \Cref{fig:WeaveThimbles} (pink).
\end{enumerate}
Weave vanishing paths will typically be considered up to compactly supported planar isotopies in the complement of the trivalent vertices of the weave, and assumed to be transversal to $\w$.
\qed
\end{definition}

\begin{example}\label{ex:weavevanishingpaths_as_segments} The $\w$-matching path starting at $v_i$, defined as the composition of the segments $\delta_{v_i}^\downarrow$ and $\delta_{v_{i-1}}^\circlearrowright,\dots,\delta_{v_{1}}^\circlearrowright$, is equivalent to the right weave vanishing path $\wp_r(p_i)$. Similarly, the $\w$-matching path starting at $v_i$ given by the composition of $\delta_{v_i}^\downarrow$ and $\delta_{v_{i-1}}^\circlearrowleft,\dots,\delta_{v_{1}}^\circlearrowleft$ gives a representative for $\wp_l(p_i)$, which will yield equivalent objects, up to a shift.\qed
\end{example}

\noindent In order to illustrate the geometric idea that in this article $\w$ is being treated as a set of branch cuts, we introduce the following notion:

\begin{definition}[Weave transport]\label{def:combinatorial_transport} Let $\wp$ be an oriented embedded segment in the plane transverse to $\w$. By definition, the $\w$-transport of an object $L\in\D$ along $\wp$ is the object
\begin{equation}\label{eq:weave_transport}
\mbox{Hol}_\w(L;\wp)\in\D
\end{equation}
obtained by starting with $L\in\D$ and following the rules of \Cref{fig:Weave_ParallelTransportModels}.(3)-(4) as one travels forward along $\wp$. Namely, if the path $\wp$ crosses a weave line decorated by $s_i\in W$ left-to-right, we apply the twist functor $T_{S_i}^{-1}$, and if $\wp$ crosses a weave line decorated by $s_i\in W$ right-to-left, we apply the inverse twist functor $T_{S_i}$.\qed
\end{definition}

\begin{center}
	\begin{figure}[h!]
		\centering
        \includegraphics[scale=1.2]{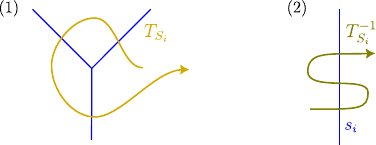}
		\caption{(1) An illustration of the fact that \Cref{fig:Weave_ParallelTransportModels}.(3) and (4) imply that circling counter-clockwise around a trivalent vertex decorated by $s_i\in W$ results in a positive twist by $S_i$, as we have $T_{S_i}=T^{-1}_{S_i} T^2_{S_i}$. (2) A depiction exemplifying the independence of the planar isotopy representative in the complement of the trivalent vertices, as we have $T^{-1}_{S_i}=T^{-1}_{S_i}T_{S_i}T^{-1}_{S_i}$.}\label{fig:Weave_ParallelTransportModels2}
	\end{figure}
\end{center}

\noindent The combinatorial rules in \Cref{fig:Weave_ParallelTransportModels}, and thus \cref{def:combinatorial_transport}, are consistent with the fact that counter-clockwise winding around a trivalent vertex leads to a positive twist, as in Picard-Lefschetz theory, cf.~\Cref{fig:Weave_ParallelTransportModels2}.(1), Similarly, they imply that performing a zig-zaging as we cross a weave line, as in \Cref{fig:Weave_ParallelTransportModels2}.(2), leads to the same transport as the straight segment in \Cref{fig:Weave_ParallelTransportModels}.(4). In particular, weave parallel transport is invariant under planar isotopies in the complement of the trivalent vertices of the weave, and thus paths $\wp$ will typically be considered up to such isotopies, see e.g.~\cref{def:weave_thimbles}. 

\noindent For instance, the spherical object $S_p$ associated to a trivalent $s_i$-vertex $p\in\w$ in \cref{def:sphericalobj} can be understood in these terms:
\begin{equation}\label{eq:sphericaltrivalent_as_holonomy}
S_{p}\simeq\mbox{Hol}_\w(S_i;\wp_r(p)),
\end{equation}

\noindent i.e.~the object $S_p\in\D$ is the $\w$-transport of $S_i$ along the right vanishing path $\wp_r(p)$. The presentation in \eqref{eq:sphericaltrivalent_as_holonomy} promotes that the objects $S_p$ are the weave analogues of vanishing cycles, which are the evaluations of the thimbles corresponding to right vanishing paths at the regular fiber at $-\infty$. The $\w$-transports $\mbox{Hol}_\w(S_i;\wp_l(p))$ along left vanishing paths are also of relevance, but they can be computed in terms of the objects in $\eqref{eq:sphericaltrivalent_as_holonomy}$, coming for the right vanishing paths, due to the following formula:

\begin{lemma}[Alternative rule for transport along left weave vanishing path]\label{lem:left_thimble_transport} Let $\w$ be a Demazure weave, $p_i\in\w$ a trivalent vertex decorated by $s_j\in W$, and $\wp_l(p_i)$ its left vanishing path. Then

\begin{equation}\label{eq:left_thimble_transport}
\on{Hol}_\w(S_j;\wp_l(p_i))\simeq T_{S_{p_1}}T_{S_{p_2}}\cdots T_{S_{p_{i-1}}}(S_{p_i}).
\end{equation}

\noindent In particular, the $\w$-transport of $S_j$ along the left vanishing path $\wp_l(p_i)$ can be computed from certain twists along the right vanishing spheres below of the vanishing sphere of $\wp_r(p_i)$.
\end{lemma}
\begin{proof}
\Cref{eq:left_thimble_transport} readily holds for right-inductive weaves. Thus the statement follows from the fact that the validity of \Cref{eq:left_thimble_transport} is preserved under weave equivalences and weave mutations.
\end{proof}



\subsection{Standard and costandard weave thimbles}\label{ssec:weave_thimbles}

In this subsection we introduce the collections of standard and costandard thimbles of a Demazure weave. They are proven to be a full exceptional collection in \cref{prop:weave_full_exceptional_collection}, and their behavior under weave equivalences and mutations is presented in \cref{lem:change_in_thimbles_under_weave_equivalence}.

\begin{definition}[Weave thimbles]\label{def:thimbles}
Let $\w$ be a Demazure weave, $v_i\in\w$ a trivalent vertex and $\F_\w$ its weave schober. By definition:
\begin{enumerate}[(1)]
    \item The standard weave thimble $\Delta_i^\w\in\glsecF$ is $\Delta_i^\w\coloneqq M_{\wp_r(v_i)}$.

    \item The costandard weave thimble $\nabla_i^\w\in\glsecF$ is $\nabla_i^\w\coloneqq M_{\wp_l(v_i)}$.\qed
\end{enumerate}
\end{definition}   

\noindent In \cref{def:thimbles} we are using the global sections from \cref{def:globalsec_from_matchingpath}, and the weave vanishing paths from \cref{def:weave_thimbles}. By \cref{ex:weavevanishingpaths_as_segments}, if we consider the $\w$-matching path $\g_i^r$, resp.~$\g_i^l$, consisting of $\delta_{v_i}^\downarrow$ composed with $\delta_{v_{i-1}}^\circlearrowright,\dots,\delta_{v_{1}}^\circlearrowright$, resp.~$\delta_{v_i}^\downarrow$ composed with $\delta_{v_{i-1}}^\circlearrowleft,\dots,\delta_{v_{1}}^\circlearrowleft$, the thimbles in \cref{def:thimbles} can be equivalently described as $\Delta_i^\w\simeq M_{\g_i^r}$ and $\nabla_i^\w\simeq M_{\g_i^l}$.\\

Given a vertex $v_i\in\rgraph_\w$, we denote by $\on{ev}_{v_i}\colon\glsecF\to \mathcal{F}_\w(v_i)$ the evaluation functor of global sections to $v_i$. By \cite[Lem.~4.11]{Chr25b}, this functor admits left and right adjoints, denoted $\on{ind}^L_{v_i}$ and $\on{ind}^R_{v_i}$ respectively, called the induction functors.  By \cite[Ex.~4.22]{Chr25b}, the weave thimbles in \cref{def:thimbles} can be expressed\footnote{The convention difference with \cite{Chr25b} has the effect that the thimbles go down on the opposite side.} as inductions:

\begin{equation}\label{eq:thimble_as_induction}
    \Delta_{i}^\w\simeq \on{ind}^L_{v_i}\circ \on{ev}_{v_i}(\Delta_i^\w)\simeq  \on{ind}^L_{v_i}((k[-1],0,0))\,
\end{equation}
and similarly
\begin{equation}\label{eq:cothimble_as_induction}
\nabla_{i}^\w\simeq \on{ind}^R_{v_i}\circ \on{ev}_{v_i}(\nabla_i^\w)\simeq  \on{ind}^R_{v_i}((k[-1],0,0))\,.
\end{equation}

\begin{remark}
 \Cref{def:thimbles}, \eqref{eq:thimble_as_induction}, and \eqref{eq:cothimble_as_induction} can be readily generalized to arbitrary $\rgraph_\w$-parametrized perverse schobers, by replacing the vanishing cycle $k\in \D^{\on{perf}}(k)$ with any chosen object in the $\infty$-category of vanishing cycles of the perverse schober at $v_i$.
 \qed
 \end{remark}

\begin{center}
	\begin{figure}[h!]
		\centering
		\includegraphics[scale=1.3]{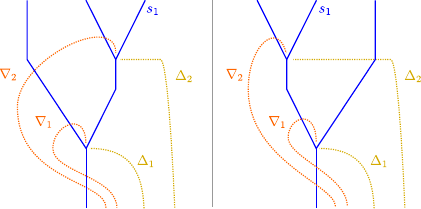}
		\caption{Two 2-stranded weaves discussed in \cref{ex:s1cube_weave_mutation}. In both cases, the standards $(\st1,\st2)$
        and costandards $(\cst1,\cst2)$ are depicted. (Left) Left inductive weave for $\s_1^3$. (Right) Right inductive weave for $\s_1^3$.}
        \label{fig:Example_WeaveMutation}
	\end{figure}
\end{center}

\noindent The adjective costandard in \cref{def:thimbles} is justified by the following fact:

\begin{proposition}[Duality for standard and costandard thimbles]\label{prop:duality_standards_costandards}
Let $\w$ be a Demazure weave.
Then the standard and costandard thimbles are dual, i.e.~
\[
\on{Mor}_{\glsecF}(\Delta^\w_i,\nabla^\w_j)\simeq \delta_{ij}\cdot k\qquad \forall i,j\in[1,m].
\]
\end{proposition}

\begin{proof}
If $i\not =j$, \eqref{eq:thimble_as_induction} and \eqref{eq:cothimble_as_induction} imply
\begin{align*}
\on{Mor}_{\glsecF}(\Delta^\w_i,\nabla^\w_j)&\simeq \on{Mor}_{\glsecF}(\on{ind}^L_{v_i}((k[-1],0,0)),\on{ind}^R_{v_j}((k[-1],0,0)))\\
& \simeq \begin{cases} \on{Mor}_{\mathcal{F}_{\w}(v_i)}(X,\on{ev}_{v_i}(\nabla^\w_j)),  & \mbox{ if }i>j\\ \on{Mor}_{\mathcal{F}_{\w}(v_j)}(\on{ev}_{v_j}(\Delta^\w_i),Y), & \mbox{ if }i<j \end{cases} \\
& \simeq 0\,
\end{align*}
where $X\in\F_{\w}(v_i)$ and $Y\in\F_{\w}(v_j)$ are the starting objects for $\Delta^\w_i$ and $\nabla^\w_j$, respectively. If $i=j$, \eqref{eq:thimble_as_induction} and \eqref{eq:cothimble_as_induction} give 
\begin{align*}
\on{Mor}_{\glsecF}(\Delta_i^\w,\nabla^\w_j)&\simeq \on{Mor}_{\glsecF}(\on{ind}^L_{v_i}((k[-1],0,0)),\on{ind}^R_{v_j}((k[-1],0,0)))\\
&\simeq  \on{Mor}_{\mathcal{F}_{\w}(v_i)}((k[-1],0,0),\on{ev}_{v_i}(\nabla^\w_j))\\
&\simeq  \on{Mor}_{\mathcal{F}_{\w}(v_i)}((k[-1],0,0),(k[-1],0,0))\\
&\simeq k\,.
\end{align*}
\end{proof}

\noindent The duality in \cref{prop:duality_standards_costandards} is in line with the classical duality between standard and costandard objects in the context of highest weight categories, see the proof of \cite[Theorem 3.11]{CPS88}, or \cite[Prop.~5.3]{braden2008galedualitykoszulduality}.\\

Let $e_i\in\Wgraph$ be an edge, and $\on{ev}_{e_i}:\glsecF\lr\D$ the evaluation functor to such an edge $e_i$. Given any standard thimble $\Delta_k^\w$, as in \cref{def:thimbles}.(1), we have

\begin{equation}\label{eq:evaluation_thimbles}
\on{ev}_{e_i}(\Delta_k^\w)\simeq\begin{cases}
0 & \mbox{ if }k<i,\\
S_{p_k} & \mbox{ if }i\leq k.
\end{cases}
\end{equation}

\noindent In particular, the evaluation $\on{ev}_{e_1}$, which we can alternatively denote $\on{ev}_{\vinf}$, always yields the spherical object $\on{ev}_{\vinf}(\Delta_k^\w)\simeq S_{p_k}$ at $p_k$. The following lemma explains how to compute the morphism objects between standard weave thimbles in terms of their evaluations:

\begin{lemma}\label{lem:morphism_objects_between_standard_thimbles}
Let $\w$ be a Demazure weave with $m$ trivalent vertices. Then for all $j\in[1,m]$, $i\in[1,j)$, the functor $\on{ev}_{e_i}$ induces an equivalence
\[
\on{Mor}_{\glsecF}(\Delta_i^\w,\Delta_j^\w)\simeq \on{Mor}_{\D}(S_{p_i},S_{p_j})\,. 
\]
In addition, $\on{Mor}_{\glsecF}(\Delta_i^\w,\Delta_i^\w)\simeq k$.
\end{lemma}

\begin{proof}
By \eqref{eq:thimble_as_induction}, $\Delta_i^\w\simeq \on{ind}^L_{v_i}((k[-1],0,0))$ and thus 
\[ \on{Mor}_{\glsecF}(\Delta_i^\w,\Delta_j^\w)\simeq \on{Mor}_{\mathcal{F}_\w(v_i)}((k[-1],0,0),\on{ev}_{v
_i}(\Delta_j^\w))\simeq\]
\[\simeq\on{Mor}_{\mathcal{F}_\w(v_i)}((k[-1],0,0),(0,S_{p_j},0))\simeq \on{Mor}_\D(S_{p_i},S_{p_j})\,.\]
The equivalence $\on{Mor}_{\glsecF}(\Delta_i^\w,\Delta_i^\w)\simeq k$ is \Cref{lem:machting_spheres_and_thimbles}.(2).
\end{proof}

\begin{example}\label{ex:s1cube_weave_mutation}
Let $\beta=\s_1^3$ and consider the two 2-stranded weaves in \Cref{fig:Example_WeaveMutation}. The combinatorial operation of exchanging these two weaves is the most fundamental instance of a weave mutation. To ease notation, we omite the superscript $\w$ from the standard and costandard objects.

\begin{enumerate}[(1)]
    \item Let $\w$ be the left inductive Demazure weave for $\beta$, as in \cref{fig:Example_WeaveMutation} (left). The evaluations at $\vinf$ of the standard thimbles $(\st1,\st2)$ are 
    \[ \evinf(\Delta_2)\simeq S_1,\quad\quad\quad \evinf(\Delta_1)\simeq S_1,\]
    and those of the costandard objects $(\cst1,\cst2)$ are
     \[ \evinf(\nabla_2)\simeq S_1[-1]\,,\quad\quad\quad \evinf(\nabla_1)\simeq S_1\,.\]
     The weave transport along the right and left vanishing paths $\wp_r$ and $\wp_l$ are illustrated in \cref{fig:Example_WeaveMutation}, with $\wp_r$ depicted in dashed yellow lines, and $\wp_l$ in dashed orange lines. By \cref{lem:morphism_objects_between_standard_thimbles}, these evaluations allow us to compute the relevant morphism groups:
     \[ \on{Mor}(\Delta_1,\Delta_2)\simeq \on{Mor}_\D(\evinf(\st1),\evinf(\st2))\simeq \on{Mor}_\D(S_1,S_1)\simeq k\oplus k[-2],\]
     \[ \on{Mor}(\nabla_2,\nabla_1)\simeq \on{Mor}_\D(\evinf(\cst2),\evinf(\cst1))\simeq \on{Mor}_\D(S_1[-1],S_1)\simeq k[1]\oplus k[-1].\]
     
     \item Let $\w$ be the right inductive Demazure weave for $\beta$, as in \cref{fig:Example_WeaveMutation} (right). In this case, see e.g.~\eqref{eq:evaluation_thimbles}, the evaluations at $\vinf$ of the standard objects $(\st1,\st2)$ are
    \[ \evinf(\Delta_2)\simeq S_1[1],\quad\quad\quad \evinf(\Delta_1)\simeq S_1.\]
    The costandard objects $(\cst1,\cst2)$ evaluate to
     \[ \evinf(\nabla_2)\simeq S_1\,,\quad\quad\quad \evinf(\nabla_1)\simeq S_1\,.\]
     \cref{lem:morphism_objects_between_standard_thimbles} implies the following computation for the morphism groups:
     \[ \on{Mor}(\Delta_1,\Delta_2)\simeq \on{Mor}_\D(\evinf(\st1),\evinf(\st2))\simeq \on{Mor}_\D(S_1,S_1[1])\simeq k[1]\oplus k[-1],\]
     \[ \on{Mor}(\nabla_2,\nabla_1)\simeq \on{Mor}_\D(\evinf(\cst2),\evinf(\cst1))\simeq \on{Mor}_\D(S_1,S_1)\simeq k\oplus k[-2]\,.\]
\end{enumerate}
\noindent In both instances (1) and (2) in this example the duality from \cref{prop:duality_standards_costandards} is readily verified.
\qed
\end{example}

\begin{center}
	\begin{figure}[h!]
		\centering
		\includegraphics[scale=1.2]{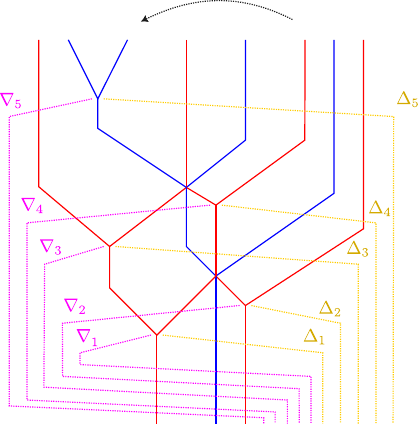}
		\caption{The weave $\w$ from \cref{ex:running_example1} and \cref{ex:running_example2}, cf.~also \cref{fig:Example_WeaveSchober1}. The associated standard objects $(\Delta_1,\ldots,\Delta_5)$, in yellow, and the costandard objects $(\nabla_1,\ldots,\nabla_5)$, in pink, constructed using \cref{subsec:matching}. The dashed black arrow on top indicates the direction of positive wrapping, which is coutner-clockwise. To ease notation, we have omitted the small hook at the beginning of the costandard objects according to \cref{fig:Weave_ParallelTransportModels}.(1).}
        \label{fig:Example_WeaveSchober_StandardsCostandards}
	\end{figure}
\end{center}

\begin{example}\label{ex:running_example2} The standards and costandards for the weave in \cref{ex:running_example1} are depicted in \cref{fig:Example_WeaveSchober_StandardsCostandards}. The spherical simples $S_{p_i}$ are computed in \cref{ex:running_example1}. In conjunction with \eqref{lem:morphism_objects_between_standard_thimbles} and \cref{lem:morphism_objects_between_standard_thimbles}, this yields the morphism groups $\mor(\Delta_i,\Delta_j)$ between the standard thimbles. For instance,
     \[ \on{Mor}(\Delta_1,\Delta_2)\simeq \on{Mor}_\D(\evinf(\st1),\evinf(\st2))\simeq \on{Mor}_\D(S_1,S_2)\simeq k[-1],\]
     \[ \on{Mor}(\Delta_3,\Delta_5)\simeq \on{Mor}_\D(\evinf(\st3),\evinf(\st5))\simeq \on{Mor}_\D(T_{S_2}^{-1}(S_1[1]),S_1[1])\simeq k.\]
\qed
 \end{example}

\begin{remark}
 \noindent Using either \Cref{def:thimbles}, or \eqref{eq:thimble_as_induction} and \eqref{eq:cothimble_as_induction}, a computation shows that the standard weave thimbles $(\Delta^\w_1,\dots,\Delta^\w_m)$ recover the costandard weave thimbles $(\nabla^\w_1,\dots,\nabla^\w_m)$ via
 $$\nabla^\w_l\simeq T_{\Delta^\w_1}^+T_{\Delta^\w_2}^+\cdots T_{\Delta^\w_{l-1}}^+(\Delta^\w_l)\,.$$
 Conversely, we also have
\begin{equation}\label{eq:standards_from_costandards}
\Delta^\w_l \simeq T_{\nabla^\w_1}^-T_{\nabla^\w_2}^-\cdots T_{\nabla^\w_{l-1}}^-(\nabla^\w_l)\,,
\end{equation}
see also \cref{rmk:exceptional_to_coexceptional} and \cite[Section 5k]{Sei08}.
\qed 
\end{remark}


\subsubsection{Exceptionality of weave thimbles and behaviour under weave operations}\label{ssec:weavethimbles_exceptional}

Let us show that the standard weave thimbles form a full exceptional collection in $\glsecF$:

\begin{proposition}[Exceptionality of weave thimbles]\label{prop:weave_full_exceptional_collection}
Let $\w$ be a Demazure weave.
    \begin{enumerate}[(1)]    
    \item The standard thimbles $(\Delta^\w_1,\dots,\Delta^\w_m)$ form a full exceptional collection in $\glsecF$.
    \item The costandard thimbles $(\nabla^\w_1,\dots,\nabla^\w_m)$ form a full coexceptional collection  in $\glsecF$. 
    \end{enumerate}
\end{proposition}

\begin{proof}
Let us denote $\Delta_i\coloneqq\Delta_i^\w$, as $\w$ is fixed throughout the argument. For Part (1), \eqref{eq:thimble_as_induction} implies
\begin{equation}\label{eq:exceptional_proof}
\on{Mor}_{\glsecF}(\Delta_i,X)\simeq \on{Mor}_{\mathcal{F}(v_i)}(\Delta_i(v_i),X(v_i))\,,\quad \forall X\in \glsecF.
\end{equation}
By considering $X=\Delta_j$ with $j\in[1,i
)$, the vanishing $\Delta_j(v_i)\simeq 0$ and \eqref{eq:exceptional_proof} imply $\on{Mor}_{\glsecF}(\Delta_i,\Delta_j)\simeq 0$. This shows the required semi-orthogonality. In conjunction with \cref{lem:morphism_objects_between_standard_thimbles} we obtain that $(\Delta_1,\dots,\Delta_m)$ is an exceptional collection.\\

Let us now show that this exceptional collection is full by an induction on the vertices of $\w$. For that, let $\w_{\leq m-1}\subset\w$ be the subweave of $\w$ containing the $m-1$ first trivalent vertices $p_1,\dots,p_{m-1}$. Then the subcategory
    \[ R\Gamma(\rgraph_{\w_{\leq m-1}},\mathcal{F}_{\w_{\leq m-1}})\simeq \on{fib}(\on{ev}_{v_m})\subset \glsecF\]
is equivalent to the full subcategory of global setions which vanish at $v_m$, and hence also at $e_m$. The induction hypothesis is that $(\Delta_1,\dots,\Delta_{m-1})$ is full. Suppose now that $\on{Mor}_{\glsecF}(\Delta_i,X)\simeq 0$ for all $i\in[1,m]$: to establish $(\Delta_1,\dots,\Delta_m)$ is full, we must show that $X\simeq 0$. Since $\mathcal{F}_{\w\leq n}(v_m)\simeq \D^{\on{perf}}(k)$ and $\Delta_m(v_m)\simeq k$, the equivalence \eqref{eq:thimble_as_induction} implies
\[
0\simeq \on{Mor}_{\glsecF}(\Delta_m,X)\simeq \on{Mor}_{\D^{\on{perf}}(k)}(k,X(v_m))\simeq X(v_m)
\]
and therefore $X\in R\Gamma(\rgraph_{\w_{\leq m-1}},\mathcal{F}_{\w_{\leq m-1}})$. Since $(\Delta_1,\dots,\Delta_{m-1})$ is full, by the induction hypothesis, we conclude find $X\simeq 0$, as required. Part (2) follows from Part (1) applied to the opposite $\infty$-category.
\end{proof}

\noindent We conclude by recording how weave thimbles behave under weave equivalence and weave mutation. By \cref{thm:weaveequivglsec}, the weave schober $\F_\w$ remains invariant under weave equivalences except for the case of a pulley move. Specifically, for a forward pulley move or a forward weave mutation from $\w$ to $\w'$, the weave schober changes according to $\F_\w\mapsto\mu_i^\#(\F_\w)$ for some $i\in[1,m]$, which itself induces an equivalence
$$\mu_i^\#:\glsecF\lr R\Gamma(\rgraph_{\w'},\F_{\w'}).$$
Similarly, for a backward pulley move or backward weave mutation from $\w$ to $\w'$, we have an equivalence
$$\mu_i^\flat:\glsecF\lr R\Gamma(\rgraph_{\w'},\F_{\w'}).$$

Therefore, by construction, the possible changes of weave thimbles under weave equivalences and weave mutations are the perverse schober analogues of Hurwitz moves for distinguished bases, see e.g.~\cite[Section (16d)]{Sei08}. Thus, if we have correctly developed the theory of weave schobers, one expects weave thimbles to either remain invariant or change by mutation of exceptional collections. This is indeed the case:

\begin{lemma}[Weave thimbles under weave equivalence and mutation]\label{lem:change_in_thimbles_under_weave_equivalence} In the notation above:

  \begin{enumerate}[$(1)$] 
\item  Let $\w\simeq \w'$ be a height-preserving weave equivalence. Then
$$\Delta_i^\w\simeq \Delta_i^{\w'}\qquad\mbox {and}\qquad \nabla_i^\w\simeq \nabla_i^{\w'},\quad \forall i\in[1,m].$$

\item Let $\w\to \w'$ be a forward pulley move or a forward weave mutation, exchanging the heights of vertices $p_{i-1}$ and $p_i$. Then $(\Delta_1^{\w'},\ldots,\Delta_m^{\w'})$ is the image under $\mu_i^\sharp$ of the forward mutation of $(\Delta_1^\w,\ldots,\Delta_m^\w)$ at $(\Delta_{i-1}^\w,\Delta_i^\w)$, i.e.~the objects
\[(\Delta_1^\w,\dots,\Delta_{i-2}^\w,T_{\Delta_{i-1}^\w}(\Delta_i^\w),\Delta_{i-1}^\w,\Delta_{i+1}^\w,\dots,\Delta_m^\w)\,,\]
where $T_{\Delta_{i-1}^\w}^{+}(\Delta_i^\w)\coloneqq \on{cof}(\on{Mor}_\C(\Delta_{i-1}^\w,\Delta_i^\w)\otimes \Delta_{i-1}^\w\longrightarrow \Delta_i^\w)$.
Similarly, the co-exceptional collection $(\nabla_1^{\w'},\dots,\nabla_m^{\w'})$ coincides with the image under $\mu_i^\sharp$ of the forward mutation of $(\nabla_1^\w,\dots,\nabla_m^\w)$ at $(\nabla_{i-1}^\w,\nabla_i^\w)$.

\item Let $\w\to \w'$ be a backward pulley move or a backward weave mutation, exchanging the heights of vertices $p_{i-1}$ and $p_i$.  Then $(\Delta_1^{\w'},\ldots,\Delta_m^{\w'})$ is the image under $\mu_i^\flat$ of the backward mutaion of $(\Delta_1^\w,\ldots,\Delta_m^\w)$ at $(\Delta_{i-1}^\w,\Delta_i^\w)$, i.e.~the objects
\[ (\Delta_1^\w,\dots,\Delta_{i-2}^\w,\Delta_{i}^\w,T_{\Delta_{i}^\w}^{-}(\Delta_{i-1}^\w),\Delta_{i+1}^\w,\dots,\Delta_m^\w)\,,\]
where $T_{\Delta_{i}^\w}^-(\Delta_{i-1}^\w)\coloneqq \on{fib}(\Delta_{i-1}^\w\longrightarrow \on{Mor}_\C(\Delta_{i-1}^\w,\Delta_{i}^\w)^*\otimes \Delta_{i}^\w)$. Similarly, the co-exceptional collection $(\nabla_1^{\w'},\dots,\nabla_m^{\w'})$ coincides with the image under $\mu_i^\flat$ of the backward mutation of $\nabla_1^\w,\dots,\nabla_m^\w$ at $(\nabla_{i-1}^\w,\nabla_i^\w)$.
\end{enumerate}
\end{lemma}

\begin{proof}
Part (1) follows directly from \cref{thm:weaveequivglsec}.(1), since $\mathcal{F}_{\w}\simeq \mathcal{F}_{\w'}$ in this case. Let us prove Part (2). For $j\in[1,m]$ with $j\not\in\{ i-1,i\}$, the diagram
\[
\begin{tikzcd}
\glsecF \arrow[rr, "\mu_i^\sharp"] \arrow[rd, "\on{ev}_{v_j}"'] &                               & R\Gamma(\rgraph_{\w'},\F_{\w'}) \arrow[ld, "\on{ev}_{v_j}"] \\
& \mathcal{F}_\w(v_j) &            
\end{tikzcd}
\]
commutes. It follows from \eqref{eq:thimble_as_induction} that 
\[\mu_i^\sharp(\Delta_j^{\w})\simeq  \mu_i^\sharp(\on{ind}^L_{v_j}((k,0,0)))\simeq \Delta_j^{\w'}\,,\quad \mbox{ if }j\not\in\{i-1,i\}.\]
The equivalences $\Delta_{i-1}^{\w'}\simeq \mu_{i}^\sharp(T_{\Delta_{i-1}}(\Delta_{i}^\w))$ and $\Delta_{i}^{\w'}\simeq \mu_i^\sharp(\Delta_{i-1}^\w)$ follow from an explicit computation of coCartesian sections, by tracing through the zig-zag of contractions from \eqref{fig:ZigZagCorrespondence_Clockwise} and using the equivalences of global sections from \Cref{prop:contraction}. Let us illustrate part of such a computation, as follows.\\

\noindent Let us suppose that $i=m$ is the top vertex, which allows us to omit the edge $e_{i+1}$ connecting to $v_i$. The computation of the equivalence $\Delta_{i}^{\w'}\simeq \mu_i^\sharp(\Delta_{i-1}^\w)$ amounts to identifying the coCartesian sections under the four equivalences of global sections from the contractions in \eqref{fig:ZigZagCorrespondence_Clockwise}. We describe below the images under the first two of these equivalences. For the first one $\pi_1^*(\mathcal{F})$, we use the perverse schober on the left side of \eqref{eq:proof_schober2}. Similarly, for the second one $(\pi_2)_*\pi_1^*(\mathcal{F})$, we use the perverse schober on the left of \eqref{eq:proof_schober3}. Let us denote the non-singular central vertex in the top middle-left ribbon graph in \eqref{fig:ZigZagCorrespondence_Clockwise} by $\tilde{v}$. We have $\pi_1^*(\mathcal{F})(\tilde{v})\simeq \on{Fun}(\Delta^1,\D)$, and thus we can denote its objects by $(\tilde{v};a\to b)$ with $a\to b$ a morphism in $\D$. Pull-push yields the following coCartesian sections, depicted here partially near $v_{i-1}$: 
\[
\begin{tikzcd}
(e_i;0)                               \\
{(v_{i-1};k[-1],0,0)} \arrow[u] \arrow[d] \\
(e_{i-1};S_{p_{i-1}})                          
\end{tikzcd} \xlongleftarrow{(\pi_1)_\ast} \begin{tikzcd}[column sep=-10]
                      & (e_i;0)                                                                                                      \\
                      & {(\tilde{v};S_{p_{i-1}}[1]\xrightarrow{\on{id}}S_{p_{i-1}}[1]))} \arrow[d] \arrow[ldd, bend right] \arrow[u] \\
                      & {(\tilde{e}_{i-1};S_{p_{i-1}}[2])}                                                                           \\
(e_{i-1};S_{p_{i-1}}) & {(v_{i-1};k[2])} \arrow[u]                                                                                  
\end{tikzcd} \xlongrightarrow{(\pi_2)_*} \begin{tikzcd}[column sep=0]
                      & {(v_i;0,S_{p_{i-1}}[1],0)} \arrow[d] \arrow[ldd, bend right] \\
                      & {(\tilde{e}_{i-1};S_{p_{i-1}}[2])}                       \\
(e_{i-1};S_{p_{i-1}}) & {(v_{i-1};k[2])} \arrow[u]                              
\end{tikzcd}
\]
The shifts that appear are inessential, as they are removed later in the zig-zag. Note however that the coCartesian section on the right takes the value $(v_i;0,S_{p_{i-1}}[1],0)$ at $v_i$, which matches the value of a standard thimble at a vertex (up to a shift), see \eqref{eq:standard_value_at_vertex}. Via such identifications, we indeed conclude in the end the required equivalence $\Delta_{i}^{\w'}\simeq \mu_i^\sharp(\Delta_{i-1}^\w)$.\\

The argument for the costandard weave thimbles in Part (2) is analogous. Finally, Part (3) follows from Part (2), as the corresponding equivalences $\mu_i^\sharp$ and $\mu_i^\flat$ between global sections are inverses of each other.
\end{proof}


\section{Proofs of Theorem \ref{thm:main2}.(2), \ref{thm:main2}.(3) and \ref{thm:main2}.(4)}\label{sec:proof_Lusztig_cycles}

The aim of this section is to prove Parts (2), (3) and (4) for Theorem \ref{thm:main2}. These results were also stated as
\Cref{thm:categorical_Lusztig_cycles_are_SMC,thm:Lusztigcycles_weaveequivalence,thm:weavemutation_SMCmutation1}, respectively.


\subsection{Categorical Lusztig cycles are a simple-minded collection}\label{ssec:LC_from_thimbles}

Let us now prove Theorem \ref{thm:main2}.(2), or cf.~\cref{thm:categorical_Lusztig_cycles_are_SMC}, showing that the categorical Lusztig cycles form a simple-minded collection. Let $\w$ be Demazure weave $\w$ with $m$ trivalent vertices, $(\Delta_1^\w,\ldots,\Delta_m^\w)$ its standard weave thimbles and $(\nabla_1^\w,\ldots,\nabla_m^\w)$ its costandard weave thimbles, as in \cref{def:thimbles}.

\begin{definition}\label{def:SMC}
\begin{enumerate}[(1)]
\item For $i\in[1,m]$, we define $\ssimp_i^\w\in\glsecF$ to be the object
 \[ \ssimp^\w_i\coloneqq T_{\Delta^\w_{1}}^{+,\geq 0}(\cdots (T_{\Delta^\w_{i-1}}^{+,\geq 0}(\Delta^\w_i)))\,.\]

\item For $i\in[1,m]$, we define $\csimp_i^\w\in\glsecF$ to be the object

\[ \csimp^\w_i\coloneqq T_{\nabla^\w_{1}}^{-,\leq 0}(\cdots (T_{\nabla^\w_{i-1}}^{-,\leq 0}(\nabla^\w_i)))\,.\]\qed
\end{enumerate}
\end{definition}

\noindent It is proven in \Cref{prop:simplesfromstandard} that $(\ssimp_1^\w,\ldots,\ssimp_m^\w)$ is a simple-minded collection in $\glsecF$, and so is $(\csimp_1^\w,\ldots,\csimp_m^\w)$. In fact, \Cref{lem:standard_costandard_SMC_coincide} implies the equivalences $\ssimp^\w_i \simeq \csimp^\w_i$ for each $i\in[1,m]$. It will nevertheless be useful for the proofs to consider both types of objects based on their description in \cref{def:SMC}. We also introduce a variant of \cref{def:SMC}, using only partially iterated semi-twists, as follows. For $i\in[1,m]$ and $j\in[i,m]$ we define
\begin{equation}\label{eq:partial_simples}
\ssimp_{i,j}^{\w}\coloneqq T^{+,\geq 0}_{\Delta_{i}^\w}(\cdots(T^{+,\geq 0}_{\Delta_{j-1}^\w}(\Delta_j^\w)))\,,
\end{equation}
where $\ssimp_{j,j}^{\w}=\Delta_j^\w$, and the case $\ssimp_{1,j}^{\w}=\ssimp_j^{\w}$ recovers \cref{def:SMC}.(1). Similarly, we define 
\begin{equation}\label{eq:partial_cosimples}
\csimp_{i,j}^{\w}\coloneqq T^{-,\leq 0}_{\nabla_{i}^\w}(\cdots(T^{-,\leq 0}_{\nabla_{j-1}^\w}(\nabla_j^\w)))\,.
\end{equation}
\noindent To ease notation, we often drop the exponent $\w$ from these notations if $\w$ can be implicitly understood by context, e.g.~we simply denote $\ssimp_j\coloneqq\ssimp_j^\w$ and  $\ssimp_{i,j}\coloneqq\ssimp_{i,j}^{\w}$ in such a case.

\begin{theorem}\label{thm:categorical_Lusztig_cycles_from_exceptionals}
Let $\w$ be a Demazure weave with $m$ trivalent vertices. Then:

\begin{enumerate}
    \item For each $i\in[1,m]$, there exists an equivalence $\simp^\w_i \simeq \ssimp^\w_i$ and, similarly, $\simp^\w_i\simeq \csimp^\w_i\,$.

    \item The categorical Lusztig cycles $\{\simp^\w_i\}$ form a simple-minded collection in $\glsecF$. 
\end{enumerate}
\end{theorem}

\begin{proof}[Proof of \Cref{thm:categorical_Lusztig_cycles_from_exceptionals}]
Let us fix $i\in[1,m]$ and omit $\w$ during the proof, so as to ease notation.\\

\noindent {\bf Part (1).A}. Let us prove the equivalence $\simp^\w_i \simeq \ssimp^\w_i$ in Part (1). First, for $j\in(i,m]$ the values
$$\simp_i(v_j)\simeq\ssimp_i(v_j)\simeq0,\quad j\in(i,m],$$
coincide at $v_j$ because they vanish at $v_j$, cf.~\cref{rmk:support_categorical_Lusztig_cycles}. Similarly, for $j\in(i,m]$ we have
$$\simp_i(e_j)\simeq\ssimp_i(e_j)\simeq0,\quad j\in(i,m].$$
At the vertex $v_i$, the values of the sections $\simp^\w_i$ and $\ssimp^\w_i$ coincide by direct computation:
$$\simp_i(v_i)\simeq \ssimp_i(v_i)\simeq (v_i;k[-1],0,0).$$
Indeed, this follows from using \Cref{def:categorical_Lusztig_cycles} for $\simp_i$, and \eqref{eq:lax_down_value_at_vertex} with $l=0$ for $\ssimp_i$, or alternatively \eqref{eq:thimble_as_induction}, as $\ssimp_i(v_i)\simeq \Delta_i(v_i)$. For their values at the edge $e_i$ we obtain
$$\simp_i(e_{i})\simeq \ssimp_i(e_{i})\simeq S_{p_i}$$
by using \Cref{def:categorical_Lusztig_cycles} for $\simp_i$, and equivalence \eqref{eq:matchingsection_down} for $\ssimp_i$, again with $l=0$.\\

\noindent The goal now is thus showing that
\begin{equation}\label{eq:equiv_three_cycles_at_v_j}
\simp_i(v_j)\simeq \ssimp_i(v_j)\qquad\mbox{ and }\qquad \simp_i(e_{j})\simeq \ssimp_i(e_{j}), \quad \forall j\in[1,i).
\end{equation}

\noindent Let us prove \eqref{eq:equiv_three_cycles_at_v_j} iteratively, starting at $j=i-1$ until we reach $j=1$. For that, fix $j\in(1,i)$, assume that \eqref{eq:equiv_three_cycles_at_v_j} holds for any $k\in(j,i)$, and let $(\beta_{j+1},{\bf a}_{j+1})$ be the weighted braid word associated to the Lusztig cycle $\g_i^\w$ at the horizontal slice directly above the vertex $v_j$. By \cref{def:categorical_Lusztig_cycles}.(1), we have
\[ \simp_i(v_j)= (v_j;\tau_{\geq 0}\on{Mor}(S_{p_j},L_{\beta_{j+1},{\bf a}_{j+1}}),L_{\beta_{j+1},{\bf a}_{j+1}},\on{ev}) \,.\]
The iterative hypothesis implies that $\ssimp_i(e_{j+1})\simeq  \simp_i(e_{j+1})\simeq L_{\beta_{j+1},{\bf a}_{j+1}}$. Since $\Delta^\w_{k}(e_j)\simeq (e_j;0)$ for $k\in[1,j)$, we obtain the equivalence $\ssimp_i(v_j)\simeq \ssimp_{j,i}(v_j)$. In order to compute $\ssimp_{j,i}(v_j)$, note that \eqref{eq:standard_value_at_vertex} implies that $\ssimp_{j+1,i}(v_j)\simeq(v_j;0,L_{\beta_{j+1},{\bf a}_{j+1}},0)$. Thus equivalence \eqref{eq:thimble_as_induction} implies
\begin{equation*}
\on{Mor}(\Delta_j,\ssimp_{j+1,i})\simeq \on{Mor}_{\mathcal{F}_\w(v_j)}((k[-1],0,0),(0,L_{\beta_{j+1},{\bf a}_{j+1}},0))\simeq \on{Mor}_\D(S_{p_j},L_{\beta_{j+1},{\bf a}_{j+1}})\,,
\end{equation*}
and therefore, together with the defining equation \eqref{eq:partial_simples} for $\ssimp_{j,i}$, we obtain
\begin{equation}\label{eq:iterative_Sij}
\ssimp_{j,i}\simeq \on{cof}(\Delta_j\otimes \tau_{\geq 0}\on{Mor}(S_{p_j},L_{\beta_{j+1},{\bf a}_{j+1}})\to \ssimp_{j+1,i})\,.
\end{equation}
\noindent Evaluating \eqref{eq:iterative_Sij} at $v_j$ leads to the equivalence
\begin{equation*}
\ssimp_{j,i}(v_j)\simeq \on{cof}((v_j;\tau_{\geq 0}\on{Mor}(S_{p_j},L_{\beta_{j+1},{\bf a}_{j+1}})[-1],0,0)\to (v_j;0,L_{\beta_{j+1},{\bf a}_{j+1}},0))\simeq \simp_i(v_j)\,.
\end{equation*}
\noindent Therefore $\ssimp_i(v_j)\simeq \ssimp_{j,i}(v_j)\simeq \simp_i(v_j)$, which proves \eqref{eq:equiv_three_cycles_at_v_j} at the vertex $v_j$. The equivalence \eqref{eq:equiv_three_cycles_at_v_j} at the edge $e_j$ follows from that, as the map $\mathcal{F}_\w(v_j\to e_{j})$ provides the required equivalence. This concludes the proof of the equivalence $\simp^\w_i \simeq \ssimp^\w_i$ in Part (1).\\

\noindent {\bf Part (1).B}. For the equivalence $\simp^\w_i \simeq \csimp^\w_i$ in Part (1), the same iterative strategy works but the computations are necessarily different. In this case, the key part of the argument is still showing
\begin{equation}\label{eq:equiv_three_cycles_at_v_j_2}
\simp_i(v_j)\simeq \csimp_i(v_j)\qquad\mbox{ and }\qquad \simp_i(e_{j})\simeq \csimp_i(e_{j}), \quad \forall j\in[1,i),
\end{equation}
as the cases $j\in[i,m]$ are identical to the previous argument. The equivalence $\csimp_i(v_j)\simeq \csimp_{j,i}(v_j)$ still holds, and now \eqref{eq:costandard_value_at_vertex} implies $\csimp_{j+1,i}(v_j)\simeq(v_j;\on{Mor}(S_{p_j},L_{\beta_j,{\bf a}_j}),L_{\beta_j,{\bf a}_j},\on{ev})$. The equivalence \eqref{eq:cothimble_as_induction} implies
\begin{equation*}
\on{Mor}(\csimp_{j+1,i},\nabla_j)\simeq \on{Mor}_{\mathcal{F}_\w(v_j)}((\on{Mor}(S_{p_j},L_{\beta,{\bf a}}),L_{\beta_{j+1},{\bf a}_{j+1}},\on{ev}),(k[-1],0,0))\simeq \on{Mor}(L_{\beta_{j+1},{\bf a}_{j+1}},S_{p_j})^\ast[-1]\,.
\end{equation*}
which, along with \eqref{eq:partial_cosimples}, leads to the equivalences
\begin{equation}\label{eq:iterative_Sijback}
\begin{aligned}\csimp_{j,i}&\simeq \on{fib}(\csimp_{j+1,i}\to \nabla_j\otimes \tau_{\leq 0}(\on{Mor}(L_{\beta_{j+1},{\bf a}_{j+1}},S_{p_j})^\ast[-1]))\\
&\simeq \on{fib}(\csimp_{j+1,i}\to \nabla_j\otimes (\tau_{\leq -1}\on{Mor}(S_{p_j},L_{\beta_{j+1},{\bf a}_{j+1}}))[1])\,,
\end{aligned}
\end{equation}
where the second equivalence uses the fact that $\D$ is $2$-Calabi--Yau. Evaluating \eqref{eq:iterative_Sijback} at $v_j$ leads to the equivalence
\begin{equation*}
\csimp_{j,i}(v_j)\simeq \on{fib}((v_j;\on{Mor}(S_{p_j},L_{\beta,{\bf a}}),L_{\beta,{\bf a}},\on{ev})\to (v_j;\tau_{\leq -1}\on{Mor}(S_{p_j},L_{\beta,{\bf a}}),0,0))\simeq \simp_i(v_j)\,.
\end{equation*}
\noindent Hence $\csimp_i(v_j)\simeq \csimp_{j,i}(v_j)\simeq \simp_i(v_j)$, which proves \eqref{eq:equiv_three_cycles_at_v_j_2} at the vertex $v_j$, and therefore also at the edges, as above. This concludes the proof of Part (1) in the statement. Part (2) follows from Part (1) and \cref{prop:simplesfromstandard}.\end{proof}


\subsection{Invariance of categorical Lusztig cycles under weave equivalences}\label{ssec:Lusztigcycles_weaveequivalence}

Let us now establish Theorem \ref{thm:main2}.(3), see also \Cref{thm:Lusztigcycles_weaveequivalence}. In order to do so, we first prove \cref{lem:restriction_induction_>=_i} and \cref{lem:two_sided_semi_twist_description_of_Lusztig_cycles}, which will be used in the argument for Theorem \ref{thm:main2}.(3).\\

Let $\w$ be a Demazure weave with $m$ trivalent vertices and $\Wgraph$ its associated graph. Given $i\in[1,m]$, consider the subgraph $\rgraph_{\w,\geq i}\subset\rgraph_\w$ form by the vertices $v_i,\dots,v_m$ and all their incident edges. In particular, the edge $e_i$ becomes external in $\rgraph_{\w,\geq i}$. Such a  subgraph specifies the corresponding exit path subcategory $\on{Exit}(\rgraph_{\w,\geq i})\subset \on{Exit}(\rgraph_{\w})$. The weave schober $\F_\w$ for $\w$, which is a functor with domain $\on{Exit}(\rgraph_{\w})$, can consequently be restricted to this subcategory: the restriction of $\mathcal{F}_\w$ along $\on{Exit}(\rgraph_{\w,\geq i})\subset \on{Exit}(\rgraph_{\w})$ will be denoted by $\mathcal{F}_{\w,\geq i}$.

\begin{lemma}\label{lem:restriction_induction_>=_i} Let $\w$ be a Demazure weave and $i\in[1,m]$. Then the following holds:
\begin{enumerate}[(1)]
    \item The restriction of coCartesian sections along the inclusion $\on{Exit}(\rgraph_{\w,\geq i})\subset \on{Exit}(\rgraph_{\w})$ defines a functor 
    \[ \on{res}_{\geq i}\colon \glsecF\to  R\Gamma(\rgraph_{\w,\geq i},\mathcal{F}_{\w,\geq i})\,.\] 
    In addition, the functor $\on{res}_{\geq i}$ admits a fully faithful left adjoint
    \[
\on{ind}^L_{\geq i}\colon R\Gamma(\rgraph_{\w,\geq i},\mathcal{F}_{\w,\geq i})\longrightarrow \glsecF\,.
\]
\item For $j\in[i,m]$, the global section $\ssimp_{i,j}^\w$ lies in the essential image of $\on{ind}^L_{\geq i}$.
\item For $j\in[i,m]$, there exist equivalences
\[ \on{res}_{\geq i}(\ssimp_j^\w)\simeq \on{res}_{\geq i}(\ssimp_{i,j}^\w),\qquad\mbox{ and }\qquad\on{res}_{\geq i}(\csimp_j^\w)\simeq \on{res}_{\geq i}(\csimp_{i,j}^\w)\,\]
in $R\Gamma(\rgraph_{\w,\geq i},\mathcal{F}_{\w,\geq i})$.

\item For $j\in[i,m]$, there exists the following equivalence in $R\Gamma(\rgraph_{\w,\geq i},\mathcal{F}_{\w,\geq i})$:
\begin{equation}\label{eq:ssimpij_equal_csimpij_abovei}
\on{res}_{\geq i}(\ssimp_{i,j}^\w)\simeq \on{res}_{\geq i}(\csimp_{i,j}^\w)
\end{equation}
\end{enumerate}

\end{lemma}

\begin{proof} Part (1) follows from inspecting the construction of $\on{ind}^L_{\geq i}$ in \cite[Lemma 4.20]{Chr25b} and using \cite[Prop.~4.18 \& Rmk.~4.19]{Chr25b}. Part (2) follows from the definition of $\ssimp_{i,j}^\w$ and the fact that each thimble $\Delta_l^\w$, with $l\in[i,m]$, lies in the essential image of $\on{ind}^L_{\geq i}$ by \eqref{eq:thimble_as_induction}. 
Part (3) follows from inspecting the definitions of $\ssimp_j^\w,\csimp_j^\w$ and the fact that $\on{res}_{\geq i}(\Delta_{l}^\w)\simeq \on{res}_{\geq i}(\nabla_{l}^\w)\simeq 0$ for $l\in[1,i)$. Finally, Part (4) follows from Part (3) and the equivalence $\ssimp_j^\w\simeq\csimp_j^\w$, established in \cref{thm:categorical_Lusztig_cycles_from_exceptionals}.
\end{proof}

Let us introduce another variant of \cref{def:SMC}, closely related to \eqref{eq:partial_simples} and \eqref{eq:partial_cosimples}. For that, consider $i\in[1,m]$ and $j\in[i+2,m]$, and define

\begin{equation}\label{eq:partial_simples_full}
\zssimp_{i+1,j}^{\w}\coloneqq T^+_{\Delta_{i+1}^\w}(\ssimp_{i+2,j}^\w)
\,.
\end{equation}
where $T^+_{\Delta_{i+1}^\w}$ denotes the full twist along ${\Delta_{i+1}^\w}$, not a semi-twist. Similarly, we define
\begin{equation}\label{eq:partial_cosimples_full}
\zcsimp_{i+1,j}^{\w}\coloneqq T^-_{\nabla_{i+1}^{\w}}(\csimp_{i+2,j}^\w)\,.
\end{equation}
with $T^-_{\nabla_{i+1}^{\w}}$ denoting the full cotwist along $\nabla_{i+1}^{\w}$. Before stating the next result, relating \cref{eq:partial_simples} to \cref{eq:partial_simples_full} and \cref{eq:partial_cosimples} to \cref{eq:partial_cosimples_full},
we recall that if $\w\to\w'$ is a backward pulley move or backwards weave mutation, exchanging the heights of the vertices $p_i$ and $p_{i+1}$, then \cref{thm:weaveequivglsec} provides an equivalence of global sections:
\[ \mu_{i+1}^\flat\colon \glsecF\stackrel{\simeq}{\lr}R\Gamma(\rgraph_{\w'},\mathcal{F}_{\w'}),\]
see also \cref{rmk:backward_pulley_schober}.

\begin{proposition}\label{lem:two_sided_semi_twist_description_of_Lusztig_cycles}
Let $\w$ be a Demazure weave with $m$ trivalent vertices, $i\in[1,m]$ and $j\in(i+1,m]$. Then the following holds:

\begin{enumerate}[(1)]
    \item In the category $\glsecF$, there exists equivalences
\begin{equation}\label{eq:Zssimp_equivalence}
\ssimp_{i+1,j}^{\w}\simeq T^{-,\leq 0}_{\Delta_{i+1}^\w}(\zssimp_{i+1,j}^{\w})
\end{equation}
and similarly
\begin{equation}\label{eq:Zcsimp_equivalence}
\csimp_{i+1,j}^{\w}\simeq T^{+,\geq 0}_{\nabla_{i+1}^\w}(\zcsimp_{i+1,j}^{\w}).
\end{equation}

\item Let $\w\to\w'$ be a backward pulley move or backwards weave mutation, exchanging the heights of the vertices $p_i,p_{i+1}$. Let 
\[ \phi\colon R\Gamma(\rgraph_{w'},\mathcal{F}_{w'})\simeq \glsecF\]
be the corresponding equivalence of global sections arising from \cref{thm:weaveequivglsec}. Then there exist equivalences in $R\Gamma(\rgraph_{\w,\geq i},\mathcal{F}_{\w,\geq i})$
\begin{equation}\label{eq:restriction_zssimp_is_zcsimp}
\on{res}_{\geq i}(\zssimp_{i+1,j}^{\w})\simeq \on{res}_{\geq i}((\mu_{i+1}^\flat)^{-1}(\zcsimp_{i+1,j}^{\w'}))\,.
\end{equation}
and
\begin{equation}\label{eq:ssimp_via_zssimp}
(\mu_{i+1}^\flat)^{-1}(\ssimp_{i,j}^{\w'})\simeq T^{-,\leq 0}_{\Delta_{i+1}^\w}T^{+,\geq 0}_{\Delta_i^\w}(\zssimp_{i+1,j}^{\w})\,.
\end{equation}
\end{enumerate}
\end{proposition}

\begin{proof} Let us start with \eqref{eq:Zssimp_equivalence} in Part (1), and again drop $\w$ from the notation for ease. By definition, the right hand side of \eqref{eq:Zssimp_equivalence} is 
$$T^{-,\leq 0}_{\Delta_{i+1}}(\zssimp_{i+1,j})=\fib(\zssimp_{i+1,j}^{\w}\lr\tau_{\leq0}(\textcolor{blue}{\Mor(\zssimp_{i+1,j},\Delta_{i+1})}^*\otimes \Delta_{i+1})).$$
To compute the morphism object highlighted in blue above, first use the defining cofiber sequence \eqref{eq:partial_simples_full}
\begin{equation}\label{eq:cofiberseq_defining_zssimp}
\Mor(\Delta_{i+1},\ssimp_{i+2,j})\otimes \Delta_{i+1} \stackrel{\on{ev}}{\lr} \ssimp_{i+2,j}\lr \zssimp_{i+1,j}=\cof(\on{ev})
\end{equation}
and apply the functor $\on{Mor}(\mhyphen,\Delta_{i+1}^\w)$ to \eqref{eq:cofiberseq_defining_zssimp}. Since $\on{Mor}(\Delta_k^\w,\Delta_{i+1}^\w)\simeq 0$ if $k\in(i+1,m]$, we obtain $\on{Mor}(\ssimp_{i+2,j},\Delta_{i+1}^\w)\simeq0$ and thus applying the functor $\on{Mor}(\mhyphen,\Delta_{i+1}^\w)$ to \eqref{eq:cofiberseq_defining_zssimp} yields the equivalence
\[ \textcolor{blue}{\on{Mor}(\zssimp_{i+1,j},\Delta_{i+1})}\simeq \on{Mor}(\on{Mor}(\Delta_{i+1},\ssimp_{i+2,j})\otimes \Delta_{i+1}[1],\Delta_{i+1})\simeq \textcolor{cyan}{\on{Mor}(\Delta_{i+1},\ssimp_{i+2,j})^*[-1]}\,,\]
where we used $\Mor(\Delta_{i+1},\Delta_{i+1})\simeq k$ in the second equivalence above. In consequence, the right hand side of \eqref{eq:Zssimp_equivalence} can be expressed as
\begin{align*}
T^{-,\leq 0}_{\Delta_{i+1}}(\zssimp_{i+1,j})& \simeq \fib(\zssimp_{i+1,j}^{\w}\lr\tau_{\leq0}((\textcolor{cyan}{\on{Mor}(\Delta_{i+1},\ssimp_{i+2,j})^*[-1]})^*)\otimes \Delta_{i+1})\\
& \simeq \fib(\zssimp_{i+1,j}^{\w}\lr\tau_{\leq0}(\on{Mor}(\Delta_{i+1},\ssimp_{i+2,j})[1])\otimes \Delta_{i+1})\\
& \simeq \fib(\zssimp_{i+1,j}^{\w}\lr\tau_{\leq-1}(\on{Mor}(\Delta_{i+1},\ssimp_{i+2,j}))[1]\otimes \Delta_{i+1})\\
& \simeq \cof(\tau_{\geq0}(\on{Mor}(\Delta_{i+1},\ssimp_{i+2,j}))\otimes \Delta_{i+1})\lr\ssimp_{i+2,j})\\
& \simeq \ssimp_{i+1,j}.
\end{align*}
This establishes the equivalence \eqref{eq:Zssimp_equivalence}. The equivalence \eqref{eq:Zcsimp_equivalence} is proven similarly.\\

\noindent For Part (2), we ease notation by omitting $\w$ and $\w'$ but denoting $\ssimp_{i,j}^{'}\coloneqq\ssimp_{i,j}^{\w'}$ and so on, i.e.~a prime in the exponent indicates the object is associated to $\w'$, instead of $\w$. We also omit the equivalence $(\mu_{i+1}^\flat)^{-1}$ in the notation, as the prime notation above already indicates where objects belong.\\

\noindent We prove the equivalence \eqref{eq:restriction_zssimp_is_zcsimp} as follows. First, \Cref{lem:change_in_thimbles_under_weave_equivalence}.(3) implies $\Delta_{i+1}\simeq\Delta_i^{'}$ and thus, by \cref{lem:restriction_induction_>=_i}.(1), $\on{res}_{\geq i}(\Delta_{i+1})\simeq  \on{res}_{\geq i}(\Delta_i^{'})$. Since $\on{res}_{\geq i}(\Delta^{'}_{i})\simeq  \on{res}_{\geq i}(\nabla_i^{'})$, we obtain
\begin{equation}\label{eq:proof_thimble_restriction_above_i}
\on{res}_{\geq i}(\Delta_{i+1})\simeq  \on{res}_{\geq i}(\nabla_i^{'}).
\end{equation}
Second, \Cref{lem:change_in_thimbles_under_weave_equivalence}.(3) also implies that $\Delta_{k}\simeq\Delta_k'$ if $k\in[i+2,m]$, and therefore
\begin{equation}\label{eq:S_to_backS_above_i}
\on{res}_{\geq i}(\ssimp_{i+2,j})\simeq \on{res}_{\geq i}(\ssimp_{i+2,j}^{'})\simeq \on{res}_{\geq i}(T^-_{\nabla_{i}^{'}}T^-_{\nabla_{i+1}^{'}}\csimp_{i+2,j}^{'})\,,
\end{equation}
where the second equivalence in \eqref{eq:S_to_backS_above_i} follows from \eqref{eq:standards_from_costandards} restricted to $\rgraph_{\w',\geq i}$. Third, the equivalence
\begin{equation}\label{eq:twist_cotwist_identity}
 \on{res}_{\geq i}(T^+_{\nabla_i^{'}}T^-_{\nabla_{i}^{'}})\simeq \mbox{id}
\end{equation}
follows from the definitions of the adjoint twist and cotwist functors. Applying these equivalences yields
\begin{align*}
\on{res}_{\geq i}(\zssimp_{i+1,j})&\stackrel{\eqref{eq:partial_simples_full}}{=} \on{res}_{\geq i}(T^+_{\Delta_{i+1}}(\ssimp_{i+2,j})) \stackrel{\eqref{eq:proof_thimble_restriction_above_i}}{\simeq} \on{res}_{\geq i}(T^+_{\nabla_i^{'}}(\ssimp_{i+2,j}))\\
& \stackrel{\eqref{eq:S_to_backS_above_i}}{\simeq} \on{res}_{\geq i}(T^+_{\nabla_i^{'}}T^-_{\nabla_{i}^{'}}T^-_{\nabla_{i+1}^{'}}(\csimp_{i+2,j}^{'}))\stackrel{\eqref{eq:twist_cotwist_identity}}{\simeq} \on{res}_{\geq i}(T^-_{\nabla_{i+1}^{'}}(\csimp_{i+2,j}^{'}))\stackrel{\eqref{eq:partial_cosimples_full}}{=} \on{res}_{\geq i}(\zcsimp_{i+1,j}^{'})\,.
\end{align*}
\noindent This establishes equivalence \eqref{eq:restriction_zssimp_is_zcsimp}. In order to show equivalence \eqref{eq:ssimp_via_zssimp}, we first note that

\begin{align*}
\on{res}_{\geq i}(\ssimp_{i,j}^{'})& \stackrel{\eqref{eq:ssimpij_equal_csimpij_abovei}}{\simeq}\on{res}_{\geq i}(\csimp_{i,j}^{'})\stackrel{\eqref{eq:Zcsimp_equivalence}}{\simeq} \on{res}_{\geq i}(T^{-,\leq 0}_{\nabla_i^{'}}T^{+,\geq 0}_{\nabla_{i+1}^{'}}(\zcsimp_{i+1,j}^{'}))\\
& \stackrel{\eqref{eq:proof_thimble_restriction_above_i}}{\simeq} \on{res}_{\geq i}(T^{-,\leq 0}_{\Delta_{i+1}}T^{+,\geq 0}_{\Delta_{i}}(\zcsimp_{i+1,j}^{'}))\stackrel{\eqref{eq:restriction_zssimp_is_zcsimp}}{\simeq} \on{res}_{\geq i}(T^{-,\leq 0}_{\Delta_{i+1}}T^{+,\geq 0}_{\Delta_{i}}(\zssimp_{i+1,j}))\,,
\end{align*}
where the second to last equivalence also used $\on{res}_{\geq i}(\nabla_{i+1}^{'})\simeq \on{res}_{\geq i}(\Delta_i)$, which is implied by \Cref{lem:change_in_thimbles_under_weave_equivalence}. Since both $\ssimp_{i,j}^{'}$ and $T^{-,\leq 0}_{\Delta_{i+1}}T^{+,\geq 0}_{\Delta_{i}}(\zssimp_{i+1,j})$ lie in the image of $\on{ind}^L_{\geq i}$, the equivalence  \eqref{eq:ssimp_via_zssimp} follows from the fact that the functor $\on{ind}^L_{\geq i}$ is fully faithful, cf.~\cref{lem:restriction_induction_>=_i}.(1).
\end{proof}

\begin{proof}[Proof of Theorem \ref{thm:main2}.$(3)$] To ease notation, we omit the equivalence of global sections of the weave schobers induced by the given weave equivalence $\w\sim\w'$, cf.~\cref{thm:weaveequivglsec}.(3).
By decomposing the given weave equivalence $\w\simeq \w'$ into height-preserving weave equivalences and pulley moves, it suffices to prove the assertion for these two types of weave equivalences. In the former case, if the weave equivalence $\w\sim\w'$ is height-preserving, \cref{lem:change_in_thimbles_under_weave_equivalence} implies that there is an identification of their standard weave thimbles $\Delta_i^{\w}\simeq \Delta_i^{\w'}$ for all $i\in[1,m]$. Hence, by construction, the objects of the simple minded collections $\ssimp_i^{\w}\simeq \ssimp_i^{\w'}$ are equivalent for all $i\in[1,m]$. \Cref{thm:categorical_Lusztig_cycles_from_exceptionals} thus implies that ${\bm L}_i^{\w}\simeq {\bm L}_i^{\w'}$ for all $i\in[1,m]$, concluding Theorem \ref{thm:main2}.$(3)$ for height-preserving equivalences. The case of a pulley move, which is more subtle, is proven as follows.\\

Consider a pulley move $\w\to \w'$ exchanging the heights of the vertices $p_i,p_{i+1}$, for some $i\in[1,m]$. Swapping $\w$ and $\w'$ if necessary, we can and do assume that the weave equivalence is a backward pulley move from $\w$ to $\w'$. As above, we ease notation by omitting $\w$ and $\w'$, keeping the prime as an exponent in the latter case. Since \cref{lem:change_in_thimbles_under_weave_equivalence} implies that $\Delta_{k}\simeq \Delta_k^{'}$ for all $k\in[1,i)$, we certainly have $\ssimp_k\simeq \ssimp_k^{'}$, and thus ${\bm L}_k\simeq {\bm L}_k^{'}$ for all $k\in[1,i)$ by \Cref{thm:categorical_Lusztig_cycles_from_exceptionals}. Our next task is showing that ${\bm L}_j\simeq {\bm L}_j^{'}$ for all $j\in[i,m]$.\\

For the cases $j\in\{i,i+1\}$, it suffices to show the equivalences
\begin{equation}\label{eq:equivalenceLi}
\ssimp_{i,i+1}\simeq \ssimp_{i,i}^{'}\qquad\mbox{ and }\qquad \ssimp_{i,i}\simeq \ssimp_{i,i+1}^{'}\,.
\end{equation}
Indeed, since $\Delta_{k}\simeq \Delta_k^{'}$ for all $k\in[1,i)$, as already used above, the equivalences in \eqref{eq:equivalenceLi} respectively imply ${\bm L}_{i+1}\simeq {\bm L}_{i+1}^{'}$ and ${\bm L}_i\simeq {\bm L}_i^{'}$. To establish the first equivalence in \eqref{eq:equivalenceLi}, we use \cref{lem:connectiveHom,lem:morphism_objects_between_standard_thimbles} to obtain the vanishing 
\begin{equation}\label{eq:vanishing_standard_hom}
\tau_{\geq 0} \on{Mor}(\Delta_i,\Delta_{i+1})\simeq \tau_{\geq 0}\on{Mor}_\D(S_{p_{i}},S_{p_{i+1}})\simeq 0.
\end{equation}
From such vanishing, it thus follows that 
\[ \ssimp_{i,i+1}\stackrel{\eqref{eq:partial_simples}}{\simeq} T^{+,\geq 0}_{\Delta_i}(\Delta_{i+1})\stackrel{\eqref{eq:vanishing_standard_hom}}{\simeq} \Delta_{i+1}\simeq \Delta_i^{'}\stackrel{\eqref{eq:partial_simples}}{\simeq} \ssimp_{i,i}^{'},\] 
where we used \cref{lem:change_in_thimbles_under_weave_equivalence} to obtain $\Delta_{i+1}\simeq \Delta_i^{'}$. For the second equivalence in \eqref{eq:equivalenceLi}, we proceed as follows. First, note that \cref{lem:change_in_thimbles_under_weave_equivalence} implies $\Delta_{i+1}^{'}\simeq T^-_{\Delta_{i+1}}(\Delta_{i})$. Second, since $\on{Mor}(\Delta_{i+1},\Delta_i^)\simeq 0$ by exceptionality, we obtain
\begin{equation}\label{eq:equivalence_for_Liplusone}
T^+_{\Delta_{i+1}}(\Delta_{i+1}^{'})\simeq T^+_{\Delta_{i+1}}T^-_{\Delta_{i+1}}(\Delta_{i})\simeq \Delta_i.
\end{equation}\\
Third, the vanishing in \eqref{eq:vanishing_standard_hom} implies that
\begin{equation*}\label{eq:vanishing_negdegrees}
\on{Mor}(\Delta_{i+1},\Delta_{i+1}^{'})\simeq \on{Mor}(\Delta_{i+1},T^-_{\Delta_{i+1}}(\Delta_{i}))\simeq \on{fib}(\underbrace{\on{Mor}(\Delta_{i+1},\Delta_i)}_{\simeq 0}\to \underbrace{\on{Mor}(\Delta_{i+1},\Delta_{i+1})}_{\simeq k}\otimes \on{Mor}(\Delta_i,\Delta_{i+1})^*)
\end{equation*}
is concentrated in positive degrees. Since $\Delta_i^{'}=\Delta_{i+1}$ by \cref{lem:change_in_thimbles_under_weave_equivalence}, this thus implies the equivalence
\begin{equation}\label{eq:semitwist_is_twist_case}
T^{+,\geq 0}_{\Delta_{i}^{'}}(\Delta_{i+1}^{'})\simeq T^{+}_{\Delta_{i}^{'}}(\Delta_{i+1}^{'}).
\end{equation}
In consequence, we obtain
\begin{equation}\label{eq:second_equivalence_Liplusone}
\ssimp_{i,i+1}^{'}\stackrel{\eqref{eq:partial_simples}}{\simeq} T^{+,\geq 0}_{\Delta_{i}^{'}}(\Delta_{i+1}^{'})\stackrel{\eqref{eq:semitwist_is_twist_case}}{\simeq} T^{+}_{\Delta_{i}^{'}}(\Delta_{i+1}^{'})\stackrel{\eqref{eq:equivalence_for_Liplusone}}{\simeq} \Delta_i\stackrel{\eqref{eq:partial_simples}}{\simeq}\ssimp_{i,i}\,.
\end{equation}
\noindent \Cref{thm:categorical_Lusztig_cycles_from_exceptionals} and \eqref{eq:equivalenceLi}, now established, thus imply ${\bm L}_i\simeq {\bm L}_i^{'}$ and ${\bm L}_{i+1}\simeq {\bm L}_{i+1}^{'}$. The remaining task is establishing that ${\bm L}_j\simeq {\bm L}_j^{'}$ for all $j\in[i+2,m]$. For that it suffices to prove the equivalence
\begin{equation}\label{eq:equivalence_for_Lj_higherj}
\ssimp_{i,j}\simeq \ssimp_{i,j}^{'},\qquad\forall j\in[i+2,m].
\end{equation}
By \eqref{eq:Zssimp_equivalence} and \eqref{eq:ssimp_via_zssimp}, equivalence \eqref{eq:equivalence_for_Lj_higherj} is itself equivalent to
\begin{equation}\label{eq:twists_of_Z_j_commute}
\framebox{$T^{+,\geq 0}_{\Delta_i}T^{-,\leq 0}_{\Delta_{i+1}}(\zssimp_{i+1,j})\simeq T^{-,\leq 0}_{\Delta_{i+1}}T^{+,\geq 0}_{\Delta_i}(\zssimp_{i+1,j})$}\qquad\forall j\in[i+2,m].\end{equation}
\noindent The remainder of the proof is to establish \eqref{eq:twists_of_Z_j_commute}. For that, we will first prove the following two equivalences
\begin{equation}\label{eq:equivalence1_proof_invarianceLC}
\tau_{\geq 0}\on{Mor}(\Delta_i,\zssimp_{i+1,j})\simeq \tau_{\geq 0} \on{Mor}(\Delta_i,T^{-,\leq 0}_{\Delta_{i+1}}(\zssimp_{i+1,j}))
\end{equation}
\begin{equation}\label{eq:equivalence2_proof_invarianceLC}
\tau_{\leq 0}\on{Mor}(\zssimp_{i+1,j},\Delta_{i+1})^*\simeq \tau_{\leq 0}\on{Mor}(T^{+,\geq 0}_{\Delta_i}(\zssimp_{i+1,j}),\Delta_{i+1})^*    
\end{equation}
To prove \eqref{eq:equivalence1_proof_invarianceLC}, consider the following cofiber sequence, obtained by applying $\on{Mor}(\Delta_i,\mhyphen)$ to the defining fiber sequence for $T^{-,\leq 0}_{\Delta_{i+1}}(\zssimp_{i+1,j})$:
\begin{equation}\label{eq:cofiberseq_equivalence1proof}
\textcolor{blue}{\on{Mor}(\Delta_i,\tau_{\leq 0}\on{Mor}(\zssimp_{i+1,j},\Delta_{i+1})^*\otimes \Delta_{i+i})[-1]}\to \on{Mor}(\Delta_i,T^{-,\leq 0}_{\Delta_{i+1}}(\zssimp_{i+1,j}))\stackrel{f}{\to} \on{Mor}(\Delta_i,\zssimp_{i+1,j})
\end{equation}
The first term is equivalent to
$$\textcolor{blue}{\on{Mor}(\Delta_i,\tau_{\leq 0}\on{Mor}(\zssimp_{i+1,j},\Delta_{i+1})^*\otimes \Delta_{i+i})[-1]}\simeq \on{Mor}(\Delta_i,\Delta_{i+1})[-1]\otimes \tau_{\leq 0}\on{Mor}(\zssimp_{i+1,j},\Delta_{i+1})^*$$
and thus \eqref{eq:vanishing_standard_hom} implies that such term is concentrated in degrees less equal than $-2$. Hence, the truncation $\tau_{\geq0}(f)$ of the morphism $f$ in \eqref{eq:cofiberseq_equivalence1proof} is an equivalence, thus establishing equivalence \eqref{eq:equivalence1_proof_invarianceLC}.\\

\noindent To prove \eqref{eq:equivalence2_proof_invarianceLC}, we similarly consider the following cofiber sequence, obtained by applying $\on{Mor}(\mhyphen,\Delta_{i+1})^*$ to the defining cofiber sequence for $T^{+,\geq 0}_{\Delta_i}(\zssimp_{i+1,j})$:
\begin{equation}\label{eq:cofiberseq_equivalence2proof}
\on{Mor}(\zssimp_{i+1,j},\Delta_{i+1})^* \stackrel{g}{\to} \on{Mor}(T^{+,\geq 0}_{\Delta_i}(\zssimp_{i+1,j}),\Delta_{i+1})^*\to \textcolor{cyan}{\on{Mor}(\tau_{\geq 0}\on{Mor}(\Delta_i,\zssimp_{i+1,j})\otimes \Delta_i,\Delta_{i+1})^*[1]}\,.
\end{equation}
The third term is equivalent to
$$\textcolor{cyan}{\on{Mor}(\tau_{\geq 0}\on{Mor}(\Delta_i,\zssimp_{i+1,j})\otimes \Delta_i,\Delta_{i+1})^*[1]}\simeq \tau_{\geq 0}\on{Mor}(\Delta_i,\zssimp_{i+1,j})\otimes \on{Mor}(\Delta_i,\Delta_{i+1})^*[1]$$
and \eqref{eq:vanishing_standard_hom} implies that it is concentrated in degrees greater equal than $2$. Therefore, the truncation $\tau_{\leq0}(g)$ of the morphism $g$ in \eqref{eq:cofiberseq_equivalence2proof} is an equivalence. This establishes \eqref{eq:equivalence2_proof_invarianceLC}.\\

Finally, let us conclude the required equivalence \eqref{eq:twists_of_Z_j_commute} from  \eqref{eq:equivalence1_proof_invarianceLC} and  \eqref{eq:equivalence2_proof_invarianceLC}. For that, consider the following commutative square, which is not biCartesian:

\begin{equation}\label{eq:square_endproofequivalence}
\begin{tikzcd}
{\tau_{\geq 0} \on{Mor}(\Delta_i,T^{-,\leq 0}_{\Delta_{i+1}}(\zssimp_{i+1,j}))\otimes \Delta_i} \arrow[d] \arrow[r] & \zssimp_{i+1,j} \arrow[d]                                                   \\
0 \arrow[r]                                                             & {\tau_{\leq 0}\on{Mor}(T^{+,\geq 0}_{\Delta_i}(\zssimp_{i+1,j}),\Delta_{i+1})^*\otimes \Delta_{i+1}}
\end{tikzcd}
\end{equation}

\noindent The equivalence \eqref{eq:equivalence1_proof_invarianceLC}, applied to the upper left corner of \eqref{eq:square_endproofequivalence}, implies that \eqref{eq:square_endproofequivalence} is equivalent to

\begin{equation}\label{eq:square1_endproofequivalence}
\begin{tikzcd}
{\tau_{\geq 0}\on{Mor}(\Delta_i,\zssimp_{i+1,j})\otimes \Delta_i} \arrow[d] \arrow[r] & \zssimp_{i+1,j} \arrow[d]                                                   \\
0 \arrow[r]                                                             & {\tau_{\leq 0}\on{Mor}(T^{+,\geq 0}_{\Delta_i}(\zssimp_{i+1,j}),\Delta_{i+1})^*\otimes \Delta_{i+1}}
\end{tikzcd}
\end{equation}

\noindent By the definition of the semi-twist along $\Delta_{i}$, the horizontal cofiber of \eqref{eq:square1_endproofequivalence} yields the coevaluation morphism
$$T_{\Delta_i}^{+,\geq,0}(\zssimp_{i+1,j})\lr {\tau_{\leq 0}\on{Mor}(T^{+,\geq 0}_{\Delta_i}(\zssimp_{i+1,j}),\Delta_{i+1})^*\otimes \Delta_{i+1}}.$$
By definition of the semi-twist along $\Delta_{i+1}$, the (vertical) fiber of such a morphism is precisely the right hand side of \eqref{eq:twists_of_Z_j_commute}. That is, the right hand side of \eqref{eq:twists_of_Z_j_commute} is equivalent to the vertical fiber of the horizontal cofiber of \eqref{eq:square1_endproofequivalence}.\\

\noindent Similarly, the equivalence \eqref{eq:equivalence2_proof_invarianceLC}, applied to the lower right corner of \eqref{eq:square_endproofequivalence}, implies that \eqref{eq:square_endproofequivalence} is equivalent to the square
\begin{equation}\label{eq:square2_endproofequivalence}
\begin{tikzcd}
{\tau_{\geq 0} \on{Mor}(\Delta_i,T^{-,\leq 0}_{\Delta_{i+1}}(\zssimp_{i+1,j}))\otimes \Delta_i} \arrow[d] \arrow[r] & \zssimp_{i+1,j} \arrow[d]                                                   \\
0 \arrow[r]                                                             & {\tau_{\leq 0}\on{Mor}(\zssimp_{i+1,j},\Delta_{i+1})^*\otimes \Delta_{i+1}}
\end{tikzcd}
\end{equation}
\noindent The vertical fiber of \eqref{eq:square2_endproofequivalence} gives the evaluation morphism
$${\tau_{\geq 0} \on{Mor}(\Delta_i,T^{-,\leq 0}_{\Delta_{i+1}}(\zssimp_{i+1,j}))\otimes \Delta_i}\lr T_{\Delta_{i+1}}^{-,\leq0}(\zssimp_{i+1,j}).$$
\noindent The (horizontal) cofiber of such a morphism is precisely the left hand side of \eqref{eq:twists_of_Z_j_commute}. That is, the left hand side of \eqref{eq:twists_of_Z_j_commute} is equivalent to the horizontal cofiber of the vertical fiber of \eqref{eq:square1_endproofequivalence}.\\

\noindent Since the horizontal cofiber of the vertical fiber of a square in a stable $\infty$-category is equivalent to the vertical fiber of the horizontal cofiber of the same square, the computations above show that both sides of \eqref{eq:twists_of_Z_j_commute} must themselves be equivalent. This establishes \eqref{eq:twists_of_Z_j_commute}, and thus the equivalence ${\bm L}_j\simeq {\bm L}_j^{'}$ between the categorical Lusztig cycles for $j\in[i+2,m]$. By considering the argument altogether, we obtain the equivalence ${\bm L}_j\simeq {\bm L}_j^{'}$ between the categorical Lusztig cycles for all $j\in[1,m]$, as required.
\end{proof}


\subsection{Tilting of categorical Lusztig cycles under weave mutation}\label{ssec:Lusztigcycles_weavemutation}

Let us now prove Theorem \ref{thm:main2}.(4), cf.~also \Cref{thm:weavemutation_SMCmutation1}. At a technical level, an important difference in the following proof of Theorem \ref{thm:main2}.(4), in contrast to the proof of Theorem \ref{thm:main2}.(3) in \cref{ssec:Lusztigcycles_weaveequivalence}, is that equivalence \eqref{eq:twists_of_Z_j_commute} does not hold for a weave mutation. In a nutshell, under a weave mutation, the difference between the two terms in \eqref{eq:twists_of_Z_j_commute} is non-trivial and, as we will now establish, leads to the simple-minded collection of categorical Lusztig cycles undergoing a mutation.\\

\begin{proof}[Proof of Theorem \ref{thm:main2}.$(4)$] Let us prove the result for a backward weave mutation $\w\to\w'$, exchanging the vertices $p_{i-1}$ and $p_{i}$. The case of a forward weave mutation follows from it. As in the proofs above, we omit explicitly writing the equivalence $\glsecF\simeq R\Gamma(\rgraph_{\w'},\mathcal{F}_{\w'})$, and the superscripts $\w,\w'$, keeping the prime as an exponent for the latter case. By \cref{thm:categorical_Lusztig_cycles_from_exceptionals}, our task is equivalent to showing that the backward mutation of the simple-minded collection $\{\ssimp_j\}$ is precisely $\{\ssimp_j^{'}\}$. Specifically, for each $j\in[1,m]$, we must show that

\[
\ssimp_j^{'}\simeq\begin{cases}
\ssimp_i[-1] & \mbox{ if }j=i\,,\\
\on{cof}(\on{Ext}^1(\ssimp_i,\ssimp_j)\otimes \ssimp_i[-1] \stackrel{\on{ev}}{\to} \ssimp_j) & \mbox{ if }j\neq i\,.\\
\end{cases}
\]

\noindent To ease notation, we abbreviate $\psi^\flat_{\ssimp_{i}}(\ssimp_j)\coloneqq\on{cof}(\on{Ext}^1(\ssimp_i,\ssimp_j)\otimes \ssimp_i[-1] \to \ssimp_j)$. As in the proofs above, the cases $j\in[1,i)$ are reasonably straightforward, the cases $j\in\{i-1,i\}$ have more content, and the challenging cases are $j\in[i+1,m]$. Let us study each of these cases.\\

{\bf Case $j\in[1,i-1)$.} In this range, \cref{lem:change_in_thimbles_under_weave_equivalence} implies that $\ssimp_j\simeq \ssimp_j^{'}$. By the local model of backward weave mutation, the Lusztig cycle $\g_i$ ends at the vertex $p_{i-1}$, and therefore $\on{ev}_{v_{l}}(\ssimp_i)\simeq 0$ for all vertices with $l\in[1,i-1)$. Hence we obtain the vanishing $\on{Ext}^1(\ssimp_{i},\ssimp_j)\simeq 0$ in this range, and thus
$$\psi^\flat_{\ssimp_{i}}(\ssimp_j)\simeq \on{cof}(0\otimes \ssimp_i[-1] \to \ssimp_j)\simeq \ssimp_j \stackrel{(\ref{lem:change_in_thimbles_under_weave_equivalence})}{\simeq}\ssimp_j^{'},\qquad\forall j\in[1,i-1),$$
as we needed to show.\\

{\bf Case $j\in\{i-1,i\}$.} Let us start with $j=i-1$. In the local model for a backward weave mutation
\begin{equation}\label{eq:tilting_simple_i}
\ssimp_{i}\simeq \on{cof}(\Delta_{i-1}\to \Delta_{i})\,.
\end{equation} Now, on the one hand, since $\on{Mor}(\Delta_i,\ssimp_{i-1})\simeq 0$, we obtain
\[
\on{Ext}^1(\ssimp_i,\ssimp_{i-1})\simeq \on{Ext}^1(\Delta_{i-1}[1],\ssimp_{i-1})\simeq  k\,,
\]
and thus
\begin{equation}\label{eq:flat_tilting_iminus1}
\psi^\flat_{\ssimp_i}(\ssimp_{i-1})\simeq \on{cof}(\ssimp_i[-1] \to \ssimp_{i-1})\,.
\end{equation}
On the other hand, \cref{lem:change_in_thimbles_under_weave_equivalence} implies $\Delta_{i-1}^{'}\simeq \Delta_i$ and so \eqref{eq:tilting_simple_i} leads to a cofiber sequence
\begin{equation}\label{eq:cofiberseq_tilting_iminus1}
\Delta_{i-1}\longrightarrow \Delta_{i-1}^{'}\longrightarrow \ssimp_i\,.
\end{equation}
Since $\on{Mor}(\Delta_j,\ssimp_i)\simeq 0$ in the range $j\in[1,i-2)$, \eqref{eq:cofiberseq_tilting_iminus1} implies the existence of a cofiber sequence
\[
\ssimp_{i-1}\longrightarrow \ssimp_{i-1}^{'}\longrightarrow\ssimp_i\,,
\]
which itself implies $\ssimp_{i-1}'\simeq \on{cof}(\ssimp_i[-1] \to \ssimp_{i-1})$. Therefore
$$\ssimp_{i-1}'\simeq\on{cof}(\ssimp_i[-1] \to \ssimp_{i-1})\stackrel{\eqref{eq:flat_tilting_iminus1}}{\simeq}\psi^\flat_{\ssimp_i}(\ssimp_{i-1}),$$
concluding the case $j=i-1$. For the case $j=i$, we proceed as follows. First, consider the two sequences of equivalences
\begin{equation}\label{eq:tilting_equivalences_res1}
\on{res}_{\geq i-1}(\Delta_{i-1})\simeq \on{res}_{\geq i-1}(\nabla_{i-1})\simeq \on{res}_{\geq i-1}(\nabla_{i}^{'})
\end{equation}
\begin{equation}\label{eq:tilting_equivalences_res2}
\on{res}_{\geq i-1}(\Delta_{i})\simeq \on{res}_{\geq i-1}(\Delta_{i-1}^{'})\simeq \on{res}_{\geq i-1}(\nabla_{i-1}^{'})
\end{equation}
which follow from \Cref{lem:restriction_induction_>=_i} and the definition of the weave thimbles, cf.~\cref{def:weave_thimbles}. Second, in the local model for a backward weave mutation we have
\begin{equation}\label{eq:tilting_equivalence_cosimplei}
\csimp_{i}^{'}\simeq \on{fib}(\nabla_{i}^{'}\to \nabla_{i-1}^{'})\,.
\end{equation}
From this, again using \Cref{lem:restriction_induction_>=_i}, we therefore obtain:
\begin{align*}
\ssimp_{i}&\simeq \on{ind}^L_{\geq i-1}\on{res}_{\geq i-1}(\ssimp_{i})\\
&\stackrel{\eqref{eq:tilting_simple_i}}{\simeq}  \on{ind}^L_{\geq i-1}\on{res}_{\geq i-1}(\on{cof}(\Delta_{i-1}\to \Delta_{i}))\\
& \stackrel{\eqref{eq:tilting_equivalences_res1}}{\simeq} \on{ind}^L_{\geq i-1}\on{res}_{\geq i-1}(\on{cof}(\nabla_{i}^{'}\to \Delta_{i}))\\
& \stackrel{\eqref{eq:tilting_equivalences_res2}}{\simeq} \on{ind}^L_{\geq i-1}\on{res}_{\geq i-1}(\on{cof}(\nabla_{i}^{'}\to \nabla_{i-1}^{'}))\\
& \simeq \on{ind}^L_{\geq i-1}\on{res}_{\geq i-1}(\on{fib}(\nabla_{i}^{'}\to \nabla_{i-1}^{'}))[1]\\
& \stackrel{\eqref{eq:tilting_equivalence_cosimplei}}{\simeq}  \on{ind}^L_{\geq i-1}\on{res}_{\geq i-1}(\csimp_{i}^{'})[1]\\
& \simeq \ssimp_{i}^{'}[1].
\end{align*}
This establishes the required result for $j\in\{i-1,i\}$.\\

{\bf Case $j\in(i,m]$.} We prove this case in two steps. The first step will be to establish the existence of a cofiber sequence
\begin{equation}\label{eq:cofiberseq_tiltingproof}
\framebox{$\ssimp_i^{\oplus r} \longrightarrow \ssimp_j \longrightarrow \ssimp_j^{'}$}
\end{equation}
for some non-negative integer $r\in\Z_{\geq0}$. The second step will be to show that
\begin{equation}\label{eq:cofiberseq_tiltingproof2}
\framebox{$\on{dim}_k\on{Ext}^1(\ssimp_i,\ssimp_j)=r$}
\end{equation}
In both cases, we will exhibit $r$ as the rank of a linear map.\\

{\bf Proof of \eqref{eq:cofiberseq_tiltingproof}.} Let us focus on the first step, i.e.~establishing the cofiber sequence \eqref{eq:cofiberseq_tiltingproof}. For that, let us first apply the functor $\on{Mor}(\Delta_{i-1},\mhyphen)$ to the defining fiber sequence for $T^{-,\leq 0}_{\Delta_{i}}(\zssimp_{i,j})$, and the functor $\on{Mor}(\mhyphen,\Delta_{i})^*$ to the defining cofiber sequence for $T^{+,\geq 0}_{\Delta_{i-1}}(\zssimp_{i,j})$. The two resulting cofiber sequences read:

\begin{equation}\label{eq:cofiberseq_tilting_proof1}
\on{Mor}(\Delta_{i-1},T^{-,\leq 0}_{\Delta_{i}}(\zssimp_{i,j}))\to \on{Mor}(\Delta_{i-1},\zssimp_{i,j})\to \textcolor{orange}{\on{Mor}(\Delta_{i-1},\tau_{\leq 0}\on{Mor}(\zssimp_{i,j},\Delta_{i})^*\otimes \Delta_{i})}
\end{equation}

\begin{equation}\label{eq:cofiberseq_tilting_proof2}
\textcolor{purple}{\on{Mor}(\tau_{\geq 0}\on{Mor}(\Delta_{i-1},\zssimp_{i,j})\otimes \Delta_{i-1},\Delta_{i})^*}\to \on{Mor}(\zssimp_{i,j},\Delta_{i})^*\to \on{Mor}(T^{+,\geq 0}_{\Delta_{i-1}}(\zssimp_{i,j}),\Delta_{i})^*\,.
\end{equation}

\noindent The third term in \eqref{eq:cofiberseq_tilting_proof1}, highlighted in orange, and the first term in \eqref{eq:cofiberseq_tilting_proof2}, highlighted in purple can be expressed as

\begin{equation}\label{eq:cofiberseq_tilting_proof_term1}
\textcolor{orange}{\on{Mor}(\Delta_{i-1},\tau_{\leq 0}\on{Mor}(\zssimp_{i,j},\Delta_{i})^*\otimes \Delta_{i})}\simeq \tau_{\leq 0}\on{Mor}(\zssimp_{i,j},\Delta_{i})^*\otimes (k\oplus k[-2])\,,
\end{equation}

\begin{equation}\label{eq:cofiberseq_tilting_proof_term2}
\textcolor{purple}{\on{Mor}(\tau_{\geq 0}\on{Mor}(\Delta_{i-1},\zssimp_{i,j})\otimes \Delta_{i-1},\Delta_{i})^*}\simeq  \tau_{\geq 0}\on{Mor}(\Delta_{i-1},\zssimp_{i,j})\otimes (k\oplus k[2])\,.
\end{equation}

\noindent Combining \eqref{eq:cofiberseq_tilting_proof1} and \eqref{eq:cofiberseq_tilting_proof_term1}, resp.~\eqref{eq:cofiberseq_tilting_proof2} and \eqref{eq:cofiberseq_tilting_proof_term2}, we obtain the two equivalences

\begin{equation}\label{eq:cofiberseq_tilting_proof_term3}
\tau_{\geq 1} \on{Mor}(\Delta_i,T^{-,\leq 0}_{\Delta_{i}}(\zssimp_{i,j}))\simeq \tau_{\geq 1}\on{Mor}(\Delta_{i-1},\zssimp_{i,j})\,,
\end{equation}

\begin{equation}\label{eq:cofiberseq_tilting_proof_term4}
\tau_{\leq -1}\on{Mor}(\zssimp_{i,j},\Delta_{i})^*\simeq \tau_{\leq -1}\on{Mor}(T^{+,\geq 0}_{\Delta_{i-1}}(\zssimp_{i,j}),\Delta_{i})^*\,.
\end{equation}
In particular, considering the cohomologies of \eqref{eq:cofiberseq_tilting_proof1} and \eqref{eq:cofiberseq_tilting_proof2} in degree 0, and using \eqref{eq:cofiberseq_tilting_proof_term3} and \eqref{eq:cofiberseq_tilting_proof_term4} above, we obtain the following two exact sequence of vector spaces:

\begin{equation}\label{eq:exactseq_tilting_proof1}
0\to \on{Ext}^0(\Delta_{i-1},T^{-,\leq 0}_{\Delta_{i}}(\zssimp_{i,j}))\to \on{Ext}^0(\Delta_{i-1},\zssimp_{i,j})\xrightarrow{\phi} \on{Ext}^0(\zssimp_{i,j},\Delta_{i})^*\,,
\end{equation}

\begin{equation}\label{eq:exactseq_tilting_proof2}
\on{Ext}^0(\Delta_{i-1},\zssimp_{i,j}) \xrightarrow{\phi} \on{Ext}^0(\zssimp_{i,j},\Delta_{i})^*\to \on{Ext}^0(T^{+,\geq 0}_{\Delta_{i-1}}(\zssimp_{i,j}),\Delta_{i})^*\to 0 \,.
\end{equation}
\noindent 
\noindent Note that the two morphisms labeled $\phi$ in \eqref{eq:exactseq_tilting_proof1} and \eqref{eq:exactseq_tilting_proof2} are the same morphism of $k$-vector spaces. Let $r\coloneqq\on{rk}(\phi)$ be the rank of such a morphism, i.e.~the dimension of its image in $\on{Ext}^0(\zssimp_{i,j},\Delta_{i})^*$. The next step to obtain the cofiber sequence \eqref{eq:cofiberseq_tiltingproof} is to consider the following diagram:

\begin{equation}\label{eq:tilting_proof_3times3_diagram}
\begin{tikzcd}
{\Delta_{i-1}^{\oplus r}[-1]} \arrow[d] \arrow[r,"f_1"]                                                        & 0 \arrow[r,"f_2"] \arrow[d]          & \Delta_{i}^{\oplus r} \arrow[d]                                                      \\
{\tau_{\geq 0}\on{Mor}(\Delta_{i-1},T^{-,\leq 0}_{\Delta_i}(\zssimp_{i,j}))\otimes \Delta_{i-1}} \arrow[r,"g_1"] \arrow[d, hook] & \zssimp_{i,j} \arrow[r,"g_2"] \arrow[d, "\on{id}"] & {\tau_{\leq 0}\on{Mor}(\zssimp_{i,j},\Delta_i)^*\otimes \Delta_{i}} \arrow[d, two heads]            \\
{ \tau_{\geq 0}\on{Mor}(\Delta_{i-1},\zssimp_{i,j})\otimes \Delta_{i-1}} \arrow[r,"h_1"]                                            & \zssimp_{i,j} \arrow[r,"h_2"]           & {\tau_{\leq 0}\on{Mor}(T^{+,\geq 0}_{\Delta_{i-1}}(\zssimp_{i,j}),\Delta_i)^*\otimes \Delta_{i}}
\end{tikzcd}
\end{equation}
\noindent The columns of \eqref{eq:tilting_proof_3times3_diagram} are cofiber sequences, whereas horizontally each row is a 3-term coherent chain complex in $R\Gamma(\Wgraph,\F_\w)$. Specifically, the first column of \eqref{eq:tilting_proof_3times3_diagram} is obtained from the cofiber sequence \eqref{eq:cofiberseq_tilting_proof1} by truncating with $\tau_{\geq0}$, tensoring with $\Delta_{i-1}$, and using $\phi$ as defined by \eqref{eq:exactseq_tilting_proof1}. The third column of \eqref{eq:tilting_proof_3times3_diagram} is obtained from the cofiber sequence \eqref{eq:cofiberseq_tilting_proof2} by truncating with $\tau_{\leq0}$, tensoring with $\Delta_{i}$, and using $\phi$ as in \eqref{eq:exactseq_tilting_proof2}. The second column should be self-explanatory. Note that by definition of the semi-twists
\begin{equation}\label{eq:tilting_proof_totalization_preterm1}
\fib(g_2)\simeq T^{-,\leq 0}_{\Delta_i}(\zssimp_{i,j})\,,
\end{equation}
\begin{equation}\label{eq:tilting_proof_totalization_preterm2}
\cof(h_1)\simeq T^{+,\geq 0}_{\Delta_{i-1}}(\zssimp_{i,j})\,.
\end{equation}

Consider the cofiber sequence obtained by taking the fibers of the maps from the second column of \eqref{eq:tilting_proof_3times3_diagram} to its third column, and then the cofibers of the natural maps from the first column of \eqref{eq:tilting_proof_3times3_diagram} to the resulting such fibers. Intrinsically, this is the functorial cofiber totalization of \eqref{eq:tilting_proof_3times3_diagram}, shifted by $[-1]$, understood as a fiber and cofibers sequence of 3-term coherent chain complexes.
Explicitly, each term of this new cofiber sequence, corresponding to rows one, two and three respectively, is computed as follows.
\begin{equation}\label{eq:tilting_proof_totalization_term1}
\cof({\Delta_{i-1}^{\oplus r}[-1]}\lr\fib(f_2))\simeq \cof({\Delta_{i-1}^{\oplus r}[-1]}\lr\Delta_{i}^{\oplus r}[-1])\stackrel{\eqref{eq:tilting_simple_i}}{\simeq} \textcolor{green}{\ssimp_i^{\oplus r}[-1]}\,,
\end{equation}
\
\begin{equation}\label{eq:tilting_proof_totalization_term2}
\begin{aligned}
\cof(\tau_{\geq 0}\on{Mor}(\Delta_{i-1},T^{-,\leq 0}_{\Delta_i}(\zssimp_{i,j}))\otimes \Delta_{i-1}\to\fib(g_2))&\stackrel{\eqref{eq:tilting_proof_totalization_preterm1}}{\simeq}\\
\stackrel{\eqref{eq:tilting_proof_totalization_preterm1}}{\simeq}\cof(\tau_{\geq 0}\on{Mor}(\Delta_{i-1},T^{-,\leq 0}_{\Delta_i}(\zssimp_{i,j}))\otimes \Delta_{i-1}\to T^{-,\leq 0}_{\Delta_i}(\zssimp_{i,j})) & \simeq \textcolor{blue}{T^{+,\geq 0}_{\Delta_{i-1}}T^{-,\leq 0}_{\Delta_{i}}(\zssimp_{i,j})} \,,
\end{aligned}
\end{equation}

\begin{equation}\label{eq:tilting_proof_totalization_term3}
\begin{aligned}
\cof(\tau_{\geq 0}\on{Mor}(\Delta_{i-1},\zssimp_{i,j})\otimes \Delta_{i-1}\to\fib(h_2))&\stackrel{\eqref{eq:tilting_proof_totalization_preterm2}}{\simeq}\\
\stackrel{\eqref{eq:tilting_proof_totalization_preterm2}}{\simeq}\cof(\tau_{\geq 0}\on{Mor}(\Delta_{i-1},\zssimp_{i,j})\otimes \Delta_{i-1}\to T^{+,\geq 0}_{\Delta_{i-1}}(\zssimp_{i,j})) &\simeq\textcolor{cyan}{T^{-,\leq 0}_{\Delta_{i}}T^{+,\geq 0}_{\Delta_{i-1}}(\zssimp_{i,j})}\,.
\end{aligned}
\end{equation}
Therefore, from the horizontal (shifted) totalization of \eqref{eq:tilting_proof_3times3_diagram}, and \eqref{eq:tilting_proof_totalization_term1}, \eqref{eq:tilting_proof_totalization_term2} and \eqref{eq:tilting_proof_totalization_term3}, we obtain a cofiber sequence

\begin{equation}\label{eq:tilting_proof_cofiberseq_almost}
\textcolor{green}{\ssimp_i^{\oplus r}[-1]}\lr \textcolor{blue}{T^{+,\geq 0}_{\Delta_{i-1}}T^{-,\leq 0}_{\Delta_{i}}(\zssimp_{i,j})}\lr\textcolor{cyan}{T^{-,\leq 0}_{\Delta_{i}}T^{+,\geq 0}_{\Delta_{i-1}}(\zssimp_{i,j})}\,.
\end{equation}

\noindent Now, equivalences \eqref{eq:Zssimp_equivalence} and \eqref{eq:ssimp_via_zssimp} in \cref{lem:two_sided_semi_twist_description_of_Lusztig_cycles} respectively imply the equivalences
\begin{equation}\label{eq:equivalence1_case_higherj}
\ssimp_{i-1,j}\simeq \textcolor{blue}{T^{+,\geq 0}_{\Delta_{i-1}}T^{-,\leq 0}_{\Delta_{i}}(\zssimp_{i,j})}\,,
\end{equation}
\begin{equation}\label{eq:equivalence2_case_higherj}
\ssimp_{i-1,j}^{'}\simeq \textcolor{cyan}{T^{-,\leq 0}_{\Delta_{i}}T^{+,\geq 0}_{\Delta_{i-1}}(\zssimp_{i,j})}\,.
\end{equation}
Therefore \eqref{eq:tilting_proof_cofiberseq_almost} is equivalent to

\begin{equation}\label{eq:tilting_proof_cofiberseq_almost2}
\textcolor{green}{\ssimp_i^{\oplus r}[-1]}\lr \textcolor{blue}{\ssimp_{i-1,j}} \lr \textcolor{cyan}{\ssimp_{i-1,j}^{'}}\,.
\end{equation}

\noindent Since $\on{Mor}(\Delta_l,\ssimp_i)\simeq 0$ for all $l\in[1,i-1)$, applying $T_{\Delta_1}^{+,\geq0}\cdots T_{\Delta_{i-2}}^{+,\geq0}$ to the cofiber sequence \eqref{eq:tilting_proof_cofiberseq_almost2} induces a cofiber sequence 
\[
(\ssimp_i)^{\oplus r}[-1] \stackrel{\varepsilon}{\longrightarrow} \ssimp_j \longrightarrow \ssimp_j^{'}\,.
\]
This thus proves the existence of the cofiber sequence \eqref{eq:cofiberseq_tiltingproof}, and concludes the first step of the proof in this case $j\in(i,m]$.\\

{\bf Proof of \eqref{eq:cofiberseq_tiltingproof2}.} Let us focus on the second step, i.e.~proving $\on{dim}_k\on{Ext}^1(\ssimp_i,\ssimp_j)=r$, cf.~\eqref{eq:cofiberseq_tiltingproof2}, where $r\in\Z_{\geq0}$ is as in the cofiber sequence \eqref{eq:cofiberseq_tiltingproof}. By the proof of \eqref{eq:cofiberseq_tiltingproof} above, such $r=\on{rk}(\phi)$ is the rank of the morphism
\begin{equation}\label{eq:tilting_proof_part2_morphismphi}
\on{Ext}^0(\Delta_{i-1},\zssimp_{i,j}) \xrightarrow{\phi} \on{Ext}^0(\zssimp_{i,j},\Delta_{i})^*
\end{equation}
featuring in \eqref{eq:exactseq_tilting_proof1} and \eqref{eq:exactseq_tilting_proof2}. In order to compute $\on{dim}_k\on{Ext}^1(\ssimp_i,\ssimp_j)$, we first study the object $\on{Mor}(\ssimp_i,\ssimp_j)$. By \eqref{eq:tilting_simple_i},
\begin{equation}\label{eq:tilting_proof_dimension_MorSiSj}
\on{Mor}(\ssimp_i,\ssimp_j)\simeq \on{fib}(\textcolor{orange}{\on{Mor}(\Delta_i,\ssimp_j)}\to \textcolor{purple}{\on{Mor}(\Delta_{i-1},\ssimp_j)})
\end{equation}
The first term in \eqref{eq:tilting_proof_dimension_MorSiSj}, highlighted in orange, can be described as
\begin{align*}
\textcolor{orange}{\on{Mor}(\Delta_i,\ssimp_j)}& \simeq \on{Mor}(\Delta_{i},\ssimp_{i,j}) \stackrel{\eqref{eq:Zssimp_equivalence}}{\simeq} \on{Mor}(\Delta_i,T^{-,\leq 0}_{\Delta_i}(\zssimp_{i,j}))
\\
& \simeq \on{fib}(\underbrace{\on{Mor}(\Delta_i,\zssimp_{i,j})}_{\simeq 0}\to \on{Mor}(\Delta_i,\tau_{\leq 0}\on{Mor}(\zssimp_{i,j},\Delta_i)^*\otimes \Delta_i))\\
& \simeq \textcolor{orange}{\tau_{\leq 0}\on{Mor}(\zssimp_{i,j},\Delta_i)^*[-1]}
\end{align*}
since $\on{Mor}(\Delta_i,\Delta_i)\simeq k$. The second term in \eqref{eq:tilting_proof_dimension_MorSiSj}, highlighted in purple, is equivalent to
\begin{align*}
    \textcolor{purple}{\on{Mor}(\Delta_{i-1},\ssimp_j)}&\simeq\on{Mor}(\Delta_{i-1},\ssimp_{i-1,j})\stackrel{\eqref{eq:equivalence1_case_higherj}}{\simeq} \on{Mor}(\Delta_{i-1},T_{\Delta_{i-1}}^{+,\geq 0}T^{-,\leq 0}_{\Delta_i}(\zssimp_{i,j}))\\
    & \simeq \on{cof}(\on{Mor}(\Delta_{i-1},\tau_{\geq 0}\on{Mor}(\Delta_{i-1},T^{-,\leq 0}_{\Delta_i}(\zssimp_{i,j}))\otimes \Delta_{i-1})\to \on{Mor}(\Delta_{i-1},T^{-,\leq 0}_{\Delta_i}(\zssimp_{i,j})))\\
    & \simeq \tau_{<0}\on{Mor}(\Delta_{i-1},T^{-,\leq 0}_{\Delta_i}(\zssimp_{i,j}))\\
    & \simeq \tau_{<0}\on{fib}(\on{Mor}(\Delta_{i-1},\zssimp_{i,j})\to \on{Mor}(\Delta_{i-1},\tau_{\leq 0}\on{Mor}(\zssimp_{i,j},\Delta_i)^*\otimes \Delta_i))\\
    &\simeq  \textcolor{purple}{\tau_{<0}\on{fib}(\on{Mor}(\Delta_{i-1},\zssimp_{i,j})\stackrel{(\Phi,\Theta)}{\lr} \tau_{\leq 0}\on{Mor}(\zssimp_{i,j},\Delta_i)^*\otimes (k\oplus k[-2]))}
\end{align*}
since $\on{Mor}(\Delta_{i-1},\Delta_i)\simeq k\oplus k[-2]$ in the local model of a backward weave mutation. By construction, in the last term of the equivalences above, the morphism
$$\Phi:\on{fib}(\on{Mor}(\Delta_{i-1},\zssimp_{i,j})\lr \tau_{\leq 0}\on{Mor}(\zssimp_{i,j},\Delta_i)^*,$$
cohomologically restricts in degree 0 to the morphism $\phi$ in \eqref{eq:tilting_proof_part2_morphismphi}, i.e.~$\on{Ext}^0(\Phi)=\phi$, while $\Theta$ denotes another morphism, to the component tensored with $k[-2]$.\\

Finally, by \eqref{eq:tilting_proof_dimension_MorSiSj} and these two sequences of equivalences above, the object $\on{Mor}(\ssimp_i,\ssimp_j)$ is the fiber of a morphism
\begin{equation}\label{eq:tilting_proof_part2_fibmorphism}
\textcolor{orange}{\tau_{\leq 0}\on{Mor}(\zssimp_{i,j},\Delta_i)^*[-1]}\to\textcolor{purple}{\tau_{<0}\on{fib}(\on{Mor}(\Delta_{i-1},\zssimp_{i,j})\stackrel{(\Phi,\Theta)}{\lr} \tau_{\leq 0}\on{Mor}(\zssimp_{i,j},\Delta_i)^*\otimes (k\oplus k[-2]))}
\end{equation}
\noindent By taking the homology of \eqref{eq:tilting_proof_part2_fibmorphism} in degree $-1$, i.e.~cohomologically in degree $1$, we thus conclude that
$$\on{dim}_k\on{Ext}^1(\ssimp_i,\ssimp_j)=\on{dim}_k(H_{-1}(\on{Mor}(\ssimp_i,\ssimp_j)))=\on{rk}(\on{Ext}^0(\Phi))=\on{rk}(\phi)=r,$$
which indeed establishes \eqref{eq:cofiberseq_tiltingproof2}.\\

{\bf Conclusion of proof for case $j\in(i,m]$.} Since $\ssimp_{j}^{'}$ is indecomposable, and \eqref{eq:cofiberseq_tiltingproof2} established that $\on{dim}_k\on{Ext}^1(\ssimp_i,\ssimp_j)=r$, the morphism $\varepsilon$ in \eqref{eq:tilting_proof_cofiberseq_almost2}, i.e.~the left morphism of \eqref{eq:cofiberseq_tiltingproof}, must agree with the evaluation morphism
$$\ssimp_i\otimes \on{Ext}^{1}(\ssimp_i,\ssimp_j)[-1]\stackrel{\on{ev}}{\longrightarrow} \ssimp_j.$$
Therefore, $\ssimp_j^{'}\simeq\on{cof}(\on{Ext}^1(\ssimp_i,\ssimp_j)\otimes \ssimp_i[-1] \stackrel{\on{ev}}{\to} \ssimp_j)$ for all $j\in(i,m]$. This concludes the case $j\in(i,m]$ and thus the proof of Theorem \ref{thm:main2}.$(4)$.
\end{proof}

\section{Silting objects and their weave cycles}\label{sec:silting_dualLusztig}

In this section we conclude our results on categorical weave calculus by studying two silting collections $\{\lsilt_i^\w\}$ and $\{\rsilt_i^\w\}$, respectively left and right dual to the simple-minded categorical Lusztig cycles $\{\LC_i\}$. The goal is to describe the objects of these silting collections in terms of both the weave $\w$ and its perverse schober $\F_\w$. In order to achieve this, we first prove a type of Ringel duality for weave schobers in \cref{ssec:two_sided_schobers}, cf.~\cref{prop:upwards_downwards_thimbles_duality}. The silting collections are then introduced in \cref{ssec:silting_collections}, and we then use such duality to establish the desired characterizations in \cref{ssec:weaverealization_silting}, cf.~\cref{thm:description_of_silting_dual_Lusztig_cycles} and \cref{cor:description_of_silting_dual_Lusztig_cycles}.


\subsection{Two-sided weave schobers}\label{ssec:two_sided_schobers}

Let $\w$ be a Demazure weave with $m$ trivalent vertices. In this subsection we first introduce three ribbon graphs $\rgraph_\w^\downarrow$, $\rgraph_\w^\updownarrow$ and $\rgraph_\w^\uparrow$, and three corresponding perverse schobers $\mathcal{F}_\w^{\downarrow}$, $\mathcal{F}_\w^{\updownarrow}$ and $\mathcal{F}_\w^{\uparrow}$.\\

First, $\rgraph_\w^{\downarrow}\coloneqq \rgraph_\w$ is defined to be the ribbon graph $\rgraph_\w$ of $\w$ introduced in \cref{sec:weave_schobers}. By definition, the ribbon graph $\rgraph_\w^\updownarrow$ is obtained from $\rgraph_\w^\downarrow$ by adding an external edge $\einftop$ incident to the top vertex $v_m$. By definition, the ribbon graph $\rgraph_\w^\uparrow$ is obtained from $\rgraph_\w^\updownarrow$ by removing the external edge $e_1$ incident to $v_1$. Second, their associated perverse schobers are defined using \Cref{not:schobers} as follows:

\begin{definition}
Let $\w$ be a Demazure weave with trivalent vertices $p_1,\dots,p_m$. Then $\mathcal{F}_\w^{\updownarrow}$, resp.~$\mathcal{F}_\w^{\uparrow}$, is the $\rgraph_\w^{\updownarrow}$-parametrized, $\rgraph_\w^{\uparrow}$-parametrized, perverse schober defined by:
$$
\mathcal{F}_\w^{\updownarrow}\coloneqq 
\begin{tikzcd}
{} \arrow[d, "\varrho_1", no head]                                        \\
\mhyphen\otimes S_{p_m} \arrow[d, "{(\varrho_2,\varrho_1)}", no head]     \\
\mhyphen\otimes S_{p_{m-1}} \arrow[d, "{(\varrho_2,\varrho_1)}", no head] \\
\vdots \arrow[d, "{(\varrho_2,\varrho_1)}", no head]                      \\
\mhyphen\otimes S_{p_1} \arrow[d, "\varrho_2", no head]                   \\
{}                                                                       
\end{tikzcd}
\qquad\qquad\qquad\qquad
\mathcal{F}_\w^{\uparrow}\coloneqq \begin{tikzcd}
{} \arrow[d, "\varrho_1", no head]                                        \\
\mhyphen\otimes S_{p_m} \arrow[d, "{(\varrho_2,\varrho_1)}", no head]     \\
\mhyphen\otimes S_{p_{m-1}} \arrow[d, "{(\varrho_2,\varrho_1)}", no head] \\
\vdots \arrow[d, "{(\varrho_2,\varrho_1)}", no head]                      \\
\mhyphen\otimes S_{p_1}                                                  
\end{tikzcd}
$$
By definition, we also set $\mathcal{F}_\w^{\downarrow}\coloneqq \mathcal{F}_\w$.\qed
\end{definition}

Let us now show that the $\infty$-category of downward global sections $\Glsec{\rgraph_\w^{\downarrow}}{\mathcal{F}_\w^{\downarrow}}$ and the $\infty$-category of upward global sections $\Glsec{\rgraph_\w^{\uparrow}}{\mathcal{F}_\w^{\uparrow}}$ both embed into the two-sided global sections $\Glsec{\rgraph_\w^{\updownarrow}}{\mathcal{F}_\w^{\updownarrow}}$. In fact, these two full subcategories of $\Glsec{\rgraph_\w^{\updownarrow}}{\mathcal{F}_\w^{\updownarrow}}$ differ by a mutation equivalence of semiorthogonal decompositions, in the sense of \cite[Def.~2.4.1]{DKSS21} and are thus equivalent as $\infty$-categories. The precise statement we prove reads as follows: 

\begin{lemma}\label{lem:upwards_downwards_duality_schobers}
\begin{enumerate}[$(1)$]
    \item The $\infty$-category of global sections $\Glsec{\rgraph_\w^{\downarrow}}{\mathcal{F}_\w^{\downarrow}}$ is equivalent to the full subcategory of $\Glsec{\rgraph_\w^{\updownarrow}}{\mathcal{F}_\w^{\updownarrow}}$ consisting of global sections whose image under $\on{ev}_\infty$ vanishes. The arising inclusion
    $$\iota^{\downarrow}\colon \Glsec{\rgraph_\w^{\downarrow}}{\mathcal{F}_\w^{\downarrow}}\lr \Glsec{\rgraph_\w^{\updownarrow}}{\mathcal{F}_\w^{\updownarrow}}$$
    admits both left and right adjoints $(\iota^{\downarrow})^L,(\iota^{\downarrow})^R$.
    \item The $\infty$-category of global sections $\Glsec{\rgraph_\w^{\uparrow}}{\mathcal{F}_\w^{\uparrow}}$ is equivalent to the full subcategory of $\Glsec{\rgraph_\w^{\updownarrow}}{\mathcal{F}_\w^{\updownarrow}}$ consisting of global sections whose image under $\on{ev}_1$ vanishes. The arising inclusion
    $$\iota^{\uparrow}\colon \Glsec{\rgraph_\w^{\uparrow}}{\mathcal{F}_\w^{\uparrow}}\lr \Glsec{\rgraph_\w^{\updownarrow}}{\mathcal{F}_\w^{\updownarrow}}$$
    admits both left and right adjoints $(\iota^{\uparrow})^L,(\iota^{\uparrow})^R$.
    \item The four composite functors 
    \[(\iota^{\downarrow})^\bullet\circ \iota^{\uparrow}\colon \Glsec{\rgraph_\w^{\uparrow}}{\mathcal{F}_\w^{\uparrow}} \longrightarrow \Glsec{\rgraph_\w^{\downarrow}}{\mathcal{F}_\w^{\downarrow}}\qquad\mbox{ and }\qquad(\iota^{\uparrow})^\bullet\circ \iota^{\downarrow}\colon \Glsec{\rgraph_\w^{\downarrow}}{\mathcal{F}_\w^{\downarrow}} \longrightarrow \Glsec{\rgraph_\w^{\uparrow}}{\mathcal{F}_\w^{\uparrow}}\]
    are equivalences of $\infty$-categories, where $\bullet\in\{L,R\}$.
\end{enumerate}
\end{lemma}

\begin{proof}
For Parts (1) \& (2), \cite[Theorem 4.6]{Chr25b} implies that the two evaluation functors
$$\on{ev}_{e_1}\colon \Glsec{\rgraph_\w^{\updownarrow}}{\mathcal{F}_\w^{\updownarrow}}\to \mathcal{F}_\w^\updownarrow(e_1)\simeq\D\,\quad\mbox{and}\quad\on{ev}_{e_\infty}\colon \Glsec{\rgraph_\w^{\updownarrow}}{\mathcal{F}_\w^{\updownarrow}}\to \mathcal{F}_\w^\updownarrow(e_\infty)\simeq\D$$
define a perverse schober parametrized by the $2$-spider. Therefore, they each admit left and right adjoints which are fully faithful. In order to show that $\Glsec{\rgraph_\w^{\downarrow}}{\mathcal{F}_\w^{\downarrow}}\simeq \on{fib}(\on{ev}_{\einftop})$ is equivalent to the full subcategory of global sections whose image under $\on{ev}_{\einftop}$ vanishes, it suffices to note that the passage to the fiber commutes with the passage to global sections, as limits commute with limits, and that the fiber of $\mathcal{F}_\w^\updownarrow(v_m)\to \mathcal{F}_\w^\updownarrow(e_\infty)$ is equivalent to $\mathcal{F}_\w^\downarrow(v_m)$. The latter follows from \cite[Lemma 4.26]{Chr22}. This concludes Part (1). A similar argument shows that $\Glsec{\rgraph_\w^{\uparrow}}{\mathcal{F}_\w^{\uparrow}}\simeq \on{fib}(\on{ev}_{e_1})$, thus concluding the proof of Part (2).\\
    
Part (3) can be proven as follows. By using the existence of the right adjoints of $\on{ev}_{e_1},\on{ev}_{\einftop}$, from Parts (1) and (2),  \cite[Prop.~2.3.2]{DKSS21} implies that $\Glsec{\rgraph_\w^{\updownarrow}}{\mathcal{F}_\w^{\updownarrow}}$ admits two semiorthogonal decompositions 
    \[ (\on{ev}_{e_1}^R(\mathcal{F}_\w^\updownarrow(e_1)),\iota^\uparrow(\Glsec{\rgraph_\w^{\uparrow}}{\mathcal{F}_\w^{\uparrow}}))\qquad\mbox{ and }\qquad(\iota^\downarrow(\Glsec{\rgraph_\w^{\downarrow}}{\mathcal{F}_\w^{\downarrow}}),\on{ev}_{e_\infty}^R(\mathcal{F}_\w^\updownarrow(e_1))\,.\] 
   The two essential images $\on{ev}_{e_1}^R(\mathcal{F}_\w^\updownarrow(e_1),\on{ev}_{e_\infty}^R(\mathcal{F}_\w^\updownarrow(e_1)$ agree by the definition of perverse schober on the $2$-spider. Since the functors $(\iota^{\uparrow})^R\circ \iota^{\downarrow}$ and $(\iota^{\downarrow})^L\circ \iota^{\uparrow}$ describe the corresponding mutation equivalences, cf.~\cite[Def.~2.4.1 \& Prop.~2.4.2]{DKSS21}, they are therefore equivalences. Similarly, by using left adjoints (instead of right adjoints as above), the corresponding argument shows that the functors $(\iota^{\downarrow})^R\circ \iota^{\uparrow},(\iota^{\uparrow})^L\circ \iota^{\downarrow}$ are also equivalences.
\end{proof}

\noindent The $\infty$-categories and functors discussed in \cref{lem:upwards_downwards_duality_schobers} are depicted in Diagram \eqref{eq:diagram_twosided}.

\begin{equation}\label{eq:diagram_twosided}
\begin{tikzcd}
    \D\simeq\mathcal{F}_\w^\updownarrow(e_1) & & \mathcal{F}_\w^\updownarrow(e_\infty)\simeq\D\\
    & \Glsec{\rgraph_\w^{\updownarrow}}{\mathcal{F}_\w^{\updownarrow}} \arrow[dr, "{(\iota^{\uparrow})^L}"', bend right] 
        \arrow[dr, "{(\iota^{\uparrow})^R}", bend left] \arrow[dl, "{(\iota^{\downarrow})^L}"', bend right] 
        \arrow[dr, "{(\iota^{\uparrow})^R}", bend left] \arrow[dl, "{(\iota^{\downarrow})^R}", bend left] \arrow[ur, "\on{ev}_{e_\infty}"] \arrow[ul, "\on{ev}_{e_1}"'] & \\
    \fib(\on{ev}_{e_\infty})\simeq\Glsec{\rgraph_\w^\downarrow}{\mathcal{F}_\w^\downarrow} \arrow[rr,"(\iota^{\uparrow})^L\circ \iota^{\downarrow}","(\iota^{\uparrow})^R\circ \iota^{\downarrow}"',bend right] \arrow[ur,hook, "\iota^{\downarrow}"] & & \Glsec{\rgraph_\w^{\uparrow}}{\mathcal{F}_\w^{\uparrow}}\simeq\fib(\on{ev}_{e_1}) \arrow[ul,hook', "\iota^{\uparrow}"']
\end{tikzcd}
\end{equation}

\begin{lemma}\label{lem:two_sided_3CY}
    The functor $(\on{ev}_{e_1},\on{ev}_{\einftop})\colon \Glsec{\rgraph_\w^{\updownarrow}}{\mathcal{F}_\w^{\updownarrow}}\to \D\times\D$ admits a relative weak right $3$-Calabi--Yau structure.
\end{lemma}

\begin{proof}
The proof is similar to the proof of \Cref{prop:rel_3_CY_str}.
\end{proof}

Let us denote $\Delta_i^{\w,\downarrow}\coloneqq \Delta_i^\w$ and $\nabla_i^{\w,\downarrow}\coloneqq \nabla_i^\w$ for the downward facing weaves thimbles in $\Glsec{\rgraph_\w^{\downarrow}}{\mathcal{F}_\w^{\downarrow}}$, as introduced in \Cref{def:thimbles}. See also Equations \eqref{eq:thimble_as_induction} and \eqref{eq:cothimble_as_induction} for their descriptions in terms of inductions from $v_i$. The corresponding upwards weave thimbles can be defined directly via induction from $v_i$, as follows:

\begin{definition}
    \begin{enumerate}[(1)] 
        \item The upwards standard thimble $\Delta_i^{\w,\uparrow}\in \Glsec{\rgraph_\w^\uparrow}{\mathcal{F}_\w^\uparrow}$ is defined as 
        \begin{equation}\label{eq:upward_cothimble}
        \Delta_i^{\w,\uparrow}\coloneqq \on{ind}^L_{v_i}((k,S_{p_i},\on{id}_{S_{p_i}}))\,.
        \end{equation}
        \item The upwards costandard thimble $\nabla_i^{\w,\uparrow}\in \Glsec{\rgraph_w^\uparrow}{\mathcal{F}_\w^\uparrow}$ is defined as 
        \[
        \nabla_i^{\w,\uparrow}\coloneqq \on{ind}^R_{v_i}((k,S_{p_i},\on{id}_{S_{p_i}}))\,.
        \]
        \qed
    \end{enumerate}
\end{definition}

The following result shows that the equivalences $(\iota^{\uparrow})^R\circ \iota^{\downarrow}$ and $(\iota^{\uparrow})^L\circ \iota^{\downarrow}$ in \cref{lem:upwards_downwards_duality_schobers} can be understood as a type of Ringel duality in the context of perverse schobers, exchanging downward standard and upward costandard thimbles, up to a shift.

\begin{proposition}[Ringel-type duality]\label{prop:upwards_downwards_thimbles_duality}
There exist equivalences in $\Glsec{\rgraph^{\uparrow}_\w}{\mathcal{F}_{\w}^{\uparrow}}$
\begin{equation}\label{eq:Ringel_duality}
((\iota^\uparrow)^{R}\circ \iota^\downarrow)(\Delta_i^{\w,\downarrow})\simeq \nabla_i^{\w,\uparrow}[-1]
\end{equation}
and
\begin{equation}\label{eq:Ringel_duality_2} ((\iota^{\uparrow})^L\circ \iota^{\downarrow})(\nabla_i^{\w,\downarrow})\simeq \Delta_i^{\w,\uparrow}[2]
\,.\end{equation}
\end{proposition}

\begin{proof}
Let us prove equivalence \eqref{eq:Ringel_duality}. Denote by 
\[ T^{\on{gl}}\coloneqq \on{fib}(\on{id}_{\Glsec{\rgraph_\w^\updownarrow}{\mathcal{F}^\updownarrow_\w}}\xrightarrow{\on{unit}}(\on{ind}^R_{\einftop},\on{ind}^R_{e_1})\circ (\on{ev}_{\einftop},\on{ev}_{e_1}))\] 
the cotwist functor $T^{\on{gl}}:\Glsec{\rgraph_\w^\updownarrow}{\mathcal{F}^\updownarrow_\w}\lr\Glsec{\rgraph_\w^\updownarrow}{\mathcal{F}^\updownarrow_\w}$ of the spherical adjunction
\[
(\on{ev}_{\einftop},\on{ev}_{e_1})\colon \Glsec{\rgraph_\w^\updownarrow}{\mathcal{F}^\updownarrow_\w}\longleftrightarrow \mathcal{F}^\updownarrow_\w(\einftop)\times \mathcal{F}^\updownarrow_\w(e_1)\noloc (\on{ind}^R_{\einftop},\on{ind}^R_{e_1}).
\]
Similarly, for each $i\in[1,m]$, denote by $T_{v_i}$ the cotwist functor of the spherical adjunction 
\[ (\varrho_2,\varrho_1)=(\mathcal{F}(v_i\to e_i),\mathcal{F}(v_i\to e_{i+1}))\colon \mathcal{F}_\w^{\updownarrow}(v_i) \longleftrightarrow \mathcal{F}_\w^{\updownarrow}(e_i) \times \mathcal{F}_\w^{\updownarrow}(e_{i+1})\noloc (\mathcal{F}(v_i\to e_i)^R,\mathcal{F}(v_i\to e_{i+1})^R)\,.\]

\noindent For these functors, \cite[Prop.~4.23.(2)]{Chr25b} implies that
\begin{equation}\label{eq:TindL_is_indTv}
T^{\on{gl}}\circ \on{ind}^L_{v_i}\\
\simeq \on{ind}^R_{v_i}\circ T_{v_i}.
\end{equation}
Now, a direct computation shows that 
\begin{equation}\label{eq:ringel_proof_Tvi_vertex}
T_{v_i}((k[-1],0,0))\simeq (k[-1],S_i[-1],\on{id}_{S_i[-1]})\,.
\end{equation}
Therefore, we conclude the equivalences:
\begin{equation}\label{eq:Tgldown_computation1}
\begin{aligned}
(T^{\on{gl}}\circ \iota^{\downarrow})(\Delta_i^{\downarrow})&\stackrel{\eqref{eq:thimble_as_induction}}{\simeq} (T^{\on{gl}}\circ \on{ind}^L_{v_i})((k[-1],0,0))\\
& \stackrel{\eqref{eq:TindL_is_indTv}}{\simeq} \on{ind}^R_{v_i}\circ T_{v_i}((k[-1],0,0))\\
& \stackrel{\eqref{eq:ringel_proof_Tvi_vertex}}{\simeq} \on{ind}^R_{v_i}((k[-1],S_i[-1],\on{id}_{S_i[-1]}))\\
& \stackrel{\eqref{eq:upward_cothimble}}{\simeq} \iota^{\uparrow}(\nabla_i^{\uparrow})[-1]\,.
\end{aligned}
\end{equation}

\noindent Let us compute $(T^{\on{gl}}\circ \iota^{\downarrow})(\Delta_i^{\downarrow})$ differently, as follows. Since we have $(\on{ev}_{\einftop},\on{ev}_{e_1})(\Delta_i^\downarrow)\simeq (0,S_i)$, we deduce the equivalence
\begin{equation}\label{eq:Tgldown_computation2}
(T^{\on{gl}}\circ \iota^{\downarrow})(\Delta_i^{\downarrow})\simeq \on{fib}(\iota^{\downarrow}(\Delta_i^{\downarrow})\to \on{ind}^R_{e_1}(S_i))\,.
\end{equation}
Since $(T^{\on{gl}}\circ \iota^{\downarrow})(\Delta_i^{\downarrow})\in \on{Im}(\iota^{\uparrow})$, by the fully faithfulness of $\iota^\uparrow$, we obtain
\begin{equation}\label{eq:ringel_proof_eq1}
(T^{\on{gl}}\circ \iota^{\downarrow})(\Delta_i^{\downarrow})\simeq (\iota^\uparrow\circ (\iota^{\uparrow})^R)(T^{\on{gl}}\circ \iota^{\downarrow})(\Delta_i^{\downarrow})\,.
\end{equation}

\noindent Now, for any upward global section $X\in \Glsec{\rgraph^{\uparrow}_\w}{\mathcal{F}^{\uparrow}_{\w}}$, we have 
\begin{equation}\label{eq:vanishing_MorX_to_upindR}
\begin{aligned} 
\on{Mor}_{\Glsec{\rgraph^{\uparrow}_\w}{\mathcal{F}^{\uparrow}_{\w}}}(X,(\iota^{\uparrow})^R(\on{ind}^R_{e_1}(S_i)))& \simeq \on{Mor}_{\Glsec{\rgraph^{\updownarrow}_\w}{\mathcal{F}^{\updownarrow}_{\w}}}(\iota^{\uparrow}(X),\on{ind}^R_{e_1}(S_i))\\
& \simeq \on{Mor}_{\mathcal{F}(e_1)}(\underbrace{\on{ev}_{e_1}(\iota^{\uparrow}(X))}_{\simeq 0},S_i)\simeq 0,
\end{aligned}
\end{equation}
where the first equivalence is the adjunction $\iota^R\dashv (\iota^{\uparrow})^R$, and $\on{ev}_{e_1}(\iota^{\uparrow}(X))\simeq0$ follows from \cref{lem:upwards_downwards_duality_schobers}.(2).
The vanishing \eqref{eq:vanishing_MorX_to_upindR} for all $X\in \Glsec{\rgraph^{\uparrow}_\w}{\mathcal{F}^{\uparrow}_{\w}}$ implies that $(\iota^{\uparrow})^R(\on{ind}^R_{e_1}(S_i))\simeq 0$. Therefore \eqref{eq:Tgldown_computation2} yields
\begin{equation}\label{eq:ringel_proof_eq2}
(\iota^\uparrow\circ (\iota^{\uparrow})^R)(T^{\on{gl}}\circ \iota^{\downarrow})(\Delta_i^{\downarrow})\simeq (\iota^\uparrow\circ (\iota^{\uparrow})^R\circ \iota^{\downarrow})(\Delta^\downarrow_i)\,.
\end{equation}
Now, combining the above equivalences yields
\begin{equation}\label{eq:ringel_proof_eq3}
(\iota^\uparrow\circ (\iota^{\uparrow})^R\circ \iota^{\downarrow})(\Delta_i)\stackrel{\eqref{eq:ringel_proof_eq2}}{\simeq}(\iota^\uparrow\circ (\iota^{\uparrow})^R)(T^{\on{gl}}\circ \iota^{\downarrow})(\Delta_i^{\downarrow})\stackrel{\eqref{eq:ringel_proof_eq1}}{\simeq}(T^{\on{gl}}\circ \iota^{\downarrow})(\Delta_i^{\downarrow})\stackrel{\eqref{eq:Tgldown_computation1}}{\simeq}
\iota^{\uparrow}(\nabla_i^{\uparrow})[-1].
\end{equation}
Since the functor $\iota^\uparrow$ is fully faithful, \eqref{eq:ringel_proof_eq3} implies the equivalence \eqref{eq:Ringel_duality}, as required. The proof of equivalence \eqref{eq:Ringel_duality_2} is analogous, using 
$T_{v_i}((k[2],S_i[2],\on{id}_{S_i[2]}))\simeq (k[-1],0,0)$ instead of \eqref{eq:ringel_proof_Tvi_vertex}.
\end{proof}

\noindent Due to \cref{lem:upwards_downwards_duality_schobers}.(3) and Proposition \ref{prop:upwards_downwards_thimbles_duality}, we often refer to either of the four functors $(\iota^{\downarrow})^\bullet\circ \iota^{\uparrow}$ and $(\iota^{\uparrow})^\bullet\circ \iota^{\downarrow}$ as dualities between the $\infty$-categories $\Glsec{\rgraph_\w^{\downarrow}}{\mathcal{F}_\w^{\downarrow}}$ and $\Glsec{\rgraph_\w^{\uparrow}}{\mathcal{F}_\w^{\uparrow}}$, $\bullet\in\{L,R\}$.

\begin{remark}
The autoequivalence
$$(\iota^\downarrow)^L\circ \iota^\uparrow \circ (\iota^{\uparrow})^L\circ \iota^{\downarrow}\colon \Glsec{\rgraph_\w^{\downarrow}}{\mathcal{F}_\w^{\downarrow}}\to \Glsec{\rgraph_\w^{\downarrow}}{\mathcal{F}_\w^{\downarrow}}$$
describes the inverse Serre functor of the relative weakly left $3$-Calabi--Yau category $\Glsec{\rgraph_\w^{\downarrow}}{\mathcal{F}_\w^{\downarrow}}$, cf.~\cref{prop:rel_3_CY_str}. This aligns with the fact that the difference of the shifts $[-1]$ and $[2]$ appearing in \Cref{prop:upwards_downwards_thimbles_duality} is indeed $-3$.\qed
\end{remark}


\subsection{Silting collections in \texorpdfstring{$\Glsec{\rgraph_\w^{\downarrow}}{\mathcal{F}_\w^{\downarrow}}$ and $\Glsec{\rgraph_\w^{\uparrow}}{\mathcal{F}_\w^{\uparrow}}$}{downwards and upwards global sections}}\label{ssec:silting_collections}

Let us now introduce the silting collections associated to the thimbles. There are four of them, as we have standard and costandard, and upwards and downwards, thimbles.

\begin{definition}[Silting collections from thimbles]\label{def:upwards_silting_objects}
Let $\w$ be a Demazure weave with $m$ trivalent vertices. For all $i\in[1,m]$, we define
\[ \sdsilt_i\coloneqq T^{-,>0}_{\Delta_{m}^{\downarrow,\w}}(\dots (T^{-,>0}_{\Delta_{i+1}^{\w,\downarrow}}(\Delta_i^{\w,\downarrow})))\in \Glsec{\rgraph^{\downarrow}_\w}{\mathcal{F}_\w^{\downarrow}}\]
\[ \cdsilt_i\coloneqq T^{+,<0}_{\nabla_{m}^{\downarrow,\w}}(\dots (T^{+,<0}_{\nabla_{i+1}^{\w,\downarrow}}(\nabla_i^{\w,\downarrow})))\in \Glsec{\rgraph^{\downarrow}_\w}{\mathcal{F}_\w^{\downarrow}} \]
and
\[ \cusilt_i\coloneqq T^{-,>0}_{\nabla_{m}^{\uparrow,\w}}(\dots (T^{-,>0}_{\nabla_{i+1}^{\w,\uparrow}}(\nabla_i^{\w,\uparrow})))\in \Glsec{\rgraph^{\uparrow}_\w}{\mathcal{F}_\w^{\uparrow}}\]
\[\susilt_i\coloneqq T^{+,<0}_{\Delta_{m}^{\uparrow,\w}}(\dots (T^{+,<0}_{\Delta_{i+1}^{\w,\uparrow}}(\Delta_i^{\w,\uparrow})))\in \Glsec{\rgraph^{\uparrow}_\w}{\mathcal{F}_\w^{\uparrow}} \]
\qed
\end{definition}

\noindent \Cref{prop:silting_from_standards} implies that $\{\sdsilt_i\}$ and $\{\cdsilt_i\}$ in \cref{def:upwards_silting_objects} form silting collections in $\Glsec{\rgraph^{\downarrow}_\w}{\mathcal{F}_\w^{\downarrow}}$. In terms of the comparison to highest weight theory, cf.~\cref{subsec:relation_w_highest_weight_theory}
, $\{\sdsilt_i\}$ and $\{\cdsilt_i\}$ are a stable $\infty$-categorical analogue of the projective and injective objects, respectively.

\begin{remark} (1) We refer to the objects $\sdsilt_i,\cdsilt_k\in\Glsec{\rgraph_{\w}}{\mathcal{F}_\w}$ in \cref{def:upwards_silting_objects} as the dual categorical Lusztig cycles of the weave $\w$, for any $i,k\in[1,m]$. Specifically $\sdsilt_i$ is a left dual categorical Lusztig cycles, and $\cdsilt_i$ a right dual categorical Lusztig cycles. Indeed, they are, respectively, left and right dual to the simple-minded collection given by the categorical Lusztig cycles $\{\LC_i\}$.\\

\noindent (2) The analogs of \cref{thm:Lusztigcycles_weaveequivalence} and \cref{thm:weavemutation_SMCmutation1} hold for these silting collections $\{\sdsilt_i\}$ and $\{\cdsilt_i\}$ in $\Glsec{\rgraph^{\downarrow}_\w}{\mathcal{F}_\w^{\downarrow}}$. Namely, they are invariant under weave equivalences, and they undergo a mutation of silting collections under a weave mutation. This is a direct consequence of Koszul duality.\qed
\end{remark}

These collections of upwards and downwards silting objects in \Cref{def:upwards_silting_objects} are related via the the Ringel-type duality from \Cref{lem:upwards_downwards_duality_schobers} as follows:

\begin{lemma}\label{lem:Ringel_duality_and_silting_objects}
Let $\w$ be a Demazure weave with $m$ trivalent vertices. Then, for all $i\in[1,m]$, there exist equivalences in $\Glsec{\rgraph_\w^{\downarrow}}{\mathcal{F}_\w^{\downarrow}}$
\begin{equation}\label{eq:Ringeldual_silting1}
((\iota^{\downarrow})^L\circ \iota^{\uparrow})(\cusilt_i)[-1]\simeq \sdsilt_i
\end{equation}
and 
\begin{equation}\label{eq:Ringeldual_silting2}
((\iota^{\downarrow})^R\circ \iota^{\uparrow})(\susilt_i)[2]\simeq \cdsilt_i\,.
\end{equation}
\end{lemma}

\begin{proof}
These equivalences are implied by \Cref{prop:upwards_downwards_thimbles_duality}. 
\end{proof}

\noindent Combining \Cref{lem:Ringel_duality_and_silting_objects} and \Cref{prop:silting_from_standards}, we find that $\{\susilt_i\}$ and $\{\cusilt_i\}$ define silting collections in $\Glsec{\rgraph^{\uparrow}_\w}{\mathcal{F}_\w^{\uparrow}}$.\\

\noindent In a highest weight abelian category, the tilting objects each admit, by definition, both a standard and a costandard flag. The following result gives an analogue of such a property for the stable $\infty$-category $\Glsec{\rgraph_\w^{\uparrow}}{\mathcal{F}_\w^{\uparrow}}$, in which the upwards thimbles form a full exceptional collection. We leave open the interesting question whether this property characterizes the upwards silting collection.

\begin{lemma}\label{lem:standard_vs_costandard_upwards_silt}
There exists an equivalence of global sections 
\[ \cusilt_i\simeq \susilt_i\,.\] 
\end{lemma}

\begin{proof}
First, since $\{\cdsilt_j\}$ are right duals to $\{\simp_i^\w\}$, we obtain
\begin{equation}\label{eq:proof_lemma_doublefiltration1}
\on{Mor}(\simp_i^\w,\susilt_j[2])\stackrel{\eqref{eq:Ringeldual_silting2}}{\simeq} \on{Mor}(\simp_i^\w,((\iota^\uparrow)^R\circ \iota^{\downarrow})(\cdsilt_j)))\simeq \on{Mor}(\simp_i^\w,\cdsilt_j)\simeq \delta_{ij}\cdot k\,.
\end{equation}
\noindent Second, by using the relative right $3$-Calabi--Yau structure from \Cref{lem:two_sided_3CY} and \cite[Prop.~5.19]{Chr22b}, we also have
\begin{equation}\label{eq:proof_lemma_doublefiltration2}
\on{Mor}(\simp_i^\w,\susilt_j[2])\simeq \on{Mor}(\susilt_j[2],\simp_i^\w)^*[-3].
\end{equation}
\noindent Equivalences \eqref{eq:proof_lemma_doublefiltration1} and \eqref{eq:proof_lemma_doublefiltration2} thus imply
\begin{equation}\label{eq:proof_lemma_doublefiltration3}
\on{Mor}(\susilt_j[-1],\simp_i^\w)\simeq \delta_{i,j}k\,.
\end{equation}
\noindent Therefore we have the equivalences 
\[ \on{Mor}(((\iota^\downarrow)^L\circ \iota^\uparrow)(\susilt_j[-1]),\simp_i^\w) \simeq  \on{Mor}(\susilt_j[-1],\simp_i^\w)\stackrel{\eqref{eq:proof_lemma_doublefiltration3}}{\simeq} \delta_{ij}\cdot k\,,\]
which show that the collection of objects $\{((\iota^\downarrow)^L\circ \iota^\uparrow)(\susilt_j[-1])\}$ is Koszul dual to the categorical Lusztig cycles $\{\simp_i^\w\}$. In consequence,
\[
 ((\iota^\downarrow)^L\circ \iota^\uparrow)(\susilt_i[-1])\simeq \sdsilt_i\stackrel{\eqref{eq:Ringeldual_silting1}}{\simeq} ((\iota^\downarrow)^L\circ \iota^\uparrow)(\cusilt_i[-1])\,,
\]
which implies that $\susilt_i\simeq \cusilt_i$, as required.
\end{proof}

\begin{remark}
The upwards silting collection $\{\cusilt_i\}$ is dual to the categorical Lusztig cycles as follows:
 \[\on{Mor}(\iota^\uparrow(\cusilt_i),\iota^\downarrow(\simp_j^\w))\simeq \delta_{ij}\cdot k[-1]\,\]
 This is implied by \Cref{lem:Ringel_duality_and_silting_objects} and \Cref{thm:Koszul_duality_for_silting_and_SMC} by using the adjunctions.\qed
\end{remark}

\subsection{Weave realization of dual categorical Lusztig cycles}\label{ssec:weaverealization_silting}

The goal of this subsection is to describe the objects in the silting collections $\{\sdsilt_i\}$ and $\{\cdsilt_i\}$ associated to $\{\Delta_i^\w\}$ in terms of the Demazure weave $\w$ itself. By \cref{lem:Ringel_duality_and_silting_objects} and \cref{lem:standard_vs_costandard_upwards_silt}, we can equivalently aim at describing the upward objects $\{\susilt_i\}$ in terms of $\w$. This is achieved in \cref{thm:description_of_silting_dual_Lusztig_cycles} and \cref{cor:description_of_silting_dual_Lusztig_cycles}, cf.~also \cref{rem:description_weave_cycle_of_silt}.

\begin{definition}[Weave cycles realization]\label{def:weave_cycle_realization}
Let $\w$ be a Demazure weave and $X\in  \Glsec{\rgraph_\w^{\updownarrow}}{\mathcal{F}_\w^{\updownarrow}}$ a global section. By definition, $X$ is said to be realized by a weave cycle $\vartheta:E(\w)\lr\Z_{\geq0}$ if:
\begin{enumerate}
    \item $\vartheta$ satisfies the tropical Lusztig propagation rules at the $4$ and $6$-valent vertices of $\w$.

    \item For any point $q\in e_i$ in the edge $e_i\in\rgraph_\w^{\updownarrow}$, the object $X(e_i)$ admits a $(\beta_q,{\bf a}_q)$-filtration, where $(\beta_q,{\bf a}_q)$ is the weighted braid word associated to $\vartheta$ at the horizontal slice of $\w$ at the height of $q$.
\end{enumerate}

\noindent By definition, if $X$ lies in $\Glsec{\rgraph_\w^{\uparrow}}{\mathcal{F}_\w^{\uparrow}}$ or $ \Glsec{\rgraph_\w^{\downarrow}}{\mathcal{F}_\w^{\downarrow}}$, a weave cycle realization of $X$ consists of a weave cycle realization of its image in $\Glsec{\rgraph_\w^{\updownarrow}}{\mathcal{F}_\w^{\updownarrow}}$.\qed
\end{definition}

\noindent A motivating instance for \cref{def:weave_cycle_realization} is that the categorical Lusztig cycle $\simp_i^\w\in \Glsec{\rgraph_\w^\downarrow}{\mathcal{F}_\w^\downarrow}$ is realized by the Lusztig cycle $\g_i:E(\w)\lr\Z_{\geq0}$ from \cref{def:Lusztig_cycles}, as is implied by our constructions in \cref{sec:lusztigcycles,sec:filtrations_from_braid_words}. In general, the downwards facing global section $\sdsilt_i$ of the dual silting collection might not admit a weave cycle realization. Nevertheless, we will now establish that there is a natural weave cycle realization of its upwards dual $\cusilt_i$:

\begin{theorem}\label{thm:description_of_silting_dual_Lusztig_cycles}
Let $\w$ be a Demazure weave with $m$ trivalent vertices. Then, for any $i\in[1,m]$:
\begin{itemize}
    \item[$(i)$] $\cusilt_i\in \Glsec{\rgraph^{\uparrow}_\w}{\mathcal{F}_\w^{\uparrow}}$ admits a realization by a weave cycle.

    \item[$(ii)$] For any edge $e_j$ of $\rgraph_\w^{\uparrow}$, the object $\cusilt_i(e_j)\in \D$ is rigid. 
\end{itemize}
\end{theorem}

\begin{proof} Let us denote
\begin{equation}\label{eq:partial_dualsiltings}
\cusilt_{j,i}\coloneqq T^{-,>0}_{\nabla^{\uparrow,\w}_j}(\cdots(T^{-,>0}_{\nabla^{\uparrow,\w}_{i+1}}(\nabla_i^{\w,\uparrow}))),\qquad i\in[1,m],\quad  j\in[j,m].
\end{equation}

\noindent Note that $\cusilt_{i,i}=\nabla_i^{\w,\uparrow}$ and $\cusilt_{m,i}=\cusilt_i$, and that \eqref{eq:partial_dualsiltings} is in line with \eqref{eq:partial_simples} and \eqref{eq:partial_cosimples}. Let us fix the index $i\in[1,m]$, we will prove both $(i)$ and $(ii)$ by an iterative argument, scanning the weave $\w$ from bottom to top. Specifically, for each $j\in[i,m]$, we will show that $\on{ev}_{e_{j+1}}(\cusilt_{j,i})\in\D$ is a rigid object, and admits a $(\beta,{\bf a})$-filtration, where the braid word $\beta$ is given by a horizontal slice of $\w$ intersecting the edge $e_{j+1}$ and the weights ${\bf a}$ will be uniquely determined, cf.~\Cref{lem:uniqueness_of_weight}. Note that  $\on{ev}_{e_{j+1}}(\cusilt_{j,i})\simeq \on{ev}_{e_{j+1}}(\cusilt_{i})\in\D$, as follows from \eqref{eq:partial_dualsiltings} and $\on{ev}_{e_{j+1}}(\nabla_l^{\uparrow,\w})\simeq 0$ for all $l\geq j+1$. Recall that we denote $\cusilt_{j,i}(e_{j+1})\coloneqq \on{ev}_{e_{j+1}}(\cusilt_{j,i})$.\\

The iterative argument starts at $j=i$. At the edge $e_{i+1}$, we have $\cusilt_{i,i}(e_{i+1})\simeq S_{p_{i}}$. The object $S_{p_{i}}$ admits a filtration for the weighted braid word arising from the slice above $p_i$ with a unique non-zero weight $1$ concentrated on the strand exiting $p_i$ in the north-east direction.\\

\noindent For the iterative step, consider $j\in[i,m)$ and suppose that $\cusilt_{j,i}(e_{j+1})$ is rigid and it admits a $(\beta,{\bf a})$-filtration as above. First, we have the following equivalences

\begin{equation}\label{eq:weaverealization_proof_equiv1}
\begin{aligned}
\on{Mor}(\cusilt_{j,i},\nabla_{j+1}^{\w,\uparrow})^*&\simeq \on{Mor}_{\mathcal{F}^{\uparrow}_\w(v_{j+1})}((0,\cusilt_{j,i}(e_{j+1}),0),(k,S_{p_{j+1}},\on{id}_{S_{p_{j+1}}}))^*\\
& \simeq \on{Mor}_\D(\cusilt_{j,i}(e_{j+1}),S_{p_{j+1}})^*\\
& \simeq \on{Mor}_\D(S_{p_{j+1}},\cusilt_{j,i}(e_{j+1}))[2]\,.
\end{aligned}
\end{equation}
In \eqref{eq:weaverealization_proof_equiv1}, the first equivalence follows from the adjunction $\on{ev}_{v_{j+1}}\dashv \on{ind}^R_{v_{j+1}}$, the second is a direct computation, and the third follows from the fact that $\D$ is $2$-Calabi--Yau. Second, we have 

\begin{equation}\label{eq:weaverealization_proof_equiv2}
\begin{aligned}
\on{ev}_{e_{j+2}}(\cusilt_{j+1,i})&\stackrel{\eqref{eq:partial_dualsiltings}}{\simeq} \on{ev}_{e_{j+2}}(T^{-,>0}_{\nabla_{j+1}^{\w,\uparrow}}(\cusilt_{j,i}))\\
& \stackrel{\eqref{eq:weaverealization_proof_equiv1}}{\simeq} \on{ev}_{e_{j+2}}(\on{fib}(\cusilt_{j,i}\to \tau_{>0}\left(\on{Mor}_\D(S_{p_{j+1}},\cusilt_{j,i}(e_{j+1}))[2]\right)\otimes \nabla^{\w,\uparrow}_{j+1}))\\
& \simeq  \on{fib}(\cusilt_{j,i}(e_{j+2})\to \tau_{>0}\left(\on{Mor}_\D(S_{p_{j+1}},\cusilt_{j,i}(e_{j+1}))[2]\right)\otimes S_{p_{j+1}})\\
& \simeq T^{-,>0}_{S_{p_{j+1}}}(\cusilt_{j,i}(e_{j+2}))\\
& \simeq T^{-,>0}_{S_{p_{j+1}}}(\cusilt_{j,i}(e_{j+1}))\,,
\end{aligned}
\end{equation}
where the last equivalence follows from $\cusilt_{j,i}(e_{j+1})\simeq \cusilt_{j,i}(e_{j+2})$. Since $\cusilt_{j,i}(e_{j+1})$ is rigid,  \Cref{lem:semi_twists_preserve_rigidity} therefore implies that $\cusilt_{j+1,i}(e_{j+2})$ is also rigid. By iterating the argument, this establishes $(i)$ in the statement.\\

Let us now focus on $(ii)$. For the iterative step at $j\in[i,m)$, we consider the trivalent vertex $p\coloneqq p_{j+1}$ of $\w$ and assume that the south edge of $p_{j+1}$ corresponds to the $b$-th crossing of the braid word $\beta$ (i.e.~the $b$-th labeled strand in the horizontal slice of the weave). Here $\beta$ is the braid word arising from the horizontal slice of $\w$ right below $p_{j+1}$. 
By the iterative hypothesis, $\cusilt_{j,i}(e_{j+1})$ admits a $(\beta,{\bf a})$-filtration, and we let $\C_\ast$ be associated coherent chain complex such that its geometric realization
\begin{equation}\label{eq:weaverealization_proof_equiv3}
\on{tot}(\C_\ast)\simeq \cusilt_{j,i}(e_{j+1})
\end{equation} corresponds to such $(\beta,{\bf a})$-filtration, cf.~\cref{sssec:comments_filteredobject}.(B). Our task is to show that
\[ \cusilt_{j+1,i}(e_{j+2})\stackrel{\eqref{eq:weaverealization_proof_equiv2}}{\simeq} T^{-,>0}_{\nabla_{j+1}^{\w,\uparrow}}(\cusilt_{j,i}(e_{j+1}))\stackrel{\eqref{eq:weaverealization_proof_equiv3}}{\simeq} T^{-,>0}_{\nabla_{j+1}^{\w,\uparrow}}(\on{tot}(\C_\ast))\simeq T^{-,>0}_{S{p_{j+1}}}(\on{tot}(\C_\ast))\]
also admits such a filtration, now corresponding to a weighted braid word arising as a horizontal slice of $\w$ right above $p$ with an appropriate choice of weights. In order to compute $T^{-,>0}_{S{p_{j+1}}}(\on{tot}(\C_\ast))$, we first consider the following morphism between two $3$-term coherent chain complexes in $\D$:
\begin{equation}\label{eq:proof_morphism_complexes}
\begin{tikzcd}
\textcolor{orange}{\on{tot}(\C_{\ast>b})[-2]} \arrow[d, "\on{coev}"] \arrow[r] & \textcolor{orange}{S_{p}^{\oplus a_b}[-1]} \arrow[d, "\alpha"] \arrow[r]                  & \textcolor{orange}{\on{tot}(\C_{\ast<b})} \arrow[d,  "\on{coev}"] \\
\textcolor{blue}{S_{p}\otimes {\tau_{>0}(\on{Mor}(\on{tot}(\C_{\ast> b}),S_{p})^*)[-2]}} \arrow[r]       & \textcolor{blue}{S_{p}^{\oplus a_b}[1]} \arrow[r] & \textcolor{blue}{S_p\otimes \tau_{>0}\left(\on{Mor}(\on{tot}(\C_{\ast<b}),S_{p})^*\right)}           
\end{tikzcd}
\end{equation}
where the left and right vertical morphisms in \eqref{eq:proof_morphism_complexes} are the canonical coevaluations maps, and the middle vertical morphism $\alpha$ is of maximal rank. By \eqref{eq:weaverealization_proof_equiv3}, the top coherent complex \textcolor{orange}{$\C_1$} of \eqref{eq:proof_morphism_complexes}, highlighted in orange, has totalization $\on{tot}(\textcolor{orange}{\C_1})\simeq\cusilt_{j,i}(e_{j+1})$. Similarly, the totalization of the bottom coherent complex \textcolor{blue}{$\C_2$} of \eqref{eq:proof_morphism_complexes}, highlighted in blue, has totalization $\on{tot}(\textcolor{blue}{\C_2})\simeq\tau_{> 0}\on{Mor}_\D(\cusilt_{j,i}(e_{j+1}),S_{p})^*\otimes S_{p}$. The totalization of the (vertical) morphism $f:\textcolor{orange}{\C_1}\lr\textcolor{blue}{\C_2}$ in \eqref{eq:proof_morphism_complexes} is the natural coevaluation morphism. Therefore, to compute $T^{-,>0}_{S{p_{j+1}}}(\on{tot}(\C_\ast))$ we must compute the fiber coherent complex $\fib(f)$ and then its totalization.\\

Instead of directly computing the fiber coherent chain complex, we first fold the left square in \eqref{eq:proof_morphism_complexes}, yielding the following morphism $g$ of coherent $3$-term chain complexes:
\begin{equation}\label{eq:proof_morphism_complexes2}
\begin{tikzcd}[column sep=small]
\textcolor{orange}{\on{tot}(\C_{\ast>b})[-2]} \arrow[d] \arrow[r] & \textcolor{blue}{S_{p}\otimes {\tau_{>0}(\on{Mor}(\on{tot}(\C_{\ast> b}),S_{p})^*)[-2]}} \oplus \textcolor{orange}{S_{p}^{\oplus a_b}[-1]} \arrow[d, "\alpha"] \arrow[r]                  & \textcolor{orange}{\on{tot}(\C_{\ast<b})} \arrow[d,  "\on{coev}"] \\
0 \arrow[r]       & \textcolor{blue}{S_{p}^{\oplus a_b}[1]} \arrow[r] & \textcolor{blue}{S_p\otimes \tau_{>0}\left(\on{Mor}(\on{tot}(\C_{\ast<b}),S_{p})^*\right)}           
\end{tikzcd}
\end{equation}

Note that $\on{tot}(\on{fib}(f))\simeq \on{tot}(\on{fib}(g))$, since folding preserves the totalization.

Let us compute the coherent complex $\fib(g)$. In \eqref{eq:proof_morphism_complexes2}, $\on{tot}(\C_{\ast<b})$ describes the part of $\cusilt_{j,i}(e_{j+1})$ arising from the weights of $\beta$ to the left of the $b$-th crossing of $\beta$, i.e.~to the left of the $b$-th weave edge for that horizontal slice. The term $\on{tot}(\C_{\ast>b})$ similarly describes the piece of $\cusilt_{j,i}(e_{j+1})$ arising from the weights to the right of the $b$-th crossing of $\beta$. Therefore, \Cref{lem:connectiveHom} implies that the truncated morphism object
\begin{equation}\label{eq:weaverealization_proof_AR}
\tau_{>0}\left(\on{Mor}(\on{tot}(\C_{\ast<b}),S_{p})^*\right)\simeq \tau_{>0}\left(\on{Mor}(S_{p},\on{tot}(\C_{\ast<b}))[2]\right)\simeq \on{Mor}(S_{p},\on{tot}(\C_{\ast<b}))[2]
\end{equation}
is equivalent to the entire morphism object, and that
\begin{equation}\label{eq:weaverealization_proof_AL}
\tau_{>0}(\on{Mor}(\on{tot}(\C_{\ast> b}),S_{p})^*)\simeq\tau_{>0}\on{Mor}(S_{p},\on{tot}(\C_{\ast> b}))[2]\simeq \on{Ext}^1(S_{p},\on{tot}(\C_{\ast>b}))[1]
\end{equation}
is concentrated in degree $1$. 

Passing to the fiber coherent $3$-term chain complex of $g$, we thus obtain the coherent $3$-term chain complex $\on{fib}(g)$ 
\begin{equation}\label{eq:fiber_complex_proof}
\on{tot}(\C_{\ast>b})[-2] \to \on{cof}\left(S_{p}^{\oplus a_b+y}[-2]\to S_{p}^{\oplus a_b}\right)\to \on{tot}(T^{-}_{S_{p}}(\C_{\ast<b}))\,,
\end{equation}
\noindent where we have denoted
\begin{equation}\label{def:y}
\framebox{$y\coloneqq \on{dim}\on{Ext}^1(S_{p},\on{tot}(\C_{\ast>b})).$}\,
\end{equation}
By grafting \eqref{eq:proof_morphism_complexes} at its middle term, the totalization of $\fib(g)$ is equivalent to the totalization of a coherent chain complex of the following form
\begin{equation}\label{eq:final_grafting}
\on{tot}(\C_{\ast>b})[-3] \to S_{p}^{\oplus a_b+y}[-2]\to S_{p}^{\oplus a_b}\to \on{tot}(T^{-}_{S_p}(\C_{\ast<b}))\,.
\end{equation}
\noindent Now, the coherent complex \eqref{eq:final_grafting} yields an associated $(\beta',{\bf a'})$-filtration for its totalization, as in \cref{sssec:comments_filteredobject}.(B), with $\beta'$ arising from the slice of $\w$ directly above $p=p_{j+1}$ and the weights being
\[
{\bf a}_i'\coloneqq \begin{cases} {\bf a}_i & i<b \\ {\bf a}_b & i=b \\ {\bf a}_b+\on{dim}\on{Ext}^1(S_{p_{j+1}},\on{tot}(\C_{\ast>b})) & i=b+1 \\ 
{\bf a}_{i-1} & i\geq b+2\,.\end{cases}
\]
Since $\on{tot}(\on{fib}(g))\simeq \on{tot}(\on{fib}(f))$, and by \eqref{eq:proof_morphism_complexes}, the totalization of $\fib(f)$ is equivalent to $T^{-,>0}_{S_{p}}(\cusilt_{j,i}(e_{j}))$, we thus obtain that the object $T^{-,>0}_{S_{p_{j+1}}}(\cusilt_{j,i}(e_{j+1}))$ admits such a $(\beta',{\bf a'})$-filtration. This establishes $(ii)$, where the weave cycle is uniquely determined by the ranks of the graded pieces of the  $(\beta',{\bf a'})$-filtration.
\end{proof}

\noindent \Cref{thm:description_of_silting_dual_Lusztig_cycles} allows us to compute the evaluation of $\sdsilt_i$ at the bottom edge $e_1\in\rgraph_\w$, as we now show that it is equivalent to the evaluation of $\cusilt_i$ at the top edge $\einftop$.

\begin{corollary}\label{cor:description_of_silting_dual_Lusztig_cycles}
Let $\w\colon \beta\to \delta(\beta)$ be a Demazure weave. Then the following holds:

\begin{enumerate}
    \item There exists an equivalence in $\D$
\begin{equation}\label{eq:cor_eval_silting_bottom}
\on{ev}_{e_1}(\sdsilt_i)\simeq \on{ev}_{\einftop}(\cusilt_i).
\end{equation}

In particular, there exists a unique weight ${\bf a}$ for the braid word $\beta$ for which the object $\on{ev}_{e_1}(\sdsilt_i)\in \D$ admits a $(\beta,{\bf a})$-filtration.

\item The evaluation $\on{ev}_{e_1}(\sdsilt_i)\in \D$ is a rigid object.
\end{enumerate} 
\end{corollary}

\begin{proof}
By the proof of \Cref{prop:upwards_downwards_thimbles_duality}, see e.g.~\eqref{eq:ringel_proof_eq3}, we have 
\[\iota^{\downarrow}(\Delta_i^\downarrow)\simeq (\iota^{\downarrow}\circ (\iota^{\downarrow})^L\circ \iota^{\uparrow})(\nabla_i^\uparrow)[-1]\simeq ((T^{\on{gl}})^{-1}\circ \iota^{\uparrow})(\nabla_i^\uparrow)[-1]\,\]
and thus also the equivalence
\begin{equation}\label{eq:corollary_proof_eq1}
\iota^{\downarrow}(\sdsilt_i)\simeq (T^{\on{gl}})^{-1}(\iota^{\uparrow}(\cusilt_i))[-1]\,.
\end{equation}
By \cite[Theorem 2.15]{Chr20}, the inverse cotwist functor $(T^{\on{gl}})^{-1}$ is equivalent to the twist of the spherical adjunction $ (\on{ind}^L_{\einftop},\on{ind}^L_{e_1})\dashv (\on{ev}_{\einftop},\on{ev}_{e_1})$. We thus have 
\begin{equation}\label{eq:corollary_proof_eq2}
\iota^{\downarrow}(\sdsilt_i)\stackrel{\eqref{eq:corollary_proof_eq1}}{\simeq} \on{cof}(\on{ind}^L_{\einftop}\on{ev}_{\einftop}(\cusilt_i)\to \iota^{\uparrow}(\cusilt_i))[-1]\,.
\end{equation}
\noindent Since by construction $\on{ev}_{e_1}(\cusilt_i)\simeq 0$, applying $\on{ev}_{e_1}$ to \eqref{eq:corollary_proof_eq2} yields
\begin{equation}\label{eq:corollary_proof_eq3}
\on{ev}_{e_1}(\sdsilt_i)\simeq \on{cof}((\on{ev}_{e_1}\circ\on{ind}^L_{\einftop}\circ\on{ev}_{\einftop})(\cusilt_i)\to 0)[-1].
\end{equation}
\noindent Finally, note that $\on{ind}^L_{\einftop}(L)$ can be glued from local sections as in \eqref{eq:standard_value_at_vertex} for all $L\in \D$, and thus we have the equivalence of functors
\begin{equation}\label{eq:corollary_proof_eq3}
\on{ev}_{e_1}\circ \on{ind}^L_{\einftop}\simeq \on{id}_{\D}.
\end{equation}
In conclusion, we deduce the required equivalence \eqref{eq:cor_eval_silting_bottom} via
\begin{equation}\label{eq:corollary_proof_eq4}
\on{ev}_{e_1}(\sdsilt_i)\stackrel{\eqref{eq:corollary_proof_eq3}}{\simeq} \on{cof}((\on{ev}_{e_1}\circ\on{ind}^L_{\einftop}\circ\on{ev}_{\einftop})(\cusilt_i)\to 0)[-1]\stackrel{\eqref{eq:corollary_proof_eq3}}{\simeq} \on{cof}(\on{ev}_{\einftop}(\cusilt_i)\to 0)[-1]\simeq \on{ev}_{\einftop}(\cusilt_i)\,.
\end{equation}
The existence of a $(\beta,{\bf a})$-filtration follows from \Cref{thm:description_of_silting_dual_Lusztig_cycles}.$(i)$, with the uniqueness of the weight ${\bf a}$ being guaranteed by \Cref{lem:uniqueness_of_weight}. The rigidity stated in Part (2) is implied by combining \Cref{thm:description_of_silting_dual_Lusztig_cycles}.$(ii)$ and \eqref{eq:cor_eval_silting_bottom}.
\end{proof}

\begin{construction}\label{rem:description_weave_cycle_of_silt}
By the proof of \cref{thm:description_of_silting_dual_Lusztig_cycles}.(i), the weave cycle $\rho_i:E(\w)\lr\Z_{\geq0}$ realizing the object $\cusilt_i$ can be described as follows, scanning $\w$ bottom-to-top:
\begin{enumerate}
   \item For each edge $e\in E(\w)$ containing a point below $p_i$, the weave cycle vanishes: $\rho_i(e)=0$.

    \item At the trivalent vertex $p_i$, the values of the weave cycle $
    \rho_i$ are
    \begin{equation}\label{eq:dualLusztigcycle_start}
    \rho_i(e_1)=0,\qquad \rho_i(e_2)=1\,,
    \end{equation}
    where $e_1$ is the north-western edge of $p_i$, and $e_2$ is the north-eastern edge. Note that by $(1)$ above, $\rho_i$ vanishes in the south edge of $p_i$.
    
    \item At each $4$-valent or $6$-valent vertex of $\w$, the weave cycle $\rho_i$ satisfies the tropical Lusztig rules.

    \item At each trivalent vertex $p_j$, with $j\in(i,m]$, with north-western edge $e_1$, north-eastern edge $e_2$ and south edge $e_s$, the weave cycle $\rho_i$ satisfies 
    \begin{equation}\label{eq:dualLusztigcycle_trivalentrule_left}
    \rho_i(e_1)=\rho_i(e_s)
    \end{equation}
    and 
    \begin{equation}\label{eq:dualLusztigcycle_trivalentrule_right}
    \rho_i(e_2)=\rho_i(e_s)+y_{ij},
    \end{equation}
    where $y_{ij}\in\Z_{\geq0}$ is the dimension of a certain $\on{Ext}^1$-group depending on $i$ and $j$, cf.~\eqref{def:y}.\qed
\end{enumerate}
\end{construction}

\noindent Recall that \cite[Section 4.5\&4.6]{CGGLSS25} introduced intersection pairing $\langle \cdot,\cdot\rangle$ on weave cycles. We conjecture that the weave cycles $\{\rho_i\}_{i\in[1,m]}$, as described in \Cref{rem:description_weave_cycle_of_silt}, are dual to the Lusztig cycles $\{\g_i\}_{i\in[1,m]}$ with respect to this intersection pairing, i.e.~

\begin{equation}\label{eq:duality}
\langle \g_i,\rho_j\rangle =\delta_{ij}\,.
\end{equation}

\begin{center}
	\begin{figure}[h!]
		\centering
        \includegraphics[scale=1]{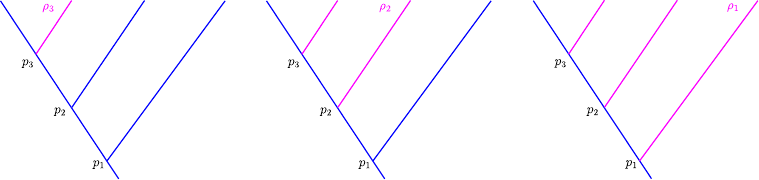}
		\caption{The weave cycles $\{\rho_i\}$ realizing the silting collection dual to the simple-minded collection of Lusztig cycles for this right-inductive weave, cf.~\cref{ex:dual_Lusztig_cycles}.(1). These weave cycles $\{\rho_i\}$ are highlighted in pink.}\label{fig:Example_DualLusztigCycles}
	\end{figure}
\end{center}

\begin{example}\label{ex:dual_Lusztig_cycles} (1) Consider the right-inductive weave $\w$ for $\beta=\s_1^4$, as in \cref{fig:Example_DualLusztigCycles}. The weave cycles $\{\rho_1,\rho_2,\rho_3\}$ realizing the silting collection $\{\cusilt_1,\cusilt_2,\cusilt_3\}$ are depicted in pink in \cref{fig:Example_DualLusztigCycles}, where a pink edge indicates weight 1, and a blue edge indicates weight 0. These weave cycles are obtained using \cref{rem:description_weave_cycle_of_silt}, as follows.\\

For $\rho_1$, we start at the trivalent vertex $p_1$, by \ref{rem:description_weave_cycle_of_silt}.(1), and assign weight 1 to the north-eastern edge of $p_1$, as in \ref{rem:description_weave_cycle_of_silt}.(2). As we scan upwards we reach the trivalent vertex $p_2$, whose assigned vanishing $2$-spherical object is $S_{p_2}\simeq S_1[1]$. To proceed according to \ref{rem:description_weave_cycle_of_silt}.(4), we use \eqref{eq:dualLusztigcycle_trivalentrule_left} and \eqref{eq:dualLusztigcycle_trivalentrule_right} to propagate upwards. From \eqref{eq:dualLusztigcycle_trivalentrule_left}, we deduce that the weight of $\rho_1$ on the north-western edge of $p_2$ is zero. To apply \eqref{eq:dualLusztigcycle_trivalentrule_right}, we must compute the value of
    $$y_{12}\coloneqq \on{dim}\on{Ext}^1(S_{p_2},\on{tot}(\C_{\ast>b}))=\on{dim}\on{Ext}^1(S_1[1],\on{tot}(\C_{\ast>b})),$$
    according to \eqref{def:y}.
    The braid word $\beta$ that defines $\C_{\ast>b}$ in this case is $\beta=\s_1^2$, as it is obtained from the horizontal slice right below $p_2$. The corresponding weights are ${\bm a}=(a_0,a_1)=(0,1)$ and hence  the coherent chain complex is
    $$\C_{\ast}=(S_1[-1]\lr 0).$$
       
    \noindent Since $b=0$ in this case, as it is the leftmost crossing of $\beta$, 
    the required totalization $\on{tot}(\C_{\ast>b})$ is
    $$\on{tot}(\C_{\ast>b})\simeq\cof(S_1[-1]\lr 0)\simeq S_1$$
    and the value of the correction term $y_{12}$ is
    \begin{equation}\label{eq:example1_weight12}
    y_{12}=\on{dim}\on{Ext}^1(S_1[1],S_1)=1,
    \end{equation}
    since $\on{Ext}^1(S_1[1],S_1)=\on{Ext}^0(S_1,S_1)=k$. By using \ref{rem:description_weave_cycle_of_silt}.(4) and \cref{eq:example1_weight12}, we therefore assign weight $y_{12}=1$ to $\rho_1$ as we move upwards. This is why the north-eastern edges of $p_1$ and $p_2$ are highlighted in pink in \cref{fig:Example_DualLusztigCycles} (right).\\

    We now continue propagating $\rho_1$ upwards until we reach $p_3$, where we must use \eqref{eq:dualLusztigcycle_trivalentrule_left} and \eqref{eq:dualLusztigcycle_trivalentrule_right} to propagate past $p_3$. As before, \eqref{eq:dualLusztigcycle_trivalentrule_left} implies that the weight of $\rho_1$ on the north-western edge of $p_3$ must be zero. From \eqref{eq:dualLusztigcycle_trivalentrule_right}, the weight on the north-eastern edge is given by
    $$y_{13}\coloneqq \on{dim}\on{Ext}^1(S_{p_3},\on{tot}(\C_{\ast>b})).$$
    The computation is similar to the above one for $y_{12}$. In this case $S_{p_3}\simeq S_1[2]$, the braid word is $\beta=\s_1^3$ and the weights are ${\bm a}=(a_0,a_1,a_2)=(0,1,1)$. Hence, the coherent chain complex $\C_{\ast}$ has the form
    \begin{equation}\label{eq:example1_coherentchaincomplex_p3}
    \C_{\ast}=(S_{\s_1}^{\oplus a_2}[-2]\lr S_{\s_1^2}^{\oplus a_1}[-1]\lr S_{\s_1^3}^{\oplus a_0})\simeq (S_1[-2]\lr S_1\lr 0)
    \end{equation}
    As above, $b=0$ also in this case. Therefore the required totalization $ \on{tot}(\C_{\ast>b})$ is equivalent to
    \begin{equation}\label{eq:example1_totalization_coherentchaincomplex_p3}
    \on{tot}(\C_{\ast>b})\simeq \cof(S_1[-1]\lr \cof(S_1\lr 0))\simeq \cof(S_1[-1]\lr S_1[1]).
    \end{equation}
    In consequence,
    \begin{equation}\label{eq:example1_y13_value}
    y_{13}\stackrel{\eqref{eq:example1_totalization_coherentchaincomplex_p3}}{=}\on{dim}\on{Ext}^1(S_1[2],\cof(S_1[-1]\lr S_1[1]))=\on{dim}\on{Ext}^1(S_1[2],S_1[1])=1,
    \end{equation}
    as the natural surjective map $\on{Ext}^1(S_1[2],S_1[1])\lr\on{Ext}^1(S_1[2],\cof(S_1[-1]\lr S_1[1]))$ from applying $\Mor(S_1[2],\mhyphen)$ to the cofiber sequence associated to \eqref{eq:example1_totalization_coherentchaincomplex_p3} is an isomorphism. This proves that $\rho_1$ as drawn in Figure \ref{fig:Example_DualLusztigCycles} is indeed the correct weave cycle realizing $\cusilt_1$. The proof of the correctness of $\rho_2$ and $\rho_3$ is much simpler, and the necessary argument and computations are already contained in the proof for $\rho_1$ above. Finally, note that the duality \eqref{eq:duality} with respect to the weave intersection pairing indeed holds in this example.\\

\begin{center}
	\begin{figure}[h!]
		\centering
		\includegraphics[scale=1]{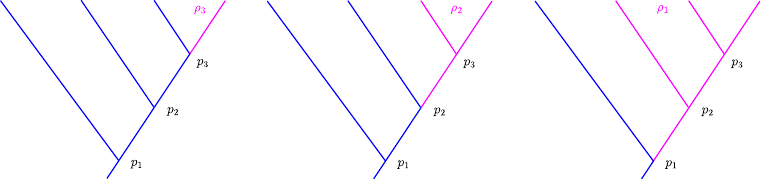}
		\caption{The dual Lusztig weave cycles $\{\rho_i\}$ from \cref{ex:dual_Lusztig_cycles}.(2), highlighted in pink.}\label{fig:Example_DualLusztigCycles2}
	\end{figure}
\end{center}

\noindent (2) Consider the left-inductive weave $\w$ for $\beta=\s_1^4$, as in \cref{fig:Example_DualLusztigCycles2}. In this case, all the correction terms $y_{12},y_{13},y_{23}$ vanish. Indeed, scanning this left-inductive weave $\w$ upwards, any trivalent vertex $p_i$ appears for the rightmost crossing of the braid word associated to the horizontal slice right below $p_i$. Thus the $b$ value in \cref{def:y} always corresponds to the index for the rightmost crossing and so $\C_{\ast>b}$ is the zero complex. Therefore the correction terms

$$y_{ij}\coloneqq \on{dim}\on{Ext}^1(S_{p_2},\on{tot}(0))=\on{dim}\on{Ext}^1(S_{p_2},0)=0$$

\noindent must always vanish. By applying  \cref{rem:description_weave_cycle_of_silt}, we obtain the weave cycles $\{\rho_1,\rho_2,\rho_3\}$ realizing the silting collection $\{\cusilt_1,\cusilt_2,\cusilt_3\}$ for this left-inductive weave $\w$. As in (1), the duality \eqref{eq:duality} with respect to the weave intersection pairing holds as well.\\

\begin{center}
	\begin{figure}[h!]
		\centering
		\includegraphics[scale=0.8]{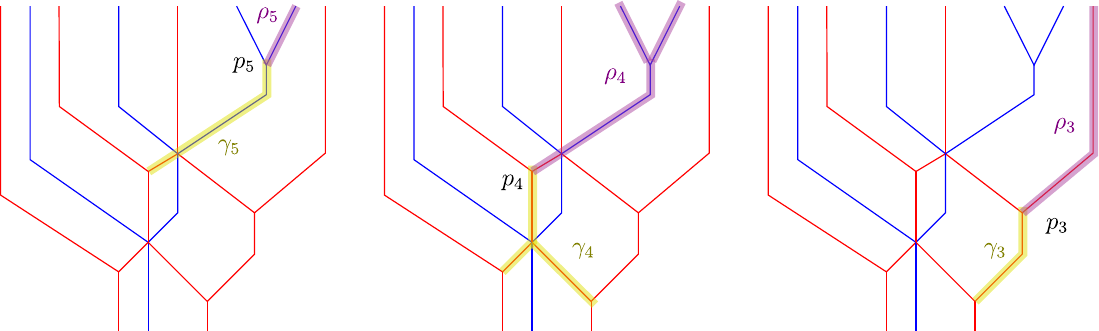}
		\caption{A few of the weave cycles $\{\rho_i\}$ from \cref{ex:dual_Lusztig_cycles}.(3), highlighted in pink, with the corresponding (absolute) Lusztig cycles $\{\g_i\}$ highlighted in yellow.}\label{fig:Example_DualLusztigCycles3}
	\end{figure}
\end{center}

\noindent (3) Consider the weave $\w$ depicted in \cref{fig:Example_DualLusztigCycles3}. The Lusztig cycles $\g_3,\g_4,\g_5$ associated to the trivalent vertices $p_3,p_4,p_5$ are highlighted in yellow. Specifically, these are all binary weave cycles, with the only values being $0$ and $1$, so we directly highlight the edges where the cycle evaluates to 1. The upwards dual weave cycles $\rho_3,\rho_4,\rho_5$
realizing the objects $\cusilt_3,\cusilt_4,\cusilt_5$ in the silting collection are highlighted in purple. These weave cycles $\rho_3,\rho_4,\rho_5$ are obtained by applying \cref{rem:description_weave_cycle_of_silt}, as follows.\\

For $\rho_5$, it suffices to apply only rules \ref{rem:description_weave_cycle_of_silt}.(1)
and \ref{rem:description_weave_cycle_of_silt}.(2), as there are no weave vertices above $p_5$. For $\rho_4$ and $\rho_3$, we use all the rules from \cref{rem:description_weave_cycle_of_silt}. In order to apply \ref{rem:description_weave_cycle_of_silt}.(4) for $\rho_4$ and $\rho_3$, we claim that the corresponding correction terms $y_{45}$, for $\rho_4$, and $y_{34}$ and $y_{35}$, for $\rho_3$, all vanish. First, the corresponding coherent chain complex $\C_{\ast>b}$ at $p_5$ for $\rho_4$ is the zero complex because $\rho_4$ has no non-zero weights to the right of $p_5$ in the horizontal slice right below $p_5$. Hence $y_{45}=0$ indeed vanishes. Second, let us argue that $y_{34}=0$. For that, we compute
\begin{equation}\label{eq:y34_1}
S_{p_4}=T_{S_2}^{-2}T_{S_1}^{-1}(S_2)\simeq T_{S_2}^{-1}(S_1)
\end{equation}
and note that the coherent chain complex  $\on{tot}(\C_{\ast})$ for $\rho_3$ at $p_4$ has the form
\begin{equation*}\label{eq:y34_2}
S_{\s_2}^{\oplus a_5}[-5]\lr S_{\s_2^2}^{\oplus a_4}[-4]\lr S_{\s_1\s_2^2}^{\oplus a_3}[-3]\lr S_{\s_2\s_1\s_2^2}^{\oplus a_2}[-2]\lr S_{\s_1\s_2\s_1\s_2^2}^{\oplus a_1}[-1]\lr S_{\s_2\s_1\s_2\s_1\s_2^2}^{\oplus a_0}
\end{equation*}
as $\beta=\s_2\s_1\s_2\s_1\s_2^2$. Since the weights for $\rho_3$ are ${\bm a}=(0,0,0,0,0,1)$, and $b=2$ for $p_4$, the totalization $\on{tot}(\C_{\ast>2})$ required to compute $y_{34}$ is \begin{equation}\label{eq:y34_3}
    \on{tot}(\C_{\ast>2})\simeq \on{tot}(S_{2}[-5]\lr 0\lr 0\lr 0\lr 0\lr 0)\simeq S_2.
    \end{equation}
The computations \eqref{eq:y34_1} and \eqref{eq:y34_3} then imply
\begin{equation}\label{eq:y34_4}
\begin{aligned}
y_{34}&\stackrel{\eqref{eq:y34_1}}{=}\on{dim}\on{Ext}^1(T_{S_2}^{-1}(S_1),\on{tot}(\C_{\ast>2}))\\
&\stackrel{\eqref{eq:y34_3}}{=}\on{dim}\on{Ext}^1(T_{S_2}^{-1}(S_1),S_2)\\
&=\on{dim}\on{Ext}^1(S_1,T_{S_2}(S_2))\\
&=\on{dim}\on{Ext}^1(S_1,S_2[-1])\\
&=\on{dim}\on{Ext}^0(S_1,S_2)=0\
\end{aligned}    
\end{equation}
Third, the vanishing $y_{35}=0$ is similar. Indeed, the spherical object at $p_5$ is $S_{p_5}\simeq T_{S_2}^{-1}(S_1)$, the braid word below $p_5$ is $\beta=(\s_2\s_1)^2\s_2\s_1\s_2$, the weights for $\rho_3$ are ${\bm a}=(0,0,0,0,0,0,1)$, and $b=5$ for $p_5$. Thus $\on{tot}(\C_{\ast>5})\simeq S_2$ in this case. In consequence \begin{equation*}\label{eq:y35}
\begin{aligned}
y_{35}=\on{dim}\on{Ext}^1(T_{S_2}^{-1}(S_1),S_2)=0.
\end{aligned}    
\end{equation*}  
This establishes the correctness of the binary weave cycles $\{\rho_3,\rho_4,\rho_5\}$ as drawn in \cref{fig:Example_DualLusztigCycles3}. Note that there are five Lusztig cycles but, since $\g_1,\g_2$ are frozen,  $\g_3,\g_4,\g_5$ are the only mutable Lusztig cycles. The corresponding intersections are
$$(\langle\g_i, \rho_j\rangle)=\begin{pmatrix}
1 & 0 & 0\\
0 & 1 & 0\\
0 & 0 & 1
\end{pmatrix},\quad \mbox{for }i,j\in[3,5],$$
with $i=j=3$ being the top left entry and $i=j=5$ the bottom right entry. The only subtle entry is $\langle \la_4,\g_5\rangle=0$, where the intersection is 0 despite the (support of the) cycles geometrically overlapping. Hence the duality \eqref{eq:duality} holds.\qed
\end{example}


\subsection{Ancillary remarks on tilting theory and exceptional collections}\label{ssec:tiltingtheory_exceptionalcollections}

Part of the results and constructions of this paper are inspired by and related to the theory of highest weight abelian categories, cf.~\cite{CPS88} or \cref{def:highest_weight_category}. This final subsection summarizes a few key aspects in the main text as they relate to this perspective.\\

\noindent Must of our results occur in the context of stable $\infty$-categories, or at least triangulated categories, in contrast to that of abelian categories. Indeed, while abelian categories suffice to tackle reduced braid words, cf.~e.g.~\eqref{eq:reduced_in_heart}, \cref{rmk:2sphericalobjects_tstructures} and \cite{GLS06,GLS08,BIRS09}, the derived setting is required in order to establish a theory for non-reduced positive braid words. Therefore, there is merit in suggesting a potential candidate for the notion of a highest weight structure in the stable context. In a nutshell, our constructions and results lead to the following potential tenet:

\begin{adjustwidth}{1cm}{1cm}
\begin{center}
{\it A stable $k$-linear $\infty$-category endowed with a full exceptional collection can be considered as the stable analog of an abelian category endowed with a highest weight structure.\\}
\end{center}
\end{adjustwidth}

\noindent This tenet is analyzed in detail in Appendix \ref{ssec:AppendixA}. First, in Appendix \ref{subsec:highest_weight_package}, we associate a simple-minded collection and a Koszul dual silting collection to a stable $\infty$-category with a full exceptional collection. These collections are the stable analogs of the simple objects and projective objects in a highest weight abelian category, respectively. Similarly to how the projective objects in a highest weight abelian category arise via universal extensions of the standard objects, the silting collection arises via universal (possibly higher) extensions between the standard objects by \Cref{prop:silting_from_standards}. We also refer to \Cref{thm:relation_with_highest_weight_abelian_cats} for a precise relationship with highest weight abelian categories.\\

\noindent Second, highest weight abelian categories admit so-called tilting objects, which are objects characterized by admitting both a standard and a costandard flag. A stable version of the tilting object, which is both a higher extension of the standard exceptional objects and a higher extension of the costandard exceptional objects, appears in \Cref{lem:standard_vs_costandard_upwards_silt}. In the theory of highest weight abelian categories, tilting objects are used to construct the Ringel duality functors, which are functors from a highest weight abelian category to its Ringel dual. Within the formalism of perverse schobers, we presented an abstract construction of two Ringel-type duality functors in \Cref{lem:upwards_downwards_duality_schobers}.(3). They intertwine the standard and costandard exceptional collections, cf.~\Cref{prop:upwards_downwards_thimbles_duality}, as well the silting object and the stable version of tilting object, cf.~\Cref{lem:Ringel_duality_and_silting_objects}. A difference with respect to classical Ringel duality is that the definitions of these functors are not module-theoretical and rather in terms of perverse schobers.\\

\noindent The following table illustrates parts of the above summary, see Appendix \ref{ssec:AppendixA} for more details.

\begin{table}[h!]
  \centering
  \label{table:summary_highestweight_stable}
\begin{tabular}{ |P{4cm}|P{5cm}|P{5cm}|}
	\hline
    
	{\bf Highest weight abelian category} & {\bf stable $k$-linear $\infty$-category} & {\bf Expression from $(\Delta_1,\ldots,\Delta_m)$}\\
	\hline

    standard objects & exceptional collection $\{\Delta_i\}$ & $\Delta_i$\\
	\hline

    \rule{0pt}{0.5cm}costandard objects & coexceptional collection $\{\nabla_i\}$ & $T_{\Delta_1}^+T_{\Delta_2}^+\cdots T_{\Delta_{i-1}}^+(\Delta_i)$\\
	\hline

    \rule{0pt}{0.5cm}simples $S_i$ & simple-minded collection $\{\lsimp_i\}$ & $T_{\Delta_{1}}^{+,\geq 0}(\cdots (T_{\Delta_{i-1}}^{+,\geq 0}(\Delta_i)))$\\
	\hline

    \rule{0pt}{0.5cm}projectives $P_i$ & silting collection $\{\lsilt_i\}$  & $T_{\Delta_{m}}^{-,>0}(\cdots (T_{\Delta_{i+1}}^{-,>0}(\Delta_i)))$\\
    \hline

    \rule{0pt}{0.5cm}injectives $I_i$ & silting collection $\{\rsilt_i\}$ & $T_{\nabla_{m}}^{+,<0}(\cdots (T_{\nabla_{i+1}}^{+,<0}(\nabla_i)))$ \\
    \hline

    \rule{0pt}{0.5cm}tilting objects $T_i$ & silting collection $\{\lsilt_i^\uparrow\}=\{\rsilt_i^\uparrow\}$ & $T^{+,<0}_{\Delta_{m}^{\uparrow}}(\dots (T^{+,<0}_{\Delta_{i+1}^{\uparrow}}(\Delta_i^{\uparrow})))$\\
    \hline
\end{tabular}
  \caption{Summary of the stable $\infty$-categorical analogues (middle) of the central objects in a highest weight abelian category (left), cf.~Appendix \ref{ssec:AppendixA}. The middle and right columns also apply to (enhanced) triangulated categories, whereas the left column is strictly restricted to abelian categories.}
\end{table}
\vspace{0.1cm}
\noindent In the above table, the index $i\in[1,m]$ always ranges from $1$ to $m$, $T_Y^{+,\geq0}$, $T_Y^{-,>0}$ and $T_Y^{+,<0}$ denote semi-twists and $T^+$ a twist, $\Delta_i^\uparrow$ denote the exceptional collection in the Ringel-type dual category, and all the exceptional and co-exceptional collections are always assumed to be full. Note also that technically we build two simple-minded collections $\{\lsimp_i\}$ and $\{\rsimp_i\}$, but we prove in \cref{lem:standard_costandard_SMC_coincide} that $\{\lsimp_i\}=\{\rsimp_i\}$.

\begin{remark}\label{rem:full_exceptional_from_schobers}
One can show that any proper $k$-linear stable $\infty$-category $\C$ with a full exceptional collection $\Delta_1,\dots,\Delta_m\in \C$ arises formally as the $\infty$-category of global sections of a perverse schober parametrized by a linear ribbon graph with one external edge, as follows. First, there is a functor $\C\to \hat{\C}$ from $\C$ into the proper $2$-Calabi--Yau completion $\hat{\C}$ of $\C$, mapping $\Delta_i$ to a $2$-spherical object $X_i$. Then, we can use the ordered collection of spherical objects $X_1,\dots,X_m\in \hat{\C}$ to define a perverse schober $\mathcal{F}$, in line with \Cref{def:weave_schobers}. The $\infty$-category of global sections of $\mathcal{F}$ is then equivalent to $\C$ in such a manner that the $i$th exceptional object $\Delta_i$ is identified with the $i$th standard exceptional object defined using left induction as in \Cref{eq:thimble_as_induction}. See also \cite{Christ_excep_coll} for more details, and note that this passage to the Calabi--Yau completion can be seen as an algebraic version of the passage to the double branched cover of a Lefschetz fibration in \cite[Section (18a)]{Sei08}.\qed
\end{remark}

Each of the notions of an exceptional, a simple-minded, or a silting collection admits a notion of mutation, cf.~\cref{subsec:recollections_on_exceptional_collections}. For the latter two, these are part of tilting and silting theory, respectively. Now, a central result in the study of the tilting theory of simple-minded collections and the silting theory of silting collections is that Koszul duality commutes with these respective notions of mutation, see \Cref{thm:silting_vs_SMC}. The relation between mutations of a given full exceptional collection and mutations of the corresponding simple-minded collection, or the corresponding silting collection, is more subtle, cf.~\Cref{thm:Lusztigcycles_weaveequivalence} and \cref{thm:weavemutation_SMCmutation1}. Based on these two results in our main text, we now present two sufficiency homological criteria ensuring that an exceptional mutation induces either a simple-minded mutation, or that it induces no mutation of the corresponding simple-minded collection.

\begin{theorem}
Let $\C$ be a proper $k$-linear stable $\infty$-category with a full exceptional collection $\Delta_1,\dots,\Delta_m$, and $i\in[1,m]$ a fixed index. Suppose that 
\[ \on{Ext}^0_\C(\Delta_i,\Delta_{i+1})\simeq k\,,\]
\[ \on{Ext}_\C^l(\Delta_j,\on{cof}(\Delta_i\to \Delta_{i+1}))\simeq 0\,,\qquad \forall l\in \mathbb{Z}_{\geq 0}, j\in [1,i)
\,,\]
as well as 
\[ \on{Ext}^{-k}_\C(\Delta_i,\Delta_{i+1})\simeq 0,\qquad \forall k\in\Z_{>0}\,.\]
Then the associated simple-minded collection $\{\lsimp_j\}_{1\leq j\leq m}$ undergoes a backwards mutation at $\lsimp_i$ under a right mutation of the exceptional collection at the pair $(\Delta_i,\Delta_{i+1})$.  
\end{theorem}

\begin{proof}
By \Cref{rem:full_exceptional_from_schobers}, we may assume that $\C$ arises as the global sections of a perverse schober on a linear graph with $m$ vertices and one external edge arising from a collection of $m$ many $2$-spherical objects, and such that the full exceptional collection arises as the standard objects obtained from left induction. Then the proof of \Cref{thm:Lusztigcycles_weaveequivalence}, see  \Cref{ssec:Lusztigcycles_weavemutation}, applies (with only cosmetic changes) to prove the assertion. 
\end{proof}

Via a similar argument, the proof of \cref{thm:weavemutation_SMCmutation1} implies the following:

\begin{theorem}
    Let $\C$ be a proper $k$-linear stable $\infty$-category with a full exceptional collection $\Delta_1,\dots,\Delta_m$, and $i\in[1,m]$ a fixed index. Suppose that
\[ \on{Ext}^{-k}_\C(\Delta_i,\Delta_{i+1})\simeq 0,\qquad \forall k\in\Z_{\geq0}.\]
Then, up to permutation, the associated simple-minded collection $\{\lsimp_i\}$ is invariant under a right mutation of the exceptional collection at $(\Delta_i,\Delta_{i+1})$.  
\end{theorem}

\appendix

\newpage 

\section{Exceptional, simple-minded, and silting collections in enhanced triangulated categories}\label{ssec:AppendixA}

The object of this appendix is to study simple-minded and silting collection in $k$-linear stable $\infty$-categories arising from full exceptional collections. Specifically, in \cref{subsec:recollections_on_exceptional_collections} we recall the notions of exceptional collections, simple-minded collections, and silting collections, as well as their corresponding mutations. In \cref{subsec:highest_weight_package}, we describe how a full exceptional collection in a $k$-linear stable $\infty$-category gives rise to Koszul dual simple-minded and silting collections. In \Cref{subsec:relation_w_highest_weight_theory}, we discuss how our constructions generalize aspects of the theory of highest weight abelian categories. We refer to \Cref{ssec:tiltingtheory_exceptionalcollections} for a discussion on how the constructions in the main part of the manuscript relate to those in this appendix and the theory of highest weight abelian categories.\\

Throughout this section, we formulate all notions and statements in the setting of $k$-linear $\infty$-categories, to keep the notation consistent with the main text. The results can equally be formulated on the level of the $k$-linear triangulated homotopy $1$-categories\footnote{For instance, the notion of a $t$-structure on a stable $\infty$-category amounts exactly to a $t$-structure on the triangulated homotopy $1$-category.} However, some of the proofs make use of the enhancement, specifically when they appeal to Koszul duality, see \Cref{thm:silting_vs_SMC}.

\subsection{Recollections}\label{subsec:recollections_on_exceptional_collections}

\noindent The objects of interest in this appendix are defined as follows.

\begin{definition}\label{def:recollection_exceptional_SMC_silting}
    Let $\C$ be a small $k$-linear stable $\infty$-category. 
    \begin{enumerate}[(1)]
        \item A finite sequence of objects $\Delta_1,\dots,\Delta_m$ in $\C$ is called an exceptional collection if 
    \begin{itemize}
        \item $\on{Mor}_{\D}(\Delta_i,\Delta_i)\simeq k$ for all $1\leq i \leq m$, and
        \item $\on{Mor}_{\D}(\Delta_i,\Delta_j)\simeq 0$ for all $1\leq j<i\leq m$.
    \end{itemize}
    The exceptional collection is said to be a full exceptional collection if in addition
    \begin{itemize}
        \item the object $\bigoplus_{i=1}^m \Delta_i$ generates $\C$ under finite limits and colimits.
    \end{itemize}
    \item A finite sequence of objects $\nabla_1,\dots,\nabla_m$ in $\C$ is called a (full) co-exceptional collection if $\nabla_m,\dots,\nabla_1$ is a (full) exceptional collection.
    \item A full exceptional collection $\Delta_1,\dots,\Delta_m$ and a full co-exceptional collection $\nabla_m,\dots,\nabla_1$ in $\C$ are said to be dual to each other if
    \[
    \on{Mor}_\C(\Delta_i,\nabla_j)\simeq \begin{cases} k & i=j, \\ 0 & \text{else}, \end{cases}\,\qquad \forall i,j\in[1,m].
    \]
    \item A finite collection of indecomposable objects $S_1,\dots,S_m$ in $\C$ is called a simple-minded collection if 
    \begin{itemize}
        \item $\on{Mor}_{\C}(S_i,S_i)$ is coconnective, i.e.~$\on{Ext}^l(S_i,S_i)\simeq 0$ for all $l<0$, for all $1\leq i \leq m$,
        \item $\on{Mor}_{\C}(S_i,S_j)[1]$ is coconnective, i.e.~$\on{Ext}^{l}(S_i,S_j)\simeq 0$ for all $l\leq 0$, for all $1\leq i,j\leq m$ with $i\not =j$, and
        \item the object $\bigoplus_{i=1}^m S_i$ generates $\C$ under finite limits and colimits as well as direct summands. 
    \end{itemize} 
    \item A finite collection of distinct indecomposable objects $P_1,\dots,P_m$ in $\C$ is called a silting collection if $P=\bigoplus_{i=1}^mP_i$ is a silting object, that is:
    \begin{itemize}
        \item $\on{Mor}_\C(P,P)$ is connective, i.e.~$\on{Ext}^i(P,P)\simeq 0$ for all $i>0$,
        \item $P$ generates $\C$ under finite limits and colimits and direct summands.
    \end{itemize}
    We note that the notion of a silting object is self-dual, i.e.~$P$ is a silting object in $\C$ if and only if it is a silting object in $\C^{\on{op}}$. 
    
    \item A simple-minded collection $S_1,\dots,S_m$ in $\C$ is said to be Koszul dual to a silting collection $P_1,\dots,P_m$ if
    \[ \on{Mor}_\C(P_i,S_i)\simeq \begin{cases} k & i=j, \\ 0 & \text{else},\,\end{cases}\qquad \forall i,j\in[1,m].\]
\qed
    \end{enumerate}
\end{definition}

\noindent For \cref{def:recollection_exceptional_SMC_silting}.(6), we remark that a simple-minded collection is uniquely determined by being Koszul dual to a given silting collection and vice versa.

\begin{remark}
\begin{enumerate}[(1)]
    \item Recall that a small $k$-linear stable $\infty$-category $\C$ is called proper if,  for all objects $X,Y\in \C$, the derived Hom object $\on{Mor}_\C(X,Y)\in \D(k)$ has finite dimensional total homology, i.e.~lies in $\D^{\on{perf}}(k)$. In this appendix, we will only study full exceptional collections in proper $k$-linear stable $\infty$-categories.
    \item One can show that if a proper $k$-linear stable $\infty$-category admits a full exceptional collection, then it is also idempotent complete and smooth.\qed
\end{enumerate}
\end{remark}

\noindent Exceptional collections, simple-minded collections and silting collections each admits a notion of mutation, as follows. (See also \cite[Section 7]{KY14}.) We begin by recalling the mutation of exceptional collections, which requires the notion of a twist.

\begin{definition}[Twists]\label{def:exceptionalmutation}
Let $\C$ be a proper $k$-linear stable $\infty$-category and $\Delta_{i-1},\Delta_i,\Delta_{i+1}\in \C$.
\begin{enumerate}[(1)]
    \item The positive twist of $\Delta_i$ by $\Delta_{i-1}$ is defined as
    \[ T_{\Delta_{i-1}}^{+}(\Delta_i)\coloneqq \on{cof}(\on{Mor}_\C(\Delta_{i-1},\Delta_i)\otimes \Delta_{i-1}\longrightarrow \Delta_i) \,,\]
    where the morphism is given by the evaluation map.\\
    
    \item The negative twist of $\Delta_i$ by $\Delta_{i+1}$ is defined as
    \[ T_{\Delta_{i+1}}^-(\Delta_i)\coloneqq \on{fib}(\Delta_i\longrightarrow \on{Mor}_\C(\Delta_{i},\Delta_{i+1})^*\otimes \Delta_{i+1}) \,,\]
    where the morphism is given by the coevaluation map.\\
\qed
\end{enumerate}
\end{definition}

\begin{lemma}\label{lem:exceptionalmutation}
Let $\C$ be a proper $k$-linear stable $\infty$-category. Let $\Delta_1,\dots,\Delta_m$ be a full exceptional collection and $\nabla_1,\dots,\nabla_m$ a full co-exceptional collection in $\C$.
\begin{enumerate}[(1)]
    \item Fix $i\in[2,m]$. Then the objects
    \[ \Delta_1,\dots,\Delta_{i-2},T_{\Delta_{i-1}}^{+}(\Delta_i),\Delta_{i-1},\Delta_{i+1},\dots,\Delta_m\]
    form an exceptional collection, called the forward or left mutation at the pair $(\Delta_{i-1},\Delta_{i})$.
      Similarly, the objects 
    \[
    \nabla_1,\dots, \nabla_{i-2},\nabla_{i}, T^+_{\nabla_i}(\nabla_{i-1}), \nabla_{i+1},\dots,\nabla_m
    \]
    form a co-exceptional collection, called the forward mutation at $(\nabla_{i-1},\nabla_i)$.
    \item  Fix $i\in[1,m-1]$. Then the objects
    \begin{equation}\label{eq:backwardsexceptionalmutation} \Delta_1,\dots,\Delta_{i-1},\Delta_{i+1},T_{\Delta_{i+1}}^{-}(\Delta_i),\Delta_{i+2},\dots,\Delta_m\end{equation}
    form an exceptional collection, called the backward or right mutation at $(\Delta_i,\Delta_{i+1})$.  
     Similarly, the objects 
    \[
    \nabla_1,\dots, \nabla_{i-1},T^-_{\nabla_{i}}(\nabla_{i+1}),\nabla_{i}, \nabla_{i+2},,\dots,\nabla_m
    \]
    form a co-exceptional collection, called the backward or right mutation at $(\nabla_i,\nabla_{i+1})$.\qed
\end{enumerate}
\end{lemma}

\begin{remark}\label{rmk:exceptional_to_coexceptional}
If $\Delta_1,\dots,\Delta_m$ is a full exceptional collection in a proper $k$-linear stable $\infty$-category $\C$, then there exists a unique dual full co-exceptional collection $\nabla_1,\dots,\nabla_m$, given by 
\begin{equation}\label{eq:costandards_from_standards}
\nabla_i=T_{\Delta_1}^+T_{\Delta_2}^+\cdots T_{\Delta_{i-1}}^+(\Delta_i)\,,
\end{equation}
see for instance \cite[Section 5k]{Sei08}.\qed
\end{remark}

For mutations of simple-minded collections, we proceed as follows. First, fixing a simple-minded collection $S_1,\dots,S_m$ in $\C$ in a $k$-linear stable $\infty$-category, it determines a bounded $t$-structure on $\C$ whose heart $H$ is a finite length abelian category, see for instance \cite[Theorem 4.4]{Sch_simple_minded}. Let us assume that the trianguled homotopy category of $\C$ is also $\on{Hom}$-finite and Krull--Schmidt. Then the object $S_i$ in the simple-minded collection
determines torsion pairs\footnote{Recall that a torsion pair $\langle F,T\rangle$ in an abelian category $A$ consists of two full subcategories $F,T\subset A$, such that $\on{Hom}_A(Y,X)=0$ for all $X\in F$ and $Y\in T$, and for every $Z\in A$, there exists a short exact sequence $X\to Z \to Y$ with $X\in F$ and $Y\in T$.} $\langle \on{add}(S_i),T\rangle$ and $\langle F,\on{add}(S_i)\rangle$, with $\on{add}(S_i)$ the additive closure of $S_i$, see \cite[Theorem 3.3]{Dug15}.\\ 

\noindent The forward tilt $H^\sharp$ of $H$ at $\on{add}(S_i)$ is defined as the abelian subcategory of $\C$ arising from the torsion pair $\langle T,\on{add}(S_i)[1]\rangle$. Similarly, the backward tilt $H^\flat$ of $H$ is defined as the abelian subcategory of $\C$ arising from the torsion pair $\langle \on{add}(S_i)[-1],F\rangle$.

\begin{definition}[Simple-minded mutation]\label{def:forward_mutation_of_simple_minded_collections}
Let $\C$ be a $k$-linear stable $\infty$-category whose triangulated homotopy category is $\on{Hom}$-finite and Krull--Schmidt. Let $S_1,\dots,S_m$ be a simple-minded collection in $\C$ describing the simple objects in a heart $H\subset \C$. The forward mutation of $S_1,\dots,S_m$ is defined as the collection of simple objects in the forward tilt $H^\sharp$, considered as objects of $\C$. The backwards mutation of $S_1,\dots,S_m$ is similarly defined as the simple objects in the backward tilt $H^\flat$.\qed
\end{definition}

\noindent The following result provides a more explicit description of the tilts at an object $S_i\in\C$ of a simple-minded collection, as in \cref{def:forward_mutation_of_simple_minded_collections}, in the case that the object $S_i$ is itself rigid.

\begin{proposition}[{$\!\!$\cite[Prop.~5.4]{KQ15}}]\label{prop:backwards_mutation_along_rigid_simple}
Let $\C$ be a $k$-linear stable $\infty$-category whose triangulated homotopy category is $\on{Hom}$-finite and Krull--Schmidt and $\{S_j\}_{1\leq j\leq m}$ a simple-minded collection in $\C$. Fix $i\in[1,m]$ and suppose that the object $S_i$ is rigid. Then the following holds: 
\begin{enumerate}[(1)]
\item The forward tilt of $\{S_j\}_{1\leq j\leq m}$ at $S_i$ is a simple-minded collection in $\C$ and it is given by
\[
\{S_i[1]\}\cup \{\psi^\sharp_{S_i}(S_j)\}_{j\in [1,m],j\not=i}
\]
where 
\[ \psi^\sharp_{S_i}(S_j)\coloneqq \on{fib}(S_j\to \on{Ext}^1(S_j,S_i)^*\otimes S_i[1])\,.\]

\item The backwards tilt of $\{S_j\}_{1\leq j\leq m}$ at $S_i$ is a simple-minded collection in $\C$ and it is given by
\[
\{S_i[-1]\}\cup \{\psi^\flat_{S_i}(S_j)\}_{j\in[1,m],j\not=i}
\]
where 
\[ \psi^\flat_{S_{i}}(S_j)\coloneqq \on{cof}(\on{Ext}^1(S_i,S_j)\otimes S_i[-1] \to S_j)\,.\]\qed
\end{enumerate}
\end{proposition}

\noindent The notion of mutation for silting collections is as follows:

\begin{definition}[{$\!\!$\cite[Def.~2.30]{AI12}}]
Let $\C$ be a $k$-linear stable $\infty$-category whose triangulated homotopy category is $\on{Hom}$-finite and Krull--Schmidt. Let $P_1,\dots,P_m\in \C$ be a silting collection.
    \begin{enumerate}[(1)]
        \item The forward mutation of $P_1,\dots,P_m$ at $P_i$ is the silting collection $$P_1,\dots,P_{i-1},\mu^+(P_i),P_{i+1},\dots,P_m$$
        where
        \[ \mu^+(P_i)=\on{cof}(P_i\to \bigoplus_{j\not= i,j=1}^m\on{Ext}^0_\C(P_i,P_j)^*\otimes P_j)\] 
        is the cofiber of a minimal left $\on{add}(\{P_j,j\not=i\})$-approximation.
        
      \item The backward mutation of $P_1,\dots,P_m$ at $P_i$ is the silting collection $$P_1,\dots,P_{i-1},\mu^-(P_i),P_{i+1},\dots,P_m$$ where
      \[ \mu^-(P_i)=\on{fib}(\bigoplus_{j\not= i,j=1}^m\on{Ext}^0_\C(P_j,P_i)\otimes P_j\to P_i)\] 
      is the fiber of a minimal left $\on{add}(\{P_j,j\not=i\})$-approximation.\qed
    \end{enumerate}
\end{definition}

\noindent The relation between these collections and their mutations above can be summarized as follows:

\begin{theorem}[Koszul-duality {$\!\!$\cite{KN_unpublished,KY14,SY19}}]\label{thm:silting_vs_SMC}
Let $A$ be a connective dg algebra. Assume that the total homology of $A$ is finite dimensional or that $A$ is smooth and that $H_i(A)$ is finite dimensional for all $i\geq 0$. Then there is bijective correspondence (called Koszul duality) between equivalence classes of 
\begin{enumerate}[(1)]
    \item silting collections in $\D^{\on{perf}}(A)$,
    \item simple-minded collections in $\D^{\on{fin}}(A)$, and
    \item bounded $t$-structures on $\D^{\on{fin}}(A)$ with length heart. 
\end{enumerate}
Further, the bijection between (1) and (2) commutes with mutations, and the corresponding simple-minded collections are Koszul dual to the corresponding silting collections.
\end{theorem}

\noindent Note that if $A$ in \Cref{thm:silting_vs_SMC} is finite dimensional and smooth, then $\D^{\on{perf}}(A)=\D^{\on{fin}}(A)$.

\subsection{From exceptional collections to simple-minded and silting collections}\label{subsec:highest_weight_package}

Let us now explain how to construct a simple-minded collection and a silting collection from a given exceptional collection. For that, we will use the following notion of semi-twists.

\begin{definition}[Semi-twists]\label{def:semitwists}
    Let $\C$ be a proper $k$-linear stable $\infty$-category and $X,Y\in \C$. 
    \begin{enumerate}[(1)]
        \item The negative inverse semi-twist of $X$ by $Y$ is defined as the fiber
    \[
    T^{-,\leq 0}_{Y}(X)\coloneqq \on{fib}(X\to \tau_{\leq 0}(\on{Mor}_\C(X,Y)^*)\otimes Y)\,,
    \]
    where $\tau_{\leq 0}(\on{Mor}_\C(X,Y)^*)\subset \on{Mor}_\C(X,Y)^*$ denotes the maximal summand concentrated in homologically negative degrees.
    \item The (strictly) positive inverse semi-twist of $X$ by $Y$ is defined as the fiber 
    \[
    T^{-,>0}_Y(X)\coloneqq \on{fib}(X\to \tau_{>0}(\on{Mor}_\C(X,Y)^*)\otimes Y)\,,
    \]
    where $\tau_{>0}(\on{Mor}_\C(X,Y)^*)\subset \on{Mor}_\C(X,Y)^*$ denotes the maximal summand concentrated in homologically strictly positive degrees.
    \item The positive semi-twist of $X$ by $Y$ is defined as the cofiber 
    \[ T^{+,\geq 0}_Y(X)\coloneqq \on{cof}(\tau_{\geq 0}(\on{Mor}_\C(Y,X))\otimes Y\to X)\,.\]
     \item The (strictly) negative semi-twist of $X$ by $Y$ is defined as the cofiber 
    \[ T^{+,< 0}_Y(X)\coloneqq \on{cof}(\tau_{< 0}(\on{Mor}_\C(Y,X))\otimes Y\to X)\]\qed
    \end{enumerate}
\end{definition}

\begin{remark}
The semi-twist functor $T^{+,\geq 0}_Y(\mhyphen)$ in \cref{def:semitwists}.(3) also appears in \cite[Appendix A]{PYK23}, referred to as the truncated twist functor.\qed
\end{remark}

\noindent Our first result constructs a simple-minded collection from an exceptional collection:

\begin{proposition}[Simple-minded from exceptional]\label{prop:simplesfromstandard}
Let $\C$ be a proper $k$-linear stable $\infty$-category.
\begin{enumerate}[(1)]
 \item Let $\Delta_1,\dots,\Delta_m\in \C$ be a full exceptional collection and define \[ \lsimp_i\coloneqq T_{\Delta_{1}}^{+,\geq 0}(\cdots (T_{\Delta_{i-1}}^{+,\geq 0}(\Delta_i)))\in \C\,.\]
Then the objects $\{\lsimp_i\}_{1\leq i\leq m}$ form a simple-minded collection in $\C$. 
\item Let $\nabla_1,\dots,\nabla_m\in \C$ be a full co-exceptional collection and define \[ \rsimp_i\coloneqq T_{\nabla_{1}}^{-,\leq 0}(\cdots (T_{\nabla_{i-1}}^{-,\leq 0}(\nabla_i)))\in \C\,.\]
Then the objects $\{\rsimp_i\}_{1\leq i\leq m}$ form a simple-minded collection in $\C$. 
\end{enumerate}
\end{proposition}

\noindent Every full exceptional collection determines a dual co-exceptional collection, and thus two associated simple-minded collections $\{\lsimp_i\}$ and $\{\rsimp_i\}$ in \cref{prop:simplesfromstandard}. We show in \Cref{lem:standard_costandard_SMC_coincide} below that these two simple-minded collections coincide.

\begin{proof}[Proof of \Cref{prop:simplesfromstandard}]
Let us prove Part (1). Part (2) follows from the same computation in the opposite $\infty$-category. For Part (1), the objects $\{\lsimp_i\}_{1\leq i\leq m}$ stably generate the exceptional collection, hence they also stably generate $\C$. We will now prove the required homological properties of the objects $\{\lsimp_i\}_{1\leq i\leq m}$, cf.~\cref{def:recollection_exceptional_SMC_silting}.(4), by an induction on the length $m$ of the exceptional collection.\\

For the base case $m=1$, we have $\lsimp_1\simeq \Delta_1$, which by exceptionality satisfies $\on{Mor}_\C(\lsimp_1,\lsimp_1)\simeq k$: the endomorphism object of $\lsimp_1$ is thus coconnective and so $\{\lsimp_1\}$ is a simple-minded collection.\\

\noindent For the induction step, from $m-1$ to $m$, we proceed as follows. For $j<i$, consider the object
\[\lsimp_{j,i}\coloneqq T_{\Delta_{j}}^{+,\geq 0}(\dots(T_{\Delta_{i-1}}^{+,\geq 0}(\Delta_i)))\,.\]
By the induction assumption, $\lsimp_1,\dots, \lsimp_{m-1}$ forms a simple-minded collection in the $k$-linear subcategory of $\C$ generated by $\Delta_1,\dots,\Delta_{m-1}$. Applying the induction assumption to $\Delta_2,\dots,\Delta_{m}$ as well, we obtain that $\lsimp_{2,2},\dots,\lsimp_{2,m}$ forms a simple-minded collection. Let $j\leq m-1$. Since $\on{Mor}_\C(\Delta_i,\Delta_1)\simeq 0$ for all $i>1$, we have
\[
\on{Mor}_\C(\lsimp_{2,m}[-1],\lsimp_{j})\simeq \on{Mor}_\C(\lsimp_{2,m}[-1],\lsimp_{2,j})\,.
\]
The first term in the fiber and cofiber sequence
\begin{equation}\label{eq:3termseq}\on{Mor}_\C(\lsimp_{2,m}[-1],\lsimp_j)\shortrightarrow \on{Mor}_\C(\lsimp_m[-1],\lsimp_{j}) \shortrightarrow 
\on{Mor}_\C(\tau_{\geq 0}(\on{Mor}_\C(\Delta_1,\lsimp_{2,m}))\otimes \Delta_1,\lsimp_j)\end{equation}
is thus coconnective. The third term is also coconnective, as follows from the fact that $\tau_{\geq 0}(\on{Mor}_\C(\Delta_1,\lsimp_{2,m}))$ is connective and the equivalences
\begin{align*} 
\on{Mor}_\C(\Delta_{1},\lsimp_j)& \simeq \on{cof}(\on{Mor}_\C(\Delta_1,\tau_{\geq 0}(\on{Mor}_\C(\Delta_{1},\lsimp_{2,j}))\otimes \Delta_{1})\to \on{Mor}_\C(\Delta_{1},\lsimp_{2,j}))\\
& \simeq \tau_{<0}(\on{Mor}_\C(\Delta_{1},\lsimp_{2,j}))\,.
\end{align*}
The coconnectivity of the outer terms in \eqref{eq:3termseq} implies that the middle term $\on{Mor}_\C(\lsimp_m[-1],\lsimp_{j})$ is coconnective. A similar computation shows that $\on{Mor}_\C(\lsimp_m,\lsimp_m)$ is also coconnective.\\

\noindent It remains to prove that $\on{Mor}_\C(\lsimp_j[-1],\lsimp_m)$ is coconnective as well for $j\in[1,m-1]$. In the case of $j=1$, using that $\lsimp_1\simeq \Delta_1$, we obtain that
\begin{align*}
\on{Mor}_\C(\lsimp_1[-1],\lsimp_m)&\simeq \on{Mor}_\C(\Delta_1[-1],\lsimp_m)\\
& \simeq \on{cof}(\on{Mor}_\C(\Delta_1[-1],\tau_{\geq 0}(\on{Mor}_\C(\Delta_1,\lsimp_{2,m}))\otimes \Delta_1)\to \on{Mor}(\Delta_1[-1],\lsimp_{2,m}))\\
&\simeq \tau_{<0}(\on{Mor}_\C(\Delta_1,\lsimp_{2,m}))[1]
\end{align*}
is indeed coconnective. Note that
\[
\on{Mor}_\C(\lsimp_{2,j}[-1],\lsimp_m)\simeq \on{Mor}_\C(\lsimp_{2,j}[-1],\lsimp_{2,m}) 
\]
is coconnective by induction. For higher $j\in[1,m-1]$, the coconnectivity of $\on{Mor}_\C(\lsimp_j[-1],\lsimp_m)$ follows from the coconnectivity of the two other terms in the following fiber and cofiber sequence:
\begin{align*}
 \on{Mor}_\C(\tau_{\geq 0}(\on{Mor}_\C(\lsimp_{2,j},\Delta_1))\otimes \Delta_1[-1],\lsimp_m)[-1] \to  \on{Mor}_\C(\lsimp_j[-1],\lsimp_m)\to \on{Mor}_\C(\lsimp_{2,j}[-1],\lsimp_m)  \,.
\end{align*}

Finally, we must show that $\lsimp_i$ is indecomposable. By contradiction, let us suppose that $\lsimp_i$ is decomposable and let $j<i$ be maximal such that $\lsimp_{j,i}$ is decomposable. Since $\lsimp_{j+1,i}$ is indecomposable by assumption, and there exists a recollement of the stable subcategory $\langle \Delta_{j},\dots,\Delta_m\rangle$ into $\langle \Delta_j\rangle$ and $\langle \Delta_{j+1},\dots,\Delta_m\rangle$, there must exist a direct sum decomposition
$$Z\oplus Z'\simeq \tau_{\geq 0}\on{Mor}(\Delta_j,\lsimp_{j+1,i})$$
such that the composite morphism
$$Z'\otimes \Delta_j\to \tau_{\geq 0}\on{Mor}(\Delta_j,\lsimp_{j+1,i})\otimes \Delta_j\to \lsimp_{j+1,i}$$
vanishes, thus giving rise to a direct summand $Z'\otimes \Delta_j[1]\subset \lsimp_{j,i}$. However, this is not possible as for any non-zero morphism $k\subset \on{Mor}(\Delta_j,\lsimp_{j+1,i})$, the arising morphism
$$k\otimes \Delta_i \to \tau_{\geq 0}\on{Mor}(\Delta_j,\lsimp_{j+1,i})\otimes \Delta_j\to \lsimp_{j+1,i}$$
is also non-zero.
\end{proof}

\noindent The analog of \cref{prop:simplesfromstandard} for silting collections is as follows.

\begin{proposition}[Silting from exceptional]\label{prop:silting_from_standards}
Let $\C$ be a proper $k$-linear stable $\infty$-category.
\begin{enumerate}[(1)]
\item Let $\Delta_1,\dots,\Delta_m$ be a full exceptional collection in $\C$. For any $i\in[1,m]$, consider the object
\[ \lsilt_i\coloneqq T_{\Delta_{m}}^{-,>0}(\cdots (T_{\Delta_{i+1}}^{-,>0}(\Delta_i)))\in \C\,.\]
Then the objects $\{\lsilt_i\}_{1\leq i\leq m}$ form a silting collection in $\C$.
\item Let $\nabla_1,\dots,\nabla_m$ be a full co-exceptional collection in $\C$. For any $i\in[1,m]$, consider the object
\[ \rsilt_i\coloneqq T_{\nabla_{m}}^{+,<0}(\cdots (T_{\nabla_{i+1}}^{+,<0}(\nabla_i)))\in \C\,.\]
Then the objects $\{\rsilt_i\}_{1\leq i\leq m}$ form a silting collection in $\C$.
\end{enumerate}
\end{proposition}
\begin{proof}
The proof follows a similar strategy as the proof of \Cref{prop:simplesfromstandard}. As in there, we only prove Part (1), since Part (2) follows from the same argument applied to the opposite $\infty$-category. For part (1), the objects $\{\lsilt_i\}_{1\leq i\leq m}$ stably generate the exceptional collection, hence they also stably generate $\C$. We will now show that $\on{Mor}(\lsilt_i,\lsilt_j)$ is connective by an induction on the length $m$ of the exceptional collection.\\

For the base case $m=1$, $\lsilt_1=\Delta_1$ is exceptional and hence a silting object. For the induction step, let us proceed as follows. For $j>i$, we denote 
\[ \lsilt_{j,i}\coloneqq  T_{\Delta_{j}}^{-,>0}(\cdots (T_{\Delta_{i+1}}^{-,>0}(\Delta_i)))\]
where $\lsilt_{i,i}=\Delta_i$. The objects $\lsilt_2,\dots,\lsilt_{m}$ form a silting collection in the stable subcategory generated by them, as do the objects $\lsilt_{m-1,1},\dots,\lsilt_{m-1,m-1}$. We must thus show that $\on{Mor}_\C(\lsilt_1,\lsilt_j)$ and $\on{Mor}_\C(\lsilt_j,\lsilt_1)$ are connective for all $j\in[1,m]$.\\

First, we observe that
\begin{align*} 
\on{Mor}_\C(\lsilt_1,\lsilt_m)&\simeq \on{Mor}_\C(\lsilt_1,\Delta_m)\\
&\simeq \on{cof}(\on{Mor}_\C(\tau_{>0}(\on{Mor}_\C(\lsilt_{m-1,1},\Delta_m)^*)\otimes \Delta_m,\Delta_m)\to \on{Mor}_\C(\lsilt_{m-1,1},\Delta_m))\\
& \simeq \tau_{\geq 0}\on{Mor}(\lsilt_{m-1,1},\Delta_m)
\end{align*}
is connective. Similarly, for $1\leq j<m$, the object 
\[
\on{Mor}_\C(\lsilt_j,\lsilt_{m-1,1})\simeq \on{Mor}_\C(\lsilt_{m-1,j},\lsilt_{m-1,1})
\]
is connective. For $j=m$, we have $\on{Mor}_\C(\lsilt_m,\lsilt_{m-1,1})\simeq 0$. Therefore, for any $j\in[1,m]$, the middle term in the exact sequence
\[  
\on{Mor}_\C(\lsilt_j,\tau_{>0}(\on{Mor}_\C(\lsilt_{m-1,1},\Delta_m)^*)[-1]\otimes \Delta_m)\to \on{Mor}_\C(\lsilt_j,\lsilt_1)\to \on{Mor}_\C(\lsilt_j,\lsilt_{m-1,1})
\]
must be connective, because the outer terms are as well. Thus $\on{Mor}_\C(\lsilt_j,\lsilt_1)$ is connective for all $j\in[1,m]$. In order to show that $\on{Mor}_\C(\lsilt_1,\lsilt_j)$  is connective, we proceed similarly. For all $j\in[2,m-1]$, we have the equivalence
\[
\on{Mor}_\C(\lsilt_1,\lsilt_{m-1,j})\simeq \on{Mor}_\C(\lsilt_{m-1,1},\lsilt_{m-1,j})\,.
\]
In the fiber and cofiber sequence 
\[
\on{Mor}_\C(\lsilt_1,\tau_{>0}(\on{Mor}_\C(\lsilt_{m-1,j},\Delta_m)^*)\otimes \Delta_m)[-1]\to \on{Mor}_\C(\lsilt_1,\lsilt_j)\to \on{Mor}_\C(\lsilt_1,\lsilt_{m-1,j})
\]
the outer terms are thus connective, and hence so is the middle term $\on{Mor}_\C(\lsilt_1,\lsilt_j)$ for any $j\in[1,m]$. This conclude the argument for connectivity. The indecomposablity of $\lsilt_{i}$ follows from a similar argument as the proof of the indecomposability of $\lsimp_i$ in \Cref{prop:simplesfromstandard}.
\end{proof}

\begin{remark}
Note that, in contrast to \cref{prop:simplesfromstandard}, cf.~\cref{lem:standard_costandard_SMC_coincide}, the two silting collections in \cref{prop:silting_from_standards} are typically different. This is a stable $\infty$-category analogue of the fact that the collections of projective and injective objects in the setting of abelian categories are typically different.\qed
\end{remark}

\noindent Let us now show in \cref{thm:Koszul_duality_for_silting_and_SMC} that the simple-minded collection in \cref{prop:simplesfromstandard}, and the silting collection in \cref{prop:silting_from_standards}, each associated with the same full exceptional collection, are Koszul dual to each other. For that, we compare the associated $t$-structures. The argument for \cref{lem:standard_costandard_SMC_coincide} will similarly use $t$-structures.\\

\noindent Note that a $t$-structure on a stable $\infty$-category is by definition a $t$-structure on the triangulated homotopy category, see for instance \cite[Def.~1.2.1.4]{HA}. All $t$-structures will be homological, meaning that the positive aisle is stable under suspension. First, we recall how to glue $t$-structure via recollement:

\begin{proposition}[{$\!\!$\cite[Thm.~1.4.10]{BBD82}}]\label{thm:gluing_t_structure_along_recollement}
Let $\C$ be a stable $\infty$-category with a recollement
\[
\begin{tikzcd}
\C_1\arrow[r, "i_*" description]  & \C  \arrow[l, "i^*"', shift right=3] \arrow[l, "i^!", shift left=3]\arrow[r, "j^*" description] & \C_2 \arrow[l, "j_*", shift left=3] \arrow[l, "j_!"', shift right=3]
\end{tikzcd}
\]
into two stable subcategories $\C_1,\C_2$ with $i^*\dashv i_*\dashv i^!$ and $j_!\dashv j^*\dashv j_*$. Suppose that $\C_1,\C_2$ are equipped with $t$-structures. Then there is a $t$-structure on $\C$ with aisles
\[
\C_{\geq 0}=\{X\in \C|j^*(X)\in (\C_2)_{\geq 0}, i^*(X)\in (\C_1)_{\geq 0}\}
\]
and
\[
\C_{\leq 0}=\{X\in \C| j^*(X)\in (\C_2)_{\leq 0}, i^!(X)\in (\C_1)_{\leq 0}\}\,.
\]
\qed
\end{proposition}

\noindent We apply \cref{thm:gluing_t_structure_along_recollement} to construct a $t$-structure $t(\Delta)$ on $\C$ from an exceptional collection $\Delta=(\Delta_1,\ldots,\Delta_m)$ in $\C$, as follows.\\

\begin{construction}\label{constr:t_structure_from_exceptional_collection}
Let $\C$ be a $k$-linear stable $\infty$-category with a full exceptional collection $\Delta_1,\dots,\Delta_m$. For any $1\leq i\leq j\leq m$, we denote by $\langle\Delta_i,\dots,\Delta_j\rangle$ the $k$-linear stable subcategory of $\C$ generated by $\Delta_i,\dots,\Delta_j$. Consider the iterated recollements
\[
\begin{tikzcd}
{\langle \Delta_1,\dots,\Delta_i\rangle }\arrow[r, "i_*" description]  & {\langle \Delta_1,\dots,\Delta_{i+1}\rangle } \arrow[l, "i^*"', shift right=3] \arrow[l, "i^!", shift left=3]\arrow[r, "j^*" description] & {\langle \Delta_{i+1}\rangle \simeq \D^{\on{perf}}(k)}\arrow[l, "j_*", shift left=3] \arrow[l, "j_!"', shift right=3]
\end{tikzcd}\,,
\]
where $i^*$ and $j_!$ are the fully faithful inclusions commuting with the inclusions into $\C$. By considering the standard $t$-structure on $\D^{\on{perf}}(k)\simeq \langle \Delta_{i+1}\rangle$, and iteratively gluing $t$-structures via \Cref{thm:gluing_t_structure_along_recollement}, we obtain a $t$-structure on $\C$ from the full exceptional collection $\Delta_1,\dots,\Delta_m$.\qed
\end{construction}

\begin{remark}
The $t$-structure from \Cref{constr:t_structure_from_exceptional_collection} coincides with the $t$-structure obtained by iteratively gluing $t$-structures in the opposite direction, i.e.~along the following recollements: 
\[
\begin{tikzcd}
{\D^{\on{perf}}(k)\simeq \langle \Delta_i\rangle }\arrow[r, "i_*" description]  & {\langle \Delta_i,\dots,\Delta_{m}\rangle } \arrow[l, "i^*"', shift right=3] \arrow[l, "i^!", shift left=3]\arrow[r, "j^*" description] & {\langle \Delta_{i+1},\dots,\Delta_m\rangle}\arrow[l, "j_*", shift left=3] \arrow[l, "j_!"', shift right=3]
\end{tikzcd}\,.
\]
\qed
\end{remark}

Let us now prove that the objects $\lsimp_i\in \C$ forming the simple-minded collection associated to a full exceptional collection $\Delta=(\Delta_1,\dots,\Delta_m)$, cf.~\Cref{prop:simplesfromstandard}, are indeed in the heart $\mathcal{H}(\Delta)\sse\C$ of the $t$-structure $t(\Delta)$ from \Cref{constr:t_structure_from_exceptional_collection}.

\begin{lemma}\label{lemma:SMC_are_the_simples_in_the_heart}
Let $\C$ be a proper $k$-linear stable $\infty$-category with a full exceptional collection $\Delta\coloneqq (\Delta_1,\dots,\Delta_m)$. Consider the associated heart $\mathcal{H}(\Delta)\sse\C$ and the associated simple-minded collection $\lsimp_1,\ldots,\lsimp_m$. Then:
\begin{enumerate}
    \item Each object $\lsimp_i\in \C$ lies in the heart $\mathcal{H}(\Delta)$.
    \item Each $\lsimp_i\in \mathcal{H}(\Delta)$ is a simple object.
\end{enumerate}
\end{lemma}

\begin{proof}
For Part (1), let us prove by a descending induction on $j\in[1,m]$ that, for all $i\in[j,m]$, the object
    \[\lsimp_{j,i}\coloneqq T_{\Delta_j}^{+,\geq 0}(\cdots(T_{\Delta_{i-1}}^{+,\geq 0}(\Delta_i))\] 
    lies in the corresponding heart $\langle \Delta_j,\dots,\Delta_m\rangle^{\heartsuit}\coloneqq \mathcal{H}(\Delta_j,\dots,\Delta_m)$. Indeed, for $j=i$, we have $\lsimp_{i,i}=\Delta_i=i_*(\Delta_i)$, and thus
    $$i^!(\lsimp_{i,i})\simeq i^*(\lsimp_{i,i})\simeq \Delta_i\in \langle \Delta_i\rangle^{\heartsuit}$$
    and $j^*(\lsimp_{i,i})\simeq 0$. By using that the morphism $\lsimp_{j+1,i}\to \lsimp_{j,i}$ describes the counit map of the adjunction $j_!\dashv j^*$ evaluated at $\lsimp_{j,i}$, we further obtain the following three inclusions:
    
    \begin{align*} j^*(\lsimp_{j,i})& \simeq \lsimp_{j+1,i}\in \langle \Delta_{j+1},\dots,\Delta_m\rangle^{\heartsuit},\\
    i^*(\lsimp_{j,i})&\simeq \tau_{\geq 0}\on{Mor}(\Delta_j,\lsimp_{j+1,i})\otimes \Delta_j[1]\in \langle \Delta_j\rangle_{\geq 0},\\
    i^!(\lsimp_{j,i})&\simeq \on{cof}(\tau_{\geq 0}\on{Mor}(\Delta_j,\lsimp_{j+1,i})\otimes \Delta_j\to \on{Mor}(\Delta_j,\lsimp_{j+1,i})\otimes \Delta_j)\\
    & \simeq \tau_{<0}\on{Mor}(\Delta_j,\lsimp_{j+1,i})\otimes \Delta_j\in \langle \Delta_j\rangle_{\leq 0}\,.
    \end{align*}
\noindent By the inductive hypothesis on $\langle \Delta_{j+1},\dots,\Delta_m\rangle^{\heartsuit}$ and the equivalences above, Part (1) thus follows.\\
    
For Part (2), we proceed by a descending induction on $j$ to show that $\lsimp_{j,i}\in \langle \Delta_j,\dots,\Delta_m\rangle^\heartsuit$ is simple. For that, given an epimorphism $\lsimp_{j,i}\to X$, we must show that it is either zero or an isomorphism. Since $j^*$ restricts to an exact functor between the hearts, $\lsimp_{j+1,i}\to j^*(X)$ is also an epimorphism. Note that $\lsimp_{j+1,i}$ is simple by the induction hypothesis and thus this later map must be either zero or an isomorphism. We treat these two cases:\\

\noindent    {\bf Case 1:} Suppose that $j^*(X)\simeq \lsimp_{j+1,i}$. Then the kernel $K$ of $\lsimp_{j,i}\to X$ lies in the image of $i_*$. Thus $i^*(K)=i^!(K)\in \langle \Delta_j\rangle^{\heartsuit}$. Since the functor $i^!$ is left exact when restricted to the hearts, we deduce that $i^!(K)\subset H_0(i^!(\lsimp_{j,i}))\simeq 0$ vanishes. Thus we obtain that the epimorphism $\lsimp_{j,i}\lr X$ is an isomorphism.\\
    
\noindent    {\bf Case 2:} Suppose that $j^*(X)\simeq0$. Since $i^*$ is right exact when restricted to the hearts, there must exist an epimorphism from the degree $0$ part of $i^*(\lsimp_{j,i})\simeq \tau_{\geq 0}\on{Mor}(\Delta_j,\lsimp_{j+1,i})\otimes \Delta_j[1]$ to the degree $0$ part of $i^*X$, and thus $\tau_{\leq 0}(i^*X)\simeq 0$. Since $i^!(X)\simeq i^*(X) \in \langle \Delta_j\rangle^\heartsuit$, we obtain $X=0$.\\

\noindent Therefore, assuming $\lsimp_{j+1,i}$ is simple, we have proven that any epimorphism $\lsimp_{j,i}\to X$ must be zero or an isomorphism, so that $\lsimp_{j,i}$ is simple as well. Hence, by induction, $\lsimp_{j,i}$ is simple for all $j\in[1,i]$. Since $\lsimp_i\simeq \lsimp_{1,i}$, this proves Part (2).
\end{proof}

\noindent \cref{constr:t_structure_from_exceptional_collection} provides a $t$-structure $t(\Delta)$ from a full exceptional collection $\Delta\coloneqq (\Delta_1,\dots,\Delta_m)$. In line with this, a silting object also gives rise to a $t$-structure:

\begin{lemma}[{$\!\!$\cite[Thm.~A.1]{BY14}}] \label{lem:t_structure_from_silting_object}
Let $\C$ be a $k$-linear stable $\infty$-category and $P\in \C$ a silting object. Then there is a $t$-structure $t(P)$ on $\C$ with aisles
    \[
    \C_{\geq 0}=\{X\in \C|\on{Ext}^i(P,X)\simeq 0\text{ for all }i\geq 1\}\,,
    \]
    \[
    \C_{\leq 0}=\{X\in \C|\on{Ext}^i(P,X)\simeq 0\text{ for all }i\leq -1\}\,.
    \]
\end{lemma}

\noindent If the silting object $P$ in \cref{lem:t_structure_from_silting_object} is associated to a full exceptional collection $\Delta\coloneqq (\Delta_1,\dots,\Delta_m)$, via \cref{prop:silting_from_standards}, the $t$-structure $t(P)$ in \cref{lem:t_structure_from_silting_object} coincides with the $t$-structure $t(\Delta)$ from \cref{constr:t_structure_from_exceptional_collection}, as we now prove.

\begin{lemma}\label{lemma:silting_t_structure_agrees_with_glued_t_structure}
Let $\C$ be a proper $k$-linear stable $\infty$-category with a full exceptional collection $\Delta=(\Delta_1,\dots,\Delta_m)$, and consider the silting object $P(\Delta)\coloneqq \bigoplus_{i=1}^m\lsilt_{i}$ associated to $\Delta_1,\dots,\Delta_m$. Then the $t$-structure $t(P(\Delta))$ coincides with the $t$-structure $t(\Delta)$.
\end{lemma}

\begin{proof}
Let us denote $\lsilt_{j,i}\coloneqq T^{-,>0}_{\Delta_j}(\cdots(T_{\Delta_{i+1}}^{-,>0}(\Delta_i)))$ for $i\leq j\in[1,m]$. Note that $\lsilt_{j,1},\dots,\lsilt_{j,j}\in \langle \Delta_1,\dots,\Delta_j\rangle$ forms a silting collection for all $j\in[1,m]$. Let us prove by induction on $m$ that the $t$-structrue induced by this silting collection coincides with the glued $t$-structure. For the base case $m=1$, $Y_{1,1}=\Delta_1$ and thus the assertion holds.\\ 

\noindent For the induction step, suppose that the statement holds for $m-1$. That is, we assume that the two $t$-structures on $\langle \Delta_1,\dots,\Delta_{m-1}\rangle$ coincide. Now, for an object $X\in \C$, we have $X\in \C_{\leq 0}$ if $j^*(X)\in \langle \Delta_m\rangle_{\leq 0}$ and $i^!(X)\in \langle \Delta_1,\dots,\Delta_{m-1}\rangle_{\leq 0}$. Since $\Delta_m=\lsilt_m$, the functor $j^*$ identifies $\lsilt_m$ with $\on{Mor}_\C(\lsilt_{m},\mhyphen)$ under the equivalence $\langle \Delta_m\rangle \simeq \D^{\on{perf}}(k)$. Now, for any $i\in[1,m-1]$, consider the fiber and cofiber sequence
    \[
    \on{Mor}_\C(\lsilt_{m-1,i},X)\to \on{Mor}_\C(\lsilt_{i},X)\to \on{Mor}(\tau_{>0}(\on{Mor}_\C(\lsilt_{m-1,i},\Delta_m)^*)\otimes \Delta_m,X)[1].
    \]
    Suppose that $X$ lies in the negative aisle of either of the two $t$-structures. Then the right term is coconnective, and the middle term is thus coconnective if and only if the left term is coconnective.
    Since $\on{Mor}_\C(\lsilt_{m-1,i},X)\in \C$ lies in the image of $i_*$, we obtain the equivalence 
    \[ \on{Mor}_\C(\lsilt_{m-1,i},X)\simeq \on{Mor}_{\langle \Delta_1,\dots,\Delta_{m-1}\rangle}(\lsilt_{m-1,i},i^!(X))\,.\] 
    Thus, $i^!(X)\in   \langle \Delta_1,\dots,\Delta_{m-1}\rangle_{\leq 0}$ if and only if $\on{Mor}_\C(\lsilt_{m-1,i},X)$ is coconnective, showing that the negative isles of the two $t$-structures on $\langle \Delta_1,\dots,\Delta_{m-1},\Delta_m\rangle$ coincide, and hence $t(\Delta)=t(P(\Delta))$. 
\end{proof}

\begin{theorem}\label{thm:Koszul_duality_for_silting_and_SMC}
Let $\C$ be a proper $k$-linear stable $\infty$-category. 
\begin{enumerate}[(1)]
    \item Let $\{\Delta_i\}_{1\leq i\leq m}$ be a full exceptional collection in $\C$, $\{\lsimp_i\}_{1\leq i\leq m}$ its associated simple-minded collection, and $\{\lsilt_{i}\}_{1\leq i\leq m}$ its associated silting collection. Then
    \[
    \on{Mor}_\C(\lsilt_i,\lsimp_i)\simeq \begin{cases} k& i=j, \\ 0& \text{else}\,,\end{cases}\qquad \forall i,j\in[1,m].
    \]
    That is, the silting collection in $\C$ associated to a full exceptional collection is Koszul dual to the simple-minded collection associated to that same exceptional collection.
    
    \item Let $\{\nabla_i\}_{1\leq i\leq m}$ be a full co-exceptional collection in $\C$, $\{\rsimp_i\}_{1\leq i\leq m}$ its associated simple-minded collection, and $\{\rsilt_{i}\}_{1\leq i\leq m}$ its associated silting collection. Then
    \[
    \on{Mor}_\C(\rsimp_i,\rsilt_j)\simeq \begin{cases} k& i=j, \\ 0& \text{else}\,,\end{cases}\qquad \forall i,j\in[1,m].
    \]
    That is, the silting collection in $\C^{\on{op}}$ associated to a full co-exceptional collection is Koszul dual to the simple-minded collection associated to that same co-exceptional collection.
\end{enumerate}
\end{theorem}

\begin{proof}
Via the bijection from \Cref{thm:silting_vs_SMC}, silting collections with $m$ indecomposable objects correspond to simple-minded collections with $m$ indecomposable objects, with the latter given by the simple objects in the heart arising from the former. Part (1) thus follows from combining \Cref{lemma:SMC_are_the_simples_in_the_heart} and \Cref{lemma:silting_t_structure_agrees_with_glued_t_structure}. Part (2) follows from the same argument applied to $\C^{\on{op}}$. 
\end{proof}

\noindent Let us finally conclude that the two simple-minded collections from \Cref{prop:simplesfromstandard}, built either from a full exceptional collection or its dual co-exceptional collection, coincide.

\begin{lemma}\label{lem:standard_costandard_SMC_coincide}
Let $\C$ be a $k$-linear stable $\infty$-category with a full exceptional collection $\Delta_1,\dots,\Delta_m$ and dual full co-exceptional collection $\nabla_1,\dots,\nabla_m$. Then the two simple-minded collections from \Cref{prop:simplesfromstandard} coincide, i.e.~ $\lsimp_i\simeq \rsimp_i$ for all $i\in [1,m]$.
\end{lemma}

\begin{proof}
    By \Cref{lemma:SMC_are_the_simples_in_the_heart}, it suffices to show that the $t$-structure $t(\Delta)$ arising from $\Delta=(\Delta_1,\dots,\Delta_m)$ of $\C$ agrees with the $t$-structure $t(\nabla)$ on $\C$ arising from $\nabla=(\nabla_1,\dots,\nabla_m)\in \C^{\on{op}}$ after passing to the opposite $\infty$-category.

\noindent Let us prove by this by induction on $m$. For the base case $m=1$, $\Delta_1=\nabla_1$, and thus the assertion holds. For the induction step, from $m-1$ to $m$, we proceed as follows, The $t$-structure $t(\nabla)$ on $\C$ arises after passing to opposite categories from gluing $t$-structures along the recollement
\[
\begin{tikzcd}
\langle \nabla_1,\dots,\nabla_{m-1}\rangle^{\on{op}} \arrow[r, "\hat{i}_*" description] & {\langle \nabla_1,\dots,\nabla_{m}\rangle^{\on{op}}} \arrow[l, "\hat{i}^!"', shift right=3] \arrow[l, "\hat{i}^*", shift left=3] \arrow[r, "\hat{j}^*" description] & {\langle \nabla_{m}\rangle^{\on{op}}} \arrow[l, "\hat{j}_*", shift left=3] \arrow[l, "\hat{j}_!"', shift right=3]
\end{tikzcd}
\]
where $\hat{i}_*$ and $\hat{j}_!$ are the fully faithful inclusions commuting with the given inclusions into $\C^{\on{op}}$. \\

Now, the functor $(\hat{i}_*)^{\on{op}}$ identifies with the functor $i_*\colon \langle \Delta_{1},\dots, \Delta_{m-1}\rangle\subset \C$, both describing the canonical inclusions of $\langle \Delta_{1},\dots, \Delta_{m-1}\rangle= \langle \nabla_{1},\dots,\nabla_{m-1}\rangle$, which coincide by the inductive hypothesis. Since passing to opposite categories exchanges left and right adjoints, we thus also have $(\hat{i}^!)^{\on{op}}\simeq i^*$.
 Furthermore, under the equivalence $\langle\nabla_m\rangle\simeq \D^{\on{perf}}(k)\simeq \langle \Delta_m\rangle$, the functor $(\hat{j}^*)^{\on{op}}$ identifies with $j^*$, since it is given by the Verdier quotient functor by $i_*$. 

Since passing to opposite categories exchanges the left and right aisles, an object $X\in \langle \nabla_1,\dots,\nabla_{m}\rangle$ lies in the positive aisle if and only if 
\[ i^!(X)\simeq (\hat{i}^!)^{\on{op}}(X)\in (\langle \nabla_1,\dots,\nabla_{m-1}\rangle^{\on{op}}_{\leq 0})^{\on{op}}= \langle \Delta_1,\dots,\Delta_{m-1}\rangle_{\geq 0}\] and 
\[ j^*(X)\simeq (\hat{j}^*)^{\on{op}}(X)\in (\langle \nabla_m\rangle^{\on{op}}_{\leq 0})^{\on{op}}=\langle \Delta_m\rangle_{\geq 0}\,.\]
Thus, the two $t$-structures coincide. 
\end{proof}

\subsection{Relation with highest weight theory}\label{subsec:relation_w_highest_weight_theory}

Highest weight categories were introduced in \cite{CPS88}. The following definition uses a slightly different formulation given in \cite{Kra17}, which itself refers to \cite{Rou08}. We will further only consider highest weight categories whose poset is a totally ordered set on $m$ elements. The full highest weight theory for more general, possibly infinite, posets can be found in \cite{BS24}.

\begin{definition}[Highest weight structure]\label{def:highest_weight_category}
    Let $A$ be a $k$-linear abelian category having only finitely many isoclasses of simple objects and finite dimensional Homs. A highest weight structure on $A$ consists of a collection of objects $\Delta_1,\dots,\Delta_m\in A$, satisfying the following conditions:
    \begin{itemize}
        \item $\on{Hom}_A(\Delta_i,\Delta_i)\simeq  k$ for all $1\leq i\leq m$. 
        \item $\on{Hom}_A(\Delta_i,\Delta_j)=0$ for all $1\leq j<i\leq m$. 
        \item For all $1\leq i\leq m$, there exists a projective cover $P_i\to \Delta_i$ whose kernel has a filtration whose associated graded consists of objects $\{\Delta_j\}_{i<j\leq m}$.
        \item The object $\bigoplus_{i=1}^m P_i$ is a projective generator of $A$, meaning that an object $X\in A$ vanishes if and only if $\on{Hom}_A(\bigoplus_{i=1}^m P_i,X)=0$.\qed
    \end{itemize}
\end{definition}

\noindent If an abelian category $A$ has a highest weight structure, then the images of the objects $\Delta_1,\dots,\Delta_m$ under the inclusion $N(A)\subset \D^b(A)$, with $N(A)$ the nerve of $A$, form a full exceptional collection, see for instance \cite[Thm.~5.2]{Kra17}. These objects are called the standard objects in $A$. Furthermore, we can label the iso-classes of simple objects in $A$ by $S_1,\dots,S_m$, and for all $i\in[1,m]$ there exists an epimorphism $\Delta_i\to S_i$ in $A$ whose kernel has composition factors in $\{S_i\}_{1\leq j<i}$.

\begin{remark}
A highest weight structure, as in \cref{def:highest_weight_category}, can equivalently be described as a quasi-hereditary structure on the endomorphism algebra of the projective generator, cf.~\cite{CPS88}.\qed
\end{remark}

The results of \Cref{subsec:highest_weight_package} in the stable setting are inspired by the notion of a highest weight abelian category. The full exceptional collection acquires the role of the standard objects, the associated simple-minded collection has the role of the simple objects, and the associated silting collection has the role of the projective objects. In the following results, we elucidate the precise relationship between the constructions from \Cref{subsec:highest_weight_package} and the theory of highest weight categories. 

\begin{lemma}\label{lem:silting_in_highest_weight_category}
Let $A$ be a $k$-linear abelian category equipped with a highest weight structure with standard objects $\Delta_1,\dots,\Delta_m\in A$. Then the silting collection $\{\lsilt_{i}\}_{1\leq i\leq m}$ in $\D^{\on{b}}(A)$ from \Cref{prop:silting_from_standards} associated with the full exceptional collection $\Delta_1,\dots,\Delta_m\in \D^{\on{b}}(A)$ describes the image of the collection of projective $A$-modules $\{P_i\}_{1\leq i\leq m}$ under the inclusion $N(A)\subset \D^{\on{b}}(A)$. 
\end{lemma}

\begin{proof}
Let us denote $\lsilt_{j,i}\coloneqq T^{-,>0}_{\Delta_j}(\cdots(T_{\Delta_{i+1}}^{-,>0}(\Delta_i)))$ for $1\leq i\leq j\leq m$. For $1\leq j\leq m$, let $A_{\leq j}\subset A$ be the abelian subcategory generated by $\Delta_1,\dots,\Delta_j$. Let us prove by induction on $j$ that $\{\lsilt_{j,i}\}_{1\leq i\leq j}$ describes the images of the projectives in $A_{\leq j}$ in $\D^b(A_{\leq j})\subset \D^b(A)$. The base case $j=1$ readily holds.\\

\noindent For the inductive step, we proceed as follows. First, the projective objects in a highest weight abelian category arise via iterated maximal extensions between the standard objects. Indeed, this follows for instance from the construction of the tilting modules described in the proof of \cite[Proposition 3.1]{Soe99} applied to the Ringel dual highest weight category. Alternatively, a direct construction can also be found in \cite[Ex.~8.1.10]{Kra22}, where these maximal extensions are called universal extensions. Now, since $A_{\leq j+1}$ is highest weight, the projective objects in $A_{\leq j+1}$ arise from the projective objects in $A_{\leq j}$ by universal extension with $\on{add}(\Delta_{j+1})$, as well as the object $\Delta_{j+1}$ itself. Since $\lsilt_{j,i}$ is projective in $A_{\leq j}$, the Statement 1 in \cite[Appendix]{DR89} implies that
\[\on{Ext}^{l}(\lsilt_{j,i},\Delta_{j+1})\simeq 0,\qquad l\in[2,m]\,.\]
Thus, $\tau_{>0}\on{Mor}(\lsilt_{j,i},\Delta_{j+1})^*$ is concentrated in degree $1$. This shows that $\lsilt_{j+1,i}=T^{-,>0}_{\Delta_{j+1}}(\lsilt_{j,i})$ corresponds to the passage to the universal extension, thus describing a projective object in $A_{\leq j+1}$.
\end{proof}

\begin{lemma}\label{lem:SMC_in_highest_weight_category}
Let $A$ be a $k$-linear abelian category equipped with a highest weight structure, and standard objects $\Delta_1,\dots,\Delta_m\in A$. Then the simple-minded collection $\{\lsimp_j\}_{1\leq j\leq m} $in $\D^{\on{b}}(A)$ from \Cref{prop:simplesfromstandard} associated with the full exceptional collection $\Delta_1,\dots,\Delta_m$ coincides with the image of the collection of simple modules $\{S_j\}_{1\leq j\leq m}$ in $A$ under the inclusion $N(A)\subset \D^{\on{b}}(A)$. 
\end{lemma}

\begin{proof}
This follows from combining \Cref{thm:silting_vs_SMC}, \Cref{thm:Koszul_duality_for_silting_and_SMC} and \Cref{lem:silting_in_highest_weight_category}. 
\end{proof}

Finally, we give a characterization of highest weight structures on hearts arising from full exceptional collections as follows. Given a proper $k$-linear stable $\infty$-category with a full exceptional collection $\Delta_1,\dots,\Delta_m$, we denote as above
$$\lsilt_{j,i}\coloneqq T^{-,>0}_{\Delta_j}(\cdots(T_{\Delta_{i+1}}^{-,>0}(\Delta_i))),\qquad 1\leq i\leq j\leq m.$$
Note that $\lsilt_{j,1},\dots,\lsilt_{j,j}$ defines a silting collection in $\langle \Delta_1,\dots,\Delta_j\rangle$ for all $j\in[1,m]$. Recall that we denote by $t(\Delta)$ the $t$-structure associated to a full exceptional collection $\Delta=(\Delta_1,\dots,\Delta_m)$ via \Cref{constr:t_structure_from_exceptional_collection}.

\begin{theorem}\label{thm:relation_with_highest_weight_abelian_cats}
Let $\C$ be a proper $k$-linear stable $\infty$-category with a full exceptional collection $\Delta=(\Delta_1,\dots,\Delta_m)$. The following are equivalent:
    \begin{enumerate}[$(1)$]
        \item  For all $j\in[1,m]$, the homology of the derived endomorphism algebra of the silting object $\bigoplus_{i=1}^j \lsilt_{j,i}$ in $\langle \Delta_1,\dots,\Delta_j\rangle$ is concentrated in degree $0$.
        \item For all $j\in[1,m]$, the silting object $\bigoplus_{i=1}^j \lsilt_{j,i}$ in $\langle \Delta_1,\dots,\Delta_j\rangle$ lies in the heart of the $t$-structure $t(\Delta_1,\ldots,\Delta_j)$.
        \item There exists a $k$-linear abelian category $A$ with a highest weight structure with standard objects $\Delta_1^A,\dots,\Delta_m^A$ and a fully faithful functor $N(A)\subset \C$, which induces an equivalence of stable $\infty$-categories $F\colon \D^b(A)\to \C$ satisfying that $F(\Delta_i^A)\simeq \Delta_i$. 
    \end{enumerate}
\end{theorem}

\begin{proof}
If the derived endomorphism algebra $R$ of a silting object is concentrated in degree $0$, then the heart of the standard $t$-structure on $\D^{\on{perf}}(R)$ contains $R$, as $R$ is then a dg module over itself in degree $0$. This proves (1)$\Rightarrow$(2).\\

For (2)$\Rightarrow$(3) we proceed as follows. Let $A$ be the heart of the $t$-structure from Lemma \ref{lem:t_structure_from_silting_object} and $A_{\leq j}$ be the heart of $\langle \Delta_1,\dots,\Delta_j\rangle$, $j\in[1,m]$. Note that $A_{\leq j}\subset A$. To establish (3), assuming (2), we verify the conditions from Definition \ref{def:highest_weight_category}. First, the silting object $\bigoplus_{i=1}^m \lsilt_{i}$ is a projective generator of $A$, which establishes the fourth condition in Definition \ref{def:highest_weight_category}. Since $\Delta_j=\lsilt_{j,j}\in A_{\leq j}$, we obtain $\Delta_j\in A$. Further, since these objects form a full exceptional collection in $\C$, we conclude that 
\[ \on{Hom}_A(\Delta_i,\Delta_j)\simeq \begin{cases} k & i=j\\ 0 & j<i \end{cases}\,.\]

\noindent This establishes the first and second conditions in Definition \ref{def:highest_weight_category}. To verify the third condition, it remains to show that the morphism $\lsilt_{i}\to \Delta_i$ is an epimorphism in $A$ whose kernel has a filtration whose associated graded consists of standard objects. Since $\lsilt_{j+1,i}\in A_{\leq j+1}$, we deduce that $j^*(\lsilt_{j+1,i})\in \langle \Delta_{j+1}\rangle^\heartsuit.$ Thus, for all $1\leq i\leq j< m$, we have a fiber and cofiber sequence in $\C$
\[
\on{Ext}^1_\C(\lsilt_{j,i},\Delta_{j+1})^*\otimes \Delta_{j+1}\to \lsilt_{j+1,i}\to \lsilt_{j,i}
\]
where $\on{Ext}^1_\C(\lsilt_{j,i},\Delta_{j+1})^*\otimes \Delta_{j+1}\simeq j_!j^*(\lsilt_{j+1,i})$. Since the fiber and cofiber sequence lies in the image of $N(A)\subset \D^b(A)$, it arises from an exact sequence in $A$. This shows that $\lsilt_{j+1,i}\to \lsilt_{j,i}$ is an epimorphism. The morphism $\lsilt_{i}\to \Delta_i$ is thus an epimorphism, as it is the composite of epimorphisms, and the kernel has the desired filtration. This proves (2)$\Rightarrow$(3).\\ 

Finally, for (3)$\Rightarrow$(1), we suppose that (3) is satisfied. Then $A_{\leq j}\subset A$ inherits a highest weight structure with standard objects $\Delta_1,\dots,\Delta_j$. By \Cref{lem:silting_in_highest_weight_category}, the silting object in $\D^b(A_{\leq j})$ then coincides with the image of the projective objects in $A_{\leq j}$ in $\D^{b}(A_{\leq j})$. Thus, the derived endomorphism algebra is concentrated in degree $0$, showing that (1) holds. 
\end{proof}

\begin{remark}
Let $A'$ be an abelian category and $\Delta_1,\dots,\Delta_m\in A'$ a collection of objects giving rise to a full exceptional collection in $\D^b(A')$. Then \cite[Theorem 5.2]{Kra17} gives a further homological criterion for when these objects describe the standard objects in a highest weight structure on the module category of the associated silting object: this is the case if and only if the exceptional collection is strictly full, which coarsely means that the higher extension groups in the exact category of objects filtered by $\Delta_1,\dots,\Delta_m$ agree with those in $\D^b(A')$. 
\end{remark}

\section{Coherent chain complexes and grafting}\label{sec:CoherentChainComplexes}

The goal of this appendix is to introduce bounded coherent chain complexes in stable $\infty$-categories, define their totalization, and describe their grafting. In the main text, we will repeatedly use grafting of coherent chain complexes, as described in \cref{lem:grafting_coherent_chain_complexes}. Let us denote $I\coloneqq N(\{1\to 0\})\simeq (\Delta^1)^{\on{op}}$.

\begin{definition}\label{def:AppB_coherent}
Let $\D$ be a stable $\infty$-category and $n\geq 0$.
    \begin{enumerate}[(1)]
        \item An $n$-cube in $\D$ is a functor $I^n\to \D$. We denote by $\on{Cube}_n(\D)\coloneqq \on{Fun}(I^n,\D)$ the $\infty$-category of $n$-cubes in $\D$.
        \item A coherent chain complex $C_\ast$ in $\D$ concentrated in degrees $0$ through $n$ consists of an $n$-cube $C_\ast\colon I^n\to \D$, satisfying that $C_{\bf J}$ is zero for all
        \[{\bf J}\in \{0,1\}^n\backslash \bigcup_{0\leq i\leq n} \{(0,\dots,\underbrace{0}_{\on{pos.~i}},1,\dots,1)\}\,.\]
        We denote by \[
        \on{Ch}_{\geq 0,\leq n}(\D)\subset \on{Cube}_n(\D)
        \]
        the full subcategory spanned by coherent chain complexes concentrated in degrees $0$ through $n$. \qed
    \end{enumerate}
\end{definition}

\begin{remark}
\cref{def:AppB_coherent} differs a priori from the standard definition of a coherent chain complex given for instance in \cite[Def.~2.3.5]{Wal22}. Let us briefly explain why the two definitions are equivalent, as follows. Let $[0,n]^{\on{op}}_\ast$ denote the nerve of $1$-category obtained from the poset $[0,n]^{\on{op}}$ (with morphisms $i\to j$ if $i>j$) by adjoining a zero object $\ast$ and consider the quotient ${[0,n]^{\on{op}}_\ast}_{/\sim} $ obtained by identifying each composite $d^2\colon i\to i-2$ with the zero morphism. The standard definition of coherent $(n+1)$-term chain complex in a stable $\infty$-category $\C$ is a pointed functor ${[0,n]^{\on{op}}_\ast}_{/\sim} \to \C$, mapping $\ast$ to $0$. We denote by $\on{Fun}^\circ({[0,n]^{\on{op}}_\ast}_{/\sim} ,\C)$ the $\infty$-category of such pointed functors. \\

\noindent Now, there is a functor $\pi\colon I^n\to {[0,n]^{\on{op}}_\ast}_{/\sim}$, mapping each ${\bf J}\in \{0,1\}^n\backslash \bigcup_{0\leq i\leq n} \{(0,\dots,\underbrace{0}_{\on{pos.~i}},1,\dots,1)\}$ to $\ast$ and $(0,\dots,\underbrace{0}_{\on{pos.~i}},1,\dots,1)$ to $n+1-i$. It can be shown that the functor $\pi$ is the localization of $\infty$-categories at all morphisms mapped under $\pi$ to $\on{id}_\ast$, for instance by using \cite[\href{https://kerodon.net/tag/05ZD}{Proposition 05ZD}]{Ker} and the observation that the full subcategory of $I^n$ generated by $\pi^{-1}(\ast)$ is weakly contractible. This thus yields the required equivalence of $\infty$-categories 
\[
 \on{Ch}_{\geq 0,\leq n}(\D) \simeq \on{Fun}^\circ({[0,n]^{\on{op}}_\ast}_{/\sim} ,\C) \,,
\]
and both definitions of coherent chain complexes are equivalent.\qed
\end{remark}

\begin{remark}\label{rem:Dold_Kan}
    Combining the two versions of the $\infty$-categorical Dold--Kan correspondence, \cite[Thm.~1.2.3.7]{HA} and \cite[Cor.~4.1.2]{Wal22}, we find that for any stable $\infty$-category $\D$ there exists an equivalence of $\infty$-categories
    \[ \on{Ch}_{\geq 0,\leq n}(\D) \simeq \on{Fun}([0,n],\D)\]
    between $(n+1)$-term coherent chain complexes and $(n+1)$-term filtered chain complexes. \qed
\end{remark}

\begin{remark}
Let $\on{Fun}^\circ(I^{n+1}\backslash \{1,\dots,1\},\D)\subset \on{Fun}(I^{n+1}\backslash \{1,\dots,1\},\D)$ be the full subcategory spanned by morphisms $C_\ast\colon I^{n+1}\backslash \{1,\dots,1\}\to \D$ satisfying that $C_{\bf J}\simeq 0$ for all \[{\bf J}\in I^{n+1}\backslash \bigcup_{0\leq i\leq n+1} \{(0,\dots,\underbrace{0}_{\on{pos.~i}},1,\dots,1)\}\,.\] 
Each object in $\on{Fun}^\circ(I^{n+1}\backslash \{1,\dots,1\},\D)$ is the right Kan extension of its restriction to $\{0\}\times I^n$, and the restriction functor 
\[
\on{Fun}^\circ(I^{n+1}\backslash \{1,\dots,1\},\D)\to \on{Ch}_{\geq 0,\leq n}(\D)
\]
is thus an equivalence of $\infty$-categories, see \cite[\href{https://kerodon.net/tag/030M}{Tag 030M}]{Ker}.\qed
\end{remark}

Let us now introduce the total cofiber of a coherent chain complex, as follows.

\begin{definition}
    We define the total fiber functor as the composite
    \[
    \on{tot-fib}\colon \on{Ch}_{\geq 0,\leq n}(\D)\simeq \on{Fun}^\circ(I^{n+1}\backslash \{1,\dots,1\},\D)\xlongrightarrow{\on{lim}}\D\,.
    \]
    The total cofiber functor is similarly defined as the composite 
    \[
    \on{tot-cof}\colon \on{Ch}_{\geq 0,\leq n}(\D)\simeq \on{Fun}^\circ(I^{n+1}\backslash \{0,\dots,0\},\D)\xlongrightarrow{\on{colim}}\D\,.
    \]
The total cofiber will be taken as the default totalization, and thus also denoted by $\on{tot}\coloneqq \on{tot-cof}$.\qed
\end{definition}

\begin{remark}\label{rem:iterated_cofiber} The total (co)fiber of a coherent chain complex is equivalent to the iterated (co)fibers in the $n$ coordinate directions of the corresponding cube, taken in any order, see e.g.~\cite[Cor.~A.13]{DJW19}. In particular, there is an equivalence $\on{tot-cof}\simeq \on{tot-fib}[n]$.
\qed
\end{remark}

\begin{lemma}\label{lem:decomposition_of_coherent_chain_complexes}
Let $\D$ be a stable $\infty$-category. 
\begin{enumerate}[$(1)$]
\item There exists a pullback diagram of stable $\infty$-categories:
\[
\begin{tikzcd}
{\on{Ch}_{\geq 0,\leq n+1}(\D)} \arrow[d] \arrow[r] \arrow[rd, "\lrcorner", phantom] & {\on{Fun}(I,\D)} \arrow[d, "\on{ev}_1"] \\
{\on{Ch}_{\geq 0,\leq n}(\D)} \arrow[r, "\on{tot-fib}"]                              & \D                                            
\end{tikzcd}
\]
\item There exists a pullback diagram of stable $\infty$-categories:
\[
\begin{tikzcd}
{\on{Ch}_{\geq 0,\leq n+1}(\D)} \arrow[d] \arrow[r] \arrow[rd, "\lrcorner", phantom] & {\on{Ch}_{\geq 0,\leq n}(\D)} \arrow[d, "\on{tot-cof}"] \\
{\on{Fun}(I,\D)} \arrow[r, "\on{ev}_0"]                                              & \D                                                     
\end{tikzcd}
\]
\item More generally, the limit of the following diagram
\[
\begin{tikzcd}
                                                   &                                                                & {\on{Ch}_{\geq 0,\leq n-i}(\D)} \arrow[d, "\on{tot-cof}"] \\
                                                   & {\on{Fun}(I,\D)} \arrow[d, "\on{ev}_1"] \arrow[r, "\on{ev}_0"] & \D                                                   \\
{\on{Ch}_{\geq 0,\leq i}(\D)} \arrow[r, "\on{tot-fib}"] & \D                                                             &                                                     
\end{tikzcd}
\]
is equivalent to $\on{Ch}_{\geq 0,\leq n+1}$.
\end{enumerate}
\end{lemma}

\begin{proof}
For Part (1), since $I^{n+1}\simeq (I^{n+1}\backslash \{1,\dots,1\})^\triangleleft$, we have 
\begin{align*}
\on{Ch}_{\geq 0,\leq n+1}(\D)& \simeq \on{Fun}(I^{n+1}\backslash \{1,\dots,1\})^\triangleleft,\D)\times_{\on{Fun}(I^{n+1}\backslash \{1,\dots,1\},\D)} \on{Fun}^\circ(I^{n+1}\backslash \{1,\dots,1\},\D)\\
& \simeq \on{Fun}(I,\D) \times_{\D}\on{Fun}(I^{n+1}\backslash \{1,\dots,1\},\D)\times_{\on{Fun}(I^{n+1}\backslash \{1,\dots,1\},\D)} \on{Fun}^\circ(I^{n+1}\backslash \{1,\dots,1\},\D)\\
&\simeq  \on{Fun}(I,\D) \times_{\D} \on{Ch}_{\geq 0,\leq n}(\D)\,,
\end{align*}
where the second equivalence uses the universal property of the limit over $I^{n+1}\backslash \{1,\dots,1\}$. This establishes Part (1). The proof of Part (2) is similar. Part (3) can be reduced to an iterated application of Parts (1) and (2).
\end{proof}

\noindent\Cref{lem:decomposition_of_coherent_chain_complexes} implies that we can ``graft'' coherent chain complexes as follows:

\begin{lemma}\label{lem:graftcplx}\label{lem:grafting_coherent_chain_complexes}
Let $m,n\geq 0$ and $i\in [0,n]$. Let $D_\ast\in \on{Ch}(\C)_{\geq 0,\leq m}$ and $C_\ast\in \on{Ch}(\C)_{\geq 0,\leq n}$, such that $C_i\simeq \on{tot}(D_\ast)$. Then there exists a canonical chain complex $\tilde{C}_\ast
\in \on{Ch}(\C)_{\geq 0,\leq n+m}$ with
\[
\tilde{C}_j=\begin{cases} 
C_j & 0\leq j <i\\
D_{j-i}& i\leq j\leq i+m\\
C_{j-m}[-m]& j>i+m
\end{cases}
\]
satisfying 
$\on{tot}(\tilde{C}_\ast)\simeq \on{tot}(C_\ast)$. The complex $\tilde{C}_*$ is said to be the grafting of $D_\ast$ onto $C_\ast$.
\end{lemma}

\begin{proof} 
There is a canonical map $\alpha\colon C_i\to \on{tot-fib}(C_{\ast\leq i-1})$, which can be considered as an element $\alpha\in \on{Fun}(I,\D)$ that satisfies $\on{ev}_{0}(\alpha)\simeq \on{tot-cof}(D_\ast)$. From the limit diagram in \Cref{lem:decomposition_of_coherent_chain_complexes}.(3) we thus obtain a coherent chain complex $\tilde{C}_{\ast\leq i+m}\in \on{Ch}_{\geq 0,\leq i+m}(\D)$ satisfying that 
\[\on{cof-tot}(\tilde{C}_{\ast\leq i+m})\simeq \on{tot-cof}(C_{\ast \leq i})\,.\] 
By \Cref{lem:decomposition_of_coherent_chain_complexes}.(3), there is a canonical morphism $\beta\colon \on{tot-cof}(C_{\ast > i})\to \on{fib-tot}(C_{\ast \leq i})$. Using the apparent equivalence $\on{tot-fib}(\tilde{C}_{\ast\leq i+m})[m]\simeq \on{tot-fib}(C_{\ast\leq i})$, the morphism $\beta$ gives rise to a morphism $\g:\on{tot-cof}(C_{\ast > j})[-m]\to \on{tot-fib}(\tilde{C}_{\ast\leq i+m})$. By using \Cref{lem:decomposition_of_coherent_chain_complexes}.(3) a last time, we use the morphism $\g$ to construct to the required coherent chain complex $\tilde{C}_\ast$ with $\on{tot}(\tilde{C}_\ast)\simeq \on{tot}(C_\ast)$.\end{proof} 

\bibliography{main} 
\bibliographystyle{alpha}

\textsc{Roger Casals: University of California Davis, Dept. of Mathematics, USA.\\}
\textit{Email address:} \texttt{casals@ucdavis.edu}
\vspace{0.2cm}

\textsc{Merlin Christ: Mathematisches Institut, Universität Bonn, Endenicher Allee 60, 53115 Bonn, Germany.\\}
\textit{Email address:} \texttt{christ@math.uni-bonn.de}

\end{document}